\documentclass[11pt]{amsart}

\usepackage{amsmath,amssymb,amsthm,tikz, nicefrac, mathrsfs,upgreek,enumitem}

\definecolor{my-blue}{rgb}{	0.05,0.1,0.5}
\definecolor{my-red}{rgb}{0.4,0.0,0.0}
\definecolor{my-green}{rgb}{0.0,0.5,0.0}
\usepackage{comment}
\usepackage[colorlinks=true, linkcolor=my-red, urlcolor=my-red, citecolor=my-red]{hyperref}

\usepackage[margin=1.1in]{geometry}

\usepackage{dsfont} 

\usepackage[utf8]{inputenc}

\usepackage{nicefrac}

\usepackage{multicol}
\PassOptionsToPackage{dvipsnames}{xcolor}
\usepackage{tikz,xcolor}

\usetikzlibrary{shapes,arrows}
\usetikzlibrary{hobby}
\usetikzlibrary{decorations.markings}
\usetikzlibrary{decorations.pathreplacing}

\usepackage[only,llbracket,rrbracket]{stmaryrd}

\def\namedlabel#1#2{\begingroup
    #2%
    \def\@currentlabel{#2}%
    \phantomsection\label{#1}\endgroup
}

\newtheorem{theorem}{\color{my-blue}{\sc Theorem}}[section]
\newtheorem{lemma}[theorem]{\color{my-blue}{\sc Lemma}}
\newtheorem{proposition}[theorem]{\color{my-blue}{\sc Proposition}}
\newtheorem{corollary}[theorem]{\color{my-blue}{\sc Corollary}}

\newtheorem{definition}[theorem]{\color{my-blue}{\sc Definition}}
\numberwithin{equation}{section}
\newtheorem{remark}[theorem]{\color{my-blue}{\sc Remark}}

\newtheorem{thmalpha}{\color{my-blue}{\sc Theorem}}

\newcommand{\be}{\begin{equation}}
\newcommand{\ee}{\end{equation}}

\providecommand{\norm}[1]{\Vert#1\Vert}

\newcommand{\eone}{\textup{e}_1}
\newcommand{\etwo}{\textup{e}_2}
\newcommand{\TV}[1]{{\big\lVert #1 \big\rVert}_{\normalfont
\text{TV}}}

\def\sl{\Delta}

\def\bE{\mathbb{E}}
\def\bP{\mathbb{P}}
\def\bQ{\mathbb{Q}}
\def\bR{\mathbb{R}}
\def\bZ{\mathbb{Z}}

\def\TF{\textup{TF}}
\def\TR{\textup{TR}}

\def\Tper{T^{\textup{per}}}
\def\gper{\gamma^{\textup{per}}}

\def\np{\tilde{n}}

\def\epb{\bar{\eta}^{\per}}

\def\WD{\textsf{W}_1}

 \def\Z{\bZ}
\def\Q{\bQ}
\def\R{\bR}

\def\LL{\mathcal{L}}
\def\N{\mathbb{N}}
\def\P{\bP}

\def \slz{\tilde{\Delta}}

\def \hLp{h^{\textup{pLPP}}}
\def \hL{h^{\textup{LPP}}}

\def \hLpt{\tilde{h}^{\textup{pLPP}}}

\def \hTASEP{h^{\textup{TASEP}}}
\def \hpTASEP{h^{\textup{pTASEP}}}

\def \hA{h^{\textup{ASEP}}}

\def \hAeps{\mathfrak{h}^{\textup{ASEP},\varepsilon}}
\def \hApeps{\mathfrak{h}^{\textup{pASEP},\varepsilon}}

\def \hLp{h^{\textup{pLPP}}}

\def \hLpeps{\mathfrak{h}^{\textup{pLPP},\varepsilon}}
\def \hfix{\mathfrak{h}}

\def \hfixp{\mathfrak{h}^{\textup{per}}}
\def \hfixpt{\tilde{\mathfrak{h}}^{\textup{per}}}
\def \hfixpa{\mathfrak{h}^{(n)}}

\def \dLeps{d^{\textup{LPP},\varepsilon}}

 \def \LL{ \mathcal{L}^{\textup{LPP}}}
 
 \def \LLpaeps{ \mathcal{L}^{\textup{LPP},\varepsilon}_n}
 \def \LLeps{ \mathcal{L}^{\textup{LPP},\varepsilon}}
 
 \def \LLpeps{ \mathcal{L}^{\textup{pLPP},\varepsilon}}

 \def \Lp{ \mathcal{L}^{\textup{per}}}
 \def \Lpb{ \bar{\mathcal{L}}^{\textup{per}}}
 \def \Lpp{ \mathcal{L}^{\textup{per}}_p}
 \def \tLp{ \tilde{\mathcal{L}}^{\textup{per}}}
 \def \Lpa{ \mathcal{L}_n}     
 \def \Lpat{ \bar{\mathcal{L}}_n}     
 \def \Lsta{ \mathcal{L}^{\textup{sta}}}
 
  \def \Pper{\mathbb{P}_{\textup{per}}}
  \def \Pfull{\mathbb{P}_{\textup{full}}} 
   
  \def \Pex{\bar{\mathbf{P}}}
  \def \Ppatch{\mathbf{P}_{\textup{\textbf{pat}}}}

 \def \Tpat{T^{(n)}} 
 \def  \Wast{\textsf{W}_{\ast}}

\def \pinp{\pi^{n,+}}
\def \pinm{\pi^{n,-}}

\def \pinpmt{\bar{\pi}^{n,\pm}}
\def \pinpt{\bar{\pi}^{n,+}}
\def \pinmt{\bar{\pi}^{n,-}}

\def \one{\normalfont\mathsf 1}

\def \zero{\normalfont\mathsf 0}
\def \onep{\normalfont\mathsf 1'}
\def \zerop{\normalfont\mathsf 0'}

\def \lbr {\llbracket}
\def \rbr {\rrbracket}
\def\modi{\textup{mod}}
\def\Aup{\normalfont\textsf{A}}
\def\Bup{\normalfont\textsf{B}}
\def\Cup{\normalfont\textsf{C}}

\def\per{\textup{per}}
\def\Ber{\textup{Ber}}
\def\full{\textup{full}}
\def\ASEP{\textup{ASEP}}
\def\pASEP{\textup{pASEP}}

\def\UC{\textup{UC}}
\def\UCp{\textup{UC}^{\textup{per}}}
\def\USC{\mathcal{C}_{\textup{up}}}

\def\USCper{\mathcal{C}^{\textup{per}}_{\textup{up}}}

\def\SetT{(\mathbb{T}\times[0,1])^2_{\uparrow}}
\def\SetTw{(\mathbb{T}\times[0,\tilde{\omega}])^2_{\uparrow}}
\def\SetR{(\mathbb{R}\times [0,1])^{2}_{\uparrow}}

\usepackage{stackengine}

\def\E{\bE}
\def\P{\bP} 

\definecolor{partcolor1}{rgb}{0.0,0.5,0.0}
\definecolor{partcolor2}{rgb}{0.0,0.5,0.0}

\definecolor{darkgreen}{rgb}{0.0,0.5,0.0}
\definecolor{darkblue}{rgb}{0.5,0.1,0.5}
\definecolor{deepblue}{rgb}{0.25,0.41,0.88}
\definecolor{nicosred}{rgb}{0.65,0.1,0.1}
\definecolor{light-gray}{gray}{0.7}
\allowdisplaybreaks[1]

\RequirePackage{datetime} 
\allowdisplaybreaks
\begin{document}
\usdate
\title[Periodic directed landscape]
{Periodic directed landscape}
\author{Amol Aggarwal}
\address{Amol Aggarwal, Stanford University, United States}
\author{Ivan Corwin}
\address{Ivan Corwin, Columbia University, United States}
\author{Dominik Schmid}
\address{Dominik Schmid, University of Augsburg, Germany}
\keywords{exclusion processes, periodic directed landscape, KPZ universality class}
\date{\today}
\begin{abstract}
We construct the periodic directed landscape, which is the conjectured scaling limit for periodic models in the Kardar--Parisi--Zhang universality class. 
We establish the convergence of periodic exponential last passage percolation to the periodic directed landscape. Moreover, we confirm  conjectures by Baik, Liao and Liu on the local structure of the periodic KPZ fixed point, and establish the convergence of periodic ASEPs to periodic KPZ fixed points, coupled according to the same periodic directed landscape. 
Our main tool to construct the periodic directed landscape (and prove convergence to it) is a technique for gluing full-space directed landscapes (and their prelimits), which is of independent interest.
\end{abstract}
\maketitle
\tableofcontents
\vspace*{-0.6cm}

\section{Introduction and main results}

\subsection{Preface}

The asymmetric simple exclusion process (ASEP) is a paradigmatic model of interacting particle systems, playing a role in non-equilibrium statistical physics akin to that of the Ising model in equilibrium statistical physics. It was introduced in the probability literature in 1970 \cite{S:InteractionMP}, though a finite-volume version already appeared two years earlier in the biology literature
\cite{https://doi.org/10.1002/bip.1968.360060102}. 
The model, a continuous time Markov process, can be described in terms of particles jumping after exponential waiting times with asymmetric rates left or right in one dimension, subject to the exclusion rule that no site can hold more than one particle (so jumps that would violate that rule are suppressed). Integrating this particle description yields a height function description; this brings the model in contact with the study of random growth models, and through a weakly asymmetric scaling limit \cite{Bertini1997} yields the continuum Kardar-Parisi-Zhang (KPZ) stochastic PDE. The model also relates, e.g. \cite{R:PDEresult} and \cite{Bertini1995}, to directed polymer and last passage percolation models, among many other types of systems.

Most of the study on ASEP has focused on questions about its long-time behavior. Initially, going back to \cite{1a062f39-83e4-32c5-854f-60bc1b9ebcc1}, this involved understanding its invariant measures and ergodicity. Starting in the early 1980s with \cite{R:PDEresult}, there was great interest in deriving hydrodynamic limits of ASEP and related interacting particle systems in which PDEs (of Hamilton-Jacobi type for the height function process) arise as law of large number scaling limits of the system when time and space and the height function are all viewed in the same scale. Once those limits were understood, attention shifted to the fluctuations around them. In the related context of the KPZ stochastic PDE, \cite{PhysRevLett.56.889} predicted universality of the fluctuations of stochastic growth models under a $3:2:1$ scaling, namely when time is scaled by $\varepsilon^{-3}$, space by $\varepsilon^{-2}$ and the height function scaled down by $\varepsilon^{1}$. In the forty years since that work, the central challenge has been to understand the precise nature of those universal fluctuations and demonstrate their universality. This has defined the study of the \emph{KPZ universality class} (for more, see reviews such as \cite{C:KPZReview, Quastel2015}).

The simplest geometry for studying ASEP and related models is 
the full-space, i.e., when particles jump on $\Z$. Especially amenable to asymptotics is the special case of ASEP called the totally asymmetric simple exclusion process (TASEP), when jumps only go in one direction. The first fluctuation results for TASEP under the above KPZ scaling date to the late 1990s and came as one-point distributions \cite{baik1999distribution, J:KPZ} or spatial multi-point distributions \cite{Praehofer2002, johansson2003discrete} for a few special types of initial data. A decade ago, this was strengthened in \cite{MQR:kpz}, which gave an exact description through transition probabilities of the limiting Markov process describing the TASEP height evolution. This was called the \emph{KPZ fixed point}. These results on the TASEP all relied on its solvability through determinantal point processes.

Returning to ASEP, this determinantal structure is mostly lost; later progress instead came through methods closely related to quantum integrable systems. Indeed, the generator of periodic ASEP can be mapped to the Hamiltonian of the periodic Heisenberg XXZ spin chain via a similarity transform \cite{PhysRevLett.68.725}, and that Hamiltonian was known to be diagonalized via the coordinate Bethe ansatz \cite{PhysRev.162.162}. The limiting one-point distribution for the ASEP under step initial data was first identified in \cite{TracyWidom2009ASEP}, nearly a decade after the TASEP case, by analyzing exact formulas derived from the coordinate Bethe ansatz (see also earlier TASEP work of \cite{Schutz1997}).

More recently, \cite{ACH:ScalingASEP,ACH:Fixed} proved the full KPZ fixed point scaling limit of ASEP under arbitrary initial data through a different framework, using the Yang--Baxter equation to map the ASEP to certain Gibbsian line ensembles \cite{Aggarwal2024} and a characterization theory \cite{AH:StrongCharacterization} of the latter. These works \cite{ACH:ScalingASEP,ACH:Fixed} also determined the scaling limit of the \emph{basic coupling} of ASEP, which builds the model simultaneously under all initial data from a space-time Poissonian field. The \emph{directed landscape}, constructed earlier in \cite{DOV:DirectedLandscape} as a limit of certain exactly solvable last passage percolation models, arises as the scaling limit of this coupling.

Already when it was introduced, it was understood that the ASEP should also be studied in other geometries, and in contact with boundary conditions, which can induce new limiting behaviors. Natural geometries include periodic (where particles hop around a ring) and half-space or interval (where there can be sources and sinks of particles at the boundary).

The above full-space methods become intricate and often no longer apply in such geometries. The most promising direction for generalizing these methods is to the half-space. For TASEP, much of the determinantal structure is replaced by Pfaffian analogs (though we will still refer to this as determinantal structure). For instance, recently the KPZ fixed point was defined as a limit of the half-space TASEP in \cite{zhang2025tasephalfspace} using this determinantal structure, and the half-space directed landscape was subsequently defined in \cite{dauvergne2026directedlandscapehalfspace} as a limit of coupled half-space TASEPs. For the half-space ASEP (which is no longer Pfaffian), only one-point distributional results under empty initial conditions have been identified \cite{BBCW:HalfspaceASEP,H:Boundary}. These works and others involving half-space Gibbsian line ensembles point to the possibility of developing the full-space program of \cite{Aggarwal2024,AH:StrongCharacterization,ACH:ScalingASEP,ACH:Fixed} for half-space models.

In finite (such as periodic) geometries, some exact Bethe ansatz formulas are known; see, for example, \cite{Liu2020} for the periodic ASEP. However, their asymptotic analysis is complicated, since they necessarily involve the Bethe roots, which are solutions to a certain family of highly nonlinear equations. The only finite-volume KPZ growth model for which these Bethe roots have been fully analyzed, leading to large-time asymptotics, is the periodic TASEP (and two close relevatives, the discrete time TASEP \cite{Liao2022} and PushASEP \cite{10.1214/24-AOP1738}), starting with the one-point distribution works of \cite{Prolhac_2015,BL:TASEPring} for special initial data. Subsequent works expanded to consider other initial data and asymptotic regimes \cite{BL:Subscale,BLS:LimitingOnePoint}  as well as  multi-point distributions 
\cite{BL:Multipoint,
BL:PeriodicGeneral, Prolhac_2020}; see 
\cite{BaikSurvey} for a survey of some of these developments. For this model, \cite{BLL:FixedPointGeneral} recently provided consistent multi-space-time limiting distributions for general initial data, thereby constructing a two-dimensional field, which they called the periodic KPZ fixed point with the given initial condition.

For non-determinantal models in finite geometries, such as the periodic ASEP, no scaling limit results had been known, even for their one-point distributions. Moreover, while the half-space ASEP should still be mappable through the Yang-Baxter equation to a line ensemble \cite{H:Boundary}, we have little reason to suspect that this would be possible for the periodic ASEP. The same holds for the ASEP on an interval, where fluctuation scaling limits are also not available; see \cite{S:MixingTASEP,ES:HighLow} for partial results in the form of tail estimates on the open TASEP height function, and \cite{BL:Cumulants} for a non-rigorous derivation of cumulant formulas for the open KPZ equation.

\begin{figure}
\centering
\includegraphics[width=.9\textwidth]{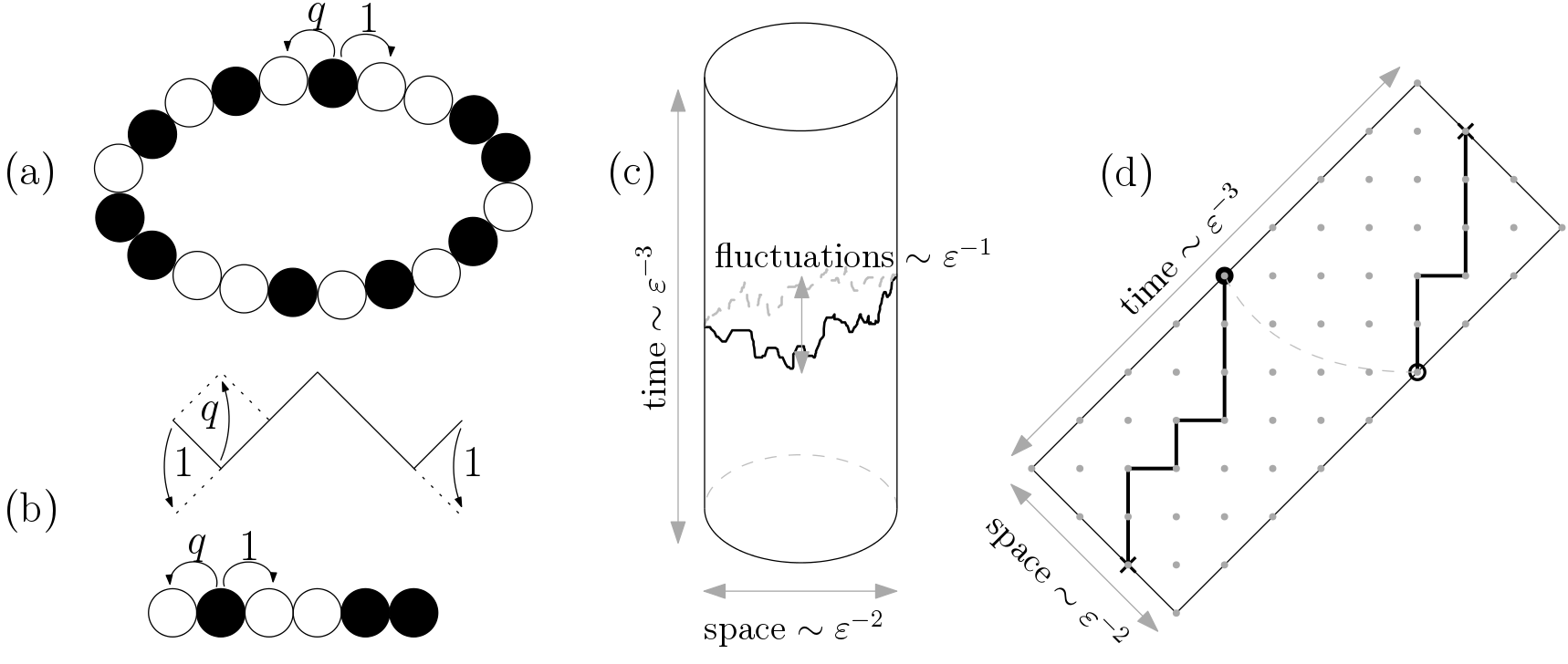}
\caption{\label{fig:intro} Models in a periodic geometry: (a) depicts periodic ASEP where particles (black) jump subject to exclusion with rates $q$ and $1$ as illustrated; (b) depicts the mapping between ASEP and a height function ($-1$ slope increments enter sites occupied by particles and $+1$ slope increments enter empty sites); (c) illustrates the characteristic KPZ scaling underwhich we show that ASEP and other models have a scaling limit (time, space and height function fluctuations live in the scales $\varepsilon^{-3}$, $\varepsilon^{-2}$ and $\varepsilon^{-1}$); (d) depicts the periodic last passage percolation model where the last passage time from the starting to ending points (both depicted by small $\times$ symbols) is the maximum over all up-right periodic paths (i.e., that wrap around the fundamental domain as depicted) of the sum of independent exponential weights (depicted by the small grey dots) along the paths. In this model, time is measured along the diagonal, space corresponds to the width of the strip, and the height function is replaced by last passage times.}
\end{figure}

The purpose of this paper is to develop a gluing technique for constructing and proving space-time scaling limits for stochastic models in periodic geometries (as depicted in Figure~\ref{fig:intro}). Periodic KPZ models have been of interest among the probability and stochastic analysis communities for several decades; see, for instance, the extensive references in the recent survey \cite{Gu2025} as well as early work such \cite{PhysRevLett.68.725} on the periodic stochastic six-vertex model and ASEP, and \cite{aa6db8ca-f0b8-3af3-a8cb-608089d8d46f} on periodic stochastically forced Hamiltonian Jacobi equations. 

The key input for us is the convergence of the full-space models to their scaling limits, which are then \emph{glued} together to construct the periodic scaling limits. Since the scaling limit of ASEP (and exponential last passage percolation) in the full-space is known, we are able to prove its convergence of its periodic version to the constructed scaling limit. Specifically, we:
\begin{itemize}[leftmargin=*]
\item Construct and uniquely characterize the \emph{periodic directed landscape};
    \item Prove convergence of exponential last passage percolation to the periodic directed landscape;
    \item Prove convergence of periodic ASEP under the basic coupling to \emph{periodic KPZ fixed points} coupled via variational problems to the periodic directed landscape. 
\end{itemize}

While the approach developed in this paper circumvents use of any exact formulas or structure for periodic models,  we believe that it is still an important problem to understand how to adapt them to the periodic setting (or other boundary conditions), and suspect that other applications (e.g., large deviations \cite{deGier2021}) will follow from that understanding.

In what follows we give an overview of these main results (in some cases with informal or less general statements than in the main body), discuss previous results that give context for our work, and finally provide the key idea behind our gluing technique.
\subsection{Overview of main results}

For the purpose of this introduction, we will consider periodic ASEP with $N$ sites and $k=N/2$ particles, assuming $N$ to be even. Section \ref{sec:ASEPconvergence} deals with arbitrary filling fraction and includes the precise statement of the below convergence theorem, as well as a precise definition of periodic ASEP and its height function. Figure \ref{fig:intro} (a) depicts an example of ASEP with $N=20$ and $k=10$. At any time, each site in this periodic domain is either occupied with exactly one particle (black), or empty (white). The particles jump independently in the clockwise direction after exponentially distributed waiting times of rate $1$, and in the counter-clockwise direction at rate $q$ for some $q \in [0,1)$. Jumps are subject to the exclusion rule, i.e., a jump is performed if and only it the target site is vacant at that time. 

The basic coupling defines the evolution of all choices of initial data simultaneously with respect to the same underlying jump clocks. Specifically, between each pair of neighboring sites, there are clockwise and counter-clockwise Poissonian clocks (i.e., Poisson point processes) of rates $1$ and $q$, respectively. When a clock rings
, if there is a particle at the site and a hole in the clockwise or counter-clockwise direction (depending on which clock rang), the particle and hole switch. This couples the evolution over all initial data to the same field of Poissonian clocks.

Periodic ASEP also induces a Markovian evolution on height functions. For the sake of this introduction and since we are currently assuming $1/2$-filling, we will define the height function as follows (in the main body, we use a slightly different convention, related by an affine transformation). Given a particle configuration, its height function $h:\Z/N\Z \to \Z$ and has $\pm 1$ increments, i.e., $h(x)-h(x-1)\in \{-1,+1\}$ where it takes value $-1$ if site $x$ is occupied by a particle and $+1$ if it is empty. At time zero, we assume $h(0)=0$ and then update the height function locally as particles in ASEP move, as depicted in Figure \ref{fig:intro} (b) (i.e., only the height at $x$ updates by $-2$ or $+2$ when the particle at $x$ moves clockwise or counter-clockwise, respectively). Owing to this, the height function contains more information than the particle configuration, namely, it tracks the net flux or winding of particles over time.

For any $M\in \N$, and for $i\in \{ 1, 2, \ldots , M \}$, let $h^{N,(i)}_0$ denote $M$ different choices of initial data and let $h^N(h^{N,(i)}_0;y,t)$ denote the height function at position $y\in \Z/N\Z$ resulting from running ASEP with initial data corresponding to $h^{N,(i)}_0$ until time $t$, subject to the basic coupling (so the evolutions over different $i$ are coupled). Assume that the ASEP initial data, at the particle level, can be sandwiched stochastically between two Bernoulli product measures with particle probabilities $1/2 + CN^{-1/2}$ and $1/2-CN^{-1/2}$ for some constant 
$C>0$. Finally, define the KPZ rescaling of the ASEP height function process, noting that $N^{-1/2}$ plays the role of $\varepsilon$ from the KPZ scaling here:
$$
\mathfrak{h}^N(h^{N,(i)}_0;Y,T):= N^{-1/2} \bigg(h^N\Big(h^{N,(i)}_0;NY,\frac{2NT}{1-q}\Big) - \frac{N^{3/2} T}{1-q}\bigg).
$$
Our first result is an informal statement of Theorem \ref{thm:ASEPcirclePDL}.

\begin{thmalpha}\label{thm:ConvergenceASEPMAIN} 
Assume $\mathfrak{h}^N(h^{N,(i)}_0;\cdot,0)$ converges weakly, jointly over $i \in \lbr M\rbr$, to some initial data $\hfix_0^{\per,(i)}(\cdot)$. Then the space-time process $\mathfrak{h}^N(h^{N,(i)}_0;\cdot,\cdot)$ likewise converges weakly, jointly over $i \in \lbr M\rbr$, to a limit $\hfixp(\hfix_0^{\per,(i)};\cdot,\cdot)$.
\end{thmalpha}

The above theorem shows that the limit of periodic ASEP is independent of $q\in [0,1)$. Thus, any limit theorem known for periodic TASEP ($q=0$) now also applies to periodic ASEP. In particular, \cite{BLL:FixedPointGeneral} (see also earlier work  \cite{BL:Subscale, BL:TASEPring, BL:Multipoint, BL:PeriodicGeneral,BaikSurvey,BLS:LimitingOnePoint}) computed the space-time multi-point distributions of the above limit for general (single) initial condition, and called the resulting stochastic process the \emph{periodic KPZ fixed point}. Those formulas now also follow for the periodic ASEP. This convergence is the periodic analog of the results proved in \cite{ACH:ScalingASEP,ACH:Fixed} for ASEP on the line.

The stochastic process $(\hfix_0^{\per},y,t)\mapsto \hfixp(\hfix_0^{\per};y,t)$ couples the periodic KPZ fixed point over multiple initial data. As in the full-space (see \cite{MQR:kpz,DOV:DirectedLandscape,NQR:Fixed}) this process is defined in terms of variational problems involving the initial data and the \emph{periodic directed landscape}, a process that we construct and characterize now. 

The following result informally summarizes several results including Propositions~\ref{pro:CauchySequence} and \ref{pro:PDLContinuity}, and parts of Theorems \ref{thm:PeriodicDirectedLandscape} and \ref{thm:ASEPcirclePDL} and Definition \ref{def:PeriodicKPZfixed}. Let $\mathbb{T}$ denote the periodic interval $[0,1]$, where $0$ and $1$ are identified.

\begin{thmalpha}\label{thm:PeriodicDirectedLandscapeMAIN}
The stochastic  process $(\hfix_0^{\per},y,t)\mapsto \hfixp(\hfix_0^{\per};y,t)$ is given via 
$$
\hfixp\big(\hfixp_0;y,t\big) := \sup_{x\in \mathbb{T}} \big(\hfixp_0(x) + \Lp(x,0;y,t)\big)
$$
where $\Lp \colon \{ (x,s;y,t)\in (\mathbb{T}\times \mathbb{R})^2 \colon s < t \}\rightarrow \R$, which we call the \emph{periodic directed landscape}, is the unique (in law) random continuous function satisfying:
\begin{enumerate}[leftmargin=*]
\item\label{eq:Property4DL} The function $\Lp$ looks locally in space and time like the full-space directed landscape $\mathcal{L}$  (see Theorem \ref{thm:PeriodicDirectedLandscape} for the precise sense in which this should hold);
\item \label{eq:Property2DL} The function $\Lp$ has independent increments, i.e., for any finite collection of disjoint intervals $(s_i,t_i)_{i\in I}$, the processes $\mathcal{L}(\cdot,s_i;\cdot,t_i)$ are independent over $i\in I$;
\item\label{eq:Property3DL} The metric composition law holds with respect to $\Lp$, i.e., almost surely for all $s<u< t$ and $x,y \in \mathbb{T}$
\begin{equation*}
\Lp(x,s;y,t) = \sup_{z \in \mathbb{T}} \big(\Lp(x,s;z,u) + \Lp(u,r;y,t)\big).
\end{equation*} 
\end{enumerate}
\end{thmalpha}

The variational representation for the periodic KPZ fixed point given above enables us to resolve several conjectures, mostly from \cite{BLL:FixedPointGeneral}, regarding its properties (see Section \ref{sec:KPZfixed}).
Later in the introduction, in Section \ref{sec:introgluing}, we describe our construction of the periodic directed landscape from gluing together patches of full-space directed landscapes. Now we describe how the periodic directed landscape arises as a limit of another model.


The periodic exponential last passage percolation (periodic LPP) is depicted in Figure \ref{fig:intro} (d), and carefully defined in Section \ref{sec:ConvergenceLPP}. Fix $N\in \mathbb{N}$ and $k\in \lbr N-1 \rbr$. Consider the strip of width (as measured by the taxi-cab distance between the left and right sides of the strip) $N$, and populate every lattice point except on the right side with independent and identically distributed exponential random variables, which we call \emph{weights}. An up-right lattice path takes unit steps up or to the right, except that we identify lattice points on the left side of the strip with those on the right side of the strip that are $N-k$ horizontal steps right and $k$ horizontal steps down. In Figure \ref{fig:intro} (d) $N=6$, $k=2$ and such a lattice path is illustrated where the black and white disks are identified. The weight of a path is the sum of the random weights over all lattice points through which the path crosses (with the convention that we exclude the endpoint from the sum), and the last passage time between a starting and ending point is the maximum over all up-right paths connecting those points of the weight of the paths. 

For periodic LPP, $N$ plays the role of the spatial width, while $k$ relates to the slope of the cylinder on which the model is defined (it is the same as the filling of periodic TASEP, which can be embedded in a different manner from the basic coupling into this model, as discussed in Remark \ref{rem:TASEPandLPP}). We identify the anti-diagonal and diagonal location of a starting/ending point with spatial and, respectively, temporal coordinates. 

Below, KPZ scaling means that $N$ and $k$ scale as $\varepsilon^{-2}$; the time horizon scales as $\varepsilon^{-3}$; and the last passage time, after centering by its law of large numbers, is scaled by $\varepsilon$. Implicit in this scaling are also explicit scaling constants (dependent on $k/N$), as well; see Section \ref{sec:ConvergenceLPP}, and Theorem \ref{thm:PeriodicLPPConvergence} in particular, for a precise statement.

\begin{thmalpha}\label{thm:ConvergenceLPPMAIN} Under KPZ scaling, the field of last passage times in periodic exponential last passage percolation converges to the periodic directed landscape.
\end{thmalpha}

\subsection{Further directions}\label{sec:furtherdirections}
We describe some further directions prompted by our work here.

In this work, we have considered ASEP and LPP and extracted scaling limits of them in periodic geometries. The fundamental input in both cases is the convergence of these models in the full-space geometry to the corresponding limits (together with some technical control on second class particles for ASEP, and on geodesics for LPP). This convergence to the KPZ fixed point and directed landscape in the full-space geometry is also known for various other models, including the KPZ stochastic PDE \cite{W:Landscape}, log-gamma polymer \cite{zhang2025convergenceloggammapolymerdirected} and stochastic six-vertex model \cite{ACH:ScalingASEP,ACH:Fixed}. In fact, the framework of \cite{ACH:ScalingASEP,ACH:Fixed}  should yield similar full-space convergence for all models in the higher-spin stochastic vertex model hierarchy \cite[Figure 1]{Aggarwal2024}. We expect that our gluing approach can be used in each of these cases (with some additional technical control required on a case-by-case basis).

Our construction of the periodic directed landscape is via gluing and thus yields a rather inexplicit object. In the full-space, the directed landscape can be constructed from the Airy line ensemble \cite{DOV:DirectedLandscape,DV:LongestSub} or from infinitely many Brownian motions \cite{dauvergne2026directedlandscapebrownianmotion}. It would be interesting if such a description were possible in the periodic setting. Returning to the full-space, there is also an (albeit somewhat inexplicit) characterization of the directed landscape in terms of certain properties (KPZ fixed point marginals, monotonicity, shift invariance and its jointly invariant measures) \cite{DZ:Characterization} and it would be interesting if such a characterization could be made in the periodic setting. Indeed, the jointly invariant measures for the periodic KPZ fixed point coupled to the periodic directed landscape are already known (called the \emph{periodic stationary horizon}) from \cite{Corwin2026}.

There are other geometries on which ASEP, LPP and related models are naturally defined. Perhaps the most studied is the interval, corresponding to open ASEP (with sources and sinks next to site $1$ and $N$) or LPP on the strip (with different rates on the left and right boundaries and where paths cannot wrap around). Given the recent works \cite{zhang2025tasephalfspace} and \cite{dauvergne2026directedlandscapehalfspace} defining the half-space KPZ fixed point and half-space directed landscape as limits of TASEP, it should be possible to implement our gluing approach to construct the interval or open version of both the KPZ fixed point or directed landscape, and show that they arise as limits from open TASEP. Note that convergence results in half-space are at the moment more limited than in full-space; the methods of \cite{ACH:ScalingASEP} have not yet been extended to show ASEP convergence in half-space. 

\begin{figure}
\centering
\includegraphics[width=.5\textwidth]{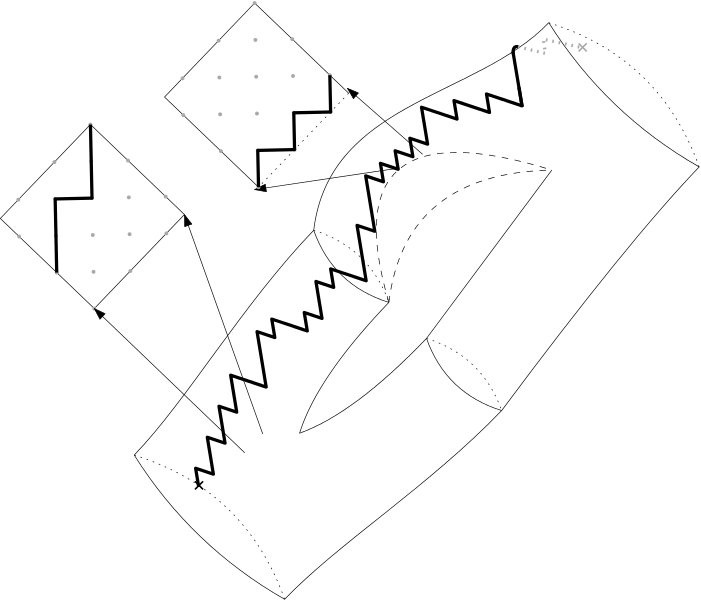}
\caption{\label{fig:intropants} An example of a time-varying domain for last passage percolation. Initially (time is measured in the diagonal direction) the domain is periodic but at some time it splits into two periodic domains and then later one of them (on the upper branch) opens into a model on a strip with different rates on the two sides, and finally the two branches join together again. The last passage time optimizes the sum of weights over all paths in this geometry.}
\end{figure}

More broadly, it should be possible to use our gluing approach to define the \emph{directed landscape on general space-time manifolds}, which would provide a natural random directed metric on them. For example, Figure \ref{fig:intropants} depicts a genus one surface with several boundary components; it describes a geometry that transitions from periodic to a pair of pants, then splits the seam on one leg and finally mends it all back together. As described in the caption, it is possible to define LPP in this geometry, though we would not expect it to admit analyzable exact formulas. Still, since it looks locally at all times like full- or half-space LPP, it should be possible through our gluing framework to extract a scaling limit and define the corresponding directed landscape. 

Our gluing approach ultimately relies on couplings and finite speed of discrepancy bounds; these properties are not exclusive to KPZ class models. For example, one could consider deterministic dynamics driven with random initial data. For the Toda lattice, \cite{aggarwal2026fluctuationstodalattice} proves a scaling limit for fluctuations on the line, and it is plausible that our gluing approach could be useful in both identifying and proving the scaling limit of the periodic Toda lattice.

Going back to early studies, e.g. \cite{1a062f39-83e4-32c5-854f-60bc1b9ebcc1}, it has been of interest to understand the convergence of ASEP and related models to their invariant measures. This can be quantified in terms of \emph{mixing times}. For TASEP on a torus of length $N$ with $k$ particles, \cite{SS:TASEPcircle} proved that its mixing time is of order $N^{2}\min(k,N-k)^{-1/2}$. They further conjectured the same statement for the ASEP (this type of conjecture has been folklore, arguably dating at least back to \cite{PhysRevLett.68.725}). We believe that our fluctuation convergence results would serve as an important step towards proving such mixing time results.

\subsection{Gluing method}\label{sec:introgluing}

In the study of PDEs and their stochastic analogs, it is often easier to first localize in time and space, and then alternate between solving in small time increments and recomposing solutions. This can be used to construct global solutions, deal with boundaries or inhomogeneities, and develop numerical approximation schemes; see, e.g., \cite{SLNH,MELENK1996289,Serre,10.1214/009117906000000115}. While there is no PDE or stochastic one in our case, we use this general idea to construct the periodic directed landscape. To show that this scheme indeed works, we will need to use probabilistic versions of finite speed of propagation estimates, in the form of geodesic fluctuation bounds or, in the case of ASEP, second-class particle fluctuation bounds.

We illustrate our gluing approach first in the context of solving the periodic heat equation (without randomness). Then, via a brief diversion in the realm of hydrodynamic limits, we will relate how we are led to implement a version of this approach for the directed landscape.

Imagine one is trying to solve the periodic heat equation $\partial_t Z = \tfrac{1}{2}\partial_x^2 Z$ for $Z(t,x)$ with $t\geq 0$ and $x\in \mathbb{T}$. Assume the initial data $Z_0$ is also periodic, so that $Z_0(0)=Z_0(1)$, and continuous (for simplicity). One way of solving this PDE is through method of images, which amounts to extending $Z_0$ periodically to be defined for $x\in \R$, and then solving the full-line heat equation for that initial data. Suppose instead that one were only allowed to solve the full-line heat equation with compactly supported initial data and for short times, as is relevant from a numerical perspective, and turns out to be very relevant for our periodic directed landscape construction.

\begin{figure}
\centering
\includegraphics[width=.9\textwidth]{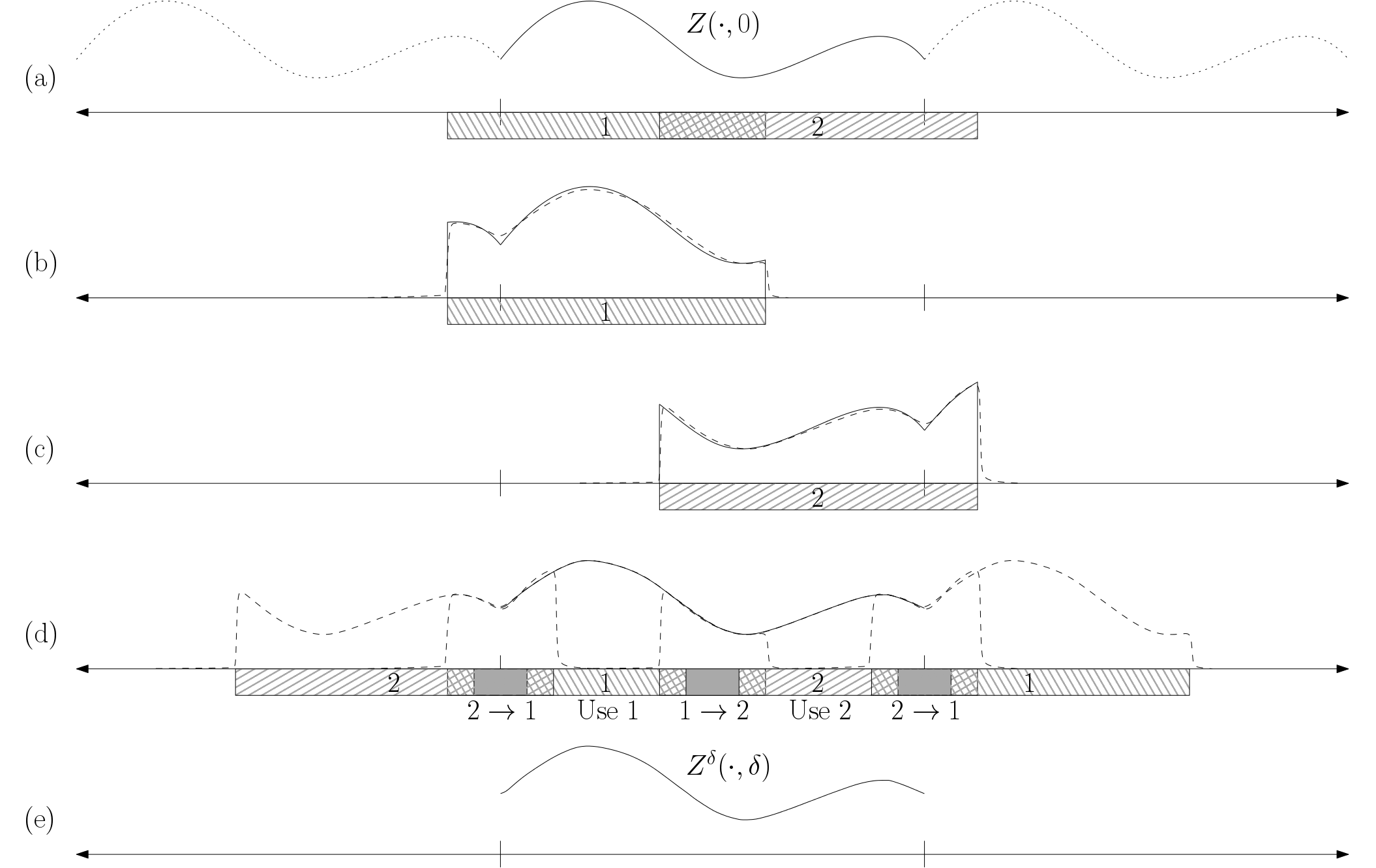}
\caption{\label{fig:introperiodicHE} An approximation scheme for the periodic heat equation.}
\end{figure}

The following scheme, illustrated in reference to Figure \ref{fig:introperiodicHE}, approximates the periodic heat equation only using the full-line heat equation with compactly supported initial data: (a) The initial data $Z(\cdot,0)=Z_0$ is periodically extended (dotted lines). (b) A new, compactly supported, initial data is created by multiplying by the indicator function of the interval labeled $1$ from $[-1/8,5/8]$. That is solved for time $\delta\ll 1$, resulting in the dashed line. (c) The same is done with respect to the interval labeled $2$ from $[3/8,9/8]$. (d) These time $\delta$ solutions are recomposed. Inside interval $1$, i.e., on $[1/16,7/16]$ we use the solution from (b) and inside the interval $2$, i.e., on $[9/16,15/16]$ we use the solution from (c). On the inside of the overlaps (in solid gray), i.e., $[0,1/16]$, $[7/16,9/16]$ and $[15/16,1]$ we linearly interpolate between the solutions from (b) and (c). (e) The result is denoted by $Z^{\delta}(\cdot,\delta)$. This becomes the initial data for step (a) and the process repeats $\delta^{-1}$ times. This defines $Z^{\delta}(x,t)$ for all times $t\in i\delta$ with $i\in \lbr \delta^{-1}\rbr$ after which we linearly interpolate in time to all $t\in [0,1]$. As $\delta\to 0$, $Z^{\delta}$ converges to the solution to the periodic heat equation with initial data $Z_0$.

The reason this scheme works is that, in a short time delta, influence in the heat equation only spreads in space diffusively, as $\delta^{1/2}$. For data of distance $x\gg \delta^{1/2}$ beyond that, the impact on the solution is of Gaussian order $e^{-c x^2/\delta}\ll 1$ (for some $c>0$). Since the recomposition step avoids using solutions in the vicinity of where they were cutoff (i.e., the solution in (b) is cutoff at $-1/8$ and $5/8$ but we only use it in the interval $[-1/16,9/16]$), even after $\delta^{-1}$ steps, the total error is bounded by $\delta^{-1}e^{-c/\delta}$ (for some other $c>0$), which goes to zero as $\delta$ does. The above Gaussian bound on influence can also be seen from the Feynman-Kac representation for the heat equation and the Gaussian distribution of Brownian motion arising in that representation.

A more involved scheme was developed by Glimm \cite{SLNH} (see also \cite[Chapter 5]{Serre}) to solve and numerically simulate hyperbolic conservation laws, such as the Burgers equation. The key input is the ability to solve the Riemann problem, i.e., when the initial data is given by two (potentially different) constants to the left and right of the origin. Splitting the domain into small spatial patches, one fits the initial data in each patch with a locally closely matched Riemann initial data. Run for short enough time, these locally mimic the true evolution and can be merged. This splitting, running and merging is repeated as above, and as the patch size and time step go to zero, the solution for general initial data is recovered. So far, we have only addressed gluing for PDEs. The next example includes randomness.

Inspired by Glimm's scheme, \cite{10.1214/009117906000000115} developed an analogous framework for proving hydrodynamic limits (laws of large numbers) for interacting particle systems such as the ASEP. The Burgers equation and related hyperbolic conservation laws arise as the deterministic limit of particle density fields for various interacting particle systems under Euler scaling (when time and space are scaled proportionally to infinity). As in Glimm's approach, one first computes this hydrodynamic limit under Riemann initial data, i.e, when the density of particles converge under Euler scaling to two constant profiles to the left and right of the origin. As with the deterministic PDE, general initial data for the particle system is decomposed locally into translates of such Riemann initial data and run for short times. The comparison between the true and patched evolution is then shown using coarse finite speed of propagation estimates, which are probabilistic analogs of the finite speed of propagation in hyperbolic conservation laws.

We now turn to the periodic directed landscape construction. This will summarize the construction of Section \ref{sec:PeriodicDL}. To approximate the periodic heat equation, we used full-line heat kernels. For the periodic directed landscape, we instead need to use full-space directed landscapes, which are random. In order to glue, we need to impose some form of periodic agreement. For a model like ASEP or LPP, this can be done at the level of the environment, i.e., making it periodic. For the directed landscape, there is no environment. Instead, we proceed as follows. 

For $\delta\ll 1$ and for $i\in \lbr \delta^{-1}\rbr$ we let $(\mathcal{L}^{i,1},\mathcal{L}^{i,2})$ be pairs of full-space directed landscapes (independent for different choices of $i$). We couple each pair $\mathcal{L}^{i,1}$ and $\mathcal{L}^{i,2}$ so as to ``agree'' with high probability on the overlap of intervals $1$ and $2$ from Figure \ref{fig:introperiodicHE}. Agreement means the following: For all times  $(i-1)\delta\leq s<t\leq i\delta$ and $x,y$ in $[3/8,5/8]$ (the overlap of intervals $1$ and $2$) we want $\mathcal{L}^{i,1}(x,s;y,t)=\mathcal{L}^{i,2}(x,s;y,t)$; and for $x,y \in [-1/8,1/8]$ we want $\mathcal{L}^{i,1}(x,s;y,t)=\mathcal{L}^{i,2}(x+1,s;y+1,t)$. This second condition means that the left part of interval $1$ and right part of interval $2$ agree in a periodic manner. Using results in  \cite{D:NonUnique,HP:Black}, we verify that this can be accomplished with probability at least $1-Ce^{-c\delta^{-2}}$ (see Lemma \ref{lem:CoupledDLs}); so, even after a union bound in $i\in \lbr \delta^{-1}\rbr$, this probability tends to $1$ as $\delta\to 0$. This ability to couple implicitly relies on rough control over the transversal fluctuations of geodesics in the directed landscape. Such control is more explicitly used repeatedly across our arguments. Indeed, it effectively replaces the Gaussian decay bounds in the periodic heat equation construction. Next, we use the above-coupled full-space directed landscapes to construct an approximation for the periodic directed landscape, and we prove these approximations form a tight family along a dyadic sequence for $\delta=2^{-n}$; the periodic directed landscape will be any limit point of this sequence (which we will confirm is unique). 

The approximation is built first for times $(i-1)\delta\leq s<t\leq i\delta$ for some  $i\in \lbr \delta^{-1}\rbr$. If $x\mod 1 \in [0,1/2)$ let $R=1$, and otherwise let $R=2$. This label $R$ refers to the choice of interval as in Figure \ref{fig:introperiodicHE}, so that interval $1$ is $[-1/8,5/8]$ and interval $2$ is $[3/8,9/8]$. For $x,y\in \R$ with $|x-y|\leq 1/16$ set $\bar{\mathcal{L}}^{\delta}(x,s;y,t) := \mathcal{L}^{i,R}(\tilde x,s;\tilde y, t)$ if the set of all geodesics for $\mathcal{L}^{i,R}$ from $(\tilde x,s)$ to $(\tilde y, t)$ stays in the rectangle with base given by interval $R$. Here $\tilde x$ is chosen to be in interval $R$ and equivalent to $x$ mod $1$, and $\tilde y = y + \tilde{x}-x$ so as to be shifted by the same amount as $x$. If some geodesics exit the rectangle, or if $|x-y|>1/16$ we set $\bar{\mathcal{L}}^{\delta}(x,s;y,t):=-\infty$. The definition is extended inductively to all $s$ and $t$ by the following metric composition property: assume $\bar{\mathcal{L}}^{\delta}$ is defined for all $t\in (s,i\delta]$ for some $i$; then, for $t\in (i\delta, (i+1)\delta]$, define 
$\bar{\mathcal{L}}^{\delta}(x,s;y,t) :=\sup_{z\in \R} \big(\bar{\mathcal{L}}^{\delta}(x,s;z,i\delta) + \bar{\mathcal{L}}^{\delta}(z,i\delta;y,t)\big)$. See Definition \ref{def:PatchedDLExt} for further details. 

The construction results in a periodic function in space, i.e., $\bar{\mathcal{L}}^{\delta}(x,s;y,t) = \bar{\mathcal{L}}^{\delta}(x+\ell,s;y+\ell,t)$ for all $\ell\in \Z$. We call this an unwrapped patched directed landscape since $x$ and $y$ are in $\R$; and define the wrapped version $\mathcal{L}^{\delta}(x,s;y,t) := \sup_{\ell\in \Z} \bar{\mathcal{L}}^{\delta}(x+\ell,s;y,t)$. Our main deduction, Proposition \ref{pro:CauchySequence} shows that for $\delta=2^{-n}$, the laws of $\mathcal{L}^{\delta}$ form a Cauchy sequence (in a suitable complete metric space on measures). The periodic directed landscape is its limit point or, rather, has the law of this limit point.
This Cauchy sequence claim underpins Theorem \ref{thm:PeriodicDirectedLandscapeMAIN}, and its proof is based on showing that $\mathcal{L}^{\delta} (x,s;y,t) = \mathcal{L}^{2\delta} (x,s;y,t)$ with high probability, whenever the magnitude of the slope of the line connecting $(x,s)$ and $(y,t)$ is not too small (at least $\delta^{1/8}$); see Lemma \ref{lem:PatchedDLsConsistently}. This once again relies on quantitative control on transversal fluctuations of geodesics, now both in the full-space directed landscape and the patched one $\mathcal{L}^{\delta}$, the latter of which are developed in Section \ref{sec:23}. A key technical difficulty is to show that the geodesics on $\mathcal{L}^{\delta}$ locally have the same properties as full space geodesics, in particular, they are not winding.

The proof of Theorem \ref{sec:furtherdirections} (convergence of periodic LPP) follows a similar route as the construction of the periodic directed landscape. However, the proof of Theorem \ref{thm:ConvergenceASEPMAIN} for ASEP is necessarily of a different nature, since the ASEP does not have a natural concept of geodesic or a variational formulation. Instead, we rely on non-standard couplings and moderate deviation bounds for second class particles (see Section \ref{sec:ASEPconvergencePeriodic}), both with origins in the recent work \cite{FS:triple}, to approximate the periodic colored ASEP by two coupled colored ASEPs on the integers in Sections \ref{sec:PatchedASEP}, \ref{sec:ConvergencePatchedASEP} and \ref{sec:ASEPfinalPart}.

We close this discussion by noting some earlier KPZ fluctuation results that involved ideas somewhat in the spirit of those discussed above. In a precursor to the construction of the directed landscape, \cite{HAMMOND_2019} introduced the 'patchwork quilt' framework using geodesic trees and used it to describe many local properties of KPZ scaling limits with general initial data.
Glimm's scheme and its lifting to hydrodynamics is also at the heart of the construction of the ASEP speed process \cite{ACG:ASEPspeed}. The idea of characterizing the directed landscape by certain locality properties, facilitated by control on transversal fluctuations, also partly underlies the work of \cite{DZ:Characterization} describing the directed landscape through KPZ fixed point marginals.

\subsection{Organization of the paper}

 In Section~\ref{sec:PeriodicDL} we construct the periodic directed landscape using suitably  patched directed landscapes. In Section~\ref{sec:ConvergenceLPP}, we prove convergence of periodic exponential last passage percolation to the periodic directed landscape. In Section~\ref{sec:KPZfixed}, we describe the periodic KPZ fixed point through the periodic directed landscape, and establish a local approximation for the former. In Section~\ref{sec:ASEPconvergence}, we verify that the periodic ASEP converges to the periodic KPZ fixed points, coupled according to the same periodic directed landscape. Moderate deviation estimates for periodic exponential last passage percolation and asymmetric simple exclusion processes are provided in the appendix.

\subsection{Notation}

We will use standard Landau notation and drop the integer brackets whenever appropriate. We will denote the discrete processes and sets in Latin letters while mostly using calligraphic notation for the continuous counterparts. The constants $c_0,C_0,c_1,C_1\dots>0$ do not depend on the system size $N$, and can change from line to line. Moreover, we will write $\lbr n \rbr := \{1,\dots,n\}$ for all $n \in \N$, and similarly $\lbr a,b \rbr= \{ \lfloor a\rfloor, \lfloor a\rfloor+1,\dots, \lfloor b\rfloor-1,\lfloor b\rfloor \}$ for $a<b$.

\subsection{Acknowledgments}

We wish to thank Yu Gu for early discussions on this subject. Amol Aggarwal was partially supported by a Packard Fellowship for Science and Engineering, a Clay Research Fellowship, by the NSF through grant DMS:1926686, and by the
IAS School of Mathematic. Ivan Corwin was partially supported by the Simons Foundation through an Investigator Award, through the W.M. Keck Foundation through a Science and Engineering Grant, and by the NSF through grants DMS:1937254, DMS:2246576,  and DMS:2348156. Much of this work was done while Dominik Schmid was a postdoctoral fellow at Columbia, funded jointly by Aggarwal's Packard Fellowship and Corwin's Simons Investigator Award. Some was also done during the 2025 Simons Symposium on Solvable Lattice Models and Interacting Particle Systems.

\section{Construction of the periodic directed landscape} \label{sec:PeriodicDL}
Our construction starts by gluing together in space and patching together in time certain restricted directed landscapes. We call the result a patched directed landscape. As the scale on which this patching occurs goes to zero (along a dyadic subsequence, we prove that the patched directed landscapes has a limit in distribution. We then show that this limit uniquely satisfies certain natural properties one would want for an ostensible periodic directed landscape; and in  Sections \ref{sec:ConvergenceLPP} and \ref{sec:ASEPconvergence} we show how periodic LPP and ASEP converge to this limit. All of this justifies our calling this limit the \emph{periodic directed landscape}.
In order to define our patched directed landscapes, we first state some properties of the full-space directed landscape.

\subsection{The full-space directed landscape}\label{sec:DLandKPZfixedFullSpace}

We recall the definition of the directed landscape as a random function on directed pairs of $\R^{4}$, and state its  basic properties. To this end, we let $\mathcal{S}^{(s)}\colon \R^2 \rightarrow \R$ denote the \textbf{Airy sheet} of scale $s$. This is a random continuous, stationary function, with its marginals explicitly given in terms of the Airy kernel; see Section~8 in~\cite{DOV:DirectedLandscape} for a precise definition via the Airy line ensemble. Define the space of directed space-time pairs
\begin{equation}
\R^{4}_{\uparrow} := \{ (x,s;y,t)\in \R^{4} \colon s<t \} , \qquad \textrm{and}\qquad [0,1]^{4}_{\uparrow} := \{ (x,s;y,t)\in [0,1]^{4} \colon s<t \}
\end{equation}
and let $\mathcal{C}(\R^{4}_{\uparrow},\R)$ (and likewise $\mathcal{C}([0,1]^{4}_{\uparrow},\R)$) denote the space of continuous functions from  $\R^{4}_{\uparrow}$ to $\R$ endowed with the uniform on compact topology. 

\begin{definition}[{\cite[Definition 10.1]{DOV:DirectedLandscape}}]\label{def:DirectedLandscape}
We denote by $\mathcal{L}=\mathcal{L}(x,s;y,t)$ the \textbf{(full-space) directed landscape}, a four parameter random function in $\mathcal{C}(\R^{4}_{\uparrow},\R)$. It satisfies:
\begin{enumerate}
\item[(i)] The marginals laws are given by the Airy sheet, i.e., for any $t,s>0$, the increment over the time interval $[t,t+s^3)$
\begin{equation}
(x,y) \mapsto \mathcal{L}(x,t;y,t+s^3)
\end{equation} is an Airy sheet $\mathcal{S}^{(s)}$ of scale $s$.
\item[(ii)] The metric composition law
\begin{equation}\label{eq:InvTriangleSup}
\mathcal{L}(x,s;y,t) = \sup_{z \in \R} \big(\mathcal{L}(x,s;z,r) + \mathcal{L}(z,r;y,t) \big)
\end{equation} holds almost surely for any $s< r < t$ and $x,y \in \R$. The supremum in \eqref{eq:InvTriangleSup} is attained. 
\item[(iii)] For all $(t_i)_{i\in \lbr \ell \rbr}$ and $(s_i)_{i \in \lbr \ell \rbr}$ with $\ell \in \N$ such that $(t_i,s_i)_{i \in \lbr \ell \rbr}$ are disjoint intervals, the random functions $
\mathcal{L}(\cdot,t_i;\cdot,s_i)
$ are independent. 
\end{enumerate}
\end{definition}

The existence and uniqueness of the full-space directed landscape is established as Theorem~10.9  in~\cite{DOV:DirectedLandscape}. Let us note that the directed landscape enjoys the scale invariance 
\begin{equation} \label{eq:ScalingDL}
\mathcal{L}(x,s;y,t) \overset{d}{=} q \mathcal{L}(q^{-2}x,q^{-3}s;q^{-2}y,q^{-3}t)
\end{equation} for all $q>0$; see Lemma 10.2(v) in \cite{DOV:DirectedLandscape}.  We also define the stationary directed landscape $\Lsta=(\Lsta(x,s;y,t))_{(x,s;y,t) \in \R^4_{\uparrow}}$ by
\begin{equation}\label{eq:StationaryDL}
    \Lsta(x,s;y,t) := \mathcal{L}(x,s;y,t) + \frac{(x-y)^2}{t-s}
\end{equation} for all $(x,s;y,t) \in \R^{4}_{\uparrow}$, and let $\lVert (x,s;y,t)\rVert_2:= \sqrt{|x-y|^2 + |s-t|^2}$ denote the $\ell_2$-distance between $(x,s)$ and $(y,t)$. 

\begin{proposition}[{\cite[Corollary~10.7]{DOV:DirectedLandscape}}]\label{pro:DLvalues} There exist constants $c,C>0$ such that the events
    \begin{equation*}
\mathsf{A}_{x,s;y,t}^z:=\left\{|\Lsta(x,s;y,t)| \leq z |t-s|^{1/3} \log^{4/3}\Big(\frac{2(2+\lVert (x,s;y,t)\rVert_2)^{\frac{3}{2}}}{t-s}\Big)\log^{2}(2+\lVert (x,s;y,t)\rVert_2) \right\} 
\end{equation*} satisfy for all  $z>0$
\begin{equation*}
    \P\left(\mathsf{A}_{x,s;y,t}^z \text{ holds for all } (x,s;y,t) \in \R^{4}_{\uparrow} \right) \geq 1 - C\exp\big(- c z^{3/2}\big) .
\end{equation*}
\end{proposition}

Moreover, we will use the following modulus of continuity result.

\begin{proposition}[{\cite[Proposition~1.6]{DOV:DirectedLandscape}}]\label{pro:ModulusOfContinuity}
Let $K \subseteq \R^{4}_{\uparrow}$ be compact. Let $\mathcal{L}=(\mathcal{L}(x,s;y,t))_{(x,s;y,t) \in \R^4_{\uparrow}}$ be a full-space directed landscape, and $\Lsta$ from \eqref{eq:StationaryDL} the respective stationary directed landscape. For all $\xi,\tau,\psi>0$, consider the events 
\begin{equation*}
\begin{split}
   \mathsf{A}_{\xi,\tau}(\psi):=  \Big\{ &\Big| \Lsta(x,s;y,t)- \Lsta(x^{\prime},s^{\prime};y^{\prime},t^{\prime}) \Big| \leq \psi \text{ for all } (x,s;y,t),(x^{\prime},s^{\prime};y^{\prime},t^{\prime}) \in K \\ &\text{ such that } |x-x^{\prime}|,|y-y^{\prime}| \leq \xi \text{ and }  |s-s^{\prime}|,|t-t^{\prime}| \leq \tau \Big\} .
\end{split}
\end{equation*} 
Then there exist constants $c,C>0$, depending only on $K$, such that for all $\xi,\tau,\phi>0$ 
\begin{equation*}
 \P\left( \mathsf{A}_{\xi,\tau}\Big( \phi\big(\tau^{1/3}\log^{2/3}(\tau^{-1}+1) + \xi^{1/2}\log^{1/2}(\xi^{-1}+1)\big)\Big) \right) \geq 1
 - C \exp(-c\phi^{3/2}) . 
\end{equation*}
\end{proposition}

Paths which achieve the supremum in \eqref{eq:InvTriangleSup} will play a key role in investigating the directed landscape. We say that a continuous path $\pi=\pi_{x,s;y,t}$ with $\pi \colon [s,t] \rightarrow \R$ is a \textbf{directed geodesic} from $u=(x,s)$ to $v=(y,t)$ in the directed landscape $\mathcal{L}$ if it satisfies
\begin{equation}\label{def:Geodesic}
\int \textup{d}\mathcal{L}\circ \pi :=  \inf_{j \in \N} \inf_{s= r_0 < r_1 < \cdots < r_{j}=t} \sum_{i=1}^{j} \mathcal{L}(\pi(r_{i-1}),r_{i-1}; \pi(r_{i}),r_i)  = \mathcal{L}(s,\pi(s);t,\pi(t)) . 
\end{equation}
We call  $\int \textup{d}\mathcal{L}\circ \pi$  the \textbf{length} of $\pi$, and write $\mathfrak{G}_{u,v}$ for the set of all geodesics from $u$ to $v$.  For $\ell \in [s,t]$ and any path $\pi$ from $u=(x,s)$ to $v=(y,t)$, not necessarily a geodesic, define $\TF_{\ell}(\pi)$, the \textbf{local transversal fluctuations} at time $\ell$,  and $\TF(\pi)$, the \textbf{transversal fluctuations}, by
\begin{equation*}
\TF_{\ell}(\pi) := \Big| \pi(\ell)-x - (y-x)\frac{\ell-s}{t-s} \Big| , \qquad \textrm{and}\qquad \TF(\pi) := \sup_{\ell \in [s,t]} \TF_{\ell}(\pi).
\end{equation*}
These measure how far $\pi$ deviates from the straight line connecting $u=(x,s)$ and $v=(y,t)$. Note that while our notation $\TF_{\ell}(\pi)$ and $\TF(\pi)$ does not explicitly include the starting and  endpoints $u$ and $v$, this is implicit in the definition.
For all pairs of sites $(u;v)\in \mathbb{R}^{4}_{\uparrow}$, define 
\begin{equation}\label{eq:PointTransversal}
\TF_{\ell}(u;v) := \sup_{\pi \in \mathfrak{G}_{u,v}}\TF_{\ell}(\pi),\qquad \textrm{and}\qquad    \TF(u;v) := \sup_{\pi \in \mathfrak{G}_{u,v}}\TF(\pi).
\end{equation} 
A geodesic $\pi^{+}_{u;v}\in\mathfrak{G}_{u,v}$ is the  \textbf{rightmost geodesic} if for all $\ell\in [s,t]$, $\pi^{+}_{u;v}(\ell) =\sup_{\pi\in \mathfrak{G}_{u,v}} \pi(\ell)$. The \textbf{leftmost geodesic} $\pi^{-}_{x,s;y,t}$ is defined similarly. The following result addresses such geodesics.
\begin{proposition}[{\cite[Theorem 12.1 and Lemma 13.2]{DOV:DirectedLandscape}}]\label{pro:Geos}
    Almost surely, for all $(u;v)\in \R^{4}_{\uparrow}$ there exists a leftmost and rightmost geodesic. For fixed $(u;v)\in \R^{4}_{\uparrow}$, there almost surely exists a unique geodesic from $u$ to $v$.
\end{proposition}
Next, we record a moderate deviations result on the local fluctuations of geodesics from \cite{GZ:fractal}; see Proposition~12.3 in~\cite{DOV:DirectedLandscape} for a similar result when $(x,s;y,t)=(0,0;0,1)$.  Recall that we denote for all $(x,s;y,t) \in \mathbb{R}^{4}_{\uparrow}$ by $\lVert (x,s;y,t)\rVert_2$ the  $\ell_2$-distance between $(x,s)$ and $(y,t)$.
\begin{proposition}[{\cite[Lemma 3.11]{GZ:fractal}}]\label{pro:ModerateDL}
There exist constants $c,C>0$, such that for all  $\phi>0$,
\begin{align*}
&\P\Bigg(\TF_{s+\theta^{-1}}(\pi) \geq \phi \theta^{-2/3}\log^{3}\big(1+\theta\lVert (x,s;y,t)\rVert_2\big),\,\,\,\forall(x,s;y,t)\in \R^{4}_{\uparrow},\pi\in \mathfrak{G}_{(x,s),(y,t)}, \theta\geq 2(t-s)^{-1}\Bigg)\\
&\leq C\exp\big({-c \phi^{9/4}}\big). 
\end{align*} By symmetry, the result holds with $s+\theta^{-1}$ replace by  $t-\theta^{-1}$. Hence,  
\begin{equation}
\begin{split}
\label{eq:PeriodicDLModerate}
&\P\Bigg(\TF(\pi) \geq \phi  |t-s|^{2/3}\log^{3}\bigg(1+\frac{\lVert (x,s;y,t)\rVert_2}{|t-s|}\bigg),\,\,\,\forall (x,s;y,t) \in \R^{4}_{\uparrow},\pi\in \mathfrak{G}_{(x,s),(y,t)}\Bigg)  \\
&\leq C\exp({-c \phi^{9/4}}).   
\end{split}
\end{equation} 
\end{proposition}

In words, Proposition \ref{pro:ModerateDL} guarantees that full-space geodesic have locally small changes in the transversal fluctuations, uniformly in the starting and ending points.  We require the following definition on the local fluctuations of geodesics. 
\begin{definition}\label{def:DeltaRegular}
     Let $\varepsilon>0$. For all $s<t$, consider the set of points
\begin{equation}\label{def:GoodIncrements}
    (z_i)_{i\in \lbr k\rbr} := \{ j\varepsilon \in [s,t] \text{ for some } j\in \N \} 
\end{equation} where $k\in \N_0$ is the number of points, and where we set $z_0=s$ and $z_{k+1}=t$, and $z_{i-1} \leq z_{i}$ for all $i\in \lbr k+1 \rbr$. 
    A geodesic $\pi\in \mathfrak{G}_{(x,s),(y,t)}$ with $(x,s;y,t)\in \R^{4}_{\uparrow}$ is $\boldsymbol{\varepsilon}$\textbf{-good} if for all $i\in \lbr k+1 \rbr$
    \begin{equation}
      \sup_{\ell \in [z_{i-1},z_i]}  \big| \pi(\ell) - \pi(z_{i-1}) \big| \leq \frac{1}{32}  . 
    \end{equation}
\end{definition}
In particular, note that every $\varepsilon$-good geodesic fluctuates by at most $\frac{1}{8}$ in a time interval of length $\varepsilon$.
We get the following consequence of Proposition \ref{pro:ModerateDL} when the  \textbf{inverse slope}
\begin{equation}\label{def:InverseSlope}
    \sl^{-1}(x,s;y,t) := \frac{y-x}{t-s} \in \R
\end{equation} between the points $(x,s)$ and $(y,t)$ with $(x,s;y,t) \in \R^{4}_{\uparrow}$ is not too large.
\begin{lemma}\label{lem:TypicallyGood}
    For each $\varepsilon>0$, define $\mathsf{G}_{\varepsilon}$ to be the event that for all $(x,s;y,t) \in [0,1]^{4}_{\uparrow}$ with $\big|\sl^{-1}(x,s;y,t)\big|\leq \varepsilon^{-\frac{1}{2}}$ and for all $\pi\in \mathfrak{G}_{(x,s),(y,t)}$, the geodesic $\pi$  is $\varepsilon\text{-good}$. There exists $\varepsilon_0>0$ and $c,C>0$ such that for all $\varepsilon\in \big(0,\varepsilon_0\big)$, 
    \begin{equation}
        \P(\mathsf{G}_{\varepsilon}) \geq  1-C\exp(-c\varepsilon^{-1/4}) . 
    \end{equation}
\end{lemma}
\begin{proof}
With the choice of $\phi=\varepsilon^{-1/6}/\log^4(\varepsilon^{-1})$, let $\mathsf{\tilde G}_\varepsilon$ denote the event on the left-hand side of \eqref{eq:PeriodicDLModerate}. The result of Proposition \ref{pro:ModerateDL} implies that there exist $\varepsilon_0>0$ and constants $c,C>0$ such that for all $\varepsilon\in (0,\varepsilon_0)$,
$\mathbb{P}(\mathsf{\tilde G}_\varepsilon)\geq 1-C\exp(-c\varepsilon^{-1/4})$. We will show below that, for small enough $\varepsilon$, $\mathsf{\tilde G}_\varepsilon\subseteq \mathsf{G}_\varepsilon$ from which the lemma we are proving is immediate. 
In what follows we will assume that the event $\mathsf{\tilde G}_\varepsilon$  holds and show that, for small enough $\varepsilon$ so too does $\mathsf{G}_\varepsilon$ hold. The event $\mathsf{G}_\varepsilon$ involves varying over all starting and ending points with suitably bounded inverse slope, and establishing control on variations over time increments of size $\varepsilon$. We will do this by considering separately the time increments at least $\varepsilon$ from $s$ or $t$ (i.e., in the bulk), and those around the starting and ending times. If $t-s$ is of order $\varepsilon$ it is possible that there are no internal increments, in which case the analysis is just around the starting and ending times.

Let us start by considering any $(x,s;y,t) \in \mathbb{R}^{4}_{\uparrow}$. Recall the notation  $(z_i)_{i \in \lbr k\rbr}$ from \eqref{def:GoodIncrements} and assume for the moment that $k\geq 2$. For all $i \in\lbr 2, k\rbr$ we will say that the interval $[z_{i-1},z_{i}]$ is in the bulk. By symmetry of our below argument, without loss of generality we can assume that $z_i \leq (s+t)/2$, i.e., we are in the first half of the time interval $[s,t]$. On the event $\mathsf{\tilde G}_\varepsilon$, it follows that for all $\varepsilon>0$ small enough
\begin{equation*}
    \begin{split}
           |z_{i-1}-\pi(z_{i-1})| &\leq (t-z_{i-1}) \varepsilon^{-\frac{1}{2}}  + (t-z_{i-1})^{\frac{2}{3}} \log^3\Big( \frac{3}{t-z_{i-1}} \Big) \phi \\
    &\leq  (t-z_{i-1}) \varepsilon^{-\frac{1}{2}}  + (t-z_{i-1})^{\frac{2}{3}} \varepsilon^{-\frac{1}{6}} \,\leq\, 2(t-z_{i-1}) \varepsilon^{-\frac{1}{2}}  , 
    \end{split}
\end{equation*}
where we use the upper bound $\big|\sl^{-1}(x,s;y,t)\big|\leq \varepsilon^{-\frac{1}{2}}$ and the bound on the transversal fluctuations by $\mathsf{\tilde G}_\varepsilon$ for the first inequality, the definition of $\phi$ for the second inequality, and the fact that $t-z_{i-1}\geq \varepsilon$ for the second and third inequality. This implies that there is $\varepsilon_0>0$ (which does not depend on $x,s,y,t$) such that for all $\pi\in \mathfrak{G}_{(x,s),(y,t)}$ and all $\varepsilon\in (0,\varepsilon_0)$
    \begin{equation}\label{eq:SlopeControl}
        \big|\sl^{-1}\big(\pi(z_{i-1}),z_{i-1};y,t\big)\big| \leq  2\varepsilon^{-\frac{1}{2}} ,  
    \end{equation}
In a similar manner, on the event $\mathsf{\tilde G}_\varepsilon$ we may control the transversal fluctuation of any geodesic $\pi'\in \mathfrak{G}_{\pi(z_{i-1},z_{i-1}),(y,t)}$ so that 
\begin{equation}\label{eq:FluctuationsGood}
        \TF_{z_{i}}(\pi^{\prime}) \leq (z_i-z_{i-1})^{\frac{2}{3}} \log^3\Big( \frac{3}{z_i-z_{i-1}} \Big) \phi \leq   \varepsilon^{1/2} 
\end{equation} with $\varepsilon\in (0,\varepsilon_0)$ small enough. Since $\pi$ restricted to the time interval $[z_{i-1},t]\in\mathfrak{G}_{(\pi(z_{i-1}),z_{i-1}),(y,t)}$ this result also applies to $\pi$. Combining the inverse slope bound \eqref{eq:SlopeControl} for $\pi$ with the transversal fluctuation bound \eqref{eq:FluctuationsGood} for $\pi$ we conclude that   
\begin{equation}\label{eq:z2k}
       \Big| \pi(z_{i}) - \pi(z_{i-1}) \Big| \leq 3\varepsilon^{1/2} 
\end{equation}
     for all $i \in \lbr 2,k\rbr$ with $(z_i)_{i \in \lbr k\rbr}$ from \eqref{def:GoodIncrements}. 
Finally,  on the event $\mathsf{\tilde G}_\varepsilon$ it follows that for all $\pi''\in \mathfrak{G}_{(\pi(z_{i-1}),z_{i-1}),(\pi(z_{i}),z_{i})}$ with \eqref{eq:z2k}, and any $\ell\in [z_{i-1},z_{i}]$, 
\begin{equation}\label{eq:z3k}
        \TF_{\ell}(\pi'') \leq \varepsilon^{1/2}.
\end{equation}
Again, this also applies to $\pi$ restricted to the time interval $[z_{i-1},z_i]$.
Putting \eqref{eq:z2k} and \eqref{eq:z3k} together, we see that on the event $\mathsf{\tilde G}_\varepsilon$, for all $(x,s;y,t) \in [0,1]^{4}_{\uparrow}$ with $\big|\sl^{-1}(x,s;y,t)\big|\leq \varepsilon^{-\frac{1}{2}}$ and for all $\pi\in \mathfrak{G}_{(x,s),(y,t)}$, the geodesic $\pi$ satisfies 
$$
\sup_{\ell \in [z_{i-1},z_i]}  \big| \pi(\ell) - \pi(z_{i-1}) \big| \leq 4\varepsilon^{1/2}\leq  \frac{1}{32}
$$
for all $i \in \lbr 2,k\rbr$, provided $\varepsilon$ is small enough so that $4\varepsilon^{1/2}\leq  \frac{1}{32}$. To show that $\pi$ is $\varepsilon\text{-good}$ it remains to control the left-hand side above for $i=1$ and $i=k+1$. This can be done in the exact same manner as for the bulk time increments. If there is no bulk, i.e. $k<2$, then the argument needs only to be applied to the starting and ending increments. Having shown that $\pi$ is $\varepsilon\text{-good}$ on the event $\mathsf{\tilde G}_\varepsilon$ we conclude that $\mathsf{\tilde G}_\varepsilon\subseteq \mathsf{G}_\varepsilon$ from which the lemma follows.
\end{proof}

\subsection{Definition of the patched directed landscape}\label{sec:PreliminaryMetricConstruction}
First some  notation. For $\mathcal{R}\subset \R^2$ let
\begin{equation}\label{def:ArrowNotation}
    (\mathcal{R})^{2}_{\uparrow} := \mathcal{R}\times \mathcal{R} \cap \R^4_{\uparrow} ,
\end{equation}
 i.e. the set of pairs of points in $\mathcal{R}$ with ordered time coordinates (thinking of the second coordinate in $\mathcal{R}$ as time. Overloading this notation slightly, we also define the space
\begin{equation}\label{def:Cylinder}
\SetT := \big\{ (x,s;y,t)\in (\mathbb{T}\times[0,1])^2 \colon s<t \big\}
\end{equation} 
where, $\mathbb{T}$ denotes the $1$-dimensional torus of length $1$, which we identify with the interval $[0,1]$, gluing the points $0$ and $1$ together.
$\SetT$ denotes the set of pairs of points $(x,s;y,t)$ on a cylinder of length and width $1$ such that the \textbf{vertical coordinates} $s$ and $t$ are ordered and the \textbf{horizontal coordinates} $x$ and $y$ are in $\mathbb{T}$. For $i,j\in \mathbb{Z}$ and $n\in \N$ define the rectangles
\begin{equation}\label{def:RectanglesDL}
\mathcal{O}^{i,j}_n:= \Big[\frac{j-1}{2}-\frac{1}{8},\frac{j-1}{2}+\frac{1}{8}\Big]\times \Big[\frac{i-1}{n},\frac{i}{n}\Big] .
\end{equation}

\begin{figure}
\centering
\begin{tikzpicture}[scale=.9]


  \draw[red, line width =1.2 pt] (4,0) -- (4,1.5) -- (14,1.5) -- (14,0) -- (4,0);
  
    \draw[red, line width =1.2 pt] (-4,4.5) -- (-4,3) -- (6,3) -- (6,4.5) -- (-4,4.5);
    
    

    \draw[darkblue, fill=darkblue!20, line width =1.2 pt] (4.05,0.05) -- (6,0.05) -- (6,1.45) -- (4.05,1.45) -- (4.05,0.05);

    \draw[darkblue, fill=darkblue!20, line width =1.2 pt] (12,0.05) -- (13.95,0.05) -- (13.95,1.45) -- (12,1.45) -- (12,0.05);

    \draw[darkblue, fill=darkblue!20, line width =1.2 pt] (-3.95,3.05) -- (-2,3.05) -- (-2,4.45) -- (-3.95,4.45) -- (-3.95,3.05);

    \draw[darkblue, fill=darkblue!20, line width =1.2 pt] (4,3.05) -- (5.95,3.05) -- (5.95,4.45) -- (4,4.45) -- (4,3.05);
    
 \draw[line width =1.5 pt, ->,dotted] (4,3.75) -- (4,0.75);    
 \draw[line width =1.5 pt, ->,dotted] (6,3.75) -- (6,0.75);       

 \draw[line width =1.5 pt, ->,dotted] (-4,3.75) -- (12,0.75);    
 \draw[line width =1.5 pt, ->,dotted] (-2,3.75) -- (14,0.75);  

\node[scale=1] (site) at (-3,3.75){$\mathcal{O}^{i,1}_n$};
\node[scale=1] (site) at (5,0.75){$\mathcal{O}^{i,2}_n$};
\node[scale=1] (site) at (5,3.75){$\mathcal{O}^{i,2}_n$};
\node[scale=1] (site) at (13,0.75){$\mathcal{O}^{i,3}_n$};

\node[scale=1] (site) at (1,3.75){$\mathcal{L}^{1}$};

\node[scale=1] (site) at (9,0.75){$\mathcal{L}^{2}$};

	\end{tikzpicture}	
\caption{\label{fig:DLpatch}  The gluing construction for the patched directed landscape uses the coupling $\mathbf{P}_{\per}$ of full-space directed landscapes $\mathcal{L}^{1}$ and $\mathcal{L}^{2}$ from Lemma~\ref{lem:CoupledDLs}. Under that coupling, $\mathcal{L}^{1}$ restricted to $\mathcal{O}^{i,1}_n$ agrees with $\mathcal{L}^{2}$ restricted to $\mathcal{O}^{i,3}_n$ and $\mathcal{L}^{1}$ and $\mathcal{L}^{2}$ restricted to $\mathcal{O}^{i,2}_n$ agrees with each other.}
 \end{figure}
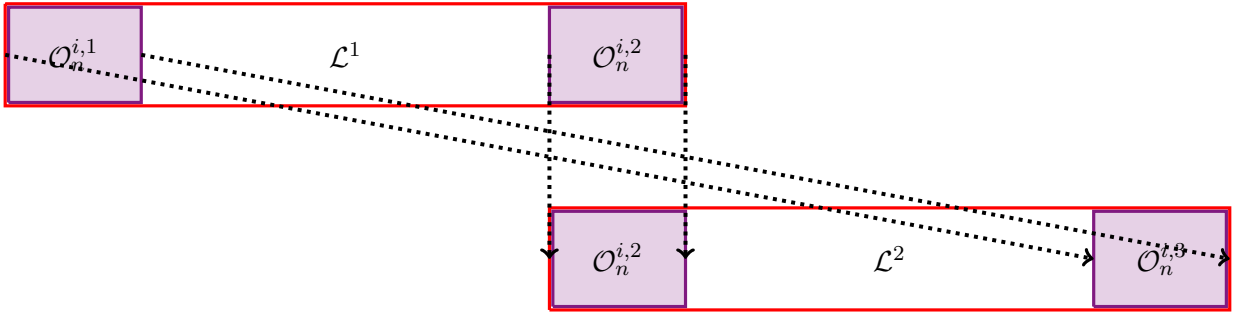

Now, recall that for probability measures $\nu,\bar{\nu}$ on a common probability space $(\Omega,\mathcal{F})$, we set
\begin{equation}\label{eq:SupCharactierzationTV}
   \TV{\nu -\bar{\nu}} := \inf\left\{ \mathbf{P}( X \neq \bar{X}) \colon  \mathbf{P} \text{ coupling of }X \sim \nu \text{ and } \bar{X} \sim \bar{\nu} \right\} = \sup_{A \in \mathcal{F}} \big(\nu(A) - \bar{\nu}(A) \big),  
\end{equation} where the second representation is known as the Kantorovich-Rubinstein duality.
The following lemma explains how to couple directed landscapes to coincide on rectangles that are short and well separated; see Figure~\ref{fig:DLpatch}. 
 The proof uses Lemma~3.2 of~\cite{HP:Black} (see also Proposition~2.6 of \cite{D:NonUnique} for a similar result), which provides an approximate independence of the directed landscape on distant rectangles.
\begin{lemma}\label{lem:CoupledDLs}
Consider two full-space directed landscapes $\mathcal{L}^{1}$, $\mathcal{L}^{2}$.
For any $\mathcal{O} \subseteq \R^2$ and $k\in \lbr 2 \rbr$, write 
$\mathcal{L}^{k}_{\mathcal{O}} := (\mathcal{L}^{k}(x,s;y,t))_{(x,s;y,t) \in (\mathcal{O})^2_{\uparrow}}$ for the restriction of $\mathcal{L}^{k}$ to $(\mathcal{O})^{2}_{\uparrow}$.
There exists $c,C>0$ such that for all $i\in \Z$ and $n\in N$ there exist a coupling $\mathbf{P}^i_n$ of $\mathcal{L}^{1}$ and $\mathcal{L}^{2}$ under which
\begin{equation}\label{eq:PeriodicCoupling}
\mathbf{P}^i_n\left( \mathcal{L}^{1}_{\mathcal{O}^{i,1}_n} =   \mathcal{L}^{2}_{\mathcal{O}^{i,3}_n} \text{ and } \mathcal{L}^{1}_{\mathcal{O}^{i,2}_n} =   \mathcal{L}^{2}_{\mathcal{O}^{i,2}_n} \right) \geq 1 - C\exp(-cn^{2}). 
\end{equation}
\end{lemma}
\begin{proof} It suffices to prove the result for $i=1$ since the coupling for other $i$ can be achieved via shifting the directed landscapes (due to their translational invariance). For the rest of the proof we will adopt the shorthand $\mathcal{O}^{j}_n=\mathcal{O}^{1,j}_n$, dropping the $i$. In preparation for showing \eqref{eq:PeriodicCoupling}, we demonstrate a related auxiliary result. Consider two independent full-space directed landscapes $\mathcal{L}^{3}$ and $\mathcal{L}^{4}$. We claim that there exists $c,C>0$ such that for all $n\in \N$, there exists couplings $\mathbf{P}'_n,\mathbf{P}''_n$ of $\mathcal{L}^{1},\mathcal{L}^{3}$, and $\mathcal{L}^{4}$ under which 
\begin{equation}\label{eq:CouplingDLs}
\begin{split}
\mathbf{P}'_n\Big( \mathcal{L}^{1}_{\mathcal{O}^{1}_n}= \mathcal{L}^{3}_{\mathcal{O}^{1}_n}  \text{ and } \mathcal{L}^{1}_{\mathcal{O}^{2}_n}= \mathcal{L}^{4}_{\mathcal{O}^{2}_n}\Big) &\geq 1 - C\exp(-cn^{2}) ,   \\ 
\mathbf{P}''_n\Big( \mathcal{L}^{1}_{\mathcal{O}^{2}_n}= \mathcal{L}^{3}_{\mathcal{O}^{2}_n}  \text{ and } \mathcal{L}^{1}_{\mathcal{O}^{3}_n}= \mathcal{L}^{4}_{\mathcal{O}^{3}_n}\Big) &\geq 1 - C\exp(-cn^{2})  .
\end{split}
\end{equation}
 We will only show the first inequality in \eqref{eq:CouplingDLs} as the second one follows analogously. By the scale invariance~\eqref{eq:ScalingDL} of the directed landscape, this is equivalent to showing that there exist $c,C>0$ such that for all $n\in \N$ there exists a coupling $\tilde{\mathbf{P}}'_n$ under which 
\begin{equation}\label{eq:ScaledBound}
    \tilde{\mathbf{P}}_n'\left( \mathcal{L}^{1}_{\tilde{\mathcal{O}}^1_n}= \mathcal{L}^{3}_{\tilde{\mathcal{O}}^1_n}  \text{ and } \mathcal{L}^{1}_{\tilde{\mathcal{O}}^2_n}= \mathcal{L}^{4}_{\tilde{\mathcal{O}}^2_n} \right) \geq 1 - C \exp(-cn^{2}) , 
\end{equation}
where we set
\begin{equation*}
    \tilde{\mathcal{O}}^1_n :=  \Big[\frac{3}{8}n^{\frac{2}{3}},\frac{5}{8}n^{\frac{2}{3}}\Big] \times  [0,1], \qquad     \tilde{\mathcal{O}}^2_n := \Big[\frac{7}{8}n^{\frac{2}{3}},\frac{9}{8}n^{\frac{2}{3}}\Big]\times [0,1]. 
\end{equation*}
Using the shift invariance of the directed landscape, the desired bound \eqref{eq:ScaledBound} follows from  Lemma~3.2 in \cite{HP:Black} (on the event $\big\{ X \leq \frac{1}{8}n^{2/3}\big\}$ in their notation).

We return now to demonstrate \eqref{eq:PeriodicCoupling}. To show the existence of the coupling $\mathbf{P}_n$ (again, we drop the $i$ superscript) we will use \eqref{eq:SupCharactierzationTV} along with the triangle inequality for total variation distance. Total variation distance requires that our measures are all defined on the same measure space. We will consider the space $\mathcal{C}\big((\mathcal{O}^{1}_n\cup \mathcal{O}^{2}_n)^2_{\uparrow}, \R\big)$ in which $\mathcal{L}^{1}$ restricted to $\mathcal{O}^{1}_n\cup \mathcal{O}^{2}_n$ takes values. By identifying spatial points modulo 1 as the same, we can treat $\mathcal{L}^{2}$ restricted to $\mathcal{O}^{3}_n\cup \mathcal{O}^{2}_n$ as taking values on the same space. Write $\nu^{1}$ for the law of the restriction of $\mathcal{L}^{1}$, and $\nu^{2}$ for the law of the restriction of $\mathcal{L}^{2}$, identified to live on the same space as $\nu^{1}$ as above. We seek to show that 
there exist $c,C>0$ such that for all $n\in \N$
\begin{equation}\label{eq:TVn12}
\TV{ \nu^{1} - \nu^{2}} \leq C \exp(-c n^{2}).    
\end{equation}
By \eqref{eq:SupCharactierzationTV}, this will imply \eqref{eq:PeriodicCoupling}.
To show \eqref{eq:TVn12}, let $\nu'$ denote the law of $\mathcal{L}^{3}$ restricted to $\mathcal{O}^{1}_n$ and $\mathcal{L}^{4}$ restricted to $\mathcal{O}^{2}_n$ where  $\mathcal{L}^{3}$ and  $\mathcal{L}^{4}$ are independent as assumed above. Likewise, let $\nu''$ denote the law of $\mathcal{L}^{4}$ restricted to $\mathcal{O}^{3}_n$ (which we treat, modulo 1, as equivalent to $\mathcal{O}^{1}_n$) and $\mathcal{L}^{3}$ restricted to $\mathcal{O}^{2}_n$. Since $\nu'$ and $\nu''$ are both product measures of restricted directed landscapes, they define the same measures and hence $\TV{ \nu' - \nu''} = 0$. From \eqref{eq:CouplingDLs} it follows that there exists $c,C>0$ such that for all $n\in \N$,
\begin{equation*}
\TV{ \nu^{1} - \nu'} \leq C \exp(-c n^{2}), \qquad 
\TV{ \nu^{2} - \nu''} \leq C \exp(-c n^{2}).
\end{equation*}
Thus, 
\begin{equation*}
\TV{ \nu^{1} - \nu^{2}} \leq  
\TV{ \nu^{1} - \nu'}
+
\TV{ \nu' - \nu''}
+
\TV{ \nu'' - \nu^{2}} \leq 2C \exp(-c n^{2}).
\end{equation*}
Replacing $2C$ by $C$ yields \eqref{eq:TVn12} as desired.
\end{proof}

Let us now briefly sketch the overall strategy in the construction of the periodic directed landscape.
Using the coupling from Lemma \ref{lem:CoupledDLs}, we will now define, for each $n\in \N$, a directed metric on $\R\times [0,1]$ built out of patching together in a periodic manner, coupled directed landscapes. We will call this the scale $n$ unwrapped patched directed landscape. Summing over integer translates of the starting point, we will likewise define the wrapped version of this object. Ultimately, we will take $n\to \infty$ and show a limit which will define our periodic directed landscape. Let us emphasize at this point that we make the selection of a sequence of couplings (namely the ones in Lemma~\ref{lem:CoupledDLs}) that we will use for the definition of the unwrapped and wrapped patched landscapes, and hence in our construction of the periodic directed landscape. After we have constructed the periodic directed landscape, our uniqueness result Theorem~\ref{thm:PeriodicDirectedLandscape} will ensures that any other choice of couplings would lead to the same limit. \\

For our construction, we require a bit of notation that we now provide. Recall \eqref{def:ArrowNotation} and set
$$\SetR := \{(x,s;y,t) \in \R^{4}_{\uparrow} \colon s,t \in [0,1] \} .$$ 
Define for all $i,j\in\Z$ and $\delta>0$ the rectangles 
\begin{equation}\label{eq:FillingRectangle}
\mathcal{R}^{i,j}_n:=\Big[\frac{1}{2}(j-1)-\frac{1}{8},\frac{j}{2}+\frac{1}{8}\Big]\times \Big[\frac{i-1}{n},\frac{i}{n}\Big]. 
\end{equation}
Note that $\mathcal{R}^{i,j}_n\cap \mathcal{R}^{i,j+1}_n = \mathcal{O}^{i,j}_n$ from \eqref{def:RectanglesDL}. Define, for all $i\in \Z$ the horizontal strips
\begin{equation}\label{eq:PartitionedRectangles}
  \mathcal{R}^{i}_n:= \bigcup_{j\in \Z}\mathcal{R}^{i,j}_n=\R\times  \Big[\frac{i-1}{n},\frac{i}{n}\Big]
\end{equation} 
and write  $\partial \mathcal{R}^{i}_n = \R \times \big\{ \frac{i}{n}\big\}$ for the upper boundary of $\mathcal{R}^{i}_n$. 

Our patched landscape on a time horizon $[0,1]$ is composed of patching together $n$ coupled (as in Lemma~\ref{lem:CoupledDLs}) pairs of directed landscapes. For each  $i\in \lbr n \rbr$, let the pair $(\mathcal{L}^{i,1},\mathcal{L}^{i,2})$ be distributed as two directed landscapes coupled via the law $\mathbf{P}^i_n$ as in \eqref{eq:PeriodicCoupling} from Lemma~\ref{lem:CoupledDLs}. For different $i$, assume these pairs are independent. We denote the resulting measure on $\big\{(\mathcal{L}^{i,1},\mathcal{L}^{i,2}):i\in  \lbr n\rbr\big\}$ by $\mathbf{P}_n$.
Define the event $\mathsf{A}_n$ that the coupling over all $i\in \lbr n\rbr$ is \textbf{successful}, i.e., 
\begin{equation}\label{eq:GoodCoupling}
\mathsf{A}_n := \bigcap_{i=1}^{n} \mathsf{A}^i_n \quad \text{ with } \quad  \mathsf{A}^i_n:=\left\{  \mathcal{L}^{i,1}_{\mathcal{O}^{1,i}_n} =   \mathcal{L}^{i,2}_{\mathcal{O}^{3,i}_n}\quad \text{ and } \quad \mathcal{L}^{i,1}_{\mathcal{O}^{2,i}_n} =   \mathcal{L}^{i,2}_{\mathcal{O}^{2,i}_n}  \right\} . 
\end{equation}
By Lemma~\ref{lem:CoupledDLs} and a union bound, there exists constants $c,C>0$ such that for each $n\in \N$
\begin{equation}\label{eq:GoodCouplingBound}
  \mathbf{P}_n( \mathsf{A}_n ) \geq 1 - C \exp\big(-cn^{2}\big).
\end{equation}

We will build the unwrapped and wrapped patched directed landscape from these coupled directed landscapes as follows. 
Define the following subset of $\R^4_{\uparrow}$ given by space-time points with distance $\frac{1}{8}$ in space and in the same $\frac{1}{n}$-sized time interval:
\begin{equation}\label{def:SlopeVerti}
\mathcal{V}_n := \Big\{ (x,s;y,t) \in \R^4_{\uparrow} \, \colon  y-x \in \Big[-\frac{1}{16}, \frac{1}{16}\Big] , \,  s,t\in \Big[\frac{i-1}{n},\frac{i}{n}\Big], \, i\in \lbr n \rbr \Big\}. 
\end{equation}
For $(x,s;y,t)\in \mathcal{V}_n$ define
\begin{equation}\label{eq:FillingRectanglePoints}
    \mathcal{R}_{x,s;y,t}:=\mathcal{R}^{i,j}_n
\end{equation}
where $j\in \Z$ is chosen such that $x\in \big[\tfrac{1}{2} (j-1),\tfrac{1}{2} j)$ and $i\in \lbr n\rbr$ is such that $s,t\in \big[\frac{i-1}{n}, \frac{i}{n}\big]$.
For $(x,s;y,t)\in \mathcal{V}_n$ and $\mathcal{R} = \mathcal{R}_{x,s;y,t}$ 
define the directed landscape associated to the rectangle $\mathcal{R}$ as
\begin{equation}\label{eq:LReqn}
    \mathcal{L}^{\mathcal{R}}(x,s;y,t) := \begin{cases} 
    \mathcal{L}^{i,1}\big(x-\frac{j-1}{2},s;y-\frac{j-1}{2},t\big) &\text{ if  } \mathcal{R}_{x,s;y,t}=\mathcal{R}^{i,j}_n \text{ for } j \text{ odd,} \\
    \mathcal{L}^{i,2}\big(x-\frac{j-2}{2},s;y-\frac{j-2}{2},t\big) &\text{ if  } \mathcal{R}_{x,s;y,t}=\mathcal{R}^{i,j}_n \text{ for } j \text{ even.}    
    \end{cases}
\end{equation}
Observe that this is $1$-periodic in the spatial variable so that $\mathcal{L}^{\mathcal{R}}(x,s;y,t)  = \mathcal{L}^{\mathcal{R}}(x+1,s;y+1,t)$. 
Finally, to any given rectangle $\mathcal{R} = \mathcal{R}^{i,j}_n$ for some $i\in \lbr n\rbr$,  $j\in \Z$ and any $(x,s;y,t)\in \R^4_{\uparrow}$, we write $\pi^{\mathcal{R},+}_{(x,s);(y,t)}$ and $\pi^{\mathcal{R},-}_{(x,s);(y,t)}$ for the rightmost and leftmost geodesics in the corresponding directed landscape $\mathcal{L}^{i,1}$ (if $j$ is odd) or $\mathcal{L}^{i,2}$ (if $j$ is even).

\begin{figure}
\centering
\begin{tikzpicture}[scale=.95]

   \draw[line width=1.2pt, dotted] (0,0) -- (15,0);
   \draw[line width=1.2pt, dotted] (0,3) -- (15,3);

  \draw[red, line width =1.2 pt] (1,0) -- (1,3) -- (9,3) -- (9,0) -- (1,0);
  
    \draw[red, line width =1.2 pt] (6,0) -- (6,3) -- (14,3) -- (14,0) -- (5,0);
    
    
    
   
\draw[line width =1 pt] (6.55+0.35-0.4,0.05) to[curve through={(6.55+0.5-0.4,0.55)..(6.55-0.4,1.05) .. (6.55+0.4-0.3,1.65) .. (6.55+0.55-0.4,2.45) }] (6.55+0.35-0.4,2.95);


\node[scale=1.2] (site) at (3.5,1.5){$\mathcal{R}_n^{i,1}$};  
\node[scale=1.2] (site) at (11.5,1.5){$\mathcal{R}_n^{i,2}$};  
\node[scale=1.2] (site) at (7.8,1.5){$\mathcal{O}_n^{i,1}$}; 

 \node[scale=0.9] (site) at (6.5,-0.25){$(x,s)$};  
  \node[scale=0.9] (site) at (6.5,3.25){$(y,t)$}; 


%
%
%
%
%
\end{tikzpicture}	
\caption{\label{fig:RectangleDL}Assignment of the different rectangles  in the construction of patched directed landscape. Note that horizontal coordinates $x,y$ are at least $\frac{1}{16}$ away from the right boundary of $\mathcal{R}_n^{i,1}$, and that we assign $\mathcal{R}_{x,s;y,t}=\mathcal{R}_n^{i,1}$ in this example. In particular, when the coupling from Lemma \ref{lem:CoupledDLs} is successful, we get that the directed landscapes $\mathcal{L}^{1}$ and $\mathcal{L}^{2}$ agree on $\mathcal{O}^{i,1}_n=\mathcal{R}_n^{i,1} \cap \mathcal{R}_n^{i,2}$.} 
\end{figure}
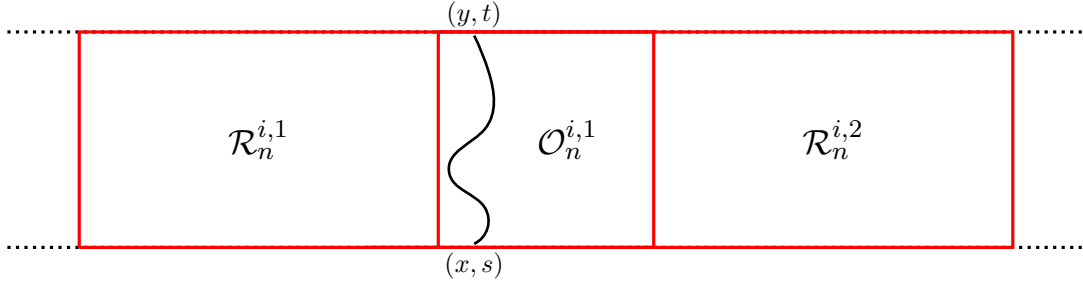

With this notation in place we may now give our construction, a visualization of which is given in Figure~\ref{fig:RectangleDL}. 
We restrict to the time interval $[0,1]$ though subsequently extend it to $\R$.
\begin{definition}\label{def:PatchedDLExt} 
Assume, as above, the directed landscapes $\big\{(\mathcal{L}^{i,1},\mathcal{L}^{i,2}):i\in  \lbr n\rbr\big\}$ are coupled via the law $\mathbf{P}_n$ described above \eqref{eq:GoodCoupling}, and recall the notation $ \mathcal{L}^{\mathcal{R}}$ from \eqref{eq:LReqn}.

We define the scale $n\in \N$ \textbf{unwrapped patched directed landscape} as the random function $\Lpat:\SetR\to\R \cup \{-\infty\}$ given as follows. For all $(x,s;y,t) \in \SetR \cap \mathcal{V}_{n}$ we set 
\begin{equation}\label{def:RestrictedRectanglesPatch}
\Lpat(x,s;y,t) := \begin{cases}  \mathcal{L}^{\mathcal{R}}(x,s;y,t) &\text{ if } \pi^{\mathcal{R},+}_{x,s;y,t}\text{ and }\pi^{\mathcal{R},-}_{x,s;y,t} \subseteq \mathcal{R}, \\
-\infty &\text{ else}. 
\end{cases}
\end{equation} 
For such $(x,s;y,t)\in \SetR \cap \mathcal{V}_{n}$ we say that a continuous path $\bar\pi:[s,t]\to \R$ is an \textbf{unwrapped inherited geodesic} for the scale $n$ unwrapped patched directed landscape if $\Lpat(x,s;y,t)\neq -\infty$ and $\bar\pi$ is a geodesic for $\mathcal{L}^{\mathcal{R}}$.

For all other $(x,s;y,t) \in \SetR \setminus \mathcal{V}_{n}$, either 
$s,t\in \big[\frac{i-1}{n},\frac{i}{n}\big]$ for some $i\in \lbr n\rbr$, or not. 
In the case when $s,t\in \big[\frac{i-1}{n},\frac{i}{n}\big]$, we set $\Lpat(x,s;y,t):=-\infty$ and say that  no unwrapped inherited geodesic exists between $(x,s)$ and $(y,t)$.
In the other case, where $s\in \big[\frac{i-1}{n},\frac{i}{n}\big]$ and $t\notin \big[\frac{i-1}{n},\frac{i}{n}\big]$ for some $i\in \lbr n\rbr$, we define $\Lpat(x,s;y,t)$ as follows. Assuming, inductively, that $\Lpa(x,s;y,t)$ is defined for $t\in \big(s,\frac{k}{n}\big]$ for some $k\geq i$ (the base case $k=i$ is above) then for $t\in \big(\frac{k}{n},\frac{k+1}{n}\big]$, 
\begin{equation}\label{eq:RecursionDL}
\Lpat(x,s;y,t) := \sup_{(\tilde{x},\tilde{s}) \in \partial \mathcal{R}^{k}_n} \Lpat(x,s;\tilde{x},\tilde{s}) +  \Lpat(\tilde{x},\tilde{s};y,t). 
\end{equation}
The unwrapped inherited geodesic is similarly inductively defined. Assume that we have specified the definition of the unwrapped inherited geodesics for $t\in \big(s,\frac{k}{n}\big]$ for some $k\geq i$ (the base case $k=i$ definition is above). Then, for $t\in \big(\frac{k}{n},\frac{k+1}{n}\big]$, we say that a path $\bar\pi:[s,t]\to\R$ is an unwrapped inherited geodesic if its restriction to $[s,\frac{k}{n}]$ and $[\frac{k}{n},t]$ are both unwrapped inherited geodesic and if the pair $\big(\bar\pi(\frac{k}{n}),\frac{k}{n}\big)$, realizes the supremum in \eqref{eq:RecursionDL}. In all cases above, if no path satisfies the condition of being an unwrapped inherited geodesic between $(x,s)$ and $(y,t)$ then we say that no such geodesic exists.

By construction, the unwrapped patched directed landscape $\Lpat$ is $1$-periodic in the horizontal coordinates, i.e., $\Lpat(x,s;y,t)=\Lpat(x+\ell,s;y+\ell,t)$ for all $\ell\in \Z$. We can now define the  scale $n\in \N$ \textbf{wrapped patched directed landscape} by setting for all $(x,s;y,t) \in \SetT$
\begin{equation}\label{eq:ExtendPatched}
\Lpa(x,s;y,t):= \sup_{\ell\in \Z}\Lpat(x+\ell,s;y,t). 
\end{equation}
A path $\pi:[s,t]\to \mathbb{T}$ on the torus is called a \textbf{wrapped inherited geodesic} if there exists an unwrapped inherited geodesic $\bar\pi:[s,t]\to \R$ such that $\bar\pi(s) = x+\ell$ where $\ell$ achieves the supremum in \eqref{eq:ExtendPatched}, and such that $\pi(\tilde s)= \bar\pi(\tilde s) \mod 1$ for all $\tilde{s}\in[s,t]$. Since we assume all inherited geodesics are continuous paths, $\bar\pi$ is uniquely determined from $\pi$ and we call it the unwrapping of $\pi$. 
We call an unwrapped inherited geodesic $\pi$ between $(x,s)$ and $(y,t)$ leftmost if for all other unwrapped inherited geodesics $\tilde\pi$, $\pi(\tilde s)\leq \pi(\tilde s)$ for all $\tilde s\in[s,t]$ (and rightmost if $\leq $ is replaced with $\geq$). We call a wrapped inherited geodesic leftmost (or rightmost) if its unwrapping is the leftmost (or rightmost) unwrapped inherited geodesic.
\end{definition}

Definition \ref{def:PatchedDLExt} introduces the notation of an inherited geodesic. It should be noted that this is not the same as if one tried to directly define a geodesic via \eqref{def:Geodesic}. An inherited geodesic is defined by the property that its restriction to each $\frac{1}{n}$-sized time interval stays in the rectangle $\mathcal{R}$ corresponding to the entry point into the interval, is a geodesic in the underlying full-space directed landscape for the rectangle $\mathcal{R}$, and that intercepts with the lines $\partial\mathcal{R}^k_n$ realize the supremum in \eqref{eq:ExtendPatched}. If we instead defined a geodesic via \eqref{def:Geodesic}, then it is possible (though very unlikely as one can see from Lemma~\ref{lem:PenaltyDistance}) that even in a given $\frac{1}{n}$-sized time interval, the geodesic could move horizontal outside of the rectangle $\mathcal{R}$ dictated by the starting point. 


\subsection{Properties of the patched directed landscape}\label{sec:23}

We now build up several properties of the patched directed landscapes with the goal of showing, in Lemma~\ref{lem:PatchedDLsConsistently}, that we can couple together scale $n$ and $2n$ patched directed landscapes so as to agree, with high probability, for all points with sufficient inverse slope $\sl^{-1}(x,s;y,t)$ from \eqref{def:InverseSlope}. To this end, we start by collecting basic properties of the unwrapped patched directed landscape in Lemma~\ref{lem:Geodesics}. We provide bounds on the lengths in the unwrapped patched directed landscape $\Lpat$ in Lemma~\ref{lem:CouplePatchedToFull} and Lemma~\ref{lem:FirstFullSpace} with sufficient inverse slope $\sl^{-1}(x,s;y,t)$ by a comparison of the inherited geodesics in $\Lpat$ with those of full-space directed landscapes. Upper bounds on the length in the unwrapped patched directed landscape for general pairs of points are established in Lemma~\ref{lem:PenaltyDistance} using a path decomposition of inherited geodesics. Lemma \ref{lem:Noncrossing} compares the unwrapped patched directed landscape to a wrapped patched directed landscape, while Lemma~\ref{lem:LocalFluctuationsPatched} addresses the local fluctuations of inherited geodesics in $\Lpat$.

For $\varepsilon>0$, define the set of all pairs of points whose inverse slope is at most $\varepsilon^{-1/2}$ in magnitude:
\begin{equation}\label{def:SlopeConti}
\mathcal{D}_{\varepsilon} = \Big\{ (x,s;y,t) \in \mathbb{R}^{4}_{\uparrow} \, \colon  \sl^{-1}(x,s;y,t) \in \big[-\varepsilon^{-\frac{1}{2}}, \varepsilon^{-\frac{1}{2}}\big]  \Big\}  . 
\end{equation} 
Moreover, similar to \eqref{def:Geodesic}, we define the \textbf{length} of a path $\pi \colon [s,t] \rightarrow \R$ (which we recall to be a continuous function from $[s,t]$ to $\R$) in the unwrapped patched directed landscape $\Lpat$ to be 
\begin{equation}\label{def:GeodesicPatched}
\int \textup{d}\Lpat\circ \pi :=  \inf_{j \in \N} \inf_{s= r_0 < r_1 < \cdots < r_{j}=t} \sum_{i=1}^{j} \Lpat(\pi(r_{i-1}),r_{i-1}; \pi(r_{i}),r_i) . 
\end{equation}
We start by addressing basic properties of the unwrapped patched directed landscape, and the existence of rightmost and leftmost inherited geodesics (recall Definition \ref{def:PatchedDLExt}) in the patched directed landscape for starting/ending point pairs in $\mathcal{V}_{n}$ and $\mathcal{D}_{n^{-1}}$.
\begin{lemma}\label{lem:Geodesics} 
  There exist $c,C>0$ such that for each $n\in \N$, there exists an event $\mathsf{B}_n$ with
 \begin{equation}\label{eq:ProbabilityBound} 
     \mathbf{P}_n(\mathsf{B}_n) \geq 1-C \exp\big(-cn^{9/8} \big) 
 \end{equation} on which:
\begin{enumerate}
    \item[(i)]\label{eq:Patched1}   The unwrapped patched directed landscape $\Lpat$ is an element of the space of upper semi-continuous functions $ \USC(\SetR,\R \cup \{-\infty\})$.
    \item[(ii)]\label{eq:Patched2}   The unwrapped patched directed landscape and its unwrapped inherited geodesics satisfy the following decomposition property: If $\bar\pi$ is an unwrapped inherited geodesic from $(x,s)$ to $(y,t)$ for some $(x,s,y,t)\in\SetR$ and $\tilde s\in (s,t)$ then $\bar\pi$ restricted to $[s,\tilde s]$ and $\bar\pi$ restricted to $[\tilde s, t]$ are both inherited geodesics, and 
\begin{equation}\label{eq:PatchedCompo}
\Lpat(x,s;y,t) = \Lpat(x,s;\bar\pi(\tilde s),\tilde s) + \Lpat(\bar\pi(\tilde s),\tilde s;y,t).
\end{equation} 
\item[(iii)]\label{eq:Patched3} For every inherited geodesic $\pi \colon [s,t] \rightarrow \R$ in the unwrapped patched directed landscape, 
\begin{equation}
  \int \textup{d}\Lpat\circ \pi  = \Lpat(s,\pi(s);t,\pi(t)) . 
\end{equation}
    \item[(iv)]\label{eq:Patched4}  For all  $(x,s;y,t)\in \SetR \cap \mathcal{V}_{n}$, we have that 
    \begin{equation}\label{eq:NonTrivialLower}
     \Lpat(x,s;y,t) > - \infty \quad \text{ for all } (x,s;y,t)\in \SetR \cap \mathcal{V}_{n} , 
 \end{equation}
and there are unique rightmost and leftmost inherited geodesics for $ \Lpat$, which we denote by $\pinpt_{x,s;y,t}$ and $\pinmt_{x,s;y,t}$, respectively.
\item[(v)]\label{eq:Patched5}  The result in \hyperref[eq:Patched2]{(ii)} holds true with $\mathcal{V}_n$ replaced by $\mathcal{D}_{n^{-1}}$ from \eqref{eq:ProbabilityBound} and/or when $\Lpat$ is replaced by $\Lpa$ (in which case we write $\pinp_{x,s;y,t}$ and $\pinm_{x,s;y,t}$ for the rightmost and leftmost wrapped inherited geodesics).
\end{enumerate}
\end{lemma}
\begin{proof}
 Recall the event $\mathsf{A}_n$ from \eqref{eq:GoodCoupling}, and that we assigned to each pair $(x,s;y,t)\in \mathcal{V}_n$ a unique rectangle $\mathcal{R}=\mathcal{R}_{x,s;y,t}$ in \eqref{eq:FillingRectanglePoints}. Proposition~\ref{pro:Geos} guarantees for all $(x,s;y,t)\in \mathcal{V}_n$ the existence of rightmost and leftmost geodesics $\pi^{\mathcal{R},-}_{x,s;y,t}$ and $\pi^{\mathcal{R},+}_{x,s;y,t}$ in $\mathcal{L}^{\mathcal{R}}$. Write $\TF^{\mathcal{R}}(\pi)$ for the transversal fluctuations of a path $\pi$ from $(x,s)$ to $(y,t)$ in the directed landscape $\mathcal{L}^{\mathcal{R}_{x,s;y,t}}$. 
Let $\mathsf{C}_{n}$ be the event that 
\begin{equation*}
 \mathsf{C}_{n}:= \left\{   \TF^{\mathcal{R}}( \pi )\leq \frac{1}{16} \, \text{ for all } \pi \in \mathfrak{G}_{(x,s),(y,t)} \ \text{and} \  (x,s;y,t) \in \SetR \cap \mathcal{V}_{n}  \right\} . 
\end{equation*}
In particular, $\mathsf{C}_{n}$ ensures that $\pi_{x,s;y,t}^{\mathcal{R}_{x,s;y,t},+}, \pi_{x,s;y,t}^{\mathcal{R}_{x,s;y,t},-} \subseteq \mathcal{R}_{x,s;y,t}$ for all $(x,s;y,t) \in \SetR \cap \mathcal{V}_{n}$ (i.e., that the first line in \eqref{def:RestrictedRectanglesPatch} holds). By using Proposition \ref{pro:ModerateDL} with $t-s$ of order $n^{-1}$ and $\phi$ of order $n^{1/2}$, it follows that there exist $c,C>0$ such that for all $n\in \N$
\begin{equation*}
    \P(\mathsf{C}_{n} ) \geq 1 - C\exp\big(-c n^{9/8} \big).
\end{equation*}
Set $\mathsf{B}_n:= \mathsf{A}_n \cap \mathsf{C}_n$ and recall the sets $\mathcal{R}^{i}_n$ from \eqref{eq:PartitionedRectangles}. Note that on the event $\mathsf{B}_n$, 
\begin{equation}\label{eq:ComparingDLsR}
  \Lpat(x,s;y,t)  =  \mathcal{L}^{\mathcal{R}}(x,s;y,t) 
\end{equation}
for all $(x,s;y,t)\in (\mathcal{R}^{i}_n)^2_{\uparrow} \cap\mathcal{V}_n$ with $i\in \lbr n \rbr$, and $\Lpat(x,s;y,t) =-\infty$ for all $(x,s;y,t) \in (\mathcal{R}^{i}_n )^2_{\uparrow} \setminus \mathcal{V}_n$. In particular, since $\mathcal{L}^{\mathcal{R}}$ is continuous, $\Lpat$ is continuous on the set  $(\mathcal{R}^{i}_n)^2_{\uparrow} \cap\mathcal{V}_n$ and upper semi-continuous on $(\mathcal{R}^{i}_n)^2_{\uparrow}$ for all $i\in \lbr n \rbr$.

For general $(x,s;y,t) \in \SetR$, item~\hyperref[eq:Patched1]{(i)} follows from the inductive variational definition of $\Lpat$ as follows. For $k\in \N_0$, define the inductive hypothesis that for all $i\in \lbr n \rbr$, $\Lpat$ is continuous and on
$(\mathcal{R}^{i}_n \cup \cdots \cup \mathcal{R}^{i+k}_n )^2_{\uparrow} \cap \mathcal{V}_n$, and $\Lpat$ is equal to $-\infty$ on $(\mathcal{R}^{i}_n \cup \cdots \cup \mathcal{R}^{i+k}_n )^2_{\uparrow} \setminus \mathcal{V}_n$. This implies upper semi-continuity. The $k=0$ base case was shown in the previous paragraph.
Given the case $k$, to show the case $k+1$ it suffices to show that $\Lpat$ is continuous on $(\mathcal{R}^{i}_n\times  \mathcal{R}^{i+k+1}_n)\cap \mathcal{V}_n$ for all $i\in \lbr n \rbr$, and $-\infty$ on $(\mathcal{R}^{i}_n\times  \mathcal{R}^{i+k+1}_n)\setminus \mathcal{V}_n$. However, this follows readily from the inductive hypothesis and the way $\Lpat$ is defined via the variational problem \eqref{eq:RecursionDL}. Indeed, for $(x,s;y,t)\in (\mathcal{R}^{i}_n\times  \mathcal{R}^{i+k+1}_n)\cap \mathcal{V}_n$, 
$$
 \Lpat(x,s;y,t)= \sup_{(\tilde{x},\tilde{s}) \in \partial \mathcal{R}^{i+k}_n} \Lpat(x,s;\tilde{x},\tilde{s}) +  \Lpat(\tilde{x},\tilde{s};y,t).
$$
Due to the inductive hypothesis, the supremum is achieved in a compact window and due to the continuity of each term on the right-hand side above, the left-hand side is likewise continuous.

For item~\hyperref[eq:Patched2]{(ii)}, consider first the special case $(x,s),(y,t) \in \mathcal{R}^{i}_n$ for some $i\in \lbr n\rbr$. 
Then the event $\mathsf{C}_n$ ensures that the two restricted paths from $(x,s)$ to $(\pi(\tilde{s}),\tilde{s})$ and from $(\pi(\tilde{s}),\tilde{s})$ to $(y,t)$, respectively, must be fully contained in the respective rectangles $\mathcal{R}^{(1)}=\mathcal{R}_{x,s;\pi(\tilde{s}),\tilde{s}}$ and $\mathcal{R}^{(2)}=\mathcal{R}_{\pi(\tilde{s}),\tilde{s};y,t}$. Moreover, on the event $\mathsf{A}_n$, we see that 
\begin{align*}
    \Lpat(x,s;\bar\pi(\tilde s),\tilde s) = \mathcal{L}^{\mathcal{R}^{(1)}}(x,s;\bar\pi(\tilde s),\tilde s)=\mathcal{L}^{\mathcal{R}^{(2)}}(x,s;\bar\pi(\tilde s),\tilde s), \\ 
    \Lpat(\bar\pi(\tilde s),\tilde s;y,t)=\mathcal{L}^{\mathcal{R}^{(1)}}(\bar\pi(\tilde s),\tilde s;y,t)=\mathcal{L}^{\mathcal{R}^{(2)}}(\bar\pi(\tilde s),\tilde s;y,t) . 
\end{align*}
Since either $\mathcal{R}_{x,s;y,t}=\mathcal{R}^{(1)}$ or $\mathcal{R}_{x,s;y,t}=\mathcal{R}^{(2)}$, we use \eqref{eq:ComparingDLsR} together with the metric composition of the full-space directed landscape to conclude \eqref{eq:PatchedCompo}. For general $(x,s),(y,t)$, let $\lbr a,b\rbr$ denote the set of all $i\in \lbr n \rbr$ such that $s < in^{-1} < t$. Then by the definition of the unwrapped patched directed landscape, 
\begin{equation}\label{eq:SplitPath}
\begin{split}
 \Lpat(x,s;y,t) &= \Lpat(x,s;\pi(an^{-1}),an^{-1})+  \sum_{i=a}^{b-1} \Lpat\Big(\pi\Big(\frac{i}{n}\Big),\frac{i}{n};\pi\Big(\frac{i+1}{n}\Big),\frac{i+1}{n}\Big) \\
 & \quad + \Lpat(\pi(bn^{-1}),bn^{-1};y,t). 
 \end{split}
\end{equation}
Now assume that $\tilde{s}\in [in^{-1},(i+1)n^{-1})$ for some $i\in \lbr a,b-1\rbr$ (the cases $\tilde{s} \in (s,an^{-1})$ and $\tilde{s} \in [bn^{-1},t)$ are similar). Then we get from the special case of item \hyperref[eq:Patched2]{(ii)} considered above that
\begin{equation}\label{eq:SplitPath2}
    \Lpat\Big(\pi\Big(\frac{i}{n}\Big),\frac{i}{n};\pi\Big(\frac{i+1}{n}\Big),\frac{i+1}{n}\Big) = \Lpat\Big(\pi\Big(\frac{i}{n}\Big),\frac{i}{n};\pi(\tilde{s}),\tilde{s}\Big) + \Lpat\Big(\pi(\tilde{s}),\tilde{s};\pi\Big(\frac{i+1}{n}\Big),\frac{i+1}{n}\Big). 
\end{equation} 
All values of $\Lpat$ on the right-hand sides of  \eqref{eq:SplitPath} and \eqref{eq:SplitPath2} are attained by  inherited geodesics corresponding to the starting and ending points. Thus, recomposing these inherited geodesics, we conclude that item \hyperref[eq:Patched2]{(ii)} holds in full generality. Item \hyperref[eq:Patched3]{(iii)} is an immediate consequence of the decomposition property \hyperref[eq:Patched2]{(ii)}.  

For item \hyperref[eq:Patched4]{(iv)}, recall that \eqref{eq:ComparingDLsR} holds on the event $\mathsf{B}_n$, and that $\pi^{\mathcal{R},-}_{x,s;y,t}$ and $\pi^{\mathcal{R},+}_{x,s;y,t}$ are the leftmost and rightmost geodesics in $\mathcal{L}^{\mathcal{R}}$. Then from the definition of $\Lpat$, we get that $\pinmt_{x,s;y,t}=\pi^{\mathcal{R},-}_{x,s;y,t}$ and $\pinpt_{x,s;y,t}=\pi^{\mathcal{R},+}_{x,s;y,t}$ are the leftmost and rightmost inherited geodesics in $\Lpat$ for all $(x,s;y,t) \in \SetR \cap \mathcal{V}_{n}$. In particular (as noted earlier), we see that \eqref{eq:NonTrivialLower} holds.

For item \hyperref[eq:Patched5]{(v)}, let now $(x,s;y,t)\in \SetR \cap \mathcal{D}_{n^{-1}}$. We have to show that \eqref{eq:NonTrivialLower} holds. To this end, we first note there exists some $n_0$ such that we can choose the constants $c,C>0$ in \eqref{eq:ProbabilityBound} so that $\mathsf{B}_n$ is the empty event when $n< n_0$. Now if $s,t \in [(i-1) n^{-1},i n^{-1}]$ for some $i\in \lbr n \rbr$, then $(x,s;y,t) \in \mathcal{D}_{n^{-1}}\subseteq\mathcal{V}_{n}$, provided that $n$ is large enough, and we conclude using \eqref{eq:NonTrivialLower} with respect to $\mathcal{V}_n$. Now let $s \in [(i-1) n^{-1},i n^{-1}]$ for some $i\in \lbr n \rbr$, and assume, inductively, that \eqref{eq:NonTrivialLower} is shown for all $(x,s;y,t)\in \SetR \cap \mathcal{D}_{n^{-1}}$ with $t\in \big(s,\frac{k}{n}\big]$ for some $k\geq i$. Let $t\in \big(\frac{k}{n},\frac{k+1}{n}\big]$. Then we see that $\Lpat(x,s;y,t)>-\infty$ as
we can find some $\tilde{z}\in \R$ such that $\big(x,s;\tilde{z},\frac{k}{n}\big),\big(\tilde{z},\frac{k}{n};y,t\big) \in \mathcal{D}_{n^{-1}}$ and use that by the induction hypothesis, $\Lpat(x,s;\tilde{z},\frac{k}{n}),\Lpat(\tilde{z},\frac{k}{n};y,t)>-\infty$. Now
let $\mathcal{I}_{x,s;y,t}$ be the set of all $z\in \R$ such that
 \begin{equation}\label{eq:VariationalGeos}
    \Lpat(x,s;y,t)= \Lpat\Big(x,s;z,\frac{k}{n}\Big) + \Lpat\Big(z,\frac{k}{n};y,t\Big) = \sup_{\tilde{x} \in \R} \Lpat\Big(x,s;\tilde{x},\frac{k}{n}\Big) + \Lpat\Big(\tilde{x},\frac{k}{n};y,t\Big) . 
 \end{equation} 
Note that $\mathcal{I}_{x,s;y,t} \subseteq [y-\frac{1}{16},y+\frac{1}{16}]$ by the definition of $\Lpat$. Moreover, we claim that $\mathcal{I}_{x,s;y,t}$ is a compact non-empty set. To see this, note that since $\Lpat$ is upper semi-continuous by item \hyperref[eq:Patched1]{(i)}, we get that the supremum in \eqref{eq:VariationalGeos} is attained, and that the level sets of $\tilde{x} \mapsto \Lpat\big(x,s;\tilde{x},\frac{k}{n}\big) + \Lpat\big(\tilde{x},\frac{k}{n};y,t\big)$ are closed. Now let $z^{+}$ and $z^{-}$ denote the largest and smallest element in $\mathcal{I}_{x,s;y,t}$, respectively. Then we obtain $\pinpt_{x,s;y,t}$ by concatenation of $\pinpt_{x,s;z^{+},k/n}$ and $\pinpt_{z^{+},k/n;y,t}$, and similarly for $\pinmt_{x,s;y,t}$. Note that by \eqref{eq:VariationalGeos}, the path $\pinpt_{x,s;y,t}$ is indeed an inherited geodesic for $\Lpat$, and hence by construction the rightmost such path. A similar argument holds for $\pinmt_{x,s;y,t}$. This gives \eqref{eq:NonTrivialLower} for all $(x,s;y,t)\in \SetR \cap \mathcal{D}_{n^{-1}}$.  
 The corresponding result for the wrapped patched directed landscape is immediate from its definition, noting that we have to take the supremum in \eqref{eq:ExtendPatched} only over finitely many values, i.e., we have that 
\begin{equation}\label{eq:ExtendPatchedProof}
\Lpa(x,s;y,t)= \sup_{\ell\in \Z}\Lpat(x+\ell,s;y,t) = \sup_{\ell\in [-n,n] \cap \Z}\Lpat(x+\ell,s;y,t)
\end{equation}
for all $(x,s;y,t) \in \mathcal{D}_{n^{-1}}$ by the definition of $\Lpat$. 
 \end{proof}

Next, our goal is to provide lower bounds on the lengths in the unwrapped patched directed landscape. To this end, we require the following lemma.

\begin{lemma}\label{lem:CouplePatchedToFull}
There exist $c,C>0$ such that for all $\delta>0$ and all $n\in \N$ with $n\geq \delta^{-2}$ the following holds. Fix any rectangle $\widetilde{\mathcal{R}} \subseteq \mathbb{R} \times [0,1]$ of side lengths $\frac{1}{4}$ by $\delta$. Then there exists a coupling $\mathbf{P}$ between a full-space directed landscape $\mathcal{L}^{\prime}$ and an unwrapped patched directed landscape $\Lpat$ such that
\begin{equation}\label{eq:patchedFullCase1}
    \mathbf{P}\left( \Lpat(x,s;y,t) = \mathcal{L}^{\prime}(x,s;y,t) \,\,\, \forall (x,s;y,t) \in (\widetilde{\mathcal{R}})^2_{\uparrow} \cap \mathcal{D}_{n^{-1}} \right) \geq 1 - C\exp(-c\delta^{-1/2}) . 
\end{equation}
The result in \eqref{eq:patchedFullCase1} holds true with $\mathcal{D}_{n^{-1}}$ replaced by $\mathcal{V}_{\delta^{-1}}$. 
\end{lemma}
\begin{proof} 
We start by constructing the full-space directed landscape $\mathcal{L}^{\prime}$ used in both statements. Recall the event $\mathsf{B}_n$ from \eqref{eq:ProbabilityBound}. In particular, we have that $\mathsf{B}_n \subseteq \mathsf{A}_n$ for $\mathsf{A}_n$ from \eqref{eq:GoodCoupling}, which is the event that the coupling between the directed landscapes $\mathcal{L}^{i,1}$ and $\mathcal{L}^{i,2}$ is successful for all $i\in \lbr n\rbr$. Recall the sets $\mathcal{R}^{i}_n$ from \eqref{eq:PartitionedRectangles}, and define for all $i \in \lbr n\rbr$ the sets
\begin{equation*}
    \widetilde{\mathcal{R}}_i := \mathcal{R}^{i}_n \cap \widetilde{\mathcal{R}} . 
\end{equation*} 
Since $\widetilde{R}$ has a horizontal side length $\frac{1}{4}$, note that there exists some $j\in \Z$ such that
\begin{equation*}
    \widetilde{R} \subseteq \bigcup_{i=1}^{n} \big(\mathcal{R}_n^{i,j} \cup \mathcal{R}_n^{i,j+1} \big) . 
\end{equation*}

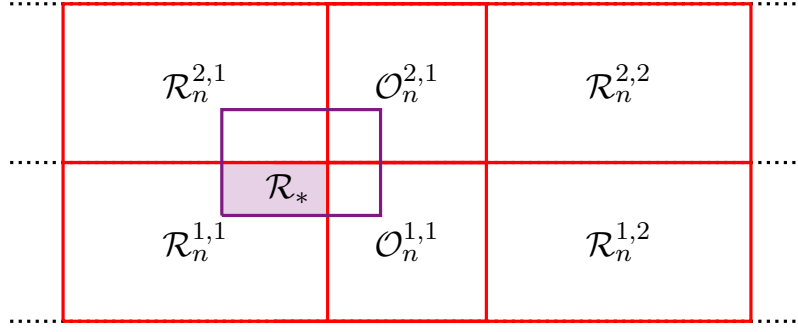
\begin{figure}
\centering
\begin{tikzpicture}[scale=.7]

   \draw[line width=1.2pt, dotted] (0,0) -- (15,0);
   \draw[line width=1.2pt, dotted] (0,3) -- (15,3);
    \draw[line width=1.2pt, dotted] (0,6) -- (15,6);

    \fill[darkblue!20, line width =1.2 pt] (4,2) -- (6,2) -- (6,3) -- (4,3) -- (4,2);

  \draw[red, line width =1.2 pt] (1,0) -- (1,3) -- (9,3) -- (9,0) -- (1,0);
  
    \draw[red, line width =1.2 pt] (6,0) -- (6,3) -- (14,3) -- (14,0) -- (5,0);

  \draw[red, line width =1.2 pt] (1,3) -- (1,6) -- (9,6) -- (9,3) -- (1,3);
  
    \draw[red, line width =1.2 pt] (6,3) -- (6,6) -- (14,6) -- (14,3) -- (5,3);

  \draw[darkblue, line width =1.2 pt] (4,2) -- (7,2) -- (7,4) -- (4,4) -- (4,2);

  \node[scale=1.2] (site) at (5.25,2.5){$\mathcal{R}_{\ast}$}; 

    
    
   


\node[scale=1.2] (site) at (3.5,1.5){$\mathcal{R}_n^{1,1}$};  
\node[scale=1.2] (site) at (11.5,1.5){$\mathcal{R}_n^{1,2}$};  
\node[scale=1.2] (site) at (7.5,1.5){$\mathcal{O}_n^{1,1}$}; 

\node[scale=1.2] (site) at (3.5,4.5){$\mathcal{R}_n^{2,1}$};  
\node[scale=1.2] (site) at (11.5,4.5){$\mathcal{R}_n^{2,2}$};  
\node[scale=1.2] (site) at (7.5,4.5){$\mathcal{O}_n^{2,1}$}; 



%
%
%
%
%
\end{tikzpicture}	
\caption{\label{fig:RectangleDLsmaller}Illustration of the subset $\mathcal{R}_{\ast}=\widetilde{\mathcal{R}}_1 \cap \mathcal{R}_n^{1,1} \setminus (\mathcal{O}^{1,1}_n \cup \mathcal{O}^{1,2}_n) \subseteq \widetilde{\mathcal{R}}$ of the rectangle $\widetilde{\mathcal{R}}$ marked in purple.} 
\end{figure}

We define $\mathcal{L}^{\prime}$ as follows: 
If $\widetilde{\mathcal{R}}_i \cap \mathcal{R}_n^{i,j} \setminus (\mathcal{O}^{i,j}_n \cup \mathcal{O}^{i,j+1}_n) \neq \emptyset$ (see Figure~\ref{fig:RectangleDLsmaller} for an illustration of this set) then take $\mathcal{L}^{\prime}$ as $\mathcal{L}^{i,1}$ on $\mathcal{R}^{i}_n$ (shifted by $\frac{j-1}{2}$) if $j$ is odd, i.e.,
\begin{equation*}
    \mathcal{L}^{\prime}(x,s;y,t)=\mathcal{L}^{i,1}\Big(x-\frac{j-1}{2},s;y-\frac{j-1}{2},t\Big)
\end{equation*} for all $(x,s;y,t) \in (\mathcal{R}^{i}_n)^2_{\uparrow}$, and take $\mathcal{L}^{\prime}=\mathcal{L}^{i,2}$ on $\mathcal{R}^{i}_n$ (shifted by $\frac{j-2}{2}$) if $j$ is even.
Similarly, if $\widetilde{\mathcal{R}}_i \cap \mathcal{R}_n^{i,j} \setminus (\mathcal{O}^{i,j}_n \cup \mathcal{O}^{i,j+1}_n) =\emptyset$, then we must have $\widetilde{R} \cap \mathcal{R}_n^{i,j+1} \neq \emptyset$, and we
take $\mathcal{L}^{\prime}=\mathcal{L}^{i,1}$ on $\mathcal{R}^{i}_n$ (shifted by $\frac{j-3}{2}$) if $j$ is even, and take $\mathcal{L}^{\prime}=\mathcal{L}^{i,2}$ on $\mathcal{R}^{i}_n$ (shifted by $\frac{j-2}{2}$) if $j$ is odd. In words, on pairs of points in $(\mathcal{R}^{i}_n)^2_{\uparrow}$, we take $\mathcal{L}^{\prime}$ as the full-space directed landscape in the construction of the unwrapped patched directed landscape $\Lpat$, which is used for a rectangle $\mathcal{R}$ with $\widetilde{R}_i \subseteq\mathcal{R}$.  
For all $(x,s;y,t) \in \R^{4}_{\uparrow}$ such that $(x,s),(y,t) \notin (\mathcal{R}^{i}_n)^2_{\uparrow}$ for some $i\in \lbr n \rbr$, we define $\mathcal{L}^{\prime}$ inductively using the metric composition \eqref{eq:InvTriangleSup}. More precisely, 
suppose that $\Lpa(x,s;y,t)$ is defined for $t\in \big(s,\frac{k}{n}\big]$ for some $k\geq i$ (with the base case $k=i$ covered above) then for $t\in \big(\frac{k}{n},\frac{k+1}{n}\big]$, 
\begin{equation}\label{eq:RecursionDLProof}
\mathcal{L}^{\prime}(x,s;y,t) := \sup_{(\tilde{x},\tilde{s}) \in \partial \mathcal{R}^{k}_n} \mathcal{L}^{\prime}(x,s;\tilde{x},\tilde{s}) +  \mathcal{L}^{\prime}(\tilde{x},\tilde{s};y,t). 
\end{equation}
Observe that under this construction, $\mathcal{L}^{\prime}$ is indeed a full-space directed landscape. 

In the remainder, for concreteness let us denote $\widetilde{\mathcal{R}}=\big[a, a+\frac{1}{4}\big] \times [b,b+\delta]$ for some $a\in \R$ and $b\geq 0$. 
We denote the geodesics from $(x,s)$ to $(y,t)$ in $\mathcal{L}^{\prime}$ by $\mathfrak{G}^{\prime}_{(x,s),(y,t)}$. Then the event
\begin{equation*}
    \mathsf{B}^{\prime}_{\delta} := \left\{ \pi^{\prime} \subseteq \Big[a -\frac{1}{8},a+\frac{3}{8}\Big] \times [b,b+\delta] \text{ for all } (x,s;y,t) \in (\widetilde{\mathcal{R}})^2_{\uparrow}\text{ and all }\pi^{\prime}\in \mathfrak{G}^{\prime}_{(x,s),(y,t)}\right\}
\end{equation*} satisfies, for some $c_0,C_0>0$ and all $\delta>0$,
\begin{equation}\label{eq:LowerBPatch}
    \P(\mathsf{B}^{\prime}_{\delta}) \geq 1- C_0\exp(-c_0\delta^{-1}).
\end{equation} This follows by bounding the transversal fluctuations on $\mathfrak{G}^{\prime}_{a,b;a,b+\delta}$ and $\mathfrak{G}^{\prime}_{a+\frac{1}{4},b;a+\frac{1}{4},b+\delta}$ using  Proposition~\ref{pro:ModerateDL}, together with the ordering of geodesics in the full-space directed landscape; see for example \cite[Lemma 2.7]{BGH:HausdorffDL}. 
Now by Lemma~\ref{lem:TypicallyGood} with $\varepsilon=n^{-1}$, we see that the event
\begin{equation*}
    \mathsf{C}^{\prime}_{n} := \left\{ \pi^{\prime}\text{ is } n^{-1}\text{-good for all } (x,s;y,t) \in (\widetilde{\mathcal{R}})^2_{\uparrow} \cap \mathcal{D}_{n^{-1}} \text{ and all }\pi^{\prime}\in \mathfrak{G}^{\prime}_{(x,s),(y,t)}\right\}
\end{equation*} satisfies, for some  $c_1,C_1>0$ and all  $n\in \N$,
\begin{equation}\label{eq:LowerCPatch}
    \P( \mathsf{C}^{\prime}_{n}) \geq 1-C_1\exp(-c_1n^{1/4}) .  
\end{equation}
Intuitively, the events $\mathsf{B}^{\prime}_{\delta}$ and  $\mathsf{C}^{\prime}_{n}$ guarantee that the geodesics in $\mathcal{L}^{\prime}$ between pairs of points in $(\widetilde{\mathcal{R}})^2_{\uparrow} \cap \mathcal{D}_{n^{-1}}$ can be taken as inherited geodesics in $\Lpat$. Now observe that on the intersection of the events $\mathsf{B}_n$ for $(\mathcal{L}^{i,1}_{n},\mathcal{L}^{i,2}_{n})_{i \in \lbr n \rbr}$ (to ensure the basic properties stated in Lemma~\ref{lem:Geodesics}) and the events $\mathsf{B}^{\prime}_{\delta}$ and $\mathsf{C}
^{\prime}_n$ for $\mathcal{L}^{\prime}$, 
we get that
\begin{equation}\label{eq:CompareCoupledPatched}
     \Lpat(x,s;y,t) = \mathcal{L}^{\prime}(x,s;y,t) \text{ for all } (x,s;y,t) \in (\widetilde{\mathcal{R}})^2_{\uparrow} \cap \mathcal{D}_{n^{-1}}.  
\end{equation} 
Similarly, on the event $\mathsf{B}_n$ for $(\mathcal{L}^{i,1}_{n},\mathcal{L}^{i,2}_{n})_{i \in \lbr n \rbr}$ and $\mathsf{B}^{\prime}_{\delta}$ for $\mathcal{L}^{\prime}$, we see that
\begin{equation}\label{eq:CompareCoupledPatched2}
     \Lpat(x,s;y,t) = \mathcal{L}^{\prime}(x,s;y,t) \text{ for all } (x,s;y,t) \in (\widetilde{\mathcal{R}})^2_{\uparrow} \cap \mathcal{V}_{\delta^{-1}} ,  
\end{equation} where we use item~\hyperref[eq:Patched2]{(ii)} from Lemma \ref{lem:Geodesics}  to split the path from $(x,s)$ to $(y,t)$ into subpaths with endpoints in $\mathcal{V}_n$. 
Combining now \eqref{eq:ProbabilityBound} from Lemma \ref{lem:Geodesics} with \eqref{eq:LowerBPatch} and \eqref{eq:LowerCPatch} for a lower bound on the probability of the events $\mathsf{B}_n$, $\mathsf{B}^{\prime}_{\delta}$ and $\mathsf{C}
^{\prime}_n$ (recalling $n\geq \delta^{-2}$), we conclude.
\end{proof}

In the following lemma, we provide explicit upper and lower  bounds on the lengths in the unwrapped patched directed landscape. 
\begin{lemma}\label{lem:FirstFullSpace}
Let $\varepsilon>0$ and $n\in \N$. Consider the rectangles 
$$\widetilde{\mathcal{R}}^{i,j}_{\varepsilon}=\Big[\frac{j-1}{16},\frac{j+1}{16}\Big]\times [(i-1)\varepsilon,i\varepsilon]$$ 
for $j\in \Z$ and $i \in \lbr \varepsilon^{-1}\rbr$, and define the events
 \begin{equation*}
     \mathsf{E}^{i,j}_{n,\varepsilon} :=\left\{ |\Lpat(x,s;y,t)|\leq  \log^{4/3}(\varepsilon^{-1}) \text{ for all } (x,s;y,t) \in (\widetilde{\mathcal{R}}^{i,j}_{\varepsilon})^2_{\uparrow}\cap \mathcal{D}_{\varepsilon} \right\} . 
\end{equation*} 
There exist $c,C>0$ such that for all $\varepsilon>0$ and all $n\in \N$ with $n \geq \varepsilon^{-2}$, 
\begin{equation}\label{eq:LowerGoal}
    \P\bigg( \mathsf{E}^{i,j}_{n,\varepsilon} \text{ holds for all } i \in \lbr \varepsilon^{-1}\rbr \text{ and } j\in \Z\bigg) \geq 1- C \varepsilon^{-1}\exp\big(-c \varepsilon^{-1/2}\big). 
\end{equation}
\end{lemma}
\begin{proof} 
By Lemma~\ref{lem:CouplePatchedToFull}, there exists for all $j \in \Z$ and $i \in \lbr \varepsilon^{-1}\rbr$ couples $\mathbf{P}_n^{i,j}$ between the unwrapped patched directed landscape $\Lpat$ on $\widetilde{\mathcal{R}}^{i,j}_{\varepsilon}$ and a full-space directed landscape $\mathcal{L}$ such that when $n \geq \varepsilon^{-2}$
\begin{equation}\label{eq:LowerAbsolutCompare}
    \mathbf{P}_n^{i,j}\left( \Lpat(x,s;y,t) = \mathcal{L}(x,s;y,t) \ \forall (x,s;y,t) \in (\widetilde{\mathcal{R}}^{i,j}_{\varepsilon})^2_{\uparrow}  \cap \mathcal{D}_{\varepsilon}\right) \geq 1-C_1\exp(-c_1\varepsilon^{-1/2})
\end{equation}
for some  $c_1,C_1>0$ that do not depend on $\varepsilon$, $n$, $i$ or $j$. For $(x,s;y,t) \in (\widetilde{\mathcal{R}}^{i,j}_{\varepsilon})^2_{\uparrow}$,  $$\frac{(x-y)^2}{t-s} \leq \varepsilon^{-1/2}|x-y|\leq 1$$ and $t-s\leq \varepsilon$. Thus, Proposition~\ref{pro:DLvalues} with $z$ of order $\varepsilon^{-1/3}$ ensures that
\begin{equation}\label{eq:LowerAbsolutValue}
       \mathbf{P}_n^{i,j}\left(|\mathcal{L}(x,s;y,t)| \leq \log^{4/3}(\varepsilon^{-1}) \ \forall (x,s;y,t) \in (\widetilde{\mathcal{R}}^{i,j}_{\varepsilon})^2_{\uparrow}  \cap \mathcal{D}_{\varepsilon}\right) \geq 1-C_2\exp(-c_2\varepsilon^{-1/2})
\end{equation} for some $c_2,C_2>0$ and all $\varepsilon>0$, where $c_2,C_2$ do not depend on $n$, $i$ or $j$. Putting  \eqref{eq:LowerAbsolutCompare} and \eqref{eq:LowerAbsolutValue} together we see that 
$$
\P( \mathsf{E}^{i,j}_{n,\varepsilon}) = \mathbf{P}_n^{i,j}( \mathsf{E}^{i,j}_{n,\varepsilon})\geq 1-C\exp(-c\varepsilon^{-1/2})
$$
for some $c,C>0$ and all $\varepsilon>0$, where $c,C$   do not depend on $n$, $i$ or $j$. Since, by the spatial periodicity of $\Lpat$, we have that $\mathsf{E}^{i,j}_{n,\varepsilon}=\mathsf{E}^{i,j+16}_{n,\varepsilon}$ 
for all $j\in \Z$, hence, we obtain \eqref{eq:LowerGoal} the above bound via a union bound over $i\in \lbr \varepsilon^{-1}\rbr$ and $j\in \lbr 16 \rbr$.
\end{proof}

Next, we provide upper bounds on the lengths in $\Lpat$. This uses a path decomposition similar to Section~3 of \cite{SS:TASEPcircle} on bounding the transversal fluctuations of geodesics in periodic last passage percolation.

\begin{lemma}\label{lem:PenaltyDistance} 
Let $\varepsilon>0$ and $n\in \N$. For all $M \in \N_0$, $i \in \lbr \varepsilon^{-1}\rbr$, and $z>0$, we define the events 
\begin{equation}\label{eq:EventAm}
    \mathsf{B}^{M,i}_{n,\varepsilon}(z) := \left\{ \sup_{ (i-1)\varepsilon\leq s < t \leq i\varepsilon,\, |x-y|\geq \frac{1}{8}M} \Lpat(x,s;y,t)  \leq -\frac{1}{260} M\varepsilon^{-1} + z\right\} . 
\end{equation}
There exist  $\varepsilon_0,c,C>0$ 
such that for all $\varepsilon \in (0,\varepsilon_0)$, all $n \geq \varepsilon^{-2}$, and all $z\in (0,\varepsilon^{-1/3}]$, 
\begin{equation}\label{eq:UniformUpperBoundDL}
    \P\Big( \mathsf{B}^{M,i}_{n,\varepsilon}(z)  \text{ holds for all } M \in \N_0 \text{ and } i \in \lbr \varepsilon^{-1} \rbr \Big) \geq 1- C\varepsilon^{-1}\exp\big(-c z^{3/2}\big) . 
\end{equation}
\end{lemma}
\begin{proof}
 We will only consider the case $i=1$, as the arguments are the same for general $i$ due to shift invariance. Intuitively, it suffices to show that with very high probability, along every path $\pi$, every horizontal transition of at least $\frac{1}{8}$ reduces the length 
 in the associated patched directed landscape by at least $\varepsilon^{-1}/260$. To formalize this, let us start by partitioning  $\R\times [0,\varepsilon]$ into  non-overlapping parts
\begin{equation}\label{eq:TargetRectangle}
 \mathcal{Q}^j_{\varepsilon} = \Big[\frac{1}{16}(j-1), \frac{1}{16}j \Big) \times [0,\varepsilon] .
\end{equation} For all $j\in \Z$, we write $\partial_+ \mathcal{Q}^j_{\varepsilon}$ and $\partial_- \mathcal{Q}^j_{\varepsilon}$ for the right and left boundary of $\mathcal{Q}^j_{\varepsilon}$,  respectively. We split the proof into two parts. As a first step, we  partition every path $\pi$ in $\R\times [0,\varepsilon]$ according to the parts spent in $\mathcal{Q}^j_{\varepsilon}$ for some $j\in \Z$. As a  second step, we bound the change of length on each of these parts separately. 
More precisely, for the first step,  note that without loss of generality, for every path $\pi$ connecting two points $(x,s)$ to $(y,t)$, we can decompose $\pi$ into finitely many parts according to the times $(t^{\pi}_{k})_{k\geq 1}$ (respectively $(s^{\pi}_{k})_{k\geq 1}$) such that we first (last) hit a given boundary $\partial_+ \mathcal{Q}^j_{\varepsilon}$ before hitting either a new boundary point in $\partial_+  \mathcal{Q}^{j-1}_{\varepsilon}$ or $\partial_+  \mathcal{Q}^{j+1}_{\varepsilon}$, or the site $(y,t)$. 
More precisely, we set $t^{\pi}_0=s$ and  
define the arrival times 
\begin{equation}\label{eq:Decomposition1}
t^{\pi}_{k} := \inf\{ \tilde{t} > t^{\pi}_{k-1} \colon (\pi(\tilde{t}),\tilde{t})  \in \partial_+\mathcal{Q}^j_{\varepsilon} \text{ with } w_{t^{\pi}_{k-1}} \notin \partial_+\mathcal{Q}^j_{\varepsilon} \text{ for some } j\in \Z\} ,
\end{equation} where the path $\pi$ first hits a new boundary $\partial_+\mathcal{Q}^j_{\varepsilon}$ as well as the departure times
\begin{equation}\label{eq:Decomposition2}
s^{\pi}_k := \sup\{ \tilde{s} \in [t^{\pi}_k,t^{\pi}_{k+1}] \colon w_{\tilde{s}} \in \partial_+\mathcal{Q}^j_{\varepsilon}\text{ with } w_{t^{\pi}_{k}} \in \partial_+\mathcal{Q}^j_{\varepsilon} \text{ for some } j\in \Z \} , 
\end{equation} where the path leaves the $\partial_+\mathcal{Q}^j_{\varepsilon}$ before hitting either $\partial_+\mathcal{Q}^{j-1}_{\varepsilon}$ or $\partial_+\mathcal{Q}^{j+1}_{\varepsilon}$.
We set $u^{\pi}_0=(x,s)$ and for $k\geq 1$, 
\begin{equation}\label{def:Upoints}
    u^{\pi}_{2k-1} := (\pi(t^{\pi}_k),t^{\pi}_k) \ \text{ as well as } \ u^{\pi}_{2k} := (\pi(s^{\pi}_k),s^{\pi}_k) . 
\end{equation}
Observe that $(u^{\pi}_k)_{k\geq 1}$ consists of $K \geq 2M$  points by our assumptions on the transversal distance between $(x,s)$ and $(y,t)$ to be at least $M/8$. Set $u^{\pi}_{K+1}=(y,t)$. A visualization of this path decomposition is given in Figure~\ref{fig:PathDecomposition}.
\begin{figure}
\centering
\begin{tikzpicture}[scale=.95]
   \draw[line width=1.2pt, dotted] (0,0) -- (11,0);
   \draw[line width=1.2pt, dotted] (0,3) -- (11,3);

  \draw[red, line width =1.2 pt] (1,0) -- (1,3) -- (10,3) -- (10,0) -- (1,0);
  
    \draw[red, line width =1.2 pt] (4,0) -- (4,3);
    \draw[red, line width =1.2 pt] (7,0) -- (7,3);

\draw[line width =1 pt] (6.5,0.05) to[curve through={(7,0.4)..(7.2,0.6)..(7,0.7)..(5.5,1.2) .. (4,1.45) .. (3.8,1.65) .. (4,1.85) ..(5.5,2.2) .. (7,2.45)..(7.3,2.75) ..(7,2.8)}] (6.5,2.95);

\node[scale=1] (site) at (2.5,1.5){$\mathcal{Q}^1_{\varepsilon}$};  
\node[scale=1] (site) at (5.5,1.5){$\mathcal{Q}^2_{\varepsilon}$};  
\node[scale=1] (site) at (8.5,1.5){$\mathcal{Q}^3_{\varepsilon}$};

 \node[scale=0.9] (site) at (6.5,-0.25){$(x,s)$};  
  \node[scale=0.9] (site) at (6.5,3.25){$(y,t)$}; 

\filldraw[deepblue,xshift=-2pt,yshift=-2pt] (6.5,-0.02) rectangle ++(4pt,4pt);

\filldraw[deepblue,xshift=-2pt,yshift=-2pt] (6.5,2.98) rectangle ++(4pt,4pt);

\filldraw[deepblue,xshift=-2pt,yshift=-2pt] (7,0.4) rectangle ++(4pt,4pt);

\filldraw[deepblue,xshift=-2pt,yshift=-2pt] (7,0.7) rectangle ++(4pt,4pt);

\filldraw[deepblue,xshift=-2pt,yshift=-2pt] (7,2.45) rectangle ++(4pt,4pt);

\filldraw[deepblue,xshift=-2pt,yshift=-2pt] (7,2.8) rectangle ++(4pt,4pt);

\filldraw[deepblue,xshift=-2pt,yshift=-2pt] (4,1.45) rectangle ++(4pt,4pt);

\filldraw[deepblue,xshift=-2pt,yshift=-2pt] (4,1.85) rectangle ++(4pt,4pt);

 \node[scale=0.7] (site) at (7.6,0.3){$(\frac{1}{8},t^{\pi}_1)$};  
 \node[scale=0.7] (site) at (7.6,0.8){$(\frac{1}{8},s^{\pi}_1)$};  

 \node[scale=0.7] (site) at (3.4,1.35){$(\frac{1}{16},t^{\pi}_2)$};  
 \node[scale=0.7] (site) at (3.4,1.95){$(\frac{1}{16},s^{\pi}_2)$}; 

 \node[scale=0.7] (site) at (7.6,2.35){$(\frac{1}{8},t^{\pi}_3)$};  
 \node[scale=0.7] (site) at (7.7,2.78){$(\frac{1}{8},s^{\pi}_3)$}; 
  


%
%
%
%
%
\end{tikzpicture}	
\caption{\label{fig:PathDecomposition}Path decomposition for $\pi$ from $(x,s)$ to $(y,t)$ with arrival times $(t^{\pi}_i)_{i \in \lbr 3 \rbr}$ and departure times $(s^{\pi}_i)_{i \in \lbr 3 \rbr}$. } 
\end{figure}
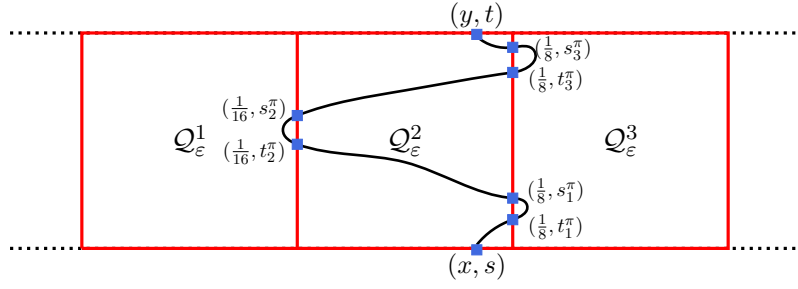

Now for the second part of the proof, let $\mathcal{L}$ be a full-space directed landscape. 
Define the event (of reasonably negative passage time between boundaries $\partial_+ \mathcal{Q}^j_{\varepsilon}$ and $\partial_- \mathcal{Q}^j_{\varepsilon}$)
\begin{equation}\label{eq:UniformModerateDL1Event}
\mathsf{A}^{\textup{full}}_{\varepsilon,j} := \Big\{\sup_{(x^{\prime},s^{\prime}) \in \partial_+ \mathcal{Q}^j_{\varepsilon}, (y^{\prime},t^{\prime}) \in \partial_- \mathcal{Q}^j_{\varepsilon}} \mathcal{L}(x^{\prime},s^{\prime};y^{\prime},t^{\prime}) \leq - \frac{1}{260}\varepsilon^{-1} - \varepsilon^{-\frac{1}{3}} \Big\} \geq 1- C_1\exp(-c_1\varepsilon^{-1}).
\end{equation}
For any $i\in \Z$, using Proposition~\ref{pro:DLvalues} with $z$ of order $\varepsilon^{-2/3}$ we find that there exist $c_1,C_1>0$  such that for all $\varepsilon>0$ 
\begin{equation}\label{eq:UniformModerateDL1}
\P\Big(\mathsf{A}^{\textup{full}}_{j,\varepsilon} \Big) \geq 1- C_1\exp(-c_1\varepsilon^{-1}) . 
\end{equation} Here, we note that $(y^{\prime}-x^{\prime})^2(t^{\prime}-s^{\prime})^{-1} \leq \varepsilon^{-1}/256$ for all $(x^{\prime},s^{\prime}) \in \partial_+ \mathcal{Q}^j_{\varepsilon}$ and $(y^{\prime},t^{\prime}) \in \partial_- \mathcal{Q}^j_{\varepsilon}$. 

We now use  Lemma~\ref{lem:CouplePatchedToFull} to transfer the result in \eqref{eq:UniformModerateDL1} about the full-space directed landscape into a result about the unwrapped patched directed landscape $\Lpat$. Define the event
\begin{equation*}
\mathsf{A}_{\varepsilon,j}^{(1)} :=\Big\{\sup_{(x^{\prime},s^{\prime}) \in \partial_+ \mathcal{Q}^j_{\varepsilon}, (y^{\prime},t^{\prime}) \in \partial_- \mathcal{Q}^j_{\varepsilon}} \Lpat(x^{\prime},s^{\prime};y^{\prime},t^{\prime}) \leq - \frac{1}{260}\varepsilon^{-1} - \varepsilon^{-1/3}\Big\} 
\end{equation*}
We invoke Lemma~\ref{lem:CouplePatchedToFull} with the parameter $\delta=\varepsilon$ and in the case where $\mathcal{V}_{\delta^{-1}}$ replaced $\mathcal{D}_{n^{-1}}$. Recalling $n\geq \varepsilon^{-2}$ by our assumptions, we see that there exist $c_2,C_2>0$ such that for all $\varepsilon>0$
\begin{equation}\label{eq:UniformModerateDL1Mod}
\P\Big(\mathsf{A}_{\varepsilon,j}^{(1)} \Big) \geq 1- C_2\exp(-c_2\varepsilon^{-1}) .  
\end{equation} 
Now, define the events (of a well-controlled upper bound)
\begin{align}\label{eq:UniformModerateDL2Event}
\begin{split}
\mathsf{A}_{\varepsilon,j,z}^{(2)} &:= \Big\{\sup_{(x^{\prime},s^{\prime}),(y^{\prime},t^{\prime}) \in \partial_+\mathcal{Q}^j_{\varepsilon}} \Lpat(x^{\prime},s^{\prime};y^{\prime},t^{\prime})  \leq  z \varepsilon^{1/3}\log^{4/3}(\varepsilon^{-1}) \Big\} , \\
\mathsf{A}_{\varepsilon,j,z}^{(3)}&:= \Big\{ \sup_{(x^{\prime},s^{\prime}) \in \mathcal{Q}^j_{\varepsilon}, (y^{\prime},t^{\prime})\in \partial_+ \mathcal{Q}^j_{\varepsilon} \cup  \partial_- \mathcal{Q}^j_{\varepsilon}}  \Lpat(x^{\prime},s^{\prime};y^{\prime},t^{\prime})  \leq  z \varepsilon^{1/3}\log^{4/3}(\varepsilon^{-1})\Big\}.
\end{split}
\end{align}
Similar to how we derived \eqref{eq:UniformModerateDL1Mod}, it follows from Proposition~\ref{pro:DLvalues} and Lemma~\ref{lem:CouplePatchedToFull} that there exist $c_3,C_3>0$  such that for all $\varepsilon>0$, and all $0<z\leq \varepsilon^{-1/3}$
\begin{equation}\label{eq:UniformModerateDL2}
 \P\Big(\mathsf{A}_{\varepsilon,j,z}^{(2)}  \Big) \geq 1- C_3\exp\big(-c_3z^{3/2} \big),\qquad 
 \P\Big( \mathsf{A}_{\varepsilon,j,z}^{(3)}\Big) \geq 1- C_3\exp\big(-c_3z^{3/2} \big).
\end{equation} 
We now assume that the event $\mathsf{A}_{\varepsilon,z}:=\bigcap_{j\in \Z}\mathsf{A}_{\varepsilon,j}^{(1)} \cap \mathsf{A}_{\varepsilon,j,z}^{(2)} \cap \mathsf{A}_{\varepsilon,j,z}^{(3)}$ holds. Since  
\begin{equation*}
\mathsf{A}_{\varepsilon,j}^{(1)} \cap \mathsf{A}_{\varepsilon,j,z}^{(2)} \cap \mathsf{A}_{\varepsilon,j,z}^{(3)} = \mathsf{A}_{\varepsilon,j+16}^{(1)} \cap \mathsf{A}_{\varepsilon,j+16,z}^{(2)} \cap \mathsf{A}_{\varepsilon,j+16,z}^{(3)}
\end{equation*} for all $j\in \Z$, by a union bound over the events in \eqref{eq:UniformModerateDL1Mod} and \eqref{eq:UniformModerateDL2} for $j\in \lbr 16 \rbr$ we conclude that there exists $c,C>0$ such that for all $\varepsilon>0$,  all $n\geq \varepsilon^{-2}$, and all $z\in (0,\varepsilon^{-1/3}]$,
$$
\P(\mathsf{A}_{\varepsilon,z})\geq 1-C\exp\big(-cz^{3/2}\big).
$$
We finish by arguing that on the event $\mathsf{A}_{\varepsilon,z}$, $\mathsf{B}^{M,i}_{n,\varepsilon}(z) $  holds for all $M \in \N_0$ and $i \in \lbr \varepsilon^{-1} \rbr$. This along with the above probability bound for $\mathsf{A}_{\varepsilon,z}$ implies \eqref{eq:UniformUpperBoundDL}. To see why  $\mathsf{B}^{M,i}_{n,\varepsilon}(z)$ holds, it suffices to observe that by Lemma \ref{lem:Geodesics}, either $\Lpat(x,s;y,t)=-\infty$ (which is certainly $\leq -M\varepsilon^{-1}/260 + z$) or there exist at least one unwrapped inherited geodesic $\bar{\pi}$ from $(x,s)$ to $(y,t)$. In that case, we can use \eqref{eq:PatchedCompo} repeatedly to decompose  
$$\Lpat(x,s;y,t)= \sum_{\ell=0}^{K} \Lpat(u_\ell;u_{\ell+1}),$$ 
where we recall the points $u_{\ell}$ from \eqref{def:Upoints}.
On the event  $\mathsf{A}_{\varepsilon,z}$, these terms alternate from being bounded above by $z\varepsilon^{1/3}\log^{4/3}(\varepsilon^{-1})$ and $-\varepsilon^{-1}/260 - \varepsilon^{-1/3}$. Taking the constants $c,C>0$ sufficiently large, we can assume that $\varepsilon^{1/3}\log^{4/3}(\varepsilon^{-1})\leq \frac{1}{3}$. Thus, since $|x-y|\geq M/8$ ensures at least $2M$ terms in the decomposition, we see that the event in $\mathsf{B}^{M,i}_{n,\varepsilon}(z)$ holds, hence completing the proof of the lemma.
\end{proof}

 Recall from \eqref{eq:ExtendPatched} the wrapped patched directed landscape $\Lpa$, which can be seen as  a function of the unwrapped landscape $\Lpat$. We have the following consequence of the previous statements on the wrapped and unwrapped patched directed landscape. 
\begin{lemma}\label{lem:Noncrossing}
  For $\delta>0$, $i \in \lbr \delta^{-1}\rbr$,  $n\in \N$ and $a\in \R$, we define the events
  \begin{equation*}
     \mathsf{A}_{n,i}^{\delta,a} := \left\{ \Lpa(x,s;y,t)=\Lpat(x,s;y,t) \,\,\, \forall (x,s;y,t) \in \big([a,a+1)\times [(i-1)\delta,i\delta]\big)^2_{\uparrow} \cap \mathcal{V}_{\delta^{-1}} \cap \mathcal{D}_{(4\delta)^2} \right\} , 
  \end{equation*} 
  where the arguments $x$ and $y$ in $\Lpa$ are understood modulo $1$.
  Then there exists  $c_1,C_1>0$ such that for all $\delta>0$, $i \in \lbr \delta^{-1}\rbr$, and  $n \geq \delta^{-2}$,  
    \begin{equation}\label{eq:GoodShort}
      \P( \mathsf{A}^{\delta,a}_{n,i} ) \geq 1- C_1 \exp\big(-c_1 \delta^{-1/2}\big) . 
  \end{equation}
Further, there exist $c_2,C_2>0$ such that for all $n\in \N$, there is a coupling $\mathbf{P}$ between a full-space directed landscape $\mathcal{L}$ and a wrapped patched directed landscape $\Lpa$ such that for all $n \geq \delta^{-2}$,
\begin{align}\label{eq:GoodShort2}
\begin{split}
    &\mathbf{P}\Big( \Lpa(x,s;y,t) = \mathcal{L}(x,s;y,t) \,\,\, \forall (x,s;y,t)\in \big([a,a+1)\times [(i-1)\delta,i\delta]\big)^2_{\uparrow} \cap \mathcal{V}_{\delta^{-1}} \cap \mathcal{D}_{(4\delta)^2}\Big) \\&\qquad\geq 1 - C_2\exp(-c_2\delta^{-1/2}) .
\end{split}
\end{align}
\end{lemma}
\begin{proof} 
We will assume $a=0$ and $i=1$ as the other cases follow by time shift-invariance. Let us first demonstrate \eqref{eq:GoodShort}. Consider the event $\mathsf{B}_{n,\delta}^{6,1}(\delta^{-1/3})$ from Lemma~\ref{lem:PenaltyDistance}. By that lemma, there exists $c_3,C_3>0$ such that for all $\delta>0$ and $n\geq \delta^{-2}$,  $\mathsf{B}_{n,\delta}^{6,1}(\delta^{-1/3})$ holds with probability at least $1-C_3 \exp(-c_3 \delta^{-1/2})$. On this event, it follows that (by taking without loss of generality $\delta$ small enough so that $-\tfrac{6}{260}\delta^{-1} + \delta^{-1/3} \leq -\tfrac{1}{45} \delta^{-1}$)
    \begin{equation}\label{eq:BoundPatchedFar1}
        \sup_{\ell \in \Z \setminus \{0\}} \Lpat(x+\ell,s;y,t) < - \frac{1}{45}\delta^{-1} \text{ for all } (x,s;y,t) \in \big([0,1)\times [0,\delta]\big)^2_{\uparrow} \cap  \mathcal{V}_{\delta^{-1}} .
    \end{equation}
In other words, the event $\mathsf{B}_{n,\delta}^{6,1}(\delta^{-1/3})$ ensures that if the sup on the left-hand side above is taken over all $\ell\in \Z$ (which defines $\Lpa$ via \eqref{eq:ExtendPatched}), then the contribution away from $\ell=0$ is bounded below by $- \frac{1}{45}\delta^{-1} $. Now we show that the $\ell=0$ term likely has a large contribution, and hence gives the value of the sup. To this end, we utilize Lemma~\ref{lem:CouplePatchedToFull} in order to couple $\Lpat$ with a full-space directed landscape while using Proposition~\ref{pro:DLvalues} with $z=\delta^{-1/3}$ to lower-bound the full-space directed landscape. Doing this, we see that there exists $c_4,C_4>0$ such that for all $\delta>0$ and $n\geq \delta^{-2}$,
\begin{equation}\label{eq:BoundPatchedFar2}
      \Lpat(x,s;y,t) > - \frac{1}{60}\delta^{-1} \text{ for all } (x,s;y,t) \in \big([0,1)\times [(i-1)\delta,i\delta]\big)^2_{\uparrow} \cap  \mathcal{V}_{\delta^{-1}} \cap \mathcal{D}_{(4\delta)^2}. 
\end{equation} Note that in our application of Proposition~\ref{pro:DLvalues}, we used the fact that $|x-y|^2(t-s)^{-1}\leq \delta^{-1}|x-y|/4\leq \delta^{-1}/64$ for all $(x,s;y,t) \in \mathcal{V}_{\delta^{-1}} \cap \mathcal{D}_{(4\delta)^2}$. Recalling again the definition of $\Lpa$ in \eqref{eq:ExtendPatched}, we see that when \eqref{eq:BoundPatchedFar1} and \eqref{eq:BoundPatchedFar2} hold, it follows that $\Lpa$ and $\Lpat$ agree on the desired range. Noting the bounds on the probabilities of these events from above, we obtain \eqref{eq:GoodShort}.
  The second statement \eqref{eq:GoodShort2} is now immediate from \eqref{eq:GoodShort} together with Lemma~\ref{lem:CouplePatchedToFull} in order to couple on $\mathcal{R}$ the unwrapped patched directed landscape and a full-space directed landscape. 
\end{proof}

Next, we control the transversal fluctuations of inherited geodesics in the unwrapped patched directed landscape. Recall from Definition~\ref{def:DeltaRegular} the notion of a geodesic to be $\varepsilon$-good.

\begin{lemma}\label{lem:LocalFluctuationsPatched}
There exist $c,C>0$ such that for all $\delta>0$ and $n\in \N$ with $n\geq \delta^{-2}$, 
\begin{equation}\label{eq:GoodBound}
 \mathsf{G}_{n,\delta} :=   \left\{ \pinpmt_{x,s;y,t} \text{ is } \delta\text{-good for all } (x,s;y,t) \in \SetR \cap \mathcal{D}_{\delta^{1/2}}\right\}  
\end{equation}
satisfies
\begin{equation}\label{eq:GoodBound2}
     \P\left(  \mathsf{G}_{n,\delta} \right) \geq 1- C\exp\big(-c\delta^{-1/8}\big) . 
\end{equation}
\end{lemma}
\begin{proof}
Suppose that the event $\mathsf{B}_n$ from Lemma~\ref{lem:Geodesics} holds, and recall $(z_i)_{i \in \lbr k \rbr}$ from \eqref{def:GoodIncrements} for a given pair of sites $(x,s;y,t) \in \mathcal{D}_{\delta^{1/2}}$. Consider any path $\pi$ from $(x,s)$ to $(y,t)$ for $(x,s;y,t) \in \mathcal{D}_{\delta^{1/2}}$, such that for some $i\in \lbr k+1\rbr$, some $x_i,y_i \in \R$ with $|x_i-y_i|> \frac{1}{8}$, and some $\ell,\ell^{\prime} \in [z_{i-1},z_i]$, $\pi(x_i)=\ell$ and $\pi(y_i)=\ell^{\prime}$. Showing \eqref{eq:GoodBound} is equivalent to showing that with probability at least $1- C\exp\big(-c\delta^{-1/8}\big)$, no such path $\pi$ can be an inherited geodesic in $\Lpat$.  
To this end, by the upper bound on the lengths in Lemma~\ref{lem:PenaltyDistance} with $\varepsilon=\delta$, $z=\delta^{-1/3}$ and $M=1$, since $n\geq \delta^{-2}$, there exists an event $\mathsf{A}^{(1)}_{\delta}$ with
\begin{equation*}
    \P\big( \mathsf{A}_{\delta}^{(1)} \big) \geq 1-C_1\exp(-c_1\delta^{-1/2})
\end{equation*}
for some  $c_1,C_1>0$, and all $\delta>0$, such that on the event $\mathsf{A}^{(1)}_{\delta}$, we get that
\begin{equation}\label{eq:RoughBound1}
    \Lpat(x_i,\ell;y_i,\ell^{\prime}) \leq - \frac{1}{280}\delta^{-1}.
\end{equation}  
\begin{figure}
\centering
\begin{tikzpicture}[scale=.95]

   \draw[line width=1.2pt, dotted] (0,0) -- (11,0);
   \draw[line width=1.2pt, dotted] (0,3) -- (11,3);

  \draw[red, line width =1.2 pt] (1,0) -- (10,0);
  \draw[red, line width =1.2 pt] (1,3) -- (10,3); 

  \draw[darkblue, dashed, line width =1.2 pt] (1,0.75) -- (10,0.75);
  \draw[darkblue, dashed, line width =1.2 pt] (1,1.5) -- (10,1.5); 
  
  \draw[red, dashed, line width =1.2 pt] (4,0) -- (4,3);
  \draw[red, dashed, line width =1.2 pt] (7,0) -- (7,3);

\draw[line width =1 pt] (6.5,0.05) to[curve through={(7,0.4)..(7.2,0.6)..(7,0.87) .. (4,1.35) .. (3.8,1.65) .. (4,1.85) ..(5.5,2.2) .. (6,2.5)}] (6.5,2.95);


 \node[scale=0.9] (site) at (6.5,-0.25){$(x,s)$};  
  \node[scale=0.9] (site) at (6.5,3.25){$(y,t)$}; 

\filldraw[deepblue,xshift=-2pt,yshift=-2pt] (6.5,-0.02) rectangle ++(4pt,4pt);

\filldraw[deepblue,xshift=-2pt,yshift=-2pt] (6.5,2.98) rectangle ++(4pt,4pt);


\filldraw[deepblue,xshift=-2pt,yshift=-2pt] (7,0.87) rectangle ++(4pt,4pt);



\filldraw[deepblue,xshift=-2pt,yshift=-2pt] (4,1.35) rectangle ++(4pt,4pt);


 \node[scale=0.7] (site) at (7.5,1){$(x_i,\ell)$};   
 \node[scale=0.7] (site) at (3.4,1.2){$(y_i,\ell^{\prime})$};



%
%
%
%
%
\end{tikzpicture}	
\caption{\label{fig:SteepPath}Path decomposition for $\pi$ on the event  $\mathsf{A}^{(2)}_{\delta}$. The dashed vertical lines in red are at least $\frac{1}{8}$ apart while the dashed horizontal lines in purple are at heights $z_{i-1}$ and $z_i$ for some $i$, respectively, so that $\ell,\ell^{\prime}\in [z_{i-1},z_i]$.  } 
\end{figure}
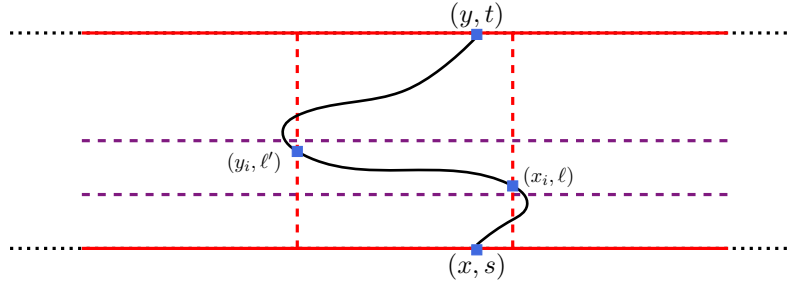
By Lemma~\ref{lem:PenaltyDistance} for $\varepsilon=\delta^{1/2}>0$, $M=0$ and $z=\delta^{-1/6}$, there exists an event $\mathsf{A}^{(2)}_{\delta}$ with
\begin{equation*}
    \P\big( \mathsf{A}^{(2)}_{\delta} \big) \geq 1-C_2\exp(-c_2\delta^{-1/4})
\end{equation*} for some $c_2,C_2>0$, and all $\delta>0$, such that on the event $\mathsf{A}^{(2)}_{\delta}$, we have (see Figure~\ref{fig:SteepPath})
\begin{equation}\label{eq:RoughBound2}
    \Lpat(x,s;x_i,\ell) + \Lpat(y_i,\ell^{\prime};y,t) \leq 2\delta^{-5/6} . 
\end{equation} Combining \eqref{eq:RoughBound1} and \eqref{eq:RoughBound2}, we find an event $\mathsf{A}_{\delta}$ of probability at least $1-C_3\exp(-c_3\delta^{-1/4})$ with some  $c_3,C_3>0$ such that on $\mathsf{A}_{\delta}$, we have from the definition of the length of a path $\pi$ in \eqref{def:GeodesicPatched} and item~\hyperref[eq:Patched3]{(iii)} of Lemma~\ref{lem:Geodesics} that
\begin{equation}\label{eq:UpperDecompo}
  \int \textup{d}\Lpat\circ \pi \leq  \Lpat(x,s;x_i,\ell)+ \Lpat(x_i,\ell ;y_i,\ell^{\prime}) + \Lpat(y_i,\ell^{\prime};y,t) \leq - \frac{1}{300}\delta^{-1} 
\end{equation} for all $\delta>0$ small enough. Since $(x,s;y,t)\in \mathcal{D}_{\delta^{1/2}}$ by our assumptions, we can find for a sequence $(u_i)_{i \in \lbr k \rbr}$ and $u_0=(x,s)$, $u_{k+1}=(y,t)$ with some $k\in \N$ such that
\begin{equation*}
   (u_{i-1};u_{i}) \in \mathcal{D}_{\delta^{1/2}} \text{ and }  u_i=(x_i^{\prime},i\delta^{1/2})
\end{equation*} for all $i \in \lbr k \rbr$, and some $(x^{\prime}_i)_{i \in \lbr k \rbr}$.
Then from Lemma \ref{lem:FirstFullSpace} with $\varepsilon=\delta^{1/2}$ and $z=\varepsilon^{-1/6}=\delta^{-1/12}$, and a union bound over $i\in \lbr \varepsilon^{-1}\rbr$, we see that 
\begin{equation}\label{eq:LowerBoundWeight}
     \Lpat(x,s;y,t) \geq  \delta^{-1/2} \min_{i \in \lbr k+1\rbr} \Lpat(u_{i-1};u_i)  \geq - C_5 \delta^{-1/2} \log^{4/3}(\delta^{-1}) 
\end{equation} with probability at least $1-C_4\exp(-c_4\delta^{-1/8})$ for some  $c_4,C_4,C_5>0$. On the event \eqref{eq:LowerBoundWeight}, we get  a lower bound on the length of the inherited geodesic from $(x,s)$ to $(y,t)$ in $\Lpat$, and on the event \eqref{eq:UpperDecompo} we get an upper bound on the length of any path $\pi$ satisfying the conditions at the beginning of this proof. Therefore, we see that on the events \eqref{eq:LowerBoundWeight} and \eqref{eq:UpperDecompo}, such a path $\pi$ can not be an inherited geodesic from $(x,s)$ to $(y,t)$, allowing us to conclude \eqref{eq:GoodBound2}. 
\end{proof}

We have now all tools in order to couple two wrapped patched directed landscapes with respect to $n$ and $2n$ for all pairs of points in $\mathcal{D}_{\varepsilon}$ with some $\varepsilon>0$ of order $n^{-1/16}$. Let us stress that though the exponent $\frac{1}{16}$ may be improved, it suffices as we ensure that any pair of points in $\SetT$ is contained in $\mathcal{D}_{\varepsilon}$ for $n$ large enough.
In the following, we slightly overload notation and write $\SetT \cap \mathcal{D}_{\varepsilon}$ for the set where we first take  $\mathcal{D}_{\varepsilon}$, and then project onto $\SetT$.

\begin{lemma}\label{lem:PatchedDLsConsistently}
There exist  $c,C>0$ and a coupling $\mathbf{P}$ of $\mathcal{L}_{n}$ under the law $\mathbf{P}_n$ and $\mathcal{L}_{2n}$ under the law $\mathbf{P}_{2n}$  such that all $\delta>0$ and $n \geq 2\delta^{-2}$, 
\begin{equation}\label{eq:CoupledDifferentn}
\mathbf{P}\left( \mathsf{D}_{n,\delta} \right) 
\geq 1- C\exp\big(-c\delta^{-\frac{1}{32}}\big),
\end{equation} 
where the event
\begin{equation*}
    \mathsf{D}_{n,\delta} := \left\{ \big(\Lpa(x,s;y,t)\big)_{(x,s;y,t) \in \SetT \cap \mathcal{D}_{\delta^{1/8}}}=\big(\mathcal{L}_{2n}(x,s;y,t)\big)_{(x,s;y,t) \in \SetT \cap \mathcal{D}_{\delta^{1/8}}}  \right\}.
\end{equation*}
\end{lemma}
\begin{proof}
We start by showing that there exists a coupling $\bar{\mathbf{P}}$ between $\Lpat$ and $\bar{\mathcal{L}}_{2n}$ such that with probability at least $1- C_0\exp\big(-c_0\delta^{-1/8}\big)$ for some $c_0,C_0>0$, 
\begin{equation}\label{eq:UnwrappedTotalCoupling}
\big(\Lpat(x,s;y,t)\big)_{(x,s;y,t) \in \SetR \cap \mathcal{D}_{\delta^{1/2}}}=\big(\bar{\mathcal{L}}_{2n}(x,s;y,t)\big)_{(x,s;y,t) \in \SetR \cap \mathcal{D}_{\delta^{1/2}}} . 
\end{equation} 
To this end, we define the sets 
\begin{equation*}
    \bar{\mathcal{R}}_{\delta}^{i,j}:=\Big[\frac{1}{16}j,\frac{1}{16}(j+4)\Big]\times \Big[\frac{i-1}{2}\delta,\frac{i+1}{2}\delta\Big] 
\end{equation*}
 for all $i\in \lbr 2\delta^{-1} \rbr$ and $j\in \Z$. Set 
\begin{equation*}
   \mathcal{\bar{R}}^{2,\delta}_{\uparrow} := \bigcup_{j\in \Z} \bigcup_{i\in \lbr 2\delta^{-1} \rbr} (\bar{\mathcal{R}}_{\delta}^{i,j}\times \bar{\mathcal{R}}_{\delta}^{i,j}) \cap \mathcal{V}_{\delta^{-1}} \cap\mathcal{D}_{(4\delta)^2} . 
\end{equation*}
 We claim that for all $i,j$ there exists a coupling $\mathbf{P}_{i,j}$ such that
\begin{equation}\label{eq:PatchiwiseCoupling}
\begin{split}
  &\mathbf{P}_{i,j}\left( \Lpat(x,s;y,t)=\bar{\mathcal{L}}_{2n}(x,s;y,t) \,\,\, \forall (x,s;y,t)\in (\bar{\mathcal{R}}_{\delta}^{i,j}\times \bar{\mathcal{R}}_{\delta}^{i,j}) \cap \mathcal{V}_{\delta^{-1}} \cap\mathcal{D}_{(4\delta)^2}  \right) \\
  &\geq 1- C_1 \exp\big(-c_1\delta^{-1/2}\big)   
  \end{split}
\end{equation} for some constants $c_1,C_1>0$ and all $\delta>0$ with $n\geq 2 \delta^{-2}$. This follows by coupling $\Lpat$ and $\bar{\mathcal{L}}_{2n}$ via Lemma~\ref{lem:Noncrossing} to the same full-space directed landscape. Note that in this way, we at the same time ensure that the couplings $(\mathbf{P}_{i,j})_{i\in [2\delta^{-1}],j\in \Z}$ are consistent, i.e., for any $i,j,i^{\prime},j^{\prime}$, we can view the couplings $\mathbf{P}_{i,j}$ and $\mathbf{P}_{i^{\prime},j^{\prime}}$ on an enlarged space such that on an event of probability at least $1- 2C_1 \exp\big(-c_1\delta^{-1/2}\big)$, all $(x,s;y,t)\in (\bar{\mathcal{R}}_{\delta}^{i,j}\times \bar{\mathcal{R}}_{\delta}^{i,j}) \cap \big(\bar{\mathcal{R}}_{\delta}^{i^{\prime},j^{\prime}}\times \bar{\mathcal{R}}_{\delta}^{i^{\prime},j^{\prime}}\big) \cap \mathcal{V}_{\delta^{-1}} \cap\mathcal{D}_{(4\delta)^2}$ are under 
$\mathbf{P}_{i,j}$ and $\mathbf{P}_{i^{\prime},j^{\prime}}$ assigned the same value $\Lpat(x,s;y,t)=\bar{\mathcal{L}}_{2n}(x,s;y,t)$. 
Thus, by \eqref{eq:PatchiwiseCoupling} and the fact that $\Lpat$ and $\bar{\mathcal{L}}_{2n}$ are $1$-periodic, there exists a coupling $\mathbf{P}$ such that  
\begin{equation}\label{eq:UnifiedCoupling}
  \mathbf{P}\left( \Lpat(x,s;y,t)=\bar{\mathcal{L}}_{2n}(x,s;y,t) \ \forall  (x,s;y,t)\in  \mathcal{\bar{R}}^{2,\delta}_{\uparrow} \right) \geq 1- C_2 \exp(-c_2\delta^{-1/2})    
\end{equation} for some  $c_2,C_2>0$ and all $\delta>0$ small enough. Recall the events $\mathsf{G}_{n,\delta}$ from \eqref{eq:GoodBound} in Lemma~\ref{lem:LocalFluctuationsPatched} which state that all inherited geodesics $\pinpmt_{(x,s;y,t)}$ for pairs of points in $\SetR\cap \mathcal{D}_{\delta^{1/2}}$ are all $\delta$-good, and that there exist $c_3,C_3>0$ such that 
\begin{equation*}
    \P(\mathsf{G}_{n,\delta}) \geq 1-C_3\exp(-c_3\delta^{-1/8})
\end{equation*}
 for all $\delta>0$ and $n\geq \delta^{-2}$. The remainder of the argument in order to obtain  \eqref{eq:UnwrappedTotalCoupling} is to reduce it to \eqref{eq:UnifiedCoupling} by using the goodness property and path decomposition. For all $(x,s;y,t) \in \big(\SetR \cap \mathcal{D}_{\delta^{1/2}} \big) \cap  \mathcal{\bar{R}}^{2,\delta}_{\uparrow}$, we conclude immediately by \eqref{eq:UnifiedCoupling}. Now for all $(x,s;y,t) \in \big(\SetR \cap \mathcal{D}_{\delta^{1/2}} \big) \setminus  \mathcal{\bar{R}}^{2,\delta}_{\uparrow}$, note that we must have $|t-s|\geq \delta/2$. Hence, on $\mathsf{G}_{n,\delta/2}$, the definition of a $\delta/2$-good geodesic ensures that for any $(x,s;y,t)$, we find a sequence of points $(z'_{\ell})_{\ell \in \lbr k'\rbr}$ with some $k'\in \N$ such that $z'_{0}=(x,s)$ and $z'_{k'+1}=(y,t)$, as well as
\begin{equation*}
    \Lpat (x,s;y, t) = \sum_{i=0}^{k'}  \Lpat
(z'_i; z'_{i+1}) ,
\end{equation*}
with the property that any pair of consecutive points $z'_{i-1}=(x'_{i-1},y'_{i-1})$ and $z'_{i}=(x'_{i},y'_{i})$ with some $i\in \lbr k'+1 \rbr$  satisfies
\begin{equation}\label{eq:DistancesBound}
    | x'_{i-1} - x'_{i} | \leq \frac{1}{8} \quad \text{and} \quad   y'_{i} - y'_{i-1} \in \Big[ \frac{\delta}{2},\delta\Big] .
\end{equation}
 Intuitively, this can be ensured by taking the points $(z_i)_{i\in \lbr k \rbr}$ in the definition of an $\delta/2$-good geodesic, and removing $z_1$ and $z_{k}$ if necessary. In particular, we see from \eqref{eq:DistancesBound} that  
\begin{equation*}
    (z'_{\ell-1};z'_{\ell}) \in (\bar{\mathcal{R}}_{\delta}^{i,j} \times \bar{\mathcal{R}}_{\delta}^{i,j}) \cap \mathcal{V}_{\delta^{-1}} \cap \mathcal{D}_{(4\delta)^2} \text{ for some } i\in \lbr 2\delta^{-1} \rbr\text{ and } j\in \Z 
\end{equation*} for all $\ell \in \lbr k\rbr$. Combining now \eqref{eq:GoodBound} and \eqref{eq:UnifiedCoupling}, we obtain \eqref{eq:UnwrappedTotalCoupling}.

We will now deduce \eqref{eq:CoupledDifferentn} for a wrapped patched directed landscape from \eqref{eq:UnwrappedTotalCoupling}, which deals with an unwrapped patched directed landscape. We claim that for some $c_4,C_4>0$, and all $\delta>0$ and $n\geq 2 \delta^{-2}$
\begin{equation}\label{eq:LocateHitting}
\begin{split}\
   \P\Big( \Lpat(x,s;y,t) &> \sup_{\ell\in \Z\colon(x+\ell,s;y,t)\notin \mathcal{D}_{\delta^{1/2}}} \Lpat(x+\ell,s;y,t) \, \forall (x,s;y,t)\in \SetR \cap\mathcal{D}_{\delta^{1/8}}   \Big) \\
    &\geq 1- C_4 \exp\Big(-c_4 \delta^{-\frac{1}{32}}\Big) .
    \end{split}
\end{equation} 
Assuming \eqref{eq:LocateHitting} holds, we see that for all $(x,s;y,t)\in \SetT \cap \mathcal{D}_{\delta^{1/8}}$, the supremum in 
\begin{equation*}
    \Lpa(x,s;y,t) := \sup_{\ell \in \Z}  \Lpat(x+\ell,s;y,t)
\end{equation*} is achieved for some $(x+\ell,s;y,t) \in \SetT \cap \mathcal{D}_{\delta^{1/2}}$, and together with \eqref{eq:UnwrappedTotalCoupling}, we conclude~\eqref{eq:CoupledDifferentn}.

It remains to show \eqref{eq:LocateHitting}. Set for all $\delta>0$
\begin{equation*}
    K_{\delta} := \left\{ (x+\ell,s;y,t)\in \SetR \cap\mathcal{D}^{\complement}_{\delta^{1/2}} \, \colon \, (x,s;y,t) \in \mathcal{D}_{\delta^{1/8}} , \ell \in \Z \right\} . 
\end{equation*}
We claim that for some $c_5,C_5>0$, and all $\delta>0$ and $n\geq 2 \delta^{-2}$
\begin{equation}\label{eq:UpperAux}
    \P\Big(  \Lpat(x,s;y,t) \leq - \frac{1}{9}\delta^{-3/16}\, \forall (x,s;y,t)\in K_{\delta}   \Big) \geq 1- C_5 \exp\Big(-c_5 \delta^{-\frac{3}{16}}\Big) .
\end{equation}
For all $(x,s;y,t) \in  K_{\delta}$ so that $|x-y| \leq \delta^{1/16}/8$, the bound is immediate from Lemma~\ref{lem:PenaltyDistance}. Similarly, when the geodesic $\pinpt_{x,s;y,t}$ satisfies  $\TF(\pinpt_{x,s;y,t})> \frac{1}{8}$, the bound follows by same arguments as in \eqref{eq:UpperDecompo}. For all $(x,s;y,t) \in \SetR \cap\mathcal{D}^{\complement}_{\delta^{1/2}}$ so that $\TF(\pinpt_{x,s;y,t})\leq \frac{1}{8}$ and $|x-y| \geq \delta^{1/16}/8$, the bound follows from Proposition~\ref{pro:ModulusOfContinuity} for $\tau=\xi=1$ and $\phi=\delta^{-1/8}$ on the full-space directed landscapes used in the definition of $\Lpat$, noting that
\begin{equation*}
    \frac{(x-y)^2}{t-s} \geq - \frac{1}{8}\delta^{-1/4}\cdot \delta^{1/16} = - \frac{1}{8}\delta^{-3/16} . 
\end{equation*}
Similarly, we claim that for some $c_6,C_6,C_7>0$, and all $\delta>0$ and $n\geq 2 \delta^{-2}$
\begin{equation}\label{eq:LowerAux}
\begin{split}
       \P\Big(  \Lpat(x,s;y,t) &\geq -C_7\delta^{-1/8}\log^{4/3}(\delta^{-1}) \, \forall (x,s;y,t)\in \SetR \cap\mathcal{D}_{\delta^{1/8}}   \Big) \\
       &\geq 1- C_6 \exp\Big(-c_6 \delta^{-\frac{1}{32}}\Big) .
       \end{split}
\end{equation}
This follows by the same arguments as in \eqref{eq:LowerBoundWeight} for a lower bound on $\Lpat(x,s;y,t)$ for all $(x,s;y,t)\in \mathcal{D}_{\delta^{1/8}}$
(replacing $\delta$ by $\delta^{1/4}$ therein). Combining \eqref{eq:UpperAux} and \eqref{eq:LowerAux}, we get \eqref{eq:LocateHitting}, allowing us to conclude.
\end{proof}

\subsection{From the patched directed landscape to the periodic directed landscape}

In the following, we utilize the wrapped patched directed landscapes $(\Lpa)_{n \in \N}$ in order to define the periodic directed landscape. To simplify notation, let
\begin{equation*}
    \USCper:=\USC(\SetT,\R \cup \{-\infty\})
\end{equation*}
be the space of upper semi-continuous functions from  $\SetT$ to $\R \cup \{-\infty\}$, equipped with the supremum norm. In particular, since $\SetT$ is compact, note that $\USCper$ is compact as well. 
We consider the space $\mathcal{M}(\USCper)$ of probability measures on $\USCper$, and equip $\mathcal{M}(\USCper)$ with the total variation distance and the Borel-$\sigma$-algebra with respect to the weak-$*$ topology. Our goal is to define the periodic directed landscape as a limit point in the space $\mathcal{M}(\USCper)$. A main challenge is that the wrapped patched directed landscapes $\Lpa$ will be $-\infty$ for all pairs of points of too large slope. Hence, we will require in the following a different topology which only puts a light weight on pairs of points with a large slopes while at the same time being weak enough so that there is still a limit point and strong enough so the limit point is unique. This will be achieved in \eqref{def:WassersteinTypeMetric} where we take a suitable geometric sum of semi-norms  in order to define our metric. 
In preparation, recall the sets $\mathcal{D}_{\delta}$ from \eqref{def:SlopeConti}, and let 
 \begin{equation}\label{def:Deltan}
     \delta_n=n^{-1/16}
 \end{equation} 
 for all $n\in \N$. Note that $(\mathcal{D}_{\delta_n})_{n\in \N}$ is an exhaustive sequence for $\R_{\uparrow}^{4}$, i.e.,
\begin{equation*}
\mathcal{D}_{\delta_n} \subseteq \mathcal{D}_{\delta_{n+1}} \text{ for all } n\in \N \ \text{ and } \  \bigcup_{n \in \N}\mathcal{D}_{\delta_{n}} = \R^{4}_{\uparrow} ,
\end{equation*} while Lemma~\ref{lem:PatchedDLsConsistently} ensures that we can couple wrapped patched directed landscapes of scales $n$ and $2n$ on the sets $\mathcal{D}_{\delta_n}$ with high probability.
Consider the space of periodic upper semi-continuous functions from $\SetT \cap \mathcal{D}_{\delta_n}$ for some $n\in \N$ to $\mathbb{R} \cup \{ -\infty\}$, i.e., 
\begin{equation*}
\USC^{(n)}:=\USCper\big(\SetT \cap \mathcal{D}_{\delta_n},\R \cup \{-\infty \} \big) .
\end{equation*} 
 Since $\SetT \cap \mathcal{D}_{\delta_n}$ is compact, note that $(\USC^{(n)},\lVert \cdot \rVert_{\infty})$ is a complete and separable metric space.
For any probability measure $\nu$ on $\USCper$ and any subset $\mathcal{D}\subset \SetT$ we define $\nu|_{\mathcal{D}}$ to be the push-forward of $\nu$ under the restriction map that takes functions $f\in \USCper$ to functions $f|_{\mathcal{D}}\in \USCper\big(\SetT \cap \mathcal{D},\R \cup \{-\infty \} \big)$ (i.e., restrict from $\SetT$ to $\mathcal{D}$).

For all probability measures $\nu,\bar{\nu}$ on $\USCper$ and for all $i\in \N$, let
\begin{equation*}
 \WD^{(i)}( \nu,\bar{\nu} ) := \TV{ \nu|_{\mathcal{D}_{\delta_i}} - \bar{\nu}|_{\mathcal{D}_{\delta_i}} }
\end{equation*} be the total variation distance between the measures $\nu,\bar{\nu}$, restricted to $\mathcal{D}_{\delta_i}$ as above, and set
\begin{equation}\label{def:WassersteinTypeMetric}
\Wast( \nu,\bar{\nu} ) := \sum_{j=1}^{\infty} 2^{-j} \WD^{(2^j)}( \nu,\bar{\nu}  ) . 
\end{equation} 
Note that for every coupling $\mathbf{P}^{(i)}$ between measures $\nu,\bar{\nu}$ on $\USCper$, we have
\begin{equation*}
 \WD^{(i)}( \nu,\bar{\nu} ) \leq \mathbf{P}^{(i)}\left( f|_{\mathcal{D}_{\delta_i}} \neq \bar{f}|_{\mathcal{D}_{\delta_i}} \text{ with } f \sim \nu, \bar{f} \sim \bar{\nu} \right)  . 
\end{equation*}
\begin{lemma}\label{lem:WassersteinTypeMetric} The space $(\mathcal{M}(\USCper),\Wast)$ with the distance function $\Wast$ defined in \eqref{def:WassersteinTypeMetric} is a complete metric space.
\end{lemma}
\begin{proof}
 Since $(\USC^{(i)},\WD^{(i)})$ are metric spaces, to verify that $\Wast$ defines a metric on $\mathcal{M}(\USCper)$, it remains to show that for any $\nu \neq \bar{\nu}$, 
\begin{equation*}
\Wast(\nu,\bar{\nu}) > 0 . 
\end{equation*}
To this end, note that since $\nu \neq \bar{\nu}$, there exists some $\varepsilon>0$ such that for any coupling $\mathbf{P}$ of $\nu$ and $\bar{\nu}$ with $f \sim \nu$ and $\bar{f} \sim \bar{\nu}$, we get that
\begin{equation*}
\mathbf{P}\left( \lVert f -\bar{f}\rVert_{\infty} \geq  \varepsilon  \right) \geq  \varepsilon .
\end{equation*} 
Now by the continuity of the measure $\mathbf{P}$ and the fact that for all $\varepsilon>0$ and $f,\bar{f} \in \USCper$
\begin{equation*}
    \bigcup_{n \in \N} \Big\{ \big\lVert f|_{{\mathcal{D}_{n^{-1}}}}-\bar{f}|_{\mathcal{D}_{n^{-1}}}\big\rVert_{\infty}\geq \frac{\varepsilon}{2}\Big\} \supseteq \{ \lVert f-\bar{f}\rVert_{\infty} \geq \varepsilon\} , 
\end{equation*}
we see that there exists some $M=M(\varepsilon)\in \N_0$ such that
\begin{equation*}
\mathbf{P}\left( \big\lVert f|_{\mathcal{D}_{\delta_M}} -\bar{f}|_{\mathcal{D}_{\delta_M}} \big\rVert_{\infty} \geq  \frac{\varepsilon}{2}  \right) \geq \frac{\varepsilon}{2} .
\end{equation*} Since $\mathbf{P}$ was an arbitrary coupling, this implies that $\Wast(\nu,\bar{\nu})>0$, which gives the desired result that $\Wast$ is indeed a metric. 

It remains to show that the space $(\mathcal{M}(\USCper),\Wast)$ is complete. Since $(\USC^{(j)},\lVert \cdot \rVert_{\infty})$ is a complete and separable metric space, it follows that $(\mathcal{M}(\USC^{(j)}),\WD^{(j)})$ is likewise a complete metric space. This can be seen by appealing to Proposition~6 in~\cite{CR:Wasserstein} and using the fact that the total variation distance can be seen as the Wasserstein distance with respect to the trivial metric $\rho(x,y) = \mathds{1}_{x \neq y}$.
We seek to show an analogous statement for the space $(\mathcal{M}(\USCper),\Wast)$. 
Let $(\nu_i)_{i \in \N}$ be a Cauchy sequence in $(\mathcal{M}(\USCper),\Wast)$ and let $\nu^{(j)}_i := \nu_i|_{\mathcal{D}_{\delta_j}}$ denote the restriction of the measure to functions on $\mathcal{D}_{\delta_j}$ as defined earlier. It is immediate that for each $j\in \N$, $\{\nu^{(j)}_i\}_{i\in \N}$ forms a Cauchy sequence in  $(\mathcal{M}(\USC^{(j)}),\WD^{(j)})$. As $(\mathcal{M}(\USC^{(j),\per}),\WD^{(j)})$ is complete for  all $j\in \N$, we see that $(\nu^{(j)}_i)_{i \in \N}$ has a unique limit point $\nu^{(j)}_{\infty}$ in $(\mathcal{M}(\USC^{(j)}),\WD^{(j)})$. Since the measures $\nu^{(j)}_i$ satisfy 
 $\nu^{(j)}_i = \nu^{(k)}_i|_{\mathcal{D}_{\delta_j}}$ for all $j \in \lbr k-1\rbr$ and $i\in \N$, 
 we can extend the measures $(\nu^{(j)}_{\infty})_{j \in \N}$ consistently to the space $\mathcal{M}(\USCper)$, and define $\nu_{\infty} \in \mathcal{M}(\USCper)$ as the projective limit of $(\nu^{(j)}_{\infty})_{j \in \N}$. 
 We claim 
\begin{equation}\label{eq:WassersteinConvergence}
\lim_{i \rightarrow \infty}\Wast(\nu_i,\nu_{\infty})= 0 . 
\end{equation} To see this, note that for any $K \in \N$, we get that
\begin{equation*}
\lim_{i \rightarrow \infty}\Wast(\nu_i,\nu_{\infty}) \leq 2^{-K}+ \lim_{i\rightarrow \infty} \sum_{j=1}^{K} \WD^{(j)}(\nu_i,\nu_{\infty}) = 2^{-K} +\sum_{j=1}^{K} \lim_{i\rightarrow \infty}\WD^{(j)}\big(\nu^{(j)}_i,\nu^{(j)}_{\infty}\big) = 2^{-K} , 
\end{equation*}
where we use the definition of $\Wast$ for the first step, the definition of $\nu_{\infty}$ as a projective limit for the second step, and the fact that $(\nu_i^{(j)})_{i\in \N}$ is a Cauchy sequence for all fixed $j\in \N$ in the last step. Since $K\in \N$ was arbitrary, we get \eqref{eq:WassersteinConvergence}, allowing us to conclude. 
\end{proof} 
 
Using that the space of Borel probability measures on $\USCper$ with $\Wast$ is a complete metric space, we have the following consequence for the patched directed landscapes. For all $n\in \N$, let $\nu_n$ denote the law of $\mathcal{L}_{n}$.

 \begin{proposition}\label{pro:CauchySequence} The laws $(\nu_{2^n})_{n\in \N}$ are a Cauchy sequence in $(\mathcal{M}(\USCper),\Wast)$. In particular, this sequence has a unique limit point in the space of upper semi-continuous periodic functions $\USCper$, endowed with the topology of uniform convergence on compact sets.  A random element $\Lp$ of $\USCper$ is called a \textbf{periodic directed landscape} if it has the law of the (unique) limit point of the sequence $(\nu_{2^{n}})_{n\in \N}$ on $(\mathcal{M}(\USCper),\Wast)$.
 \end{proposition} 

\begin{proof}
This is achieved by noting that if the landscapes agree on all slopes up to a given value, then we can couple the different wrapped patched landscapes dyadically until we achieve the desired precision. More precisely, let $\ell \geq n$ for some $\ell,n\in \N$.  Then we see that 
\begin{equation}\label{eq:CauchyProperty2}
\begin{split}
\Wast\left( \nu_{2^n} , \nu_{2^\ell}\right) &\leq \sum_{i=n}^{\ell-1} \Wast\left( \nu_{2^i} , \nu_{2^{i+1}}\right) \leq \sum_{i=n}^{\ell-1} \Big( 2^{-i+1} + \sum_{j=1}^{i+1} 2^{-j}\WD^{(j)}(\nu_{2^i} , \nu_{2^{i+1}} )  \Big) \\
&\leq 2^{-n+2} + \sum_{i=n}^{\ell-1}\sum_{j=1}^{i+1} 2^{-j} C \exp(-ci^{-1/16})\leq C_0 \exp(-c_0n^{-1/16})
\end{split}
\end{equation} for some $c_0,C_0,c,C>0$, using the triangle inequality for $\Wast$ for the first step, the definition of $\Wast$ for the second step, and Lemma~\ref{lem:PatchedDLsConsistently} for the third step (with $c,C>0$ therein). 
In particular, $(\mathcal{L}_{2^n})_{n\in \N}$ is a Cauchy sequence. Since the space $(\mathcal{M}(\USCper),\Wast)$ is complete by Lemma~\ref{lem:WassersteinTypeMetric}, the sequence of patched directed landscapes has a unique limit point in $(\mathcal{M}(\USCper),\Wast)$. 
\end{proof}


As with the full-space directed landscape in \eqref{def:Geodesic},  we define geodesics in the periodic directed landscapes as the paths $\pi:[s,t]\to \mathbb{T}$ for some $0\leq s<t\leq 1$ which satisfy \eqref{def:Geodesic} with $\mathcal{L}$ replaced by $\Lp$.
Note that the periodic directed landscape is by construction defined on the space of periodic upper semi-continuous functions $\USCper$. This can be refined as follows.
\begin{proposition}\label{pro:PDLContinuity}
The periodic directed landscape $\Lp$ is almost surely a continuous function.
\end{proposition}
The proof of Proposition~\ref{pro:PDLContinuity} is deferred to the end of Section~\ref{sec:ConvergenceLPP}, where we use the convergence of periodic last passage percolation to the periodic directed landscape in order to give a short argument for the continuity. 
Next, we have the following characterization of the periodic directed landscape involving the full-space directed landscape. This is similar to the description of the periodic directed landscape given in the introduction as Theorem~\ref{thm:PeriodicDirectedLandscapeMAIN}.

\begin{theorem}\label{thm:PeriodicDirectedLandscape} 
The periodic directed landscape
$\Lp \colon \{ (x,s;y,t)\in (\mathbb{T}\times [0,1])^2_{\uparrow}\}\rightarrow \R$ satisfies:
\begin{enumerate}
\item[\namedlabel{eq:Cond1}{(i)}] (Locally full-space directed landscape) Let $\mathcal{L}$ denote a full-space directed landscape. For any $r\in (0,1)$ there  exist $c,C>0$, such that for any $\varepsilon\in (0,1]$ and any rectangle $\mathcal{R}_{\varepsilon}=[a_{\varepsilon},a_{\varepsilon}+r]\times [b_{\varepsilon},b_{\varepsilon}+\varepsilon]$ with $a_{\varepsilon}\in \mathbb{T}$ and $b_{\varepsilon} \in [0,1-\varepsilon]$, we have that
\begin{equation}\label{eq:LocalApproximationDL}
\TV{ (\Lp(x,s;y,t))_{(x,s;y,t) \in (\mathcal{R}_{\varepsilon})^2_{\uparrow}} - (\mathcal{L}(x,s;y,t))_{(x,s;y,t) \in (\mathcal{R}_{\varepsilon})^2_{\uparrow} }} \leq C \exp( - c\varepsilon^{-1/2} ).
\end{equation} 
\item[\namedlabel{eq:Cond3}{(ii)}] (Independent increments) For all $(t_i)_{i\in \lbr \ell \rbr}$ and $(s_i)_{i \in \lbr \ell \rbr}$ with $\ell \in \N$ such that $\big((t_i,s_i)\big)_{i \in \lbr \ell \rbr}$ are disjoint subintervals of $[0,1]$, we have that the random functions $(t_i,s_i) \mapsto \Lp(\cdot,t_i;\cdot,s_i)$ are independent.  
\item[\namedlabel{eq:Cond4}{(iii)}] (Metric composition) Almost surely, for any $0 \leq s< r < t \leq 1$ and $x,y \in \mathbb{T}$
\begin{equation}\label{eq:InvTriangleSupPDL}
\Lp(x,s;y,t) = \sup_{z \in \mathbb{T}} \Big(\Lp(x,s;z,r) + \Lp(z,r;y,t) \Big).
\end{equation}
\end{enumerate}
Moreover, if we replace \hyperref[eq:Cond1]{(i)} by (the weaker assumption)
\begin{enumerate}
\item[\namedlabel{eq:Cond1prime}{(i')}] For some $r'\in (0,1)$ and some $\alpha>2$, there exists $C'>0$ such that for all $\varepsilon\in (0,1]$ and any rectangle $\mathcal{R}'_{\varepsilon}=[a'_{\varepsilon},a'_{\varepsilon}+r']\times [b'_{\varepsilon},b'_{\varepsilon}+\varepsilon]$ with $a'_{\varepsilon}\in \mathbb{T}$ and $b'_{\varepsilon} \in [0,1-\varepsilon]$, we have that
\begin{equation}\label{eq:LocalApproximationDLPrime}
\TV{ (\Lp(x,s;y,t))_{(x,s;y,t) \in (\mathcal{R}'_{\varepsilon})^2_{\uparrow}} - (\mathcal{L}(x,s;y,t))_{(x,s;y,t) \in (\mathcal{R}'_{\varepsilon})^2_{\uparrow} }} \leq C'\varepsilon^{\alpha} .
\end{equation} 
\end{enumerate}
then any random field $\Lp$ supported on the space of continuous functions that satisfies properties \hyperref[eq:Cond1prime]{(i')}, \hyperref[eq:Cond3]{(ii)}, and \hyperref[eq:Cond4]{(iii)} must be a periodic directed landscape.
\end{theorem}

\begin{proof}
Write in the following $\nu^{\per}$ for the law of $\Lp$, and $\nu^{\per}_A$ for its projection onto a set $A \subseteq \SetT$. Similarly, write $\nu^{\full}$ for the law of a full-space directed landscape $\mathcal{L}$, and $\nu^{\full}_B$ for its projection onto $B \subseteq \SetR$. 
We start by showing that $\Lp$ satisfies property \hyperref[eq:Cond1]{(i)}. We claim that it suffices to show the existence of $C_1,c_1>0$ such that  
\begin{equation*}
g(\delta,\varepsilon):=\TV{ \nu^{\per}_{\mathcal{D}_{\delta}\cap (R_{\varepsilon})^2_{\uparrow}} - \nu^{\full}_{\mathcal{D}_{\delta}\cap (R_{\varepsilon})^2_{\uparrow}} } 
\end{equation*} 
satisfies for any $\varepsilon>0$ and all $\delta>0$ sufficiently small
\begin{equation}\label{eq:SmallCoupling}
    g(\delta,\varepsilon)  \leq C_1 \exp\big( - c_1\varepsilon^{-\frac{1}{2}} \big) . 
\end{equation}
Since $\Lp$ and $\mathcal{L}$ are random variables on two Polish spaces, note that by Ulam's tightness theorem, the space of couplings between $\Lp$ and $\mathcal{L}$ is tight. 
By the coupling representation of the total variation distance, we find a family of probability measures $\mathbf{P}_N^{\varepsilon}$ such that
 \begin{equation*}
     g(N^{-1},\varepsilon) = 1-\mathbf{P}_N^{\varepsilon}\left(  \nu^{\per}_{\mathcal{D}_{N^{-1}}\cap (R_{\varepsilon})^{2}_{\uparrow}}=\nu^{\full}_{\mathcal{D}_{N^{-1}}\cap (R_{\varepsilon})^{2}_{\uparrow}}\right)
 \end{equation*} for all $N\in \N$. Then the sequence of couplings $(\mathbf{P}_N^{\varepsilon})_{N \in \N}$ is tight, and any limit point $\mathbf{P}_{\varepsilon}$ of $(\mathbf{P}_N^{\varepsilon})_{N \in \N}$ is a coupling with the property that
 \begin{equation}
     \mathbf{P}_{\varepsilon}\Big( \nu^{\per}_{  (R_{\varepsilon})^{2}_{\uparrow}}=\nu^{\full}_{  (R_{\varepsilon})^{2}_{\uparrow}} \Big) \geq 1-C_1 \exp\big( - c_1\varepsilon^{-\frac{1}{2}} \big)
 \end{equation} using that the landscapes $\Lp$ and $\mathcal{L}$ must match on all of $(R_{\varepsilon})^{2}_{\uparrow} = \bigcup_{N=1}^{\infty} D_{N^{-1}} \cap (R_{\varepsilon})^{2}_{\uparrow}$.  

Now we prove that \eqref{eq:SmallCoupling} holds. To this end, we write $\nu_n$ for the law of a wrapped patched directed landscape $\Lpa$ and $\bar{\nu}_n$ for the law of an unwrapped patched directed landscape $\Lpat$, respectively. From  Lemma~\ref{lem:LocalFluctuationsPatched} in order to bound the probability that all geodesics for pairs of points in $(R_{\varepsilon})^{2}_{\uparrow} \cap \mathcal{D}_{\delta}$ can be decomposed into sub-paths with endpoints in $ \mathcal{V}_{\delta^{-1}}$, and Lemma~\ref{lem:Noncrossing} for a coupling between $\Lpa$ and $\Lpat$ on these sub-paths, we find $c_2,C_2>0$ such that
\begin{equation}
    \TV{ \nu_n|_{\mathcal{D}_{\delta}\cap (R_{\varepsilon})^{2}_{\uparrow}} - \bar{\nu}_n|_{\mathcal{D}_{\delta}\cap (R_{\varepsilon})^{2}_{\uparrow}}} \leq C_2 \exp(-c_2 \varepsilon^{-1/2})
\end{equation} for all $\varepsilon>0$, and all $\delta=\delta(\varepsilon)>0$ sufficiently small. 
Together with Lemma~\ref{lem:CoupledDLs} in order to couple $\Lpat$ to a   full-space directed landscape on $R_{\varepsilon}$, we find  $c_3,C_3>0$ such that 
\begin{equation*}
\begin{split}
\TV{ \nu_{n}|_{\mathcal{D}_{\delta}\cap (R_{\varepsilon})^{2}_{\uparrow}} - \nu^{\full}_{\mathcal{D}_{\delta}\cap (R_{\varepsilon})^{2}_{\uparrow}}} &\leq \TV{ \bar{\nu}_{n}|_{\mathcal{D}_{\delta}\cap (R_{\varepsilon})^{2}_{\uparrow}} - \nu^{\full}_{\mathcal{D}_{\delta}\cap (R_{\varepsilon})^{2}_{\uparrow}}} + C_3 \exp\big( - c_3\varepsilon^{-\frac{1}{2}} \big)  \\
&\leq 2 C_3 \exp\big( - c_3\varepsilon^{-\frac{1}{2}} \big) 
\end{split}
\end{equation*} 
for any $\varepsilon,\delta>0$, and all $n=n(\delta,\varepsilon)$ large enough. Recalling $\Lp$ is the limit of patched directed landscapes, 
this implies \eqref{eq:SmallCoupling}, and hence property~\hyperref[eq:Cond3]{(i)}, by choosing $n$ large enough.

To show the second property~\hyperref[eq:Cond3]{(ii)}, it suffices by the almost sure continuity of the periodic directed landscape (Proposition \ref{pro:PDLContinuity}) to address the case when the closures of the intervals $(t_i, s_i)$ are mutually disjoint. This follows by construction of the underlying patched directed landscapes for all choices of rational dyadic numbers 
\begin{equation}\label{def:Dyadics}
\mathcal{S}_{\textup{dyadic}} := \left\{ z \in \Q \colon z2^k \in \Z \text{ for some } k \in \N \right\} ,
\end{equation} i.e.,
$s,t,s_i,t_i \in \mathcal{S}_{\textup{dyadic}}$ for all $i \in \lbr \ell \rbr$, respectively.
By Kolmogorov's extension theorem and the fact that $\Lp$ is continuous by Proposition~\ref{pro:PDLContinuity}, $\Lp$ is the unique extension from the rational dyadics to $\SetT$. This yields the independence in \hyperref[eq:Cond3]{(ii)} for general intervals $(t_i,s_i)$ with mutually disjoint closures (and thus for all mutually disjoint open intervals, as mentioned above).

The third property \hyperref[eq:Cond4]{(iii)} follows for all rational dyadic numbers from the existence of leftmost and rightmost geodesics in Lemma \ref{lem:Geodesics} for the wrapped patched directed landscapes in the construction of $\Lp$, and Lemma \ref{lem:PatchedDLsConsistently} to couple the underlying patched directed landscapes $\Lpa$ for different values of $n$, together with a Borel--Cantelli argument. 
For the metric composition~\hyperref[eq:Cond4]{(iii)} with general values of $(x,s;y,t) \in \SetT$, we follow similar arguments as in Lemma~10.8 of \cite{DOV:DirectedLandscape}. Consider families of sites $(x_n,s_n)_{n \in \N}$ and $(y_n,t_n)_{n \in \N}$ with 
\begin{equation*}
(x_n,s_n), (y_n,t_n) \in \mathcal{S}_{\textup{dyadic}} \times \mathcal{S}_{\textup{dyadic}} 
\end{equation*} for all $n\in \N$, and
\begin{equation*}
(x_n,s_n;y_n,t_n) \rightarrow (x,s;y,t)
\end{equation*} as $n\rightarrow \infty$. 
Moreover, let $(r_n)_{n \in \N}$ with $s_n \leq r_n \leq t_n$ be  such that
\begin{equation*}
    \lim_{n \rightarrow \infty} r_n = r \in [s,t] . 
\end{equation*}
Then by Lemma~\ref{lem:Geodesics} on the existence of a rightmost geodesics in the patched directed landscape, together with Lemma~\ref{lem:PatchedDLsConsistently} in order to couple the patched and the periodic directed landscape, there exists some $(z_n)_{n \in \N}$ in the rational dyadics such that
\begin{equation}\label{eq:PreLimConti}
    \Lp(z_n,r_n) = \sup_{z^{\prime}_n \in \mathbb{T}} \big(\Lp(x_n,s_n;z^{\prime}_n,r_n) + \Lp(z^{\prime}_n,r_n;y_n,t_n) \big)
\end{equation}  for all $n\in \N$. 
Now for any subsequential limit $z\in \mathbb{T}$ of $(z_n)_{n \in \N}$, 
the continuity of the periodic directed landscape $\Lp$ from Proposition~\ref{pro:PDLContinuity} guarantees that
\begin{equation*}
    \Lp(z_n,r_n) \rightarrow \Lp(z,r)  
\end{equation*} 
as $n \rightarrow \infty$. Hence, together with \eqref{eq:PreLimConti}, we get that 
\begin{equation*}
    \Lp(z,\tau)  \geq  \sup_{z^{\prime} \in \mathbb{T}}\big(\Lp(x,s;z^{\prime},r) + \Lp(z^{\prime},r;y,t) \big).  
\end{equation*} The reverse inequality in order to get the metric composition \hyperref[eq:Cond4]{(iii)} is immediate by the continuity of $\Lp$, and the metric composition over the dyadic rational numbers.

In order to show the uniqueness of the law of $\Lp$ under the weaker conditions \hyperref[eq:Cond1prime]{(i')}, \hyperref[eq:Cond3]{(ii)}, and \hyperref[eq:Cond4]{(iii)}, suppose that $\P(\Lp \in \, \cdot \, )$ and $\P(\tLp \in \, \cdot \, )$ denote two different measures on $\USCper(\SetT,\R)$ satisfying conditions \hyperref[eq:Cond1prime]{(i')}, \hyperref[eq:Cond3]{(ii)}, and \hyperref[eq:Cond4]{(iii)}. 
Recall the sets $\mathcal{R}^{i}_{n}$ from \eqref{eq:PartitionedRectangles} and sets $\mathcal{V}_n$ from \eqref{def:SlopeVerti} for all $n\in \N$. 
We set in the following $n=\delta^{-1}=r^{-1}$. Note that we can exhaust $\mathcal{R}^{1}_{n}$ by $\lceil r^{-1} n\rceil$ many $r \times n^{-1}$ rectangles. We claim that for all $n\in \N$, there exists a coupling $\mathbf{P}_n$ such that for some $c_4,C_4>0$, 
\begin{equation}\label{eq:FullPeriodicCoupling}
\begin{split}
    \mathbf{P}_n\Big(  \Lp(x,s;y,t)= \tLp(x,s;y,t)&=\Lpa(x,s;y,t) \ \forall (x,s;y,t) \in (\mathcal{R}^{1}_{n})^2_{\uparrow} \cap \mathcal{V}_{n} \cap \mathcal{D}_{16n^{-2}} \Big) \\
    &\geq 1- \lceil r^{-1} n  \rceil n^{-\alpha} + C_4 \exp(-c_4 n^{1/2})
\end{split}
\end{equation}
Indeed, this follows from property \hyperref[eq:Cond1prime]{(i')}, noting that $r \leq r'$ for all $n$ large enough, to couple $ \Lp$ and $\tLp$ with full-space directed landscapes, together with Lemma~\ref{lem:CoupledDLs} to couple the respective full-space directed landscapes with an unwrapped patched directed landscape $\Lpat$, and Lemma~\ref{lem:Noncrossing} to show that the coupling also holds with respect to a wrapped patched directed landscape.
Similarly, we get that for any $i \in \lbr n \rbr$, and $n$ large enough
\begin{equation}\label{eq:SmallRingAgreement}
\mathbf{P}_n\big( \Lp(x,s;y,t) = \tLp(x,s;y,t)\ \forall (x,s;y,t) \in (\mathcal{R}^{i}_{n})^2_{\uparrow} \cap \mathcal{V}_{n} \cap \mathcal{D}_{16n^{-2}} \big) \geq 1 - 2\lceil r^{-1} n  \rceil n^{-\alpha} .
\end{equation} 
Write $f^{+}_{x,s;y,t}(z)$ for the right-most maximizer in the metric composition~\hyperref[eq:Cond4]{(iii)} at location $z\in [x,y]$ for points $(x,s;y,t)$ in $\Lp$, and similarly $\tilde{f}^{+}_{x,s;y,t}(z)$ for the right-most maximizer in the metric composition~\hyperref[eq:Cond4]{(iii)} at location $z\in [x,y]$ for points $(x,s;y,t)$ in $\tLp$.
By iterating Lemma \ref{lem:PatchedDLsConsistently}, we see that Lemma \ref{lem:LocalFluctuationsPatched} continues to holds with respect to the periodic directed landscapes $\Lp$ and $\tLp$, i.e., the events
\begin{equation*}
\begin{split}
    \mathsf{C}_{\delta} &:= \left\{ \sup_{i\in \Z \colon i\delta,(i+1)\delta \in [x,y]} \sup_{z\in [0,\delta]}\big|f^{+}_{x,s;y,t}(i\delta+z)-f^{+}_{x,s;y,t}(i\delta)\big|\leq \frac{1}{32} \forall (x,s;y,t) \in \SetT \cap \mathcal{D}_{\delta^{1/8}}   \right\} \\
    \tilde{\mathsf{C}}_{\delta} &:= \left\{ \sup_{i\in \Z \colon i\delta,(i+1)\delta \in [x,y]}\sup_{z\in [0,\delta]}\big|\tilde{f}^{+}_{x,s;y,t}(i\delta+z)-\tilde{f}^{+}_{x,s;y,t}(i\delta)\big|\leq \frac{1}{32} \forall (x,s;y,t) \in \SetT \cap \mathcal{D}_{\delta^{1/8}}   \right\}
\end{split}
\end{equation*} 
hold with probability at least $1-C_5\exp(-c_5 \delta^{-1/32})$ for some $c_5,C_5>0$ and all $\delta>0$ small enough. 
Together with \eqref{eq:SmallRingAgreement}, recalling that $n=\delta^{-1}=r^{-1}$, and the independence \hyperref[eq:Cond3]{(ii)} for $\Lp$ and $\tLp$, there exists a coupling $\tilde{\mathbf{P}}$ such that for some $C_6>0$, and all $\delta>0$ small enough, 
\begin{equation*}
\tilde{\mathbf{P}}\Big( \Lp(x,s;y,t) = \tLp(x,s;y,t) \ \forall (x,s;y,t) \in \SetT \cap \mathcal{D}_{\delta^{1/8}}\Big) \geq 1 - C_6\delta^{\alpha-2} .
\end{equation*} As $\alpha>2$ by our assumptions, and $\delta>0$ was arbitrary, we conclude. 
\end{proof}

\begin{remark}\label{rem:ExtendedPDL}
Note that using the independence, and Kolmogorov's extension theorem,  we can extend (via repeated use of metric composition) the periodic directed landscape to arbitrary time horizons on the space  $(\mathbb{T} \times \R)^2_{\uparrow}=\left\{ (x,s;y,t) \in (\mathbb{T} \times \R)^2 \colon s<t \right\}$. 
\end{remark}
For the sake of completeness, we can define the periodic directed landscape also for general periodic length $p$ using the $1:2:3$ rescaling of the underlying directed landscapes in the definition of the wrapped patched directed landscape $(\Lpa)_{n\in \N}$. 
\begin{definition}\label{def:PeriodicExtension}
 For $p>0$, let $\mathbb{T}^{(p)}$ denote the $1$-dimensional torus of length $p$, and recall that
 \begin{equation*}
     (\mathbb{T}^{(p)} \times \R)^{2}_{\uparrow} = \left\{ (x,s;y,t) \in \big(\mathbb{T}^{(p)} \times \R\big)^2 \colon s<t \right\} . 
 \end{equation*}
Then we define the $\boldsymbol{p}$\textbf{-periodic directed landscape} as
\begin{equation*}
   \Lpp(x,s;y,t) := p^{-1/2}\Lp_1(px,p^{3/2}s;py,tp^{3/2})
\end{equation*} for all $(x,s;y,t) \in (\mathbb{T}^{(p)} \times \R)^{2}_{\uparrow}$, where $\Lp_1$ is a $1$-periodic directed landscape on  $(\mathbb{T} \times \R)^2_{\uparrow}$.
\end{definition}

This section has given us the existence and uniqueness of the periodic directed landscape. Our focus will now turn to establishing convergence of periodic models to this limit. We will work with two models for which the corresponding convergence is known in full-space.

\section{Convergence of periodic exponential last passage percolation}\label{sec:ConvergenceLPP}

In the following, we collect basic properties of exponential last passage percolation on the two-dimensional integer lattice in periodic as well as  independently and identically distributed (i.i.d.)~environments.  Similarly to the continuum, we approximate  periodic last passage percolation using only restricted lattice paths.
 More precisely, we will show that last passage percolation in periodic environments can be well-approximated by coupled last passage percolation on subsets of i.i.d.~last passage percolation.

\subsection{Basic definitions in last passage percolation}

For $v=(v_1,v_2) \in \Z^2$, we set $|v|:= v_1+v_2$. We let $\succeq$ denote the component-wise ordering on $\Z^2$, i.e., 
\begin{equation}
v \succeq u \quad \Leftrightarrow \quad v_1\geq u_1 \text{ and } v_2\geq u_2 . 
\end{equation}
 Let $\Pi_{u,v}$ be the set of all lattice paths $\pi$ connecting $u$ to $v$, i.e., the set of all up-right paths on $\Z^2$ given by $\pi = (z_i)_{0,\dots,|v-u|}$ such that
\begin{equation*}
 z_0=u, \ z_{|v-u|}=v \text{ and } z_{i}-z_{i-1} \in \{ \eone,\etwo\} \text{ for all } i \in \lbr |u-v| \rbr , 
\end{equation*} where $\eone:=(1,0)$ and $\etwo := (0,1)$, and we recall that $\lbr n\rbr:= \{ 1,\dots, n\}$ for all $n\in \N$. For every $u\in \Z^2$, we assign a random variable $\omega_u$. Then for a given collection of random variables $(\omega_z)_{z \in \Z^2}$ and sites $u,v \in \Z^2$ with $v \succeq u$, we define the \textbf{(exponential) last passage time}
\begin{equation}\label{def:LPPdef}
T_{u,v}:=  \sup_{\pi \in \Pi_{u,v}} \sum_{z \in \pi \setminus \{v\}} \omega_z . 
\end{equation} between $u$ and $v$, with the convention that $T_{u,v}=-\infty$ whenever $v \nsucceq u$. 
Note that in contrast to the standard definition, we will exclude the endpoint $v$ in \eqref{def:LPPdef}, as this allows us to concatenate last passage times, i.e., we have that almost surely
\begin{equation}
\sup_{w \in \Z^2\colon v \succeq w \succeq u} T_{u,w} + T_{w,v} = T_{u,v} . 
\end{equation} With a slight abuse of notation, we set
\begin{equation}
T_{u,v} = T_{(\lfloor u_1 \rfloor,\lfloor u_2 \rfloor),(\lfloor v_1 \rfloor,\lfloor v_2 \rfloor)}
\end{equation} for all $u=(u_1,u_2)\in \R^2$ and $v=(v_1,v_2) \in \R^2$, and we write
\begin{equation}
    T_{A,B} := \sup_{u \in A} \sup_{v \in B} T_{u,v}
\end{equation} for all $A,B \subseteq \Z^2$.
A lattice path $\gamma_{u,v} \in \Pi_{u,v}$ maximizing the right-hand in side in \eqref{def:LPPdef} is called a \textbf{geodesic}. Note that for every geodesic $\gamma_{u,v}$ and $w,w^{\prime} \in \gamma_{u,v}$ with $v \succeq w^{\prime} \succeq w \succeq u$, the path $\gamma_{u,v}$ restricted between sites $w$ and $w^{\prime}$ is again a geodesic.

We will be interested in two particular families of random variables $(\omega_v)_{v \in \Z^2}$, and their corresponding last passage times. First, when  $(\omega_v)_{v \in \Z^2}$ is a family of independent Exponential-$1$-distributed random variables, we refer to $(\omega_v)_{v \in \Z^2}$ as an \textbf{i.i.d.}~or \textbf{full-space  environment}. Second, let $N$ and $k=k(N) \in \lbr N-1 \rbr$ be fixed. When 
 \begin{equation}\label{def:PeriodicEnvironment}
\omega_{(v_1,v_2)}=\omega_{(v_1+N-k,v_2-k)}
\end{equation} for all $(v_1,v_2)\in \Z^2$, and the random variables are independent, otherwise, we say that  $(\omega_v)_{v \in \Z^2}$ is an $\mathbf{(N,k)}$\textbf{-periodic environment}. We simply refer to a \textbf{periodic environment} when $k$ and $N$ are clear from the context. Intuitively, while we allow $k$ to be essentially arbitrary, up to certain boundedness condition such as \eqref{eq:NonDegeneratedEnvironment} to prove the most general result, the proofs are all very similar if $k = N/2$. We will in the following denote the law of last passage times in an $(N,k)$-periodic environment by $\Pper^{N,k}$, and by $\Pfull$ for an i.i.d.~environment. Again, whenever the parameters $N,k$ are clear from the context, we simply write $\Pper$.
Let us note that for both choices of environments (full-space and periodic), a  geodesic $\gamma_{u,v}$ is almost surely unique for all $u,v \in \Z^2$ with $v \succeq u$. In the case of periodic environments, we in addition define a certain set-to-point last passage time.

\begin{definition}\label{def:LastPassageTimes}
Let $(\omega_v)_{v \in \Z^2}$ be an $(N,k)$-periodic environment for some $N \in \N$ and with $k=k(N)\in \lbr N -1\rbr$. 
We define the \textbf{periodic last passage time} $\Tper_{u,v}$ between $u$ and $v$ as
\begin{equation}\label{def:LastPassageTimeIID}
\Tper_{u,v}:=  \sup_{i \in \Z} T_{u+i(N-k,-k),v} . 
\end{equation}
\end{definition}
We refer to $\gper_{u,v}$ as the corresponding \textbf{periodic geodesic} between $u$ and $v$, and note that  $\gper_{u,v}$ is unique almost surely. 
Intuitively, the periodic last passage time $\Tper_{u,v}$ agrees with the definition of a last passage times defined on a cylindrical environment, where we identify every vertex $u\in \Z^2$ with its \textbf{periodic translates} (see Figure \ref{fig:intro}(d))
\begin{equation}\label{def:PeriodicTranslates}
\TR(u) := \left\{  w \in \Z^2 \, \colon  w= u+(i(N-k),-ik) \text{ for some } i \in \Z \right\} . 
\end{equation}
In order to study (periodic) last passage times, it will be useful to consider certain level sets in the environment. To this end, we let for all $x \in \Z$
\begin{equation}
\Lambda^{(x)} := \big\{ (z_i)_{i \in \Z} \, \colon z_0 =(x,x) \wedge  z_{i}-z_{i-1} \in \{ \eone,-\etwo\} \, \forall i \in \Z  \big\} 
\end{equation} denote the set of all down-right lattice paths through $\mathbf{x}:=(x,x) \in \Z^2$. Similarly, we define 
\begin{equation}
\Lambda^{(x)}_{N,k} := \{ \lambda \in \Lambda^{(x)} \, \colon \lambda  = \lambda  + (N-k,-k) \}
\end{equation} as the set of $(N,k)$-periodic down-right lattice paths through $x$. We simply write
\begin{equation}\label{def:LambdaNk}
\Lambda= \bigcup_{x \in \Z}\Lambda^{(x)} \ \text{ and } \ \Lambda_{N,k}= \bigcup_{x \in \Z}\Lambda^{(x)}_{N,k} . 
\end{equation}
For a given down-right path $\lambda \in \Lambda$, we define the corresponding last passage times with respect to $\lambda$ as
\begin{equation}
T_{\lambda ,v} := \sup_{u \in \lambda}\sup_{\pi \in \Pi_{u,v}} \sum_{w \in \pi \setminus \{ v\} } \omega_{w} , 
\end{equation} and similarly by $\Tper_{\lambda ,v}$ when $\lambda  \in \Lambda_{N,k}$. Intuitively, the last passage times $T_{\lambda,v}$ corresponds to setting the lengths on all sites $u \notin \lambda$ with $u \preceq v$ for some $v\in \lambda$ to $-\infty$, while taking the supremum over all sites $\Z^2$ as possible starting positions. Observe that in the special case where the initial path is given by the \textbf{step initial data}
\begin{equation}
\lambda^{(x,y)}_{\textup{step}} := \{ (z,y) \colon z \geq x \} \cup  \{ (x,w) \colon w \geq y \} 
\end{equation} for some $x\in \Z$, and the environment is i.i.d., we get that almost surely 
\begin{equation}
T_{\lambda ^{(x,y)}_{\textup{step}},v} = T_{(x,y),v}
\end{equation} for all $v \succeq (x,y)$, and we have that
\begin{equation}\label{eq:VariationalDiscrete}
 T_{\lambda,v} = \sup_{(x,y)\in \lambda}  T_{\lambda ^{(x,y)}_{\textup{step}},v} . 
\end{equation}
 Similarly, when the environment is $(N,k)$-periodic, for $x,y\in \Z$, we denote by $
\lambda^{(x,y)}_{\per,N}$ the \textbf{periodic step initial data} parametrization $\lambda^{(x,y)}_{\per,N}=(z_i)_{i \in \Z}$ such that $z_0=(x,y)$ and
\begin{equation}
\Tper_{\lambda^{(x,y)}_{\per,N},v} = \Tper_{(x,y);v}
\end{equation} almost surely for all $v \succeq (x,y)$, as well as
\begin{equation}\label{eq:VariationalDiscretePeriodic}
 \Tper_{\lambda_{N}^{\per},v} = \sup_{(x,y)\in \lambda_{N}^{\per}}  \Tper_{\lambda^{(x,y)}_{\per,N},v} . 
\end{equation}
\begin{definition}\label{def:HeightLPP}
For a down-right path $\lambda=(z_i)_{i\in \Z} \in \Lambda$ and a function $h_0 \colon \lambda \rightarrow \R \cup \{-\infty\}$ we define the \textbf{last passage percolation height function} $\hL$ with initial condition $h_0$ on $\lambda$ by
\begin{equation}\label{eq:LPPRecursion}
    \hL_{\lambda}(h_0;x,y) := \sup_{u \in \lambda} \big( h_{0}(u) + T_{u,(x,y)}\big) 
\end{equation} for all $x,y \in \Z$ with $(x,y) \succeq u $ for some $u\in \lambda$. . The \textbf{periodic height function} $\hLp_{\lambda}$ is defined as in \eqref{eq:LPPRecursion}, but with respect to the periodic last passage times $\Tper$ from Definition~\ref{def:LastPassageTimes} and a down-right path $\lambda=\lambda^{\per}_N \in \Lambda_{N,k}$. 
\end{definition}
While we will work mostly with the periodic last passage times and the LPP height function, let us also record the following notion of the TASEP height function, allowing to identify (periodic) last passage percolation and the (periodic) totally asymmetric simple exclusion process -- the $q=0$ special case of the (periodic) ASEP, which we will study in Section~\ref{sec:ASEPconvergence}.
\begin{definition}\label{def:HeightTASEP}
For the \textbf{TASEP height function} $(\hTASEP(\lambda;y,t))_{y\in \Z,t\geq 0}$, assume that $(0,0) \in \lambda$ for a fixed down-right path $\lambda \in \Lambda$. For all $t\geq 0$, consider the down-right path
\begin{equation}\label{eq:CornerExtremalBasic}
\lambda_t= (z^{t}_i)_{i \in \Z} := \{ w \in \Z^2 \, \colon \,  T_{\lambda,w} \leq t < T_{\lambda,w+(1,1)} \} 
\end{equation} such that $z^{t}_{0}=(x_t,x_t)$ for some $x_t \in \Z$. For $y \in \Z$ and $t > 0$, set $h_t(0)=x_t$ and recursively
\begin{equation}\label{eq:CornerExtremalBasic2}
h_t(i) - h_t(i-1) = \begin{cases} -1 & \text{ if } z^{t}_i-z^{t}_{i-1} = \eone , \\
 0 & \text{ if } z^{t}_i-z^{t}_{i-1} = -\etwo . 
\end{cases}
\end{equation}
 The \textbf{periodic TASEP height function} $\hpTASEP$ is defined as in \eqref{eq:CornerExtremalBasic} and \eqref{eq:CornerExtremalBasic2}, with respect to the periodic last passage times $\Tper$ from Definition~\ref{def:LastPassageTimes} and a down-right path $\lambda^{\per}_N \in \Lambda_{N,k}$.
\end{definition}

\begin{remark}\label{rem:TASEPandLPP}
The totally asymmetric simple exclusion process (TASEP) on $\Z$ can be described using exponential last passage percolation in i.i.d.~environment~\cite{R:PDEresult}. Similarly, the evolution of a TASEP with $k$ particles on a circle of length $N$ can be equivalently described using last passage percolation in an $(N,k)$-periodic environment, see for example~\cite{SS:TASEPcircle}. 
Let us mention that we will not discuss exclusion processes until Section~\ref{sec:KPZfixed} where we discuss the periodic KPZ fixed point. In Section~\ref{sec:ASEPconvergence}, we focus on the basic coupling between them instead of their correspondence to last passage percolation; see Definition~\ref{def:BasicCouplingASEP}.
\end{remark}

\subsection{Convergence to the directed landscape}

In order to capture the scaling limit of exponential last passage times in a given direction, we record the following result from \cite{DV:LongestSub} on the convergence of exponential last passage percolation on the full-space to the directed landscape. Fix some constant $\rho>0$ and define the vectors 
\begin{equation}\label{def:VectorsDL}
    \mathbf{u}:=(\rho^{1/2},-1) \quad \text{ and } \quad \mathbf{v}:=(\rho,1). 
\end{equation}
 We define the rescaled last passage times 
\begin{equation}\label{def:LPPmetric}
\dLeps(x,s; y,t) := T_{x \varepsilon^{-2}\mathbf{u} + s  \varepsilon^{-3}\mathbf{v},y \varepsilon^{-2}\mathbf{u}+t \varepsilon^{-3}\mathbf{v}} 
\end{equation}
for all $(x,s;y,t) \in \R^{4}_{\uparrow}$ such that  $(x\varepsilon^{-2}\mathbf{u} + s\varepsilon^{-3}\mathbf{v};y\varepsilon^{-2}\mathbf{u}+t\varepsilon^{-3}\mathbf{v}) \in \Z^2 \times \Z^2$. Let
\begin{equation}\label{def:MeanLPP}
M(\rho,\varepsilon,s,t) :=  (t-s) \varepsilon^{-3}(1+ \sqrt{\rho})^2
\end{equation} denote the first order of the mean of $\dLeps(x,s; y,t)$.  Moreover, set
\begin{equation}\label{def:VarianceLPP}
V= V(\rho,\varepsilon) := \frac{1}{\sqrt{2}}(1+\sqrt{\rho})^{\frac{3}{2}} \rho^{-1/4} \varepsilon^{-1} 
\end{equation} as a rescaling parameter for the last passage times in direction $\rho$. 
We define the \textbf{(full-space) LPP landscape} $\LLeps$ as an element of  $\mathcal{C}(\R^{4}_{\uparrow},\R)$ with
\begin{equation}\label{def:LPPSheet}
\LLeps(x,s;y,t) :=  \frac{1}{V(\rho,\varepsilon)} \big(\dLeps(x,s;y,t)- M(\rho,\varepsilon,s,t) \big) , 
\end{equation}
whenever $(x\varepsilon^{-2}\mathbf{u} + s\varepsilon^{-3}\mathbf{v};y\varepsilon^{-2}\mathbf{u}+t\varepsilon^{-3}\mathbf{v}) \in \Z^2 \times \Z^2$, and linear interpolation,  otherwise, and record the following convergence result. 
%
\begin{theorem}[{\cite[Theorem 1.7 and Remark 1.10]{DV:LongestSub}}]\label{thm:LPPtoDL}
On the full-space $\R^{4}_{\uparrow}$, we have convergence of the LPP landscape $\LLeps$  to the directed landscape $\mathcal{L}$  (Definition \ref{def:DirectedLandscape}), i.e., as $\varepsilon \rightarrow 0$ 
\begin{equation}
\big(\LLeps(x,s\tilde{\omega}^{-1};y,t\tilde{\omega}^{-1})\big)_{(x,s;y,t) \in \R^{4}_{\uparrow}} \rightarrow \big(\mathcal{L}(x,s;y, t)\big)_{(x,s;y,t) \in \R^{4}_{\uparrow}}
\end{equation} in distribution, where we have $\tilde{\omega}=  2^{3/2} (1 + \sqrt{\rho})^{-1/2} \rho^{1/4}$.
\end{theorem}
\begin{remark}
In the notation of \cite{DV:LongestSub}, we have that $(u_1,u_2)=(\rho^{1/2},1)$ and $(v_1,-v_2)=(\rho,-1)$, as well as
\begin{equation}
\alpha = (1+\sqrt{\rho})^2, \quad \beta=1+\frac{1}{\sqrt{\rho}}, \quad  \chi=(1+\sqrt{\rho})^{\frac{4}{3}}\rho^{-\frac{1}{6}}, \quad \tau=2(1+\sqrt{\rho})^{\frac{2}{3}}\rho^{\frac{2}{3}}
\end{equation} for the unscaled parameters and
\begin{equation}
\tilde{\alpha} = \alpha, \quad \tilde{\eta}=\frac{1}{2}(1+\sqrt{\rho})^{\frac{1}{3}}\rho^{-\frac{1}{6}}, \quad \tilde{\gamma}=0, \quad \tilde{\chi}=\frac{1}{\sqrt{2}}\rho^{-\frac{1}{4}}(1+\sqrt{\rho})^{\frac{3}{2}}, \quad \tilde{\omega}=2^{\frac{3}{2}}(1+\sqrt{\rho})^{-\frac{1}{2}}\rho^{\frac{1}{4}}
\end{equation} for the rescaled parameters. We refer to Remark~1.9 in \cite{DV:LongestSub} for an explanation of  scaling parameters. The scaling relations can also be checked for example from 
Theorem~1.7 in \cite{RRV:Beta} for the one point scaling of $\dLeps(0,0; y,t)$ to a Tracy--Widom $\textup{GUE}$ in order to obtain the above choice of parameters; see also Theorem~1.6 in \cite{J:KPZ}, where this result first appeared. In particular, note that $M$ from \eqref{def:MeanLPP} is of order $\varepsilon^{-3}$ while $V$ from \eqref{def:VarianceLPP} is of order $\varepsilon^{-1}$. 
\end{remark}

\subsection{Moderate deviations for last passage percolation}\label{sec:ModerateLPP}

We collect now some basic estimates on last passage percolation in i.i.d.~and periodic environments. Here, we largely use $\tilde{n}$ instead of $n$ in this section to avoid confusion with the parameter $n$ from Section \ref{sec:PeriodicDL}.
While moderate deviation result on last passage times in i.i.d.~environments are well-known, similar results for last passage time in periodic environments were only recently established in~\cite{SS:TASEPcircle}. We start with the following moderate deviation result on last passage times in i.i.d.~environments, which can be found in \cite{BHS:Binfinite}; see also \cite{LR:BetaEnsembles}.
 
\begin{theorem}[{\cite[formula $(2)$ and $(3)$]{BHS:Binfinite}}] \label{thm:ScalingLPPiid} Let $(\omega_v)_{v \in \Z^2}$ be an i.i.d.~environment. There exist constants $c,C,m_0>0$ such that for all $x>0$ and $v=(\np,m\np)$ with $m \in \big( \frac{m_0}{\np},1 \big]$, the last passage times satisfy
\begin{equation}\label{eq:StatementSteepUpper}
\Pfull \left( \big| T_{(0,0),v} -  (1+\sqrt{m})^2 \np \big| \geq  x  \np^{1/3}(m+m^{-1})^{-1/6}   \right) \leq C \exp\left(-c x \right) .
\end{equation} 
\end{theorem} 
Let us mention that the exponents in Theorem \ref{thm:ScalingLPPiid}  are not optimal (and can be improved for constant $m$), but suffice for our purposes. Next, consider the \textbf{transversal fluctuations} $\TF(\pi)$ of a lattice path $\pi$, i.e. for $\pi$ from $u$ to $v=u+(n,m n)$ with some $m \in (0,1)$ and $n\in \N$, we set
\begin{equation}\label{def:TransversalFluctuations}
\TF(\pi) := \max_{i,j \in \Z}\{ |j-m i| \colon u+(i,j) \in \pi \} .
\end{equation}
The following result states moderate deviations for the transversal fluctuations of geodesics in i.i.d.~environments.
\begin{lemma}[{\cite[Lemma 3.5]{SS:TASEPcircle}}] \label{lem:TransversalFluctuations} Recall that $\gamma_{(0,0),(\np,m\np)}$ denotes the geodesic from $(0,0)$ to $(\np,m\np)$ for some $m\in (0,1)$ and $\np\in \N$. There exist constants $m_0,C,c>0$  such that for all $m \in \big[\frac{m_0}{\tilde{n}},1 \big]$, $x>0$, and $n\in \N$
\begin{equation}\label{eq:TransversalIID}
\Pfull\big(\TF(\gamma_{(0,0),(\np,m \np)}) \geq x m^{2/3}\np^{2/3}\big) \leq C\exp(-cx) \, .
\end{equation} 
\end{lemma}
Again, the exponent in \eqref{eq:TransversalIID} is believed not to be optimal, but suffices for our arguments. While the above results are stated for i.i.d.~environments, similar estimates also hold for the transversal fluctuation in periodic last passage percolation. In the following, we will assume for $(N,k)$-periodic environments that
\begin{equation}\label{eq:NonDegeneratedEnvironment}
\mathfrak{a} \leq \liminf_{N \rightarrow \infty} \frac{k}{N} \leq \limsup_{N \rightarrow \infty} \frac{k}{N} \leq 1-\mathfrak{a} , 
\end{equation} with some fixed constant $\mathfrak{a}>0$. We have the following result on the moderate deviations. 

 \begin{lemma}\label{lem:GeodesicsFullPeriodicAgree} Consider a family of $(N,k)$-periodic environments such that \eqref{eq:NonDegeneratedEnvironment} holds. For all $m_0>0$,  there exist constants $\theta_0,c,C>0$, depending only on $\mathfrak{a}$ and $m_0$, such that for all sites $v=(\np,m\np)$ with some $\theta\geq \theta_0$, $m=m(\np) \in [m_0,m_0^{-1}]$, and $\np=\np(\theta) \in \N$, which satisfy
 \begin{equation}\label{eq:CorrectPeriodicScale}
 \np \leq \theta^{-1}N^{3/2} , 
 \end{equation} and all $x \in [1,\theta^{2/3}]$, we get that for all $N$ large enough
\begin{equation}
\Pper\Big( \TF\big(\gper_{(0,0),v}\big) \geq x \np^{2/3} \Big) \leq C \exp(-c x) .
\end{equation} 
\end{lemma} 

We defer the proof of Lemma \ref{lem:GeodesicsFullPeriodicAgree} to Appendix \ref{sec:AppendixLPP}. Let us note that a similar moderate deviation result was already shown as Lemma~4.4 of \cite{SS:TASEPcircle} in the special case where we set
\begin{equation}
\label{def:FixedSlope}    
m=m_{\ast}:=\frac{k^2}{(N-k)^{2}} . 
\end{equation} However, we require a more general set of slopes in the sequel as well as moderate deviations for the transversal fluctuations of (periodic) geodesics on smaller scales. To this end, we define for all lattice paths $\pi$ from $u$ to $u+(\np,m\np)$, and all $L \in \lbr \np \rbr$ the \textbf{local transversal fluctuations} \begin{equation}\label{def:TransversalFluctuationsLocal}
\TF_L(\pi) := \max\{ |j-m^{-1} L| \colon u+(j,L) \in \pi \} .
\end{equation}
The following lemma can be found as Theorem~3 in~\cite{BSS:Coalescence} for i.i.d.~environments, while we will also require a similar statement for  periodic environments. 

\begin{lemma}\label{lem:ModerateLocal}
For all $m_0>0$, there exist constants $\theta_0,c_1,C_1,L_0>0$, depending only on $\mathfrak{a}$ and $m_0$, such that for all $m=m(n) \in [m_0,m_0^{-1}]$, $\theta \geq \theta_0$, for all $x  \leq \theta$, and $L=L(\np) \in \lbr \np \rbr$ with $m L \geq L_0$, 
\begin{equation}\label{eq:FullLocalModerate}
\Pfull\Big( \TF_L\big(\gamma_{(0,0),(\np,m\np)}\big) \geq x L^{2/3} \Big) \leq C_1\exp({-c_1 x^{2}}) 
\end{equation} for all $n$ large enough. Similarly, consider a family of $(N,k)$-periodic environments such that \eqref{eq:NonDegeneratedEnvironment} and \eqref{eq:CorrectPeriodicScale} hold with respect to some $\mathfrak{a}>0$. Then we get that for all $ x  \leq \theta^{1/3}$, and $L \in \lbr \np \rbr$ with $L \geq L_0$ that we find constants $c_2,C_2>0$ so that
\begin{equation}\label{eq:PeriodicLocalModerate}
\Pper\Big( \TF_L\big(\gper_{(0,0),(\np,m\np)}\big) \geq x L^{2/3} \Big) \leq C_2\exp({-c_2 x^{2}}) 
\end{equation} for all $N$ and $\np=\np(N)$ large enough.
\end{lemma} 

Again, we defer the proof of Lemma~\ref{lem:ModerateLocal} to Appendix~\ref{sec:AppendixLPP}. 
As a consequence, similar to Lemma~\ref{lem:TypicallyGood} for the directed landscape, we can bound the local periodic transversal fluctuations of last passage times. To this end, we set for all $\vartheta \in \N$ 
\begin{equation}\label{eq:VerticalLPPdistance}
\begin{split}
    \gamma^{\per,+}_{u,v}(\vartheta) &:= \sup\big\{ \ell \in \N \colon (\ell,\vartheta) \in \gper_{u,v} \big\} ,  \\
    \gamma^{\per,-}_{u,v}(\vartheta) &:= \inf\{ \ell \in \N \colon (\ell,\vartheta) \in \gper_{u,v} \} ,
    \end{split}
\end{equation} and simply write $\gamma^{\per,\pm}_{u,v}$ if we allow for either of the two choices. 
Recall from \eqref{def:FixedSlope} that $m_{\ast}=k^{2}(N-k)^{-2}$. 
\begin{definition}\label{def:GoodLPP}
    Let $\delta>0$ and $N\in \N$. For fixed $u=(u_1,u_2)\in \Z^2$ and $v=(v_1,v_2)\in \Z^2$ with $v \succeq u$, consider the set of lattice points
\begin{equation}\label{def:GoodIncrementsDiscrete}
    (z_i)_{i\in \lbr \ell\rbr} := \{ i \delta N^{3/2}  \in [u_2,v_2] \text{ for some } i\in \N \} 
\end{equation} for some $\ell\in \N_0$, with the convention that $z_0=u_2$ and $z_{\ell+1}=v_2$, and $z_{i-1} \leq z_{i}$ for all $i\in \lbr \ell+1 \rbr$. 
    We say $\gper_{u,v}$ is $\boldsymbol{(\delta,N)}$\textbf{-good} if for all $i\in \lbr \ell +1 \rbr$ 
\begin{equation}
  \sup_{\varphi \in [z_{i-1},z_i]} \Big| \gamma^{\per,\pm}_{u,v}( \varphi ) - \gamma^{\per,\pm}_{u,v}(z_{i-1})  - m_{\ast}^{-1} (\varphi-z_{i-1}) N^{3/2}   \Big| \leq \frac{N}{32}.
\end{equation}
\end{definition}
 The next result states that with high probability, all lattice paths in a given domain are $(\delta,N)$-good when the sites belong to the set
\begin{equation}\label{def:SlopeDn}
D_{\delta}^N := \left\{ (v;w) \in \Z^{4}_{\uparrow} \colon  \,   m^{-1}_{\ast} - \delta^{-1/2} N^{-1/2}  \leq \frac{w_1-v_1}{w_2-v_2} \leq m^{-1}_{\ast} + \delta^{-1/2} N^{-1/2} \right\} . 
\end{equation} 
 Again, the proof uses standard arguments from last passage percolation, together with Lemma~\ref{lem:ModerateLocal}, and is deferred to Appendix~\ref{sec:AppendixLPP}. To state the result, we consider the set of sites
\begin{equation*}
    \mathbb{X}^{N}_{\delta} := \left\{   u=(u_1,u_2) \in \Z^2  \, \colon \, u_1 - m_{\ast}^{-1}u_2 \in [-\delta^{-1}N,\delta^{-1}N] \, \wedge \,  u_2 \in  \lbr N^{3/2} \rbr\right\} 
\end{equation*} for all $\delta>0$ and $N\in \N$.
\begin{proposition}\label{pro:TypicallyGoodLPP}
    There exist constants $C,c>0$ such that for all $\delta>0$ and $N$ large enough
    \begin{equation*}
        \P\left(  \gper_{u,v} \text{ is } (\delta,N)\text{-good for all } (u_1,u_2;v_1,v_2) \in D_{\delta}^N \cap \big(\mathbb{X}^{N}_{\delta} \times \mathbb{X}^{N}_{\delta}\big)\right) \geq 1- C\exp\big(-c \delta^{-1/6}\big).
    \end{equation*}
\end{proposition}

\subsection{Approximation of the last passage times}\label{sec:Approximation}

Next, we aim to approximate last passage percolation in periodic environments by full-space last passage percolation. This is achieved by an approximation of the i.i.d.~and periodic last passage times using refined sequences of overlapping parallelograms. The environment inside a parallelogram will be coupled with an i.i.d.~environment. Here, we utilize the control over last passage times and geodesics from the previous section. We will now formalize this construction. For all $i \in \Z$ and $j\in \N_0$, let
\begin{align*}
    r_{4}(i,j) &:= \Big(\Big\lfloor \Big( \frac{j-1}{2} - \frac{1}{8} \Big)(N-k)  + (i-1) N^{3/2}n^{-1} \Big\rfloor, \Big\lfloor - \Big( \frac{j-1}{2} - \frac{1}{8} \Big)k  +   m_{\ast} (i-1) N^{3/2}n^{-1} \Big\rfloor\Big) , \\
    r_{1}(i,j) &:= \Big(\Big\lfloor \Big( \frac{j}{2} + \frac{1}{8} \Big)(N-k) \Big\rfloor+\lfloor (i-1) N^{3/2}n^{-1} \rfloor,  \Big\lfloor- \Big( \frac{j}{2} + \frac{1}{8} \Big)k + m_{\ast} (i-1) N^{3/2}n^{-1} \Big\rfloor\Big) , 
\end{align*} as well as
\begin{align*}
    r_{2}(i,j) &:= r_1 + \big(\lfloor N^{3/2}n^{-1} \rfloor , \lfloor  m_{\ast} N^{3/2}n^{-1} \rfloor\big) ,  \\
    r_{3}(i,j) &:= r_4 + \big(\lfloor N^{3/2}n^{-1} \rfloor , \lfloor  m_{\ast} N^{3/2}n^{-1} \rfloor\big) . 
\end{align*}
For all $i \in \N_0$ and $j\in \Z$, similar to \eqref{eq:FillingRectangle}, we denote by $R_{n}^{i,j}$ the rectangle spanned by the points $(r_{\ell}(i,j))_{\ell \in\lbr 4\rbr}$ in counter-clockwise order, with $r_{1}(i,j)$ denoting the bottom right corner.
\begin{figure}
\centering
\begin{tikzpicture}[scale=.95]

   \draw[line width=1.2pt, dotted] (0,0) -- (10,0);
   \draw[line width=1.2pt, dotted] (0,3) -- (10,3);

  \draw[red, line width =1.2 pt] (1,0) -- (1,3) -- (9,3) -- (9,0) -- (1,0);
  
    
    
    
    \draw[fill=darkblue!20, line width =1.2 pt] (1.05,0.05) -- (1.05,2.95) -- (2.95,2.95) -- (2.95,0.05) -- (1.05,0.05);

        \draw[fill=darkblue!20, line width =1.2 pt] (7.05,0.05) -- (7.05,2.95) -- (8.95,2.95) -- (8.95,0.05) -- (7.05,0.05);
   
\draw[line width =1 pt] (7.85+0.35-0.4,0.05) to[curve through={(7.85+0.5-0.4,0.55)..(7.58-0.4,1.05) .. (7.85+0.4-0.3,1.65) .. (7.85+0.55-0.4,2.45) }] (7.85+0.35-0.4,2.95);


\node[scale=0.8] (site) at (5,1.5){$\bar{R}_n^{i,j}$};  
\node[scale=0.8] (site) at (2,1.5){$R_n^{i,j} \setminus \bar{R}_n^{i,j}$};  

\node[scale=0.6] (site) at (9.3,-0.25){$r_{1}(i,j)$};  
\node[scale=0.6] (site) at (9.3,3.25){$r_{2}(i,j)$};  
\node[scale=0.6] (site) at (1,3.25){$r_{3}(i,j)$};  
\node[scale=0.6] (site) at (1,-0.25){$r_{4}(i,j)$};  


\node[scale=0.6] (site) at (6.8,-0.25){$\bar{r}_{1}(i,j)$};  
\node[scale=0.6] (site) at (6.8,3.25){$\bar{r}_{2}(i,j)$};  
\node[scale=0.6] (site) at (3.1,3.25){$\bar{r}_{3}(i,j)$};  
\node[scale=0.6] (site) at (3.1,-0.25){$\bar{r}_{4}(i,j)$}; 

\node[scale=0.6] (site) at (7.85,-0.25){$u$};  
\node[scale=0.6] (site) at (7.85,3.25){$v$}; 
%
%
%
%
%

	\end{tikzpicture}	
\caption{\label{fig:RectangleComposition1} Rectangles $R_{n}^{i,j}$ and $\bar{R}_{n}^{i,j}$ used in the construction of patched last passage percolation for $\rho=1$, rotated by $\pi/4$. We denote by $r_\ell(i,j)$ and $\bar{r}_\ell(i,j)$ for $\ell \in \lbr 4\rbr$ the corners of the respective rectangles. Note that the path $\pi_{u,v}$ is an element of $V_{n}^N$ by bounding the distance between $u$ and $v$ perpendicular to the  $(1,1)$-direction.}
 \end{figure}
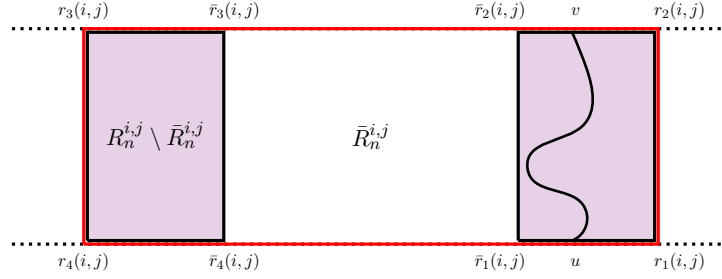
 Moreover, we set
\begin{equation}
\bar{r}_\ell(i,j) := r_\ell(i,j) + \frac{1}{16}(1-2 \mathds{1}_{\ell \in \{1,2\}}) (N-k,-k)
\end{equation} for all $\ell \in \lbr 4\rbr$, i.e., we shrink the parallelograms $R_{n}^{i,j}$ to the left and right. We denote the corresponding rectangle with corners $(\bar{r}_{\ell}(i,j))_{\ell \in \lbr 4\rbr}$ by $\bar{R}_{n}^{i,j}$. 
A visualization of this construction is given in Figure~\ref{fig:RectangleComposition1}.
 By construction, for all $n$ and $N$ large enough, the parallelograms satisfy 
\begin{equation}\label{eq:CloseRectangle}
R_{n}^{i,j} \cap R_{n}^{k,\ell} \neq \emptyset \ \text{ if and only if } \  |i-k| \leq  1  \text{ and } |j-\ell| \leq 1. 
\end{equation} 
We set in the following
\begin{equation}
\label{eq:RectanglesDiscrete}  R_{n}^{i}:=\bigcup_{j \in \Z} R_{n}^{i,j} . 
\end{equation}
and define the \textbf{restricted last passage times} onto the sets $R_{n}^{i,j}$ as
 \begin{equation}
 T^{i,j}_{u,v} := \sup_{\pi \in \Pi_{u,v}} \sum_{w \in \pi \setminus \{ v\}} \omega_{w} \mathds{1}_{ \{ w \in R_{n}^{i,j} \} }
 \end{equation} for all $u,v \in  R_{n}^{i,j}$. In the following, our goal is to replicate the i.i.d.~and periodic last passage times using only the restricted last passage times $T^{i,j}_{u,v}$ with some $i\in \lbr n  \rbr$ and $j\in\Z$. To this end, recall $m_{\ast}=k^2(N-k)^{-2}$, and set in the following
\begin{equation}
\Z^{4}_{\uparrow} := \left\{ (v,w) \in \Z^2 \times \Z^2 \colon w \succeq v\right\} .
\end{equation} We define the sets
\begin{equation*}
V_n^{N} := \left\{ (v;w) \in \Z^{4}_{\uparrow} \,  \colon w = v + (\lfloor -s(N-k)\rfloor,\lfloor sk\rfloor) + (\lfloor tN^{\frac{3}{2}}\rfloor,\lfloor tm_{\ast}N^{\frac{3}{2}}\rfloor)\, ,t>0, s\in \Big[-\frac{1}{16},\frac{1}{16}\Big]\right\}. 
\end{equation*}  
Fix $i\in \N$. Then for any pair of points $(u;v)\in V_{n}^N$ with $u,v \in R_n^{i}$, we have that $u,v \in R_n^{i,\mathfrak{j}}$ for some $\mathfrak{j} \in \Z$. Whenever $\mathfrak{j}=\mathfrak{j}(u,v)$ is unique, then let $M = \mathfrak{j}\mod 2$ and define
\begin{equation}
\Tpat_{u,v} := T^{i,M}_{u,v} . 
\end{equation} 
Otherwise, recalling \eqref{eq:CloseRectangle}, suppose that there exists some $\mathfrak{j}\in \Z$ such that $u,v \in R_{n}^{i,\mathfrak{j}} \cap R_n^{i,\mathfrak{j}+1}$. If in addition $(u;v) \in V_{n}^N$ holds, we get that either $u,v \in \bar{R}_{n}^{i,\mathfrak{j}}$ or $u,v \in \bar{R}_{n}^{i,\mathfrak{j}+1}$. In this case, let $M = \mathfrak{j} \mod 2$, and set for all $u,v\in \bigcup_{j\in \N} R_n^{j}$
\begin{equation}\label{def:PatchedShort}
\Tpat_{u,v} := \begin{cases}
T^{i,M}_{u,v} & \text{ if } u,v\in \bar{R}_{n}^{i,\mathfrak{j}} \cap V^{N}_n ,  \\
T^{i,M+1}_{u,v} & \text{ if } u,v\in \big(\bar{R}_{n}^{i,\mathfrak{j}+1} \setminus  \bar{R}_{n}^{i,\mathfrak{j}}\big) \cap V^{N}_n , \\
-\infty & \text{ otherwise.}
\end{cases}  
\end{equation}
Intuitively, for every pair $(u,v)$ sufficiently far from the boundary, such that it is fully contained in one of $n$ many patches of last passage percolation, we let $\Tpat_{u,v}$ agree with the restricted last passage times between $u$ and $v$ from the respective patch. We make the following observation on comparing the times $(\Tpat_{u,v})_{u,v \in \Z^2}$ to full-space and periodic last passage times. \begin{figure}
\centering
\begin{tikzpicture}[scale=.95]

   \draw[line width=1.2pt, dotted] (0,0) -- (15,0);
   \draw[line width=1.2pt, dotted] (0,3) -- (15,3);

  \draw[red, line width =1.2 pt] (1,0) -- (1,3) -- (9,3) -- (9,0) -- (1,0);
  
    \draw[red, line width =1.2 pt] (6,0) -- (6,3) -- (14,3) -- (14,0) -- (5,0);
    
    \draw[darkblue, line width =1.2 pt] (6.05,0.05) -- (8.1,0.05) -- (8.1,2.95) -- (6.05,2.95) -- (6.05,0.05);
    
    \draw[darkblue, line width =1.2 pt] (6.55+0.35,0.05) -- (8.6+0.35,0.05) -- (8.6+0.35,2.95) -- (6.55+0.35,2.95) -- (6.55+0.35,0.05);
    
    \draw[fill=darkblue!20, line width =1.2 pt] (6.55+0.35,0.05) -- (8.1,0.05) -- (8.1,2.95) -- (6.55+0.35,2.95) -- (6.55+0.35,0.05);
   
\draw[line width =1 pt] (6.55+0.35-0.4,0.05) to[curve through={(6.55+0.5-0.4,0.55)..(6.55-0.4,1.05) .. (6.55+0.4-0.3,1.65) .. (6.55+0.55-0.4,2.45) }] (6.55+0.35-0.4,2.95);

\draw[line width =1 pt] (8.1+0.4,0.05) to[curve through={(8.1-0.2+0.4,0.55)..(8.1+0.4,1.05) .. (8.1+0.3+0.4,1.65) .. (8.1+0.25+0.4,2.45)..(8.05+0.4,2.7) }] (8.1+0.4,2.95);

\node (site) at (3.5,1.5){$R_n^{i,j}$};  
\node (site) at (11.5,1.5){$R_n^{i,j+1}$};  

\node[scale=0.6] (site) at (1,3.25){$r_{3}(i,j)$};  
\node[scale=0.6] (site) at (1,-0.25){$r_{4}(i,j)$};  

\node[scale=0.6] (site) at (14.3,-0.25){$r_{1}(i,j+1)$};  
\node[scale=0.6] (site) at (14.3,3.25){$r_{2}(i,j+1)$};  

\node[scale=0.6] (site) at (8.5,-0.25){$\tilde{r}_{1}(i,j)$};  
\node[scale=0.6] (site) at (8.5,3.25){$\tilde{r}_{2}(i,j)$};  
\node[scale=0.6] (site) at (7.0,3.25){$\tilde{r}_{3}(i,j+1)$};  
\node[scale=0.6] (site) at (7.0,-0.25){$\tilde{r}_{4}(i,j+1)$}; 
%
%
%
%
%

	\end{tikzpicture}	
\caption{\label{fig:RectangleComposition2} Assignment of the different rectangles used in the construction of patched last passage percolation for $\rho=1$, rotated by $\pi/4$. In Lemma \ref{lem:LPPCoupleSmallSlaps}, we argue that whenever the geodesics $\pi_{\tilde{r}_{1}(i,j),\tilde{r}_{2}(i,j)}$ and $\pi_{\tilde{r}_{4}(i,j+1),\tilde{r}_{3}(i,j+1)}$ do not leave the rectangles $R_{n}^{i,j} \setminus \bar{R}_{n}^{i,j}$ and $R_{n}^{i,j+1} \setminus \bar{R}_{n}^{i,j+1}$, respectively, then the last passage times can be coupled to agree on the light purple shaded part $\bar{R}_{n}^{i,j} \cap \bar{R}_{n}^{i,j+1}$.}
 \end{figure}
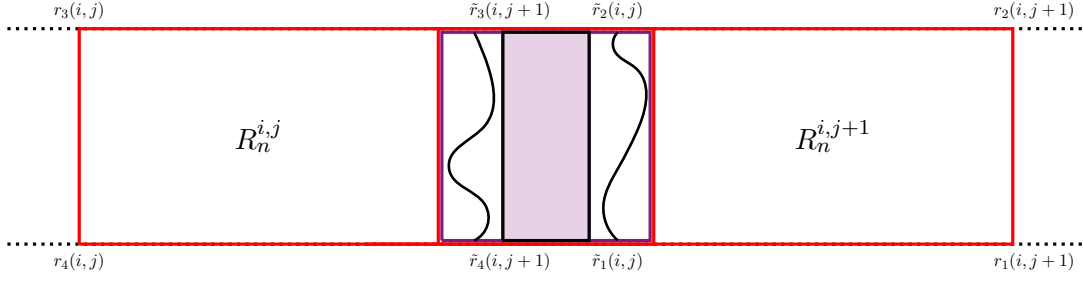

\begin{lemma}\label{lem:LPPCoupleSmallSlaps}
Let $(\omega_v)_{v \in \Z^2}$ be an i.i.d.~environment. Fix $j\in \N$. Then there exist some constants $c,C>0$ and a coupling $\mathbf{P}_{\textup{full}}$ between the i.i.d. LPP environment $\Pfull$ and the family of environments  used to determine the last passage times $(\Tpat_{u,v})_{u,v \in \Z^2}$ such that for fixed $i\in \N$,
\begin{equation}\label{eq:LPPCoupleSmallSlaps1}
 \mathbf{P}_{\textup{full}}\Big( T_{u,v}= \Tpat_{u,v} \, \forall (u;v) \in (R_n^{i,j}\times R_n^{i,j}) \cap V^{N}_n, j\in \lbr n -1 \rbr \Big) \geq 1 - C \exp\big(-cn^{1/2}\big) 
\end{equation} for all $n\in \N$, and $N$ sufficiently large. 
Similarly, let $(\omega_v)_{v \in \Z^2}$ be an $(N,k)$-periodic environment and assume that \eqref{eq:NonDegeneratedEnvironment} holds. Then there exists a coupling $\mathbf{P}_{\per}$ between  $\big(\Tper_{u,v}\big)_{u,v \in \Z^2}$ and $\big(\Tpat_{u,v}\big)_{u,v \in \Z^2}$  such that for fixed $i\in \N$
\begin{equation}\label{eq:LPPCoupleSmallSlaps2}
\mathbf{P}_{\per}\Big( \Tper_{u,v}= \Tpat_{u,v} \, \forall (u;v) \in (R_n^{i}\times R_n^{i})  \cap V^{N}_n   \Big) \geq 1 - C \exp\big(-cn^{1/2}\big) 
\end{equation} for all $n\in \N$, and $N$ sufficiently large. 
\end{lemma}


\begin{proof}
First consider the case that $(\omega_v)_{v \in \Z^2}$ is an  $(N,k)$-periodic environment. 
Let $\gamma^{(n)}_{u,v}$ and $\gper_{u,v}$ denote the geodesics between $u$ and $v$ corresponding to the patched and periodic last passage times $(\Tpat_{u,v})_{u,v \in \Z^2}$ and $(\Tper_{u,v})_{u,v \in \Z^2}$, respectively. Now let the environments $(\omega^{(n)}_v)_{v \in R_{n}^{i,j}}$ used in the construction of $(\Tpat_{u,v})_{u,v \in \Z^2}$ agree with the periodic environment $(\omega_v)_{v \in \Z^2}$, and note that $(\omega^{(n)}_v)_{v \in R_{n}^{i,j}}$ are i.i.d.~for given $i,j$ and $n$. Recall the points $r_{\ell}(i,j)$ and $\bar{r}_{\ell}(i,j)$ defined for all $i,j$ and $\ell \in \lbr 4 \rbr$. We define
\begin{equation}
\tilde{r}_{\ell}(i,j) := \frac{1}{2} r_{\ell}(i,j) + \frac{1}{2} \bar{r}_{\ell}(i,j)
\end{equation} as the midpoint between $r_{\ell}(i,j)$ and $\bar{r}_{\ell}(i,j)$; see Figure \ref{fig:RectangleComposition2}.
Recalling the moderate deviations estimates on the transversal fluctuations from Lemma~\ref{lem:GeodesicsFullPeriodicAgree}, and using the ordering of geodesics, 
\begin{align*}
\big\{ \Tper_{u,v}= \Tpat_{u,v} \, \forall u,v \in (\bar{R}_n^{i,j}\times \bar{R}_n^{i,j}) \cap V^{N}_n \big\} &\subseteq \Big\{ \TF(\gamma^{(n)}_{\bar{r}_1(i,j-1),\bar{r}_2(i,j-1)}), \TF(\gper_{\bar{r}_1(i,j-1),\bar{r}_2(i,j-1)}) \leq \frac{N}{32} \Big\} \\
&\cap \Big\{ \TF(\gamma^{(n)}_{\bar{r}_4(i,j+1),\bar{r}_3(i,j+1)}), \TF(\gper_{\bar{r}_4(i,j+1),\bar{r}_3(i,j+1)}) \leq \frac{N}{32} \Big\} .
\end{align*} 
Hence, using Lemma \ref{lem:GeodesicsFullPeriodicAgree} and a union bound over $j \in \lbr n\rbr$, and hence $j\in \Z$, we get that
\begin{equation}
\mathbf{P}_{\per}\big( T_{u,v}= \Tpat_{u,v} \, \forall u,v \in R_n^{i} \, \colon (u,v)\in V^{N}_n \big) \geq 1- C_1n\exp\big(-c_1n^{1/2}\big)  
\end{equation} for some constants $c_1,C_1>0$, and all $n\in \N$, provided that $N$ is sufficiently large. This shows~\eqref{eq:LPPCoupleSmallSlaps2}. Let us now turn to the case \eqref{eq:LPPCoupleSmallSlaps1} of i.i.d.~environments. This follows by the same arguments as \eqref{eq:LPPCoupleSmallSlaps2}, replacing the geodesic $\gper$ by $\gamma$, and noting that by Lemma \ref{lem:TransversalFluctuations}
together with a union bound over $j\in  \lbr n \rbr$, we get that
\begin{equation*}
\mathbf{P}_{\textup{full}}\big( T_{u,v}=\Tpat_{u,v} \, \forall u,v \in R_n^{i,j}, i \in \lbr n-1\rbr \, \colon (u,v)\in V^{N}_n\big) \geq 1- C_2 n\exp\big(-c_2n^{1/2}\big)
\end{equation*} for some constants $c_2,C_2>0$, and all $n\in \N$, with $N$ large enough.
\end{proof}

In words, Lemma \ref{lem:LPPCoupleSmallSlaps} ensures that the restricted last passage times $\Tpat_{u,v}$ behave as full-space and periodic last passage times, provided that $u$ and $v$ satisfy suitable slope conditions, and are not too far apart. In particular, this allows us to couple periodic and full-space last passage percolation for patches of size $N$ by $n^{-1} N^{3/2}$ for sufficiently large $n$ and $N$.
 As a consequence of the previous estimates, we record the following bound on the local
 fluctuations of geodesics in periodic last passage percolation, which will be a key input in  establishing continuity of the periodic directed landscape. Define the set
\begin{equation*}
    \bar{B}^N_{\delta}(u) := \left\{   u=(u_1 + z +m_{\ast}^{-1}\tau,u_2 + \tau ) \in \Z^2  \, \colon \, z \in [- \delta N,\delta N] \, \wedge \, \tau \in  [ - \delta N^{3/2} , \delta N^{3/2} ]\right\} 
\end{equation*} for all $u\in \Z^2$, and all $\delta>0$ and $N\in \N$.
 
 \begin{lemma}\label{lem:ConsistentLPPLocally} For $\delta>0$ and $N\in \N$, let $(u;v) = (u_1,u_2;v_1,v_2) \in D_{\delta^{1/4}}^N \cap \big(\mathbb{X}^{N}_{\delta^{1/4}} \times \mathbb{X}^{N}_{\delta^{1/4}}\big)$ be such that $v_2-u_2 \geq 3\delta^{1/2}N^{3/2}$.  There exist $c,C>0$, which do not depend on $\delta>0$ or $N\in \N$, such that 
the events
\begin{equation*}
    \hat{D}_{N,\delta}^{u,v} := \Big\{\gamma^{\per}_{u',v}(u_2+z) \in [u_1-5\delta^{3/5}N,u_1+5\delta^{3/5}N] \, \forall \, u' \in \bar{B}^N_{\delta}(u) , z\in [\delta N^{3/2},3\delta N^{3/2}] \Big\}
\end{equation*}
 satisfy for all $\delta>0$ and all $N$ large enough
 \begin{equation}\label{eq:Ball1Local}
\mathbb{P}\big(\hat{D}_{N,\delta}^{u,v} \big) \geq 1 - C\exp(-c\delta^{-1/24}) .
\end{equation}  Similarly, we get that the events
\begin{equation*}
    \check{D}_{N,\delta}^{u,v} := \Big\{\gamma^{\per}_{u,v'}(v_2-z) \in [v_1-5\delta^{3/5}N,v_1+5\delta^{3/5}N] \, \forall \, v' \in \bar{B}^N_{\delta}(v) , z\in [\delta N^{3/2},3\delta N^{3/2}] \Big\}
\end{equation*}
satisfy for all $\delta>0$ and all $N$ large enough
 \begin{equation}\label{eq:Ball2Local}
\mathbb{P}\big(\check{D}_{N,\delta}^{u,v} \big) \geq 1 - C\exp(-c\delta^{-1/24}) .
\end{equation}
 \end{lemma}
 The proof of Lemma~\ref{lem:ConsistentLPPLocally} is deferred to Appendix~\ref{sec:AppendixLPP}. Let us now return to establishing the convergence of periodic last passage percolation to the periodic direct landscape. 
We argue in the sequel that we can extend the definition of the last passage times from short time scales $N^{3/2}n^{-1}$ to larger time scales of order $N^{3/2}$. More precisely, similar to Definition \ref{def:PatchedDLExt} of the patched directed landscape, we extend the approximation of periodic last passage times via a  variational characterization. To this end, we introduce some additional notation. For all $n\in \N$,
 \begin{equation}\label{def:BigRectangle}
 R_{n} := \bigcup_{i\in \lbr n  \rbr, j\in \Z} R_{n}^{i,j} =  \bigcup_{i\in \lbr n \rbr} R_{n}^{i}
 \end{equation}
for the union of the parallelograms in \eqref{eq:RectanglesDiscrete}, and
\begin{equation}
(R_n)^2_{\uparrow} := \Big\{ (u;v) \in \Z^{4}_{\uparrow} \colon u,v \in R_{n} \Big\} 
\end{equation}
for the directed set of pairs of points in $R_n$. Furthermore, let 
\begin{equation}
\partial \bar{R}_n^{i,j} := S(\bar{r}_{2}^{i,j},\bar{r}_{3}^{i,j})
\end{equation} be the upper boundary of the rectangles $\bar{R}_n^{i,j}$, where we let $S(u,v)$ with 
\begin{equation}\label{def:Segment}
S(v,w) := \big\{ u \in \Z^2 \colon u= \big(\lfloor v_1 s + w_1 (1-s) \rfloor , \lfloor v_2 s + w_2 (1-s) \rfloor \big) \text{ for some } s\in [0,1] \big\}  
\end{equation}
denote the discrete segment joining $u$ and $v$. Similarly,  we write  $\partial \bar{R}_n^{i}= \bigcup_{j\in \Z}\partial \bar{R}_n^{i,j}$.
Note that the last passage times  $T^{(n)}_{u,v}$ are with high probability well-defined for all pairs of sites in $(R^{j}_n \times R^{j}_n) \cap V^{N}_n$ due to Lemma \ref{lem:LPPCoupleSmallSlaps}, and can be extended to $(R_n)^2_{\uparrow}$ as follows.

\begin{definition}\label{def:ExtendedPeriodicLPT}
Define the \textbf{patched last passage times} $\Tpat$ on $(R_n)^2_{\uparrow}$ of length $N^{3/2}$ via a maximization procedure over the last passage times on the restricted rectangles. More precisely, we take $\Tpat_{u,v}$ when $u,v \in R_{n}^{i}$ for some $i\in \lbr n\rbr$ as in \eqref{def:PatchedShort}. For all $u,v \in R_n$ such that $u \in R_n^{i}$ and $v \in R_n^{j}$ with  $1\leq i<j \leq n$, we proceed by induction as in Definition~\ref{def:PatchedDLExt}. When $j=i+1$, 
\begin{equation}
\Tpat_{u,v} :=   \sup_{w\in \partial \bar{R}_n^{i}} \big(\Tpat_{u,w} +  \Tpat_{w,v}\big).   
\end{equation} 
When $k=j-i>1$, suppose that the patched last passage times have been defined for all $i^{\prime},j^{\prime}<k$, and set
\begin{equation}
\Tpat_{u,v} :=   \sup_{w\in \partial \bar{R}_n^{j-1}} \big(\Tpat_{u,w} +  \Tpat_{w,v}\big).   
\end{equation} 
\end{definition}

\subsection{Convergence to periodic directed landscape}\label{sec:ApproximationLargePeriodic} 
We have the following strategy to use the patched last passage times in order to show convergence of the periodic last passage percolation to the periodic directed landscape. Our first step is to couple periodic last passage percolation with the patched last passage times. In a second step, we establish the convergence of the patched last passage times to patched directed landscapes. In a third step, we use the convergence of patched directed landscapes to the periodic directed landscape and conclude. 

Recall the patched last passage times $\Tpat$ from Definition \ref{def:ExtendedPeriodicLPT}, and $\mathbf{u},\mathbf{v}$ from \eqref{def:VectorsDL}. Further, recall $M$ and $V$ from \eqref{def:MeanLPP} and \eqref{def:VarianceLPP}, respectively. As for last passage percolation on the full-space in \eqref{def:LPPSheet}, we define now the \textbf{patched LPP landscape} by a suitable rescaling, i.e.,  set 
\begin{equation}\label{def:PatchedRenormalized}
\LLpaeps(x,s;y,t):=  \frac{1}{V(\rho,\varepsilon)} \Big(\Tpat_{x \varepsilon^{-2}\mathbf{u} + s  \varepsilon^{-3}\mathbf{v},y \varepsilon^{-2}\mathbf{u}+t \varepsilon^{-3}\mathbf{v}}- M(\rho,\varepsilon,s,t) \Big)  
\end{equation} for all $ (x,s;y,t) \in \R^4_{\uparrow} \cap \mathcal{V}_{n}$ such that $x \varepsilon^{-2}\mathbf{u} + s  \varepsilon^{-3}\mathbf{v},y \varepsilon^{-2}\mathbf{u}+t \varepsilon^{-3}\mathbf{v} \in \Z^2$, and linear interpolation otherwise. In order to match with periodic last passage percolation, we will set 
\begin{equation}\label{eq:EpsilonScale}
\varepsilon = \varepsilon(k,N) := k^{-1/2} \quad \text{ and } \quad  \rho=m_{\ast}^{-1}=\frac{(N-k)^2}{k^2} , 
\end{equation} i.e.,  we get that for all $j\in \Z$
\begin{equation}
\LLpaeps(x+j,s;y+j,t) = \LLpaeps(x,s;y,t) .
\end{equation} In particular, $\LLpaeps$ yields a random function on the cylinder $\SetT$. We remark that this definition can also be directly extended to $(\mathbb{T}\times [0,T])^{2}_{\uparrow}$ for any fixed $T>0$.
Next, we define a corresponding periodic LPP landscape with respect to periodic last passage percolation. To this end, let $N\in \N$ and $k=k(N) \in \lbr N\rbr$ be such that \eqref{eq:NonDegeneratedEnvironment} holds. 
We  define the \textbf{periodic LPP landscape} as the random function on $\SetT$ with
\begin{equation}\label{def:periodicLPPSheet}
\LLpeps(x,s;y,t) :=  \frac{1}{V(\rho,\varepsilon)} \Big(\Tper_{x\varepsilon^{-2}\mathbf{u} + s\varepsilon^{-3}\mathbf{v},y\varepsilon^{-2}\mathbf{u}+t\varepsilon^{-3}\mathbf{v}}- M(\rho,\varepsilon,s,t) \Big) 
\end{equation} with respect to the periodic last passage times $\Tper$ when $x\varepsilon^{-2}\mathbf{u} + s\varepsilon^{-3}\mathbf{v},y\varepsilon^{-2}\mathbf{u}+t\varepsilon^{-3}\mathbf{v} \in \Z^2$, and linear interpolation otherwise. Recall the definition of $\mathcal{D}_{\varepsilon}$ from \eqref{def:SlopeConti}.
We have the following coupling between the periodic and the patched LPP landscapes. 
\begin{lemma}\label{lem:LPPtoDLRestricted} Fix some $c_0>0$. 
There exists a coupling $\mathbf{P}_{N,n}$ between the patched and periodic LPP landscapes, and constants $c,C>0$, such that for all $n\in \N$, we find some $N_0=N_0(n)$ such that for all $N \geq N_0$, 
\begin{equation*}
\mathbf{P}_{N,n}\Big( \LLpeps(x,s;y,t)=\LLpaeps(x,s;y,t) \, \forall (x,s;y,t) \in \mathcal{D}_{c_0n^{-1}} \cap \SetT \Big) \geq 1 - C\exp(-cn^{\frac{1}{6}}) .
\end{equation*} 
\end{lemma}
\begin{proof} 
Note that by construction, we have that
\begin{equation}\label{eq:UpperComp}
    \LLpaeps(x,s;y,t)  \leq \LLpeps(x,s;y,t) 
\end{equation} almost surely for all $(x,s;y,t) \in \mathcal{D}_{n^{-1}} \cap \SetT$. By construction, the equality holds when the geodesic $\gper_{u,v}$, mapping to $(x,s;y,t)$ after rescaling of $u$ and $v$, is $(n^{-1},N)$-good.
Hence, by Proposition ~\ref{pro:TypicallyGoodLPP}, there exists an event of probability at least $1-C\exp(-cn^{1/6})$ with some constants $c,C>0$ such that for all $(x,s;y,t) \in \mathcal{D}_{c_0n^{-1}} \cap \SetT$ for all $c_0>0$, 
\begin{equation}
     \LLpaeps(x,s;y,t)  \geq \LLpeps(x,s;y,t) . 
\end{equation}
Together with \eqref{eq:UpperComp}, this gives the desired result. 
\end{proof}

Next, we establish the convergence of the patched LPP landscape $\LLpaeps$ from Definition~\ref{def:ExtendedPeriodicLPT} to the patched directed landscapes $\Lpa$ from Definition \ref{def:PatchedDLExt}. Recall that $\tilde{\omega}=2^{\frac{3}{2}}(1+\sqrt{\rho})^{-\frac{1}{2}}\rho^{\frac{1}{4}}$.

\begin{lemma}\label{lem:TruncatedConvergence}
For all $n\in \N$,  the patched last passage times $\LLpaeps$ converge  to the patched directed landscape $\Lpa$ as $\varepsilon \rightarrow 0$, i.e., there exist events $(\mathcal{A}_n)_{n \in \N}$ with $\mathcal{A}_n=\mathcal{A}_n(N)$ such that 
\begin{equation}
    \P(\mathcal{A}_n) \geq 1-C\exp\big(-cn^{1/6}\big)
\end{equation}
 with some constants $c,C>0$ and all $n\in \N$, while we have that as $\varepsilon \rightarrow 0$
\begin{equation*}
\P\Big(\big(\LLpaeps\big)_{(x,s\tilde{\omega}^{-1};y,t\tilde{\omega}^{-1}) \in \SetTw \cap \mathcal{D}_{n^{-1}}} \in \cdot \, \Big| \, \mathcal{A}_n \Big) \rightarrow \P\Big( (\Lpa)_{(x,s;y,t ) \in \SetT \cap \mathcal{D}_{n^{-1}}} \in \cdot \, \Big| \, \mathcal{A}_n \Big) . 
\end{equation*} 
\end{lemma}
\begin{proof}
We will in the following only consider the case that $\tilde{\omega}=1$ as a similar argument applies for general $\tilde{\omega}$. Recall the rectangles $\mathcal{R}_n^{i,j}$ from \eqref{eq:FillingRectangle}, and that $\mathcal{R}_n^{i}= \mathcal{R}_n^{i,1} \cup \mathcal{R}_n^{i,2}$. Note that by the definition of the patched last passage times and the patched directed landscape, and the change of $D_{n^{-1}}^N$ to $\mathcal{D}_{n^{-1}}$, it suffices to show that as $\varepsilon\rightarrow 0$
\begin{equation}
(\LLpaeps)_{(x,s;y,t) \in \big(\mathcal{R}_n^{1}\times \mathcal{R}_n^{1} \big) \cap \mathcal{D}_{n^{-1}} } \rightarrow (\Lpa)_{(x,s;y,t) \in \big(\mathcal{R}_n^{1} \times \mathcal{R}_n^{1}\big) \cap \mathcal{D}_{n^{-1}}} 
\end{equation} in distribution conditioning on the events $\mathcal{A}_n$ (which we will exactly describe later on). In order to simplify notation, consider the projections 
\begin{align*}
\LLeps_{n,(1)}&=(\LLpaeps(x,s;y,t))_{(x,s;y,t) \in \SetT \cap (\mathcal{R}_n^{1,1} \times \mathcal{R}_n^{1,1}) \cap \mathcal{D}_{n^{-1}}} \\
\LLeps_{n,(2)}&=(\LLpaeps(x,s;y,t))_{(x,s;y,t) \in \SetT \cap (\mathcal{R}_n^{1,2} \times \mathcal{R}_n^{1,2}) \cap \mathcal{D}_{n^{-1}}}  . 
\end{align*} Then by \eqref{eq:LPPCoupleSmallSlaps1} in Lemma \ref{lem:LPPCoupleSmallSlaps}, we can couple $\LLeps_{n,(1)}$ and $\LLeps_{n,(2)}$ with two full-space LPP sheets $\LL_{(j)}$ on $\SetT \cap (\mathcal{R}_n^{1,j} \times \mathcal{R}_n^{1,j} ) \cap \mathcal{D}_{n^{-1}}$ for $j=1,2$
such that for some constants $c_1,C_1>0$, all $n\in \N$, and all $\varepsilon>0$ small enough
\begin{equation}\label{eq:CouplingSmaller}
\mathbf{P}\left( \LLeps_{n,(1)} = \LL_{(1)} \text{ and } \LLeps_{n,(2)} = \LL_{(2)} \right) \geq 1-C_1\exp(-c_1n^{1/6}) . 
\end{equation}
From Theorem \ref{thm:LPPtoDL} for the convergence of full-space last passage percolation to the full-space directed landscape, we see that on the event in \eqref{eq:CouplingSmaller}, as $\varepsilon \rightarrow 0$,
\begin{equation}\label{eq:ConvergencePatched}
\LLeps_{n,(j)} \rightarrow (\mathcal{L}_{(i)}(x,s;y,t))_{(x,s;y,t) \in \SetT \cap (\mathcal{R}_n^{1,j} \times \mathcal{R}_n^{1,j} ) \cap \mathcal{D}_{n^{-1}}}
\end{equation} for $j=1,2$, where $\mathcal{L}_{(j)}$ are full-space directed landscapes. Since $\LLeps_{n,(1)}$ and $\LLeps_{n,(2)}$ agree on the intersection of the underlying rectangles for all $\varepsilon>0$, the same holds for the limiting objects, i.e., on the event in \eqref{eq:ConvergencePatched}, we get that with probability at least $1-C_2\exp(-c_2n^{1/6})$ for some constants $c_2,C_2>0$, 
\begin{equation}\label{eq:CoupledLimitObjects}
 \mathcal{L}_{(1)}(x,s;y,t)=\mathcal{L}_{(2)}(x,s;y,t) \text{ for all } (x,s;y,t) \in \SetT \cap \big(\mathcal{R}_n^{1,1}\cap \mathcal{R}_n^{1,2}\big)^2 \cap \mathcal{D}_{n^{-1}} .
\end{equation} Let $\tilde{\mathcal{A}}_n(N)$ be the intersections of the events in \eqref{eq:CouplingSmaller} and \eqref{eq:CoupledLimitObjects} and recall Definition~\ref{def:PatchedDLExt} for the patched directed landscape and the geodesics $ \pi^{\mathcal{R}}_{x,s;y,t}$ with respect to $\LLeps_{n,(1)}$ and $\LLeps_{n,(2)}$. Then on the event $\tilde{\mathcal{A}}_n(N)$, we see that for some constants $c_3,C_3>0$
\begin{equation}\label{eq:RestrictionEvent}
    \P\Big( \pi^{\mathcal{R}}_{x,s;y,t} \subseteq \mathcal{R} \, \forall (x,s),(y,t)\in \mathcal{R}^{i}_n \cap \mathcal{V}_{n} , i \in \lbr n \rbr\Big)  \geq 1 - C_3 \exp(-c_3 n^{-1/4})
\end{equation} by Proposition \ref{pro:ModerateDL}. 
 Let $\mathcal{A}_n(N)$ be the intersection of the event $\tilde{\mathcal{A}}_n(N)$ and the event in \eqref{eq:RestrictionEvent}. Then on $\mathcal{A}_n(N)$, 
 the random function defined as in Definition \ref{def:PatchedDLExt} with respect to $\mathcal{L}_{(1)}$ and $\mathcal{L}_{(2)}$ is a patched directed landscape on $(\mathbb{T} \times [0,n^{-1}])^2_{\uparrow}$, allowing us to conclude. 
\end{proof}

We have the following convergence of periodic last passage percolation, which formalizes the convergence to the periodic directed landscape stated in Theorem~\ref{thm:ConvergenceLPPMAIN}. 

\begin{theorem}\label{thm:PeriodicLPPConvergence}  
Recall the definition of the periodic LPP landscape $\LLpeps$ with respect to the periodic last passage percolation $(T_{u,v}^{\per})_{u,v\in \Z^2}$ from \eqref{def:periodicLPPSheet}. Let $\tilde{\omega}=2^{\frac{3}{2}}(1+\sqrt{\rho})^{-\frac{1}{2}}\rho^{\frac{1}{4}}$. Then 
\begin{equation}
 (\LLpeps(x,s\tilde{\omega}^{-1};y,t\tilde{\omega}^{-1}))_{(x,s;y,t) \in \SetTw} \rightarrow (\Lp(x,s;y,t))_{(x,s;y,t) \in \SetT}
\end{equation} in the sense of uniform convergence on compact sets,  where $\Lp$ is a periodic directed landscape.
\end{theorem}
\begin{proof}
Assume again that $\tilde{\omega}=1$ as the arguments are similar for general $\tilde{\omega}$. 
Recall that we set $\varepsilon=k^{-1/2}$ with $k=k(N)$ under assumption \eqref{eq:NonDegeneratedEnvironment}. Let $\E^{\per}_N$, $\E^{\per}_{N,n}$, $\bar{\E}_n$ and $\E^{\per}$ denote the expectation under the laws of $\LLpeps$, $\LLpaeps$, $\Lpa$ and $\Lp$, respectively. 
Then we have to verify that for any test function $f$ supported on some compact set $K \subseteq \SetT$, we get that the periodic LPP landscape  converges to the periodic directed landscape, i.e.,
\begin{equation}
\lim_{N \rightarrow \infty}\E^{\per}_{N}[f] = \E^{\per}[f] . 
\end{equation}
We claim that without loss of generality, we can assume that
$K=\mathcal{D}_{M^{-1}}$ for some $M\in \N$. To see this, note that the inverse slope 
$(x,s;y,t) \mapsto \sl^{-1}(x,s;y,t)=(t-s)/(x-y)$ is a continuous function on $\SetT$, and thus obtains its maximum $\sl_{\max}$ and minimum $\sl_{\min}$ on any compact set $K \subseteq \SetT$. Hence, we see that $K \subseteq \mathcal{D}_{M^{-1}}$ with $M=\max(|\sl_{\max}|,|\sl_{\min}|)$. Recall the coupling $\mathbf{P}_{N,n}$ from Lemma~\ref{lem:LPPtoDLRestricted} between $\LLpeps$ and $\LLpaeps$. Since $\lVert f \rVert_{\infty}<\infty$, we see that for all $\delta>0$, there exists some $N_0,M_0>0$ such that for all $N \geq N_0$, and all $n\geq \max(M,M_0)$, 
\begin{equation}\label{eq:ApproxLPP1}
\big| \E^{\per}_{N}[f] - \E^{\per}_{N,n}[f] \big| \leq \lVert f \rVert_{\infty} \mathbf{P}_{N,n}\left( \LLpaeps \neq \LLpaeps \text{ on } \mathcal{D}_{M^{-1}}\right) \leq \frac{\delta}{3} , 
\end{equation}
where we use Lemma \ref{lem:LPPtoDLRestricted} for the last step. 
Next, we get that for all $\delta>0$, there exists some $N_1 \in \N$ such that for all $N \geq N_1$, and $n \in \N$
\begin{equation}\label{eq:ApproxLPP2}
\Big| \E^{\per}_{N,n}[f] - \bar{\E}_n[f] \Big| \leq \frac{\delta}{3} . 
\end{equation}
This is due to the convergence of the restricted last passage times in the uniform on compact topology to the patched directed landscape by Lemma \ref{lem:LPPtoDLRestricted}. 
Finally, since $\Lp$ is the unique limit point of $(\Lpa)_{n\in \N}$ with respect to $\Wast$, we get that for all $\delta>0$, there exists some $M_1\in \N$ such that for all $n \geq M_1$, 
\begin{equation}
    \Wast(\Lpa,\Lp) \leq \frac{\delta}{3M \lVert f \rVert_{\infty}} . 
\end{equation}
Using the definition of $\Wast$, there exists  for all $n\geq M_1$ a coupling $\mathbf{P}_{n}$ such that 
\begin{equation}
\mathbf{P}_{n}\left( \Lpa = \Lp \text{ on } \mathcal{D}_{M^{-1}} \right) \geq 1- \frac{\delta}{3 \lVert f \rVert_{\infty}} . 
\end{equation}
In particular, for any test function $f$, any $\delta>0$, and $n\geq M_1$, we get that 
\begin{equation}\label{eq:ApproxLPP3}
\big|\bar{\E}_n[f] - \E^{\per}[f] \big| \leq \frac{\delta}{3} . 
\end{equation} 
Combining the approximation steps \eqref{eq:ApproxLPP1}, \eqref{eq:ApproxLPP2} and \eqref{eq:ApproxLPP3}, we see that 
\begin{equation}
  \big| \E^{\per}_{N}[f] - \E^{\per}[f] \big| \leq \big| \E^{\per}_{N}[f] - \E^{\per}_{N,n}[f] \big| + \Big| \E^{\per}_{N,n}[f] - \bar{\E}_n[f] \Big| +\big|\bar{\E}_n[f] - \E^{\per}[f] \big| \leq \delta
\end{equation}
for every $\delta>0$, provided that $n\geq \max(M,M_0,M_1)$ and $N \geq \max(N_0,N_1)$.
Since $\delta>0$ was arbitrary, we conclude. 
\end{proof}   

\begin{remark}\label{rem:GeneralPeriod} 
In the same way as above, we obtain the $p$-periodic directed landscape as a scaling limit under a different renormalization of the periodic last passage times. Moreover, as in Remark~\ref{rem:ExtendedPDL}, the convergence can be extended to arbitrary finite time horizons $[0,T]$. 
\end{remark}

We have now all tools in order to deduce the continuity of the periodic directed landscape.

\begin{proof}[Proof of Proposition~\ref{pro:PDLContinuity}] Recall from Theorem~\ref{thm:PeriodicLPPConvergence} the convergence of periodic last passage percolation to the periodic directed landscape, and let us stress that Theorem~\ref{thm:PeriodicLPPConvergence} is derived without using Proposition~\ref{pro:PDLContinuity} or Theorem~\ref{thm:PeriodicDirectedLandscape}, respectively. Let $B_{\delta}(z,\tau)$ denote for all $(z,\tau)\in \SetT$ the $\ell_2$-ball of radius $\delta>0$ around $(z,\tau)$. Due to Theorem~\ref{thm:PeriodicLPPConvergence}, it suffices for the continuity of $\Lp$ to show that for all $(x,s;y,t)\in (\mathbb{T} \times (0,1))^{2}_{\uparrow}$ and $\varepsilon'>0$, we find some $\delta>0$ such that 
\begin{equation}
 \P\left( \sup_{(x',s') \in B_{\delta}(x,s)}\sup_{(y',t') \in B_{\delta}(y,t)} \big| \LLpeps(x',s';y',t') -  \LLpeps(x,s;y,t) \big| \leq \varepsilon'\right) \geq 1- \varepsilon' 
\end{equation} for all $\varepsilon>0$ small enough.  
Using the triangle inequality and symmetry, it suffices to show that for any $\varepsilon'>0$, we find some $\delta>0$ such that for  all $\varepsilon >0$ small enough
\begin{equation}\label{eq:OneSidedConti}
    \P\left( \sup_{(x',s') \in B_{\delta}(x,s)} \big| \LLpeps(x',s';y,t) -  \LLpeps(x,s;y,t) \big| \leq \varepsilon' \right) \geq 1- \varepsilon' . 
\end{equation}
For all $(x,s;y,t)\in \SetT$, note that $\LLpeps(x,s;y,t)$ satisfies the metric composition law
\begin{equation}\label{eq:MetricLPPLandscape}
   \LLpeps(x,s;y,t) = \sup_{z \in \mathbb{T}} \big( \LLpeps(x,s;z,\tau) + \LLpeps(z,\tau;y,t)\big) , 
\end{equation}
for all $\tau \in [s,t]$, and we write
\begin{equation}
    \Sigma^{(\varepsilon)}_{x,s,y,t}(\tau) := \left\{ z\in \mathbb{T} \, \colon \, \LLpeps(x,s;z,\tau) + \LLpeps(z,\tau;y,t) =  \LLpeps(x,s;y,t)\right\}  
\end{equation} for the set of $z\in \mathbb{T}$ which achieve the supremum in \eqref{eq:MetricLPPLandscape}. Note that $\Sigma^{(\varepsilon)}_{x,s,y,t}(\tau) \neq \emptyset$ by the continuity of $\LLpeps$ and compactness.
By Lemma~\ref{lem:ConsistentLPPLocally}, we find some constant $C_{\ast}>0$ so that for all $\delta>0$, and all $\varepsilon>0$ small enough
\begin{equation}\label{eq:LocateExit}
    \P\left( \Sigma^{(\varepsilon)}_{x',s',y,t}(x+2\delta) \in \big[x-C_{\ast}\delta^{\frac{3}{5}},x+C_{\ast}\delta^{\frac{3}{5}}\big] \, \forall (x',s') \in B_{\delta}(x,s)\right) \geq 1- C_1 \exp(-c_1\delta^{-\frac{1}{24}}) , 
\end{equation} where $c_1,C_1>0$ do not depend on $\delta,\varepsilon>0$.
By Lemma~\ref{lem:LPPCoupleSmallSlaps}, we find a coupling $\mathbf{P}$ between periodic last passage percolation and full-space last passage percolation so that for $c_2,C_2>0$ and all $\delta>0$, and $\varepsilon=\varepsilon(\delta)>0$ sufficiently small, 
\begin{equation*}
\begin{split}
        \mathbf{P}&\left( \LLpeps(x',s';z',\tau)=\LLeps(x',s';z',\tau)  \, \forall (x',s';z',\tau) \in ([x-C_{\ast}\delta^{\frac{3}{5}},x-C_{\ast}\delta^{\frac{3}{5}}]\times[s-\delta,s+2\delta])^2_{\uparrow} \right) \\
        &\geq 1 - C_2 \exp(-c_2\delta^{-1/6}) . 
\end{split}
\end{equation*}
Hence, applying the estimates from Proposition~\ref{pro:DLvalues} (with $t-s$ of order $\delta$, $|x-y|$ of order at most $\delta^{3/5}$, and $z$ of order $\delta^{-1/9}$) stated the full-space directed landscape to the periodic directed landscape $\Lp$, we see that the events
\begin{equation*}
\begin{split}
        &A_{x,s;y,t}^{(\varepsilon),\delta} \\
        &:= \left\{ \sup_{z\in [x-C_{\ast}\delta^{3/5},x+C_{\ast}\delta^{3/5}]} \!\!\!\!\!\!\!\!\!\!\!\!\!\!\!\!\!\!|\LLpeps(x,s;z,s+2\delta) - \LLpeps(x',s';z,s+2\delta) | \leq \delta^{\frac{1}{6}} \,\, \forall (x',s') \in B_{\delta}(x,s)  \right\} 
\end{split}
\end{equation*}
satisfy for some $c_3,C_3>0$, all $\delta>0$, and $\varepsilon=\varepsilon(\delta)>0$ sufficiently small, 
\begin{equation}\label{eq:ModulusPeriodic}
     \P\left( A_{x,s;y,t}^{(\varepsilon),\delta}  \right) \geq 1- C_3 \exp(-c_3\delta^{-1/6}) .
\end{equation}
Let $f^{+}_{x,s;z,s+2\delta}$ denote the largest element of $\Sigma^{(\varepsilon)}_{x,s,y,t}(x+2\delta)$.
On the events in \eqref{eq:LocateExit} and \eqref{eq:ModulusPeriodic},
\begin{equation*}
\begin{split}
     \LLpeps(x,s;y,t) &= \LLpeps(x,s;f^{+}_{x,s;z,s+2\delta},s+2\delta) + \LLpeps(f^{+}_{x,s;z,s+2\delta},s+2\delta;y,t)  \\
     &\leq  \LLpeps(x',s';f^{+}_{x,s;z,s+2\delta},s+2\delta) + \LLpeps(f^{+}_{x,s;z,s+2\delta},s+2\delta;y,t) + \delta^{\frac{1}{6}} \\
     &\leq \LLpeps(x',s';y,t) + \delta^{\frac{1}{6}}
\end{split}
\end{equation*} for all $(x',s') \in B_{\delta}(x,s)$, using the metric composition \eqref{eq:MetricLPPLandscape} for the first equality and the last inequality.
Swapping the roles of $(x',s')$ and $(x,s)$, on the events in \eqref{eq:LocateExit} and \eqref{eq:ModulusPeriodic},
\begin{equation*}
     \LLpeps(x',s';y,t) \leq \LLpeps(x,s;y,t) + \delta^{\frac{1}{6}}
\end{equation*} holds for all $(x',s') \in B_{\delta}(x,s)$. In total, we obtain \eqref{eq:OneSidedConti}, allowing us to conclude.
\end{proof}

\section{The periodic KPZ fixed point coupled to the periodic directed landscape}\label{sec:KPZfixed}

The full-space KPZ fixed point is a Markov process that was introduced in  \cite{MQR:kpz} as a limit of TASEP on the line, and is described by exact formulas for its transition probabilities. It was subsequently identified in Corollary  4.2 of \cite{NQR:Fixed} with certain variational problem marginals of the full-space directed landscape. That formulation has several uses -- it allows one to prove convergence to the KPZ fixed point using some of the softer (less formula-based) approaches relevant to proving directed landscape convergence (e.g. \cite{W:Landscape,ACH:ScalingASEP,ACH:Fixed}); leads to soft proofs of natural properties of the fixed point; and provides a coupling of KPZ fixed points over all initial data which arises as the limit of natural couplings.

Recently, \cite{BLL:FixedPointGeneral} proved that the periodic TASEP admits a scaling limit at the level of consistent families of formulas for its limiting finite dimensional distributions in space and time. In parallel with the full-space, they called this limiting stochastic process the \emph{periodic KPZ fixed point}. The existence of a continuous version of this process, and its identification as a Markov process is not clear from the finite dimensional distributions (unlike with the transition probability formulas in \cite{MQR:kpz}). There are other natural properties of this process that are not evident from the formulas, and are conjectured in \cite{BLL:FixedPointGeneral} and earlier work such as \cite{L:Multi}.

In this section we first (in Section \ref{sec:PeriodicKPZfixed}) identify the periodic KPZ fixed point with certain variational problem marginals of our newly constructed periodic directed landscape, and then use this (in Section \ref{sec:localstruct}) to prove the above mentioned basic and conjectured properties. In particular, we resolve Conjecture (1.20) in \cite{BLL:FixedPointGeneral} on a variational characterization of the periodic KPZ fixed point in terms of the periodic directed landscape (see Theorem~\ref{thm:PeriodicKPZFixed}), Conjecture (1.16) in \cite{BLL:FixedPointGeneral} on the local structure of the periodic KPZ fixed point with narrow wedge initial data (see Proposition~\ref{pro:LargePeriodLimit})  and  Conjecture $(5)$ in~\cite{BL:PinchedUp} on the approximation of the periodic KPZ fixed point by the full space KPZ fixed when the period length diverges (see Corollary~\ref{cor:KPZfixedWedge}). Note, the equation numbers in \cite{BLL:FixedPointGeneral} are based on the first arXiv posting of that paper.

This first result (identification of periodic KPZ fixed point with periodic directed landscape) is given here as Theorem~\ref{thm:PeriodicKPZFixed}  and proved by embedding TASEP as a variational problem in the periodic exponential LPP model (see Remark~\ref{rem:TASEPandLPP}), and using the convergence of that model to the periodic directed landscape (see Theorem~\ref{thm:PeriodicLPPConvergence}). LPP yields one coupling of TASEPs with different initial data. Another, different coupling is via the basic coupling and the joint convergence under that coupling is shown in Section \ref{sec:ASEPconvergence}, not just for TASEP but for ASEP.

\subsection{Variational characterization of the periodic KPZ fixed point}\label{sec:PeriodicKPZfixed} 
We start by briefly recalling the definition of the full-space KPZ fixed point as this will be useful both by analogy and for specific comparison results to the periodic setting. Then, we summarize the main result of \cite{BLL:FixedPointGeneral} and provide, as Theorem \ref{thm:PeriodicKPZFixed}, the variational characterization of the periodic KPZ fixed point via the periodic directed landscape. 
Though the proper approach would be to define the full-space KPZ fixed point via its transition probabilities and as a TASEP scaling limit, and then to quote the following variational formula as a theorem, we will be satisfied here with simply defining the fixed point via the directed landscape. In the periodic setting we proceed  properly.

\begin{definition}\label{def:KPZfixedFullSpace} Let $\UC$ denote the space of upper semi-continuous functions $f \colon \R \rightarrow \R \cup \{ -\infty\}$ such that for some $x\in \R$, $f(x)\neq -\infty$ and such that there exists some constants $m,C>0$ so that for all $x\in \R$
$f(x) \leq m |x| + C$. For any $\mathfrak{h}_0 \in \UC$ and $\mathcal{L}$ an independent (of $\mathfrak{h}_0$) full-space directed landscape, we define the \textbf{full-space KPZ fixed point} with initial data $\mathfrak{h}_0$ as the random function $(y,t) \mapsto \hfix(\mathfrak{h}_0;y,t)$, where
\begin{equation}\label{eq:VariationFixed}
\hfix(\mathfrak{h}_0;y,t) := \sup_{x\in \R} \big(\mathfrak{h}_0(x) + \mathcal{L}(x,0;y,t)\big).
\end{equation} For multiple initial data this provides a coupling of full-space KPZ fixed points.
\end{definition}

Turning to the periodic setting, we start by recalling the main result of \cite{BLL:FixedPointGeneral} which defines the periodic KPZ fixed point as the limit of periodic TASEP and identifies its finite dimensional distributions via explicit formulas. 

As in Definition 1.1 of \cite{BLL:FixedPointGeneral}, we let $\UCp$ denote the space of $1$-periodic upper semi-continuous functions $f \colon \mathbb{T} \rightarrow \R \cup \{ -\infty\}$ such that $f(y) \neq -\infty$ for some $y \in \mathbb{T}$, and we equip this space with the topology of local Hausdorff convergence on the hypographs. 
Recall the height function $\hpTASEP$ of the periodic TASEP from Definition \ref{def:HeightTASEP} (along with some other notation like $\Lambda_{N,k}$ used below). Note, it is defined there already in terms of the coupling to periodic exponential LPP though noted in Remark \ref{rem:TASEPandLPP} that this is equivalent to the more standard particle process description of the model.

\begin{theorem}[{\cite[Theorem 1.2 and Theorem 1.5]{BLL:FixedPointGeneral}}]\label{thm:BaikPKPZ}
Suppose that $\lambda_N \in \Lambda_{N,k}$ with $k=\lfloor \theta N\rfloor$ for some $\theta \in (0,1)$ is a family of TASEP initial conditions so that for some $\hfixp_0 \in \UCp$
\begin{equation}\label{eq:ConvergenceIniPeriodic}
  \frac{-\hpTASEP(\lambda_N;\alpha N,0) + \theta\alpha N}{-\sqrt{2\theta(1-\theta)}N^{1/2}}   \rightarrow \hfixp_0(\alpha)
\end{equation} as $N\to \infty$ where the convergence is as a process in $\alpha$ in the space $\UCp$. Then for all $M \in \N$, 
and all $\tau=(\tau_i)_{i \in \lbr M \rbr}$, $\alpha=(\alpha_i)_{i \in \lbr M \rbr}$ and $\beta=(\beta_i)_{i \in \lbr M \rbr}$, 
 \begin{equation}\label{eq:BaikLiuCharacterization}
 \lim_{N \rightarrow \infty} \P\left( \bigcap_{i=1}^{M} \left\{ \frac{-\hpTASEP(\lambda_N;y_i,t_i)+\theta y_i-\theta(1-\theta)t_i}{ -\sqrt{2\theta(1-\theta)}N^{1/2}}\leq \beta_i \right\} \right) = F_{\hfix_0}^M(\beta,\alpha,\tau) , 
 \end{equation}
 where we set
 \begin{equation}\label{eq:PointShift}
     t_i := \frac{\tau_i N^{3/2}}{\sqrt{\theta(1-\theta)}}, \quad y_i := \alpha_i N + (1-2\theta)t_i , 
 \end{equation}
 and $(F_{\hfix_0}^M)$ has an explicit formula given in Definition~2.1 of \cite{BLL:FixedPointGeneral}. Moreover, the family $(F_{\hfix_0}^M)$ defines a random field $(\alpha,\tau) \mapsto \mathcal{H}^{\per}_1(\hfixp_0;\alpha,\tau)$, called the \textbf{periodic KPZ fixed point}. 
\end{theorem}

 \begin{remark}\label{rem:FactorTwo}
To match our Definition~\ref{def:HeightTASEP} of the height function $\hpTASEP(\lambda_N;\alpha N,0)$ to the definition of the height function $H^{\textup{pTASEP}}$ used in \cite{BLL:FixedPointGeneral}, we have to apply the substitution $-h^{\textup{pTASEP}}(\lambda_N;0,y) = 2 H^{\textup{pTASEP}}(y,0) + y$ for all $y\in \Z$. Moreover, let us point out that due to the change of the spatial coordinate by a factor of $2$ later on in Theorem \ref{thm:PeriodicKPZFixed}, we impose in the normalization in \eqref{eq:ConvergenceIniPeriodic} and \eqref{eq:BaikLiuCharacterization} an additional factor of $\sqrt{2}$. 
 \end{remark}
Note that in a similar way, one can obtain the $p$-periodic KPZ fixed point $\mathcal{H}^{\per}_p$ on the space of $p$-periodic upper semi-continuous functions with some $p>0$ by the scaling relation 
 \begin{equation*}
 \mathcal{H}^{\per}_p(\gamma,\tau) := p^{1/2}  \mathcal{H}^{\per}_1(p^{-1}\gamma,p^{-3/2}\tau) . 
 \end{equation*}

We now define the candidate variational formulation of the periodic KPZ fixed point (analogously to Definition~\ref{def:KPZfixedFullSpace} in the full-space setting). Its identification with the periodic KPZ fixed point is given as Theorem \ref{thm:PeriodicKPZFixed}. Recall, from earlier, the periodic directed landscape (see Proposition~\ref{pro:CauchySequence} and Remark~\ref{rem:ExtendedPDL}).

\begin{definition}\label{def:PeriodicKPZfixed}
    Let $\hfixp_0 \in \UCp$ and $\Lp$ be a periodic directed landscape on $(\mathbb{T} \times \R)^{2}_{\uparrow}$. Then define the random function $(y,t) \mapsto \hfixp(\hfixp_0;y,t)$ for all $y\in \mathbb{T}$ and $t >0$ by
\begin{equation}\label{eq:VariationFixedPeriodic}
 \hfixp\big(\hfixp_0;y,t\big) := \sup_{x\in \mathbb{T}} \big( \hfixp_0(x) + \Lp(x,0;y,t)\big).
\end{equation}
\end{definition}

The next statement shows that the above limit is indeed a variational characterization of the periodic KPZ fixed point. 

\begin{theorem}\label{thm:PeriodicKPZFixed} 
The function $(y,t) \mapsto \hfixp(\hfixp_0;y,t)$ in Definition~\ref{def:PeriodicKPZfixed} is (up to scaling the spatial coordinate) a periodic KPZ fixed point from Theorem ~\ref{thm:BaikPKPZ}  with initial data $\hfixp_0 \in \UCp$, i.e. 
\begin{equation}
      \hfixp\big(\hfixp_0;y,\tau\big) \overset{(d)}{=} \mathcal{H}^{\per}_2(\hfixpt_0;2y,\tau)
\end{equation}
jointly over $y\in \mathbb{T}$ and $\tau\geq 0$.
Here, we set $\hfixpt_0(x) := \hfixp_0(2x)$ for all $x\in \mathbb{T}$. 
\end{theorem}

\begin{remark}
Theorem~\ref{thm:PeriodicKPZFixed} resolves (1.20) in \cite{BLL:FixedPointGeneral}, where a variational characterization of the periodic KPZ fixed point in terms of the periodic directed landscape was conjectured. Let us emphasize that as noted in the footnote at page 5 of \cite{BLL:FixedPointGeneral}, the above descriptions of the periodic KPZ point differ in the spatial coordinate by a factor of $2$. This is consistent with Theorem \ref{thm:PeriodicDirectedLandscape} stating that the periodic direct landscape can be locally approximated by a full-space directed landscape, and equation (1.16) in \cite{BLL:FixedPointGeneral} where the local convergence of the periodic KPZ fixed point to the full-space KPZ fixed point is conjectured.
\end{remark} 
\begin{remark}
    In view of Theorem~\ref{thm:PeriodicKPZFixed}, we see that the variation characterization of the periodic KPZ fixed point from Definition~\ref{def:PeriodicKPZfixed} together with property \hyperref[eq:Cond4]{(iii)} from Theorem~\ref{thm:PeriodicDirectedLandscape} readily implies that the periodic KPZ fixed point satisfies the Markov property. Moreover, from the continuity of the periodic directed landscape in Proposition~\ref{pro:PDLContinuity}, we get that the periodic KPZ fixed point is supported on the space of continuous functions. 
\end{remark}

As a first step towards proving Theorem~\ref{thm:PeriodicKPZFixed}, we use Theorem~\ref{thm:PeriodicLPPConvergence} in order to show convergence of periodic exponential LPP to the limit in \eqref{eq:VariationFixedPeriodic}, i.e. the convergence of periodic last passage percolation to the periodic KPZ fixed point. For given $\hfixp_0 \in \UCp$, we assign the values 
\begin{equation}
    h_0^{N}(u_1,-u_2) := \hfixp_0(u_2/(N-k))
\end{equation} for all $(u_1,u_2)\in \lambda_N^{\ast}$. 
Consider the tilted height function $\hLpt$ 
\begin{equation}
    \hLpt_N(\hfixp_0;y,t) := \hLp_{\lambda_N^{\ast}}(h_0^{N};\mathbf{u}y+\mathbf{v}t)
\end{equation} for all $y\in \R$ and $t\geq 0$, where we recall the height functions $\hLp_{\lambda_N^{\ast}}$ of $(N,k)$-periodic last passage percolation from Definition~\ref{def:HeightLPP}. Note that
\begin{equation}
    \hLpt_{N}(\hfixp_0;y+k,t) = \hLpt_{N}(\hfixp_0;y,t)
\end{equation} for all $y\in \R$ and $t\geq 0$.  Recall the functions $M=M(\rho,\varepsilon,s,t)$ from \eqref{def:MeanLPP} and $V=V(\rho,\varepsilon)$ from \eqref{def:VarianceLPP}, as well as that $\tilde{\omega}=2^{\frac{3}{2}}(1+\sqrt{\rho})^{-\frac{1}{2}}\rho^{\frac{1}{4}}$.
Let  $\varepsilon=k^{-1/2}$ for all $k=k(N) \in \N$, and define the \textbf{rescaled periodic LPP height function}
\begin{equation}\label{def:RescaledHeightLPP}
\hLpeps(\hfixp_0;x,t) := \frac{1}{V(\rho,\varepsilon)} \left( \hLp_{N}(\hfixp_0;\varepsilon^{-2}x,\varepsilon^{-3}t) -M(\rho,\varepsilon,t)\right)  
\end{equation}
with linear interpolation if $(\varepsilon^{-2}x,\gamma\varepsilon^{-3}t) \notin \Z^2$. 

\begin{proposition}\label{pro:LPPgeneralData} Fix $T>0$. Assume that \eqref{eq:NonDegeneratedEnvironment} holds
 and let $\hfixp_0 \in \UCp$. Then the limit
\begin{equation}
\tilde{\mathfrak{h}}\big(\hfixp_0;y,t\big) := \lim_{\varepsilon \rightarrow 0} \hLpeps(\hfixp_0;y,t\tilde{\omega}^{-1})
\end{equation} exists jointly in $y\in \mathbb{T}$ and $t\in [0,T]$, where $(y,t) \mapsto \tilde{\mathfrak{h}}\big(\hfixp_0;y,t\big)$ satisfies \eqref{eq:VariationFixedPeriodic}.
\end{proposition}
\begin{proof}
Recall from Theorem~\ref{thm:PeriodicLPPConvergence} that the periodic LPP landscape (defined in \eqref{def:periodicLPPSheet}) converges to the periodic directed landscape under assumption~\eqref{eq:NonDegeneratedEnvironment}. 
From Definition~\ref{def:HeightLPP} of the periodic LPP height function in terms of exponential last passage percolation, we see that the rescaled periodic LPP height function from \eqref{def:RescaledHeightLPP} with  initial conditions
\begin{equation}\label{def:DeltaInitialConditions}
    \hfixp_{(x)}(z) := \begin{cases} 0 & \text{ if } z=x + j \text{ for some $j\in \Z$} , \\
    - \infty & \text{ otherwise } z \neq x, 
    \end{cases}
\end{equation} must satisfy for any $x\in \R$, 
\begin{equation}\label{eq:ConvergeToPDL}
   \Lp(x,0;y,t) = \lim_{\varepsilon \rightarrow 0} \hLpeps(\hfixp_{(x)};y,t\tilde{\omega}^{-1}) . 
\end{equation} For general initial data $\hfixp_0 \in \UCp$, the variational characterization for the right-hand side of \eqref{eq:ConvergeToPDL} follows by noting that
\begin{equation*}
\begin{split}
    \lim_{\varepsilon \rightarrow 0} \hLpeps(\hfixp;y,t) &= \lim_{\varepsilon \rightarrow 0} \sup_{z \in \R} \big(\hfixp(z) +\hLpeps(\hfixp_{(z)};y,t) \big) \\ 
    &= \sup_{z \in \R}  \lim_{\varepsilon \rightarrow 0}\big(\hfixp(z) +\hLpeps(\hfixp_{(z)};y,t) \big) \\
    & = \sup_{z \in \R} \big(\hfixp(z) + \Lp(z,0;y,t) \big) , 
    \end{split}
\end{equation*}
 using the variational characterization of periodic LPP in \eqref{eq:VariationalDiscretePeriodic}, as well as \eqref{eq:ConvergeToPDL} and the continuity of $\Lp$ from Proposition \ref{pro:PDLContinuity} to exchange the limit and taking the supremum.
\end{proof}

\begin{figure}
\begin{tikzpicture}[scale=0.37]

\fill[gray!20] 	 (1,8) -- ++(2,0) -- ++(0,-1) -- ++(1,0) -- ++(0,-1) -- ++(1,0) -- ++(0,-1) -- ++(2,0)-- ++(0,-1) -- ++(1,0)-- ++(0,-1) -- ++(1,0)-- ++(0,-1) -- ++(2,0)-- ++(0,-1) -- ++(1,0)--++(-11,0)--++(0,8);

\foreach \x in{1,...,10}{
	\draw[gray!50,thin](1,\x) to (18,\x);   }

\foreach \x in{1,...,18}{
	\draw[gray!50,thin](\x,1) to (\x,10);   }
	
\node[shape=circle,scale=0.9,draw] (H0) at (0*1.2,-1){} ;

\foreach \x in{1,...,16}{
 	\node[shape=circle,scale=0.9,draw] (H\x) at (\x*1.2,-1){} ;  }

    \node at (3*1.2,-2) (l1){$1$} ;
    \node at (9*1.2,-2) (lN){$N$} ;

    \draw(-0.7*1.2,-1) to (H0);
	\draw(H0) to (H1);  
	\draw(H1) to (H2);  
	\draw(H2) to (H3);  
	\draw(H3) to (H4);  
	\draw(H4) to (H5);  
	\draw(H5) to (H6);  
	\draw(H6) to (H7);  
	\draw(H7) to (H8);  
	\draw(H8) to (H9);  
	\draw(H9) to (H10);
	\draw(H10) to (H11);  
	\draw(H11) to (H12);  
	\draw(H12) to (H13);  
	\draw(H13) to (H14);  
	\draw(H14) to (H15);				
	\draw(H15) to (H16);  
	\draw(H16) to (16.7*1.2,-1);	
	
	\draw[densely dotted, line width=1.5pt](2.5*1.2,-0.5) to (2.5*1.2,-1.5);  
	\draw[densely dotted, line width=1.5pt](9.5*1.2,-0.5) to (9.5*1.2,-1.5);   
	\draw[densely dotted, line width=1.5pt](16.5*1.2,-0.5) to (16.5*1.2,-1.5);

\draw[red,line width =2pt] 	 (1,8) -- ++(2,0) -- ++(0,-1) -- ++(1,0) -- ++(0,-1) -- ++(1,0) -- ++(0,-1) -- ++(2,0)-- ++(0,-1) -- ++(1,0)-- ++(0,-1) -- ++(1,0)-- ++(0,-1) -- ++(2,0)-- ++(0,-1) -- ++(1,0);


\node[shape=circle,fill=red,scale=0.8] at (H2){} ;

\node[shape=circle,fill=red,scale=0.8] at (H4){} ;
\node[shape=circle,fill=red,scale=0.8] at (H6){} ;
\node[shape=circle,fill=red,scale=0.8] at (H9){} ;

\node[shape=circle,fill=red,scale=0.8] at (H11){} ;
\node[shape=circle,fill=red,scale=0.8] at (H13){} ;
\node[shape=circle,fill=red,scale=0.8] at (H16){} ;

\fill[purple!10] (24+1,8) -- ++(2,0) -- ++(0,-1) -- ++(1,0) -- ++(0,-1) -- ++(1,0) -- ++(0,-1) -- ++(2,0)-- ++(0,-1) -- ++(1,0)-- ++(0,-1) -- ++(1,0)-- ++(0,-1) -- ++(2,0)-- ++(0,-1) -- ++(6,0) -- ++(0,2) -- ++(-2,0) -- ++(0,1) -- ++(-2,0) -- ++(0,2)-- ++(-2,0) -- ++(0,1) -- ++(-2,0) -- ++(0,2) -- ++(-2,0) -- ++(0,1) -- ++(-6,0) -- ++(0,-2) ;

\fill[gray!20] 	 (1+24,8) -- ++(2,0) -- ++(0,-1) -- ++(1,0) -- ++(0,-1) -- ++(1,0) -- ++(0,-1) -- ++(2,0)-- ++(0,-1) -- ++(1,0)-- ++(0,-1) -- ++(1,0)-- ++(0,-1) -- ++(2,0)-- ++(0,-1) -- ++(1,0)--++(-11,0)--++(0,8);

\foreach \x in{1,...,10}{
	\draw[gray!50,thin](1+24,\x) to (18+24,\x);   }

\foreach \x in{1,...,18}{
	\draw[gray!50,thin](\x+24,1) to (\x+24,10);   }
	
\node[shape=circle,scale=0.9,draw] (I0) at (24+0*1.2,-1){} ;

\foreach \x in{1,...,16}{
 	\node[shape=circle,scale=0.9,draw] (I\x) at (24+\x*1.2,-1){} ;  }

    \node at (24+3*1.2,-2) (j1){$1$} ;
    \node at (24+9*1.2,-2) (jN){$N$} ;

     \draw(24-0.7*1.2,-1) to (I0);
	 \draw(I0) to (I1);  
	 \draw(I1) to (I2);  
	 \draw(I2) to (I3);  
	 \draw(I3) to (I4);  
	\draw(I4) to (I5);  
	\draw(I5) to (I6);  
	\draw(I6) to (I7);  
	\draw(I7) to (I8);  
	\draw(I8) to (I9);  
	\draw(I9) to (I10);
	\draw(I10) to (I11);  
	\draw(I11) to (I12);  
	\draw(I12) to (I13);  
	\draw(I13) to (I14);  
	\draw(I14) to (I15);				
	\draw(I15) to (I16);  
	\draw(I16) to (24+16.7*1.2,-1);	
	
	\draw[densely dotted, line width=1.5pt](24+2.5*1.2,-0.5) to (24+2.5*1.2,-1.5);  
	\draw[densely dotted, line width=1.5pt](24+9.5*1.2,-0.5) to (24+9.5*1.2,-1.5);   
	\draw[densely dotted, line width=1.5pt](24+16.5*1.2,-0.5) to (24+16.5*1.2,-1.5);

\draw[black!80,line width =1pt] 	 (24+1,8) -- ++(2,0) -- ++(0,-1) -- ++(1,0) -- ++(0,-1) -- ++(1,0) -- ++(0,-1) -- ++(2,0)-- ++(0,-1) -- ++(1,0)-- ++(0,-1) -- ++(1,0)-- ++(0,-1) -- ++(2,0)-- ++(0,-1) -- ++(1,0);

\draw[red,line width =2pt] 	 (24+5,8+2) -- ++(2,0) -- ++(0,-1) -- ++(2,0) -- ++(0,-2) -- ++(2,0) -- ++(0,-1) -- ++(2,0) -- ++(0,-2)-- ++(2,0) -- ++(0,-1) -- ++(2,0) -- ++(0,-2);

\node[shape=circle,fill=red,scale=0.8] at (I0){} ;
\node[shape=circle,fill=red,scale=0.8] at (I1){} ;

\node[shape=circle,fill=red,scale=0.8] at (I4){} ;
\node[shape=circle,fill=red,scale=0.8] at (I7){} ;
\node[shape=circle,fill=red,scale=0.8] at (I8){} ;

\node[shape=circle,fill=red,scale=0.8] at (I11){} ;
\node[shape=circle,fill=red,scale=0.8] at (I14){} ;
\node[shape=circle,fill=red,scale=0.8] at (I15){} ;

\end{tikzpicture}
\caption{\label{fig:LPPPicture} Visualization of the correspondence between the periodic TASEP and periodic last passage percolation at time $0$ (left) and at some time $T>0$ (right). At the right-hand side, all shaded sites below the red curve have a last passage time of at most $T$ to the black curve, while all sites above the red curve have a last passage time larger than $T$ to the black curve. The red curve's vertical / horizontal increments map to particles / holes depicted below.}
\end{figure}
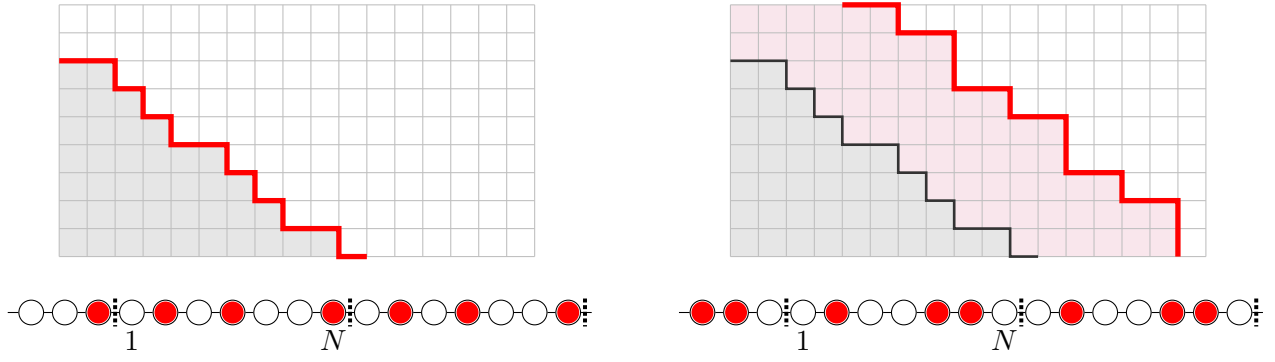

In preparation for Theorem~\ref{thm:PeriodicKPZFixed}, 
recall $\rho=k^{2}(N-k)^{-2}$ from \eqref{eq:EpsilonScale}, and the vectors $ \mathbf{u}:=(\rho^{1/2},-1)$ and $\mathbf{v}:=(\rho,1)$ from \eqref{def:VectorsDL}. In the proof, we will rely on the correspondence between periodic last passage percolation and the periodic TASEP from Definition~\ref{def:HeightTASEP} (see Figure~\ref{fig:LPPPicture}). More precisely, this correspondence allows us to use Theorem~\ref{thm:LPPtoDL} on the convergence of the periodic last passage percolation to the directed landscape to approximate the level sets of last passage percolation for given initial data up to order $o(\varepsilon^{-1})$ error terms. We then convert this into the convergence of the height function of the periodic TASEP to the variational characterization of the periodic KPZ fixed point from Definition~\ref{def:PeriodicKPZfixed}.

\begin{proof}[Proof of Theorem \ref{thm:PeriodicKPZFixed}]
Consider for all $u_{\varepsilon}, v_{\varepsilon}\in \Z^2$ the points $(x^{u}_{\varepsilon},s^{u}_{\varepsilon})$ and $(y^{v}_{\varepsilon},t^{v}_{\varepsilon})$ such that
\begin{equation}
    u_{\varepsilon} = x^{u}_{\varepsilon}\varepsilon^{-2}\mathbf{u} + s^{u}_{\varepsilon} \varepsilon^{-3} \mathbf{v} , \qquad v_{\varepsilon} = y^{v}_{\varepsilon}\varepsilon^{-2}\mathbf{u} + t^{v}_{\varepsilon} \varepsilon^{-3} \mathbf{v}  .
\end{equation} 
 Theorem~\ref{thm:PeriodicLPPConvergence} and the continuity of the limit $\Lp$ from Proposition~\ref{pro:PDLContinuity} ensure that 
for any fixed $T>0$ and $\varepsilon'>0$, we find some $\delta'>0$ so that 
\begin{equation}\label{eq:CloseApprox}
 \varepsilon^{1/2}\Big| \Tper_{u_{\varepsilon},v_{\varepsilon}} - M(\rho,\varepsilon,s_{\varepsilon},t_{\varepsilon}) - V(\rho,\varepsilon)\Lp(x,s\tilde{\omega};y,t\tilde{\omega}) \big| \leq \varepsilon'
\end{equation} for all $(u_{\varepsilon},v_{\varepsilon}) \in \Z^2$ with
\begin{equation}
  \varepsilon^{3/2}  \lVert u_{\varepsilon} - (x\varepsilon^{-2}\mathbf{u} + s\varepsilon^{-3} \mathbf{v}) \rVert_2  \leq \delta', \qquad \textrm{and}\qquad
\varepsilon^{3/2}  \lVert v_{\varepsilon} - (y\varepsilon^{-2}\mathbf{u} + t\varepsilon^{-3} \mathbf{v}) \rVert_2  \leq \delta' 
\end{equation} with respect to some $(x,s;y,t)\in (\mathbb{T}\times[-T,T])^{2}_{\uparrow}$, and $\varepsilon>0$ small enough. As in the proof of Theorem \ref{thm:PeriodicLPPConvergence}, let $\varepsilon=k^{-1/2}$ for $k=k(N)$. 
Recall from Definition \ref{def:HeightTASEP} that we have 
\begin{equation*}
  \hpTASEP(\lambda_N;Y+1,t) - \hpTASEP(\lambda_N;Y,t)=-\eta^{\per}_t(Y+1).
\end{equation*} for all $Y \in \Z$, and $t\geq 0$, where $\eta^{\per}_t(Y)$ denotes the particle configuration of a periodic TASEP at time $t\geq 0$ and location $Y \textup{ mod } N$. 
For step initial data at $X$ with
\begin{equation*}
 \hpTASEP(\lambda_N;X,0) =-S,
\end{equation*}
we consider the corresponding down-right path $\lambda^{(X+S,S)}_{\per,N}$ passing through $(X+S,S)$. Similar, set
$Z=-\hpTASEP(\lambda_N;Y,t)$. As pointed out in Remark~\ref{rem:TASEPandLPP}, the coupling between the periodic TASEP and periodic last passage percolation ensures that for all $X,S,Y \in \Z$ and $T\geq 0$, 
\begin{equation*}
    h^\per_{X,S}(Y,T): =\hpTASEP\big(\lambda^{(X+S,S)}_{\per,N};Y,T\big) =-Z 
\end{equation*}
if and only if 
\begin{equation*}
    \Tper_{(X+S,S),(Y+Z,Z)}\leq T , 
\quad \text{and} \quad
    \Tper_{(X+S,S),(Y+Z+1,Z+1)}>T.
\end{equation*}
In the following, we define
\begin{equation*}
    T_\varepsilon
    :=
    \frac{\tau}{\theta^2\sqrt{(1-\theta)}}\varepsilon^{-3},
    \qquad
    Y_\varepsilon
    :=
    \frac{y}{\theta}\varepsilon^{-2}
    +
    (1-2\theta)T_\varepsilon.
\end{equation*}
for fixed $\tau>0$ and $y\in \R$. Moreover, for fixed $\alpha \in \R$, we set $X_\varepsilon:=\frac{\alpha}{\theta}\varepsilon^{-2}$. Then by assumption \eqref{eq:ConvergenceIniPeriodic}, recalling that $k=\theta N$, we can find for all $\varepsilon>0$ some  
\begin{equation*}
  S_\varepsilon= 
   - \alpha\varepsilon^{-2}
    -
    2^{-1/2}\sqrt{1-\theta}\,
    \hfixp_0(\alpha)\varepsilon^{-1}
    +
    o(\varepsilon^{-1})
\end{equation*} 
so that $-  h^\per_{X_\varepsilon,S_\varepsilon}(X_\varepsilon,0) = S_{\varepsilon}$.
A direct conversion of the point $(X_\varepsilon+S_\varepsilon,S_\varepsilon)$ into the basis $(\mathbf u,\mathbf v)$, recalling $\rho = (1-\theta)^2\theta^{-2}$, gives
\begin{equation}\label{eq:CorrectStartDecomp}
    (X_\varepsilon+S_\varepsilon,S_\varepsilon)
    =
    \left(
    \alpha\varepsilon^{-2}
    +
    \frac{1-2\theta}{\sqrt{2(1-\theta)}}
    \hfixp_0(\alpha)\varepsilon^{-1}
    \right)\mathbf u
    -
    \frac{\theta}{\sqrt{2(1-\theta)}}
    \hfixp_0(\alpha)\varepsilon^{-1}
    \mathbf v
    +
    o(\varepsilon^{-1}),
\end{equation}
where the error term $o(\varepsilon^{-1})$ is with respect to both components. Now fix $r\in\mathbb R$ and define
\begin{equation*}
    Z_\varepsilon(r)
    :=
    -\theta Y_\varepsilon
    +
    \theta(1-\theta)T_\varepsilon
    -
    2^{-1/2}\sqrt{1-\theta}\,
    r\varepsilon^{-1}.
\end{equation*}
Converting the point $(Y_\varepsilon+Z_\varepsilon(r),Z_\varepsilon(r))$ into the basis $\mathbf{u}$ and $\mathbf{v}$ gives
\begin{equation*}
\begin{split}
    (Y_\varepsilon+Z_\varepsilon(r),Z_\varepsilon(r))
    &=
    \left(
    y\varepsilon^{-2}
    +
    \frac{1-2\theta}{\sqrt{2(1-\theta)}}
    r\varepsilon^{-1}
    \right)\mathbf u 
    +
    \left(
    \theta^2T_\varepsilon
    -
    \frac{\theta}{\sqrt{2(1-\theta)}}
    r\varepsilon^{-1}
    \right)\mathbf v . 
\end{split}
\end{equation*}
Therefore the difference between the $\mathbf v$-coordinates of
$(Y_\varepsilon+Z_\varepsilon(r),Z_\varepsilon(r))$ and $(X_\varepsilon+S_\varepsilon,S_\varepsilon)$ is
\begin{equation*}
    \theta^2T_\varepsilon
    -
    \frac{\theta}{\sqrt{2(1-\theta)}}
    \bigl(r-\hfixp_0(\alpha)\bigr)\varepsilon^{-1}
    +
    o(\varepsilon^{-1}) . 
\end{equation*}
Recall $M$ and $V$ from \eqref{def:MeanLPP} and \eqref{def:VarianceLPP}, the quantity $\tilde{\omega}=2^{3/2}(1-\theta)^{1/2}$, and that $\varepsilon=\theta^{-1/2} N^{-1/2}$. Using that by definition of $T_{\varepsilon}$
\begin{equation*}
   \tilde{\omega} \theta^2 T_{\varepsilon} =  2^{\frac{3}{2}} \tau \varepsilon^{-3} , 
\end{equation*}
the approximation in \eqref{eq:CloseApprox} yields 
\begin{equation*}
\Tper_{(X_\varepsilon+S_\varepsilon,S_\varepsilon),(Y_\varepsilon+Z_\varepsilon(r),Z_\varepsilon(r))}
    = T_\varepsilon  -
    \frac{1}{\theta\sqrt{2(1-\theta)}}
    \left(
    r-\hfixp_0(\alpha)-\Lp(\alpha,0;y,2^{3/2}\tau)
    \right)\varepsilon^{-1}
    +
    o(\varepsilon^{-1}).
\end{equation*}
Now let $\delta>0$. If
\begin{equation*}
    r>
    \hfixp_0(\alpha)+\Lp(\alpha,0;y,2^{3/2}\tau)+\delta,
\end{equation*}
then we see that
\begin{equation*}
\Tper_{(X_\varepsilon+S_\varepsilon,S_\varepsilon),(Y_\varepsilon+Z_\varepsilon(r),Z_\varepsilon(r))}<T_\varepsilon
\end{equation*}
for all sufficiently small $\varepsilon>0$. If
\begin{equation*}
    r<
    \hfixp_0(\alpha)+\Lp(\alpha,0;y,2^{3/2}\tau)-\delta,
\end{equation*}
then 
\begin{equation*}
\Tper_{(X_\varepsilon+S_\varepsilon,S_\varepsilon),(Y_\varepsilon+Z_\varepsilon(r),Z_\varepsilon(r))}>T_\varepsilon
\end{equation*}
for all sufficiently small  $\varepsilon>0$. Again, recalling that $\varepsilon^{-2}=\theta N$, this ensures that as $\varepsilon \rightarrow 0$ 
\begin{equation}\label{eq:ConvergenceFinalBit}
    \frac{
    h^\per_{X_{\varepsilon},S_{\varepsilon}}(Y_\varepsilon,T_\varepsilon)
    -
    \theta Y_\varepsilon
    +
    \theta(1-\theta)T_\varepsilon
    }{
    \sqrt{2(1-\theta)}\,\varepsilon^{-1}
    }
    \rightarrow \hfixp(\alpha;y,2^{3/2}\tau) := 
    \hfixp_0(\alpha)+\Lp(\alpha,0;y,2^{3/2}\tau).
\end{equation}
uniformly in $\alpha \in \R$, giving the desired result for step initial data. For general initial data $\lambda_N$ which satisfies \eqref{eq:ConvergenceIniPeriodic}, note that 
\begin{equation*}
   \hpTASEP(\lambda_N;Y,T) = \sup_{X \in \R } h^{\per}_{X,- \hpTASEP(\lambda_N;X,0)}(Y,T) .
\end{equation*}
as well as that for all $y\in\R$ and $t \geq 0$ 
\begin{equation*}
    \hfixp(\hfixp_0;y,t) = \sup_{\alpha \in \R} \hfixp(\alpha;y,t) . 
\end{equation*} 
Taking the supremum over $\alpha$ on the right-hand side of \eqref{eq:ConvergenceFinalBit}, and using that for all $y \in \R$ and $t\geq 0$, the $1\colon2\colon 3$ scaling
\begin{equation*}
    \hfixp(\hfixp_0;y,t) \overset{(d)}{=} 2^{1/2} \hfixp(\bar{\hfix}^{\per}_0;y/2,2^{-3/2}t) 
\end{equation*} with  $\bar{\hfix}^{\per}_0(x):=2^{-1/2}\hfixp_0(x / 2)$ for all $x\in \R$, which holds by a change of variables in the convergence of \eqref{eq:ConvergenceFinalBit}, we obtain the desired result.
\end{proof}

\subsection{Local structure of the periodic KPZ fixed point}\label{sec:localstruct}

In the following, we investigate the local structure of the periodic KPZ fixed point. Similar to property \hyperref[eq:Cond1]{(i)} in Theorem \ref{thm:PeriodicDirectedLandscape}, where we approximate the periodic directed landscape locally by a full-space directed landscape, we  provide a local approximation of the periodic KPZ fixed point $\hfixp$ by a full-space KPZ fixed point $\hfix$. To this end, we denote by
\begin{equation}\label{def:Support}
    \textup{supp}(\hfix_0) = \{ z \in \R \, \colon  \hfix_0(z)> - \infty \} 
\end{equation} the support of the initial data $\hfix_0$ and set
\begin{equation}\label{def:ExtendedSupport}
    \textup{ext}_{\delta}(\hfix_0) = \big\{ z + \theta \text{ for some } z \in \textup{supp}(\hfix_0) , \theta \in [-\delta^{1/2},\delta^{1/2}] \big\} 
\end{equation} as the $\delta$-extended support of $\hfix_0$. 
\begin{proposition}\label{pro:LocalPeriodicFixedPoint}
Fix initial data $\hfix_0 \in \UC$ and $\hfixp_0 \in \UCp$ such that
\begin{equation}\label{eq:PeriodicallyExtendedConf}
\hfix_0(x) = \hfixp_0(x)
\end{equation} for all $x\in \R$, where with a slight abuse of notation, we extend $\hfixp$ periodically from $\mathbb{T}$ to $\R$. 
Fix some $a\in \R$. 
Then there exist constants $c,C>0$, and a coupling $\mathbf{P}$ between the periodic KPZ fixed point $\hfixp$ and the KPZ fixed point $\hfix$, started from $\hfixp_0$ and $\hfix_0$, respectively, such that 
\begin{equation}\label{eq:CoupleEventPeriodic}
    \mathcal{A}_{\varepsilon} := \left\{ \hfixp\big(\hfixp_0;y,\varepsilon\big) = \hfix(\hfix_0;y,\varepsilon) \ \forall y \in \Big[a-\frac{1}{4},a+\frac{1}{4}\Big] \cap  \textup{ext}_{\varepsilon}(\hfix_0)\right\}
\end{equation}
satisfies for all $\varepsilon>0$
\begin{equation}\label{eq:CouplePeriod}
\mathbf{P}\left( \mathcal{A}_{\varepsilon} \right) \geq 1- C\exp(-c\varepsilon^{-1/2}) .
\end{equation} Moreover, we get that for all $y\in \mathbb{T}\setminus\textup{supp}(\hfix_0)$ 
\begin{equation}\label{eq:EscapePeriodic}
   \lim_{\varepsilon \rightarrow 0} \hfix\big(\hfix_0;y,\varepsilon\big) =  \lim_{\varepsilon \rightarrow 0} \hfixp\big(\hfixp_0;y,\varepsilon\big) = -\infty ,
\end{equation} in the sense of almost sure convergence in supremum-norm of the function in $y$.
\end{proposition}  Let us note that the choice of the constant $\frac{1}{4}$ in \eqref{eq:CoupleEventPeriodic} is not optimal; any value strictly less than $\frac{1}{2}$ would suffice. In order to show Proposition~\ref{pro:LocalPeriodicFixedPoint}, recall the moderate deviations from Proposition~\ref{pro:ModulusOfContinuity} on the full-space directed landscape, and the set $\mathcal{D}_{\delta}$ from \eqref{def:SlopeConti}. We require the following estimates for the periodic directed landscape similar to Lemma~\ref{lem:FirstFullSpace} and Lemma~\ref{lem:PenaltyDistance}.

\begin{lemma}\label{lem:ModeratePeriodic}
Let $\delta>0$. For all $M \in \N_0$, $i \in \lbr \delta^{-1}\rbr$, and $z>0$, we define the events 
\begin{align*}\label{eq:EventAm}
    \mathcal{B}_{\delta}^{(1)} &:= \bigcap_{i\in \lbr \delta^{-1} \rbr}\bigg\{ \sup_{ |x-y|> \frac{1}{8}, \, (i-1)\delta\leq s < t \leq i\delta} \Lp(x,s;y,t)  \leq - \frac{1}{300}\delta^{-1} \bigg\} , \\
    \mathcal{B}_{\delta}^{(2)} &:= \bigcap_{i\in \lbr \delta^{-1} \rbr}\bigg\{  \inf_{ (x,s;y,t) \in  \mathcal{D}_{\delta} \colon \, (i-1)\delta\leq s < t \leq i\delta} \ \Lp(x,s;y,t)  \geq - \delta^{-1/2} \bigg\} . 
\end{align*}
There exist constants $c,C>0$ 
such that for all $\delta>0$
\begin{equation}\label{eq:UniformUpperBoundDLNew}
    \Pper\Big( \mathcal{B}_{\delta}^{(1)}  \cap  \mathcal{B}_{\delta}^{(2)} \Big) \geq 1- C\delta^{-1}\exp\big(-c \delta^{-1/2}\big) . 
\end{equation}
\end{lemma}
\begin{proof}
We start with a bound on the probability of the event $\mathcal{B}_{\delta}^{(1)}$. 
In the following, we will only consider $i=1$ as the other cases are similar. Recall the rectangles $(\mathcal{Q}_j)_{j \in \Z}$ from \eqref{eq:TargetRectangle}. 
As in \eqref{eq:FullPeriodicCoupling}, there exists a coupling $\mathbf{P}$ between a periodic directed landscape $\Lp$ and full-space directed landscapes $\mathcal{L}^{(j)}$ such that for every $j \in \lbr \delta^{-1}\rbr$, the events
\begin{equation}\label{eq:CouplingPeriodicFull}
    \mathcal{C}_j := \left\{ \Lp(x,s;y,t)=\mathcal{L}^{(j)}(x,s;y,t) \ \forall (x,s;y,t) \in (\mathcal{Q}_j \times \mathcal{Q}_j) \, \colon s<t\right\} , 
\end{equation}
 satisfy for some constants $c_1,C_1>0$ and all $\delta>0$
\begin{equation}\label{eq:PeriodicPatched}
  \mathbf{P}_n(\mathcal{C}_j \ \forall j \in \lbr 16 \rbr) \geq 1- C_1 \exp(-c_1\delta^{-1/2}) .
\end{equation} 
Recall the points $(u_k)_{k\in \lbr K \rbr \cup \{0\}}$ for some $K\in \N$  with $(s_k)_{k\geq 1}$ and $(t_k)_{k\geq 1}$ from \eqref{eq:Decomposition1} and \eqref{eq:Decomposition2} in the path decomposition from $u$ to $v$ for $u=(x,s)$ and $v=(y,t)$ with some $0\leq s<t\leq\delta$ and $|x-y|>\frac{1}{8}$. Note that in this case, as $u_0=u$ and $u_K=v$, we must have $K\geq 3$. The remaining part of the proof to bound the lengths $\Lp(u_{i-1},u_{i})$ for all $i \in \lbr K \rbr$, and hence the probability of the event $\mathcal{B}_{\delta}^{(1)}$, follows mutatis muntandis as  Lemma~\ref{lem:PenaltyDistance}, using the bounds  \eqref{eq:UniformModerateDL1} and \eqref{eq:UniformModerateDL2}.  
For the second event $\mathcal{B}_{\delta}^{(2)}$, we again only consider the case $i=1$. Observe that the straight line connecting $u$ to $v$ for $(u;v) \in \mathcal{D}_{\delta}$ intersects at most two rectangles $(\mathcal{Q}_j)_{j \in \lbr 16 \rbr}$ for all $\delta>0$ small enough. Assume without loss of generality that $u,v \in \mathcal{Q}_1$ as the case of two rectangles is similar. 
Now using the coupling \eqref{eq:CouplingPeriodicFull} to the full-space directed landscape, and that by our assumptions $u=(x,s)$ and $v=(y,t)$ satisfy $
    \frac{|y-x|^2}{|t-s|} \leq 1 $, Proposition~\ref{pro:DLvalues} with $z=\delta^{-1/3}$ guarantees that 
\begin{equation*}
    \P\left( \Lp(u;v) \geq \delta^{-1/2} \ \forall (u;v)\in (\mathcal{Q}_1 \times \mathcal{Q}_1) \cap \mathcal{D}_{\delta}\right) \geq  1 - C_2 \exp(-c_2 \delta^{-1/2}) 
\end{equation*}
for some constants $c_2,C_2>0$ and all $\delta>0$. A union bound over the different rectangles $(\mathcal{Q}_j)_{j \in \lbr 16 \rbr}$ gives the desired lower bound on the probability of the event $\mathcal{B}_{\delta}^{(2)}$. 
\end{proof}

\begin{proof}[Proof of Proposition~\ref{pro:LocalPeriodicFixedPoint}]
Without loss of generality, we assume that $a=\frac{1}{2}$. 
For fixed initial data $\hfix_0 \in \UC$, we define the set of maximizers in \eqref{eq:VariationFixed} for the KPZ fixed point $\hfix$ at time $t>0$ and location $y \in \R$ as
\begin{equation}\label{def:LocalExtremes}
    \begin{split}
        \Sigma(y,t) := \left\{ z\in \R \, \colon \,  \hfix(\mathfrak{h}_0;y,t) =   \hfix_0(z) + \mathcal{L}(z,0;y,t) \right\} .
    \end{split}
\end{equation}
Similarly, for initial data $\hfixp_0 \in \UCp$, we define the set of maximizers in \eqref{eq:VariationFixedPeriodic} for the KPZ fixed point at time $t>0$ and location $y \in \mathbb{T}$ as
\begin{equation}\label{def:LocalExtremesPeriodic}
    \begin{split}
        \Sigma^{\per}(y,t) := \left\{ z\in \mathbb{T} \, \colon \,  \hfixp(\mathfrak{h}^{\per}_0;y,t) =   \hfixp_0(z) + \Lpb(z,0;y,t) \right\} . 
    \end{split}
\end{equation}
Note that $\Sigma(y,t),\Sigma^{\per}(y,t) \neq \emptyset$ by compactness and continuity of $\Lp$.
For fixed $\delta>0$, define
\begin{align*}
\mathcal{B}_{\delta} &:= \left\{ \Sigma(y,t) \subseteq \Big[ \frac{1}{8}, \frac{7}{8} \Big]\ \forall y \in \Big[\frac{1}{4},\frac{3}{4}\Big] \cap \textup{ext}_{\delta}(\hfix_0) \right\} ,  \\
\mathcal{B}^{\per}_{\delta} &:= \left\{  \Sigma^{\per}(y,t) \subseteq \Big[ \frac{1}{8}, \frac{7}{8} \Big]  \ \forall y \in \Big[\frac{1}{4},\frac{3}{4}\Big] \cap \textup{ext}_{\delta}(\hfix_0)  \right\} . 
\end{align*}
By Proposition \ref{pro:ModulusOfContinuity} to get an upper bound on $\mathcal{L}(x,0;y,\delta)$ when $|x-y|>\frac{1}{8}$, and a lower bound when $x\in  \textup{ext}_{\delta}(\hfix_0)$, there exist $c_0,C_0>0$, depending only on $\hfix_0$, such that for all $\delta>0$
\begin{equation}\label{eq:LocateStart1} 
       \P(\mathcal{B}_{\delta}) \geq 1- C_0\exp(-c_0\delta^{-1/2}) .  
\end{equation}
Similarly, by Lemma~\ref{lem:ModeratePeriodic} to get an upper bound on $\Lp(x,0;y,\delta)$ whenever $|x-y|>\frac{1}{8}$, and a lower bound when $x\in  \textup{ext}_{\delta}(\hfixp_0)$, there exist $c_1,C_1>0$, depending only on $\hfixp_0$, such that for all $\delta>0$
\begin{equation}\label{eq:LocateStart2} 
     \Pper(\mathcal{B}^{\per}_{\delta}) \geq 1- C_1\exp(-c_1\delta^{-1/2}) .  
\end{equation}
By Theorem \ref{thm:PeriodicDirectedLandscape}\hyperref[eq:Cond1]{(i)}, as well as \eqref{eq:LocateStart1} and \eqref{eq:LocateStart2} with respect to $\hfixp_0 \equiv 0$, there exists a coupling $\mathbf{P}$ between the periodic directed landscape $\Lp$ and a full-space directed landscape $\mathcal{L}$ such that
\begin{equation}\label{eq:CoupleLandscapes}
    \mathcal{C}_{\delta} :=
     \left\{  \mathcal{L}(x,s;y,t)=\Lp(x,s;y,t) \ \forall (x,s;y,t) \in \left(\Big[ \frac{1}{4}, \frac{3}{4} \Big] \times [0,\delta]\right)^2_{\uparrow} \right\} 
\end{equation}
 satisfies with some $c_2,C_2>0$ and all $\delta>0$
\begin{equation}
        \mathbf{P}\left( \mathcal{C}_{\delta} \right) \geq  1-C_2\exp(-c_2\delta^{-1/2}) . 
\end{equation}
Whenever $\mathcal{B}_{\delta}$ occurs for $\mathcal{L}$ and $\hfix_0$, as well as $\mathcal{B}^{\per}_{\delta}$ for $\Lp$ and $\hfixp_0$,  then on the event $\mathcal{C}_{\delta}$, 
\begin{equation}
\hfixp\big(\hfixp_0;y,\delta\big) = \hfix(\hfix_0;y,\delta) 
\end{equation} for all $y\in \big[ \frac{1}{4}, \frac{3}{4} \big] \cap \textup{ext}_{\delta}(\hfix_0)$. This yields  \eqref{eq:CouplePeriod}. For  \eqref{eq:EscapePeriodic}, note that $\hfix(\hfix_0;y,\varepsilon) \rightarrow -\infty$ as $\varepsilon \rightarrow 0$ by Proposition~\ref{pro:ModulusOfContinuity}. For $\hfixp$, the desired estimate follows from an upper bound by  Lemma~\ref{lem:ModeratePeriodic} on $\Lp(x,0;y,\delta)$ whenever $|x-y|\geq \frac{1}{8}$, and property \hyperref[eq:Cond1]{(i)} in Theorem~\ref{thm:PeriodicDirectedLandscape} to couple $\Lp(x,0;y,\delta)$ and $\mathcal{L}(x,0;y,\delta)$ for all  $|x-y|< \frac{1}{8}$.
\end{proof}

As a consequence, we approximate for narrow wedge initial conditions the periodic KPZ fixed point by the full-space KPZ fixed point. 
This confirms Conjecture $(5)$ in~\cite{BL:PinchedUp}.
\begin{corollary}\label{cor:KPZfixedWedge} Fix some $T>0$. Consider the narrow wedge initial data $\hfix^{\vee}_0$, i.e.,
\begin{equation}
    \hfix^{\vee}_0(y) := \begin{cases}
        0 & \text{ if } z =0, \\
        -\infty & \text{ otherwise.}
    \end{cases}
\end{equation}
As $\varepsilon \rightarrow 0$, we have convergence as a process in $x\in \R$ and $t\geq 0$ of
   \begin{equation}
   \varepsilon^{-\frac{1}{3}} \hfixp(\hfix^{\vee}_0;y\varepsilon^{\frac{2}{3}};t\varepsilon)\rightarrow   \hfix(\hfix^{\vee}_0;y;t).
\end{equation}
\end{corollary}
\begin{proof}
Since $y \varepsilon^{2/3} \in \textup{ext}_{\varepsilon}(\hfix^{\vee}_0)$ for all $\varepsilon$ sufficiently small, the result follows from Proposition~\ref{pro:LocalPeriodicFixedPoint}, 
and the fact that the full-space KPZ fixed point satisfies  the scaling relation
\begin{equation}
    \varepsilon^{-\frac{1}{3}} \hfix\big(\hfix^{\vee}_0;y\varepsilon^{\frac{2}{3}};t\varepsilon\big) \overset{(d)}{=} \hfix\big(\hfix^{\vee}_0;y;t\big) 
\end{equation} for any $y\in \R$, and all $t>0$ and $\varepsilon>0$; see Theorem~4.5 in~\cite{MQR:kpz}.
\end{proof}

Next, we study the $p$-periodic KPZ fixed point, taking the period $p$ length to infinity, and approximate it locally by a full-space KPZ fixed point. To this end, fix some $\hfix_0 \in \UC$ and define its \textbf{periodized cutoff} as
\begin{equation}
    \hfix^{\textup{cutoff}}_p(z) := \hfix_0\left( z - p\Big\lfloor \frac{z}{p}+\frac{1}{2} \Big\rfloor\right)
\end{equation} for all $z \in \R$. Note that $\hfix^{\textup{cutoff}}_p \in \UC_p$.
The following verifies Conjecture~(1.16) from~\cite{BLL:FixedPointGeneral}. 
\begin{proposition}\label{pro:LargePeriodLimit}
    Let $\hfix_0 \in \UC$, and $\hfix^{\textup{cutoff}}_p$ its corresponding periodized cutoff.  Let $\hfixp_p$ be a $p$-periodic KPZ fixed point. As $p\to \infty$, we have convergence as a process in $y\in \R$ and $t\geq 0$ of
    \begin{equation}\label{eq:CutoffRescale}
        \hfixp_p(\hfix^{\textup{cutoff}}_p;y,t) \rightarrow \hfix(\hfix_0;y,t).
    \end{equation}
\end{proposition}
\begin{proof}
    We first argue that we can couple the $p$-periodic directed landscape with high probability to agree on $R_{p,t} :=\big[ \frac{p}{4},\frac{p}{4}\big] \times [0,t]$, i.e., we find $c_1,C_1>0$, depending only on $t$ such that for all $p>0$ large enough, there exist a coupling $\mathbf{P}_p$ with
    \begin{equation}\label{eq:TwoDLsCoupling}
        \mathbf{P}_p\left(  \Lp_p(x,s;y,t)=\mathcal{L}(x,s;y,t) \ \forall (x,s;y,t) \in R_{p,t} \right) \geq 1- C_1\exp(-c_1p^{3/4}) . 
    \end{equation} 
    This follows immediately from Definition~\ref{def:PeriodicExtension} and property \hyperref[eq:Cond1]{(i)} of Theorem~\ref{thm:PeriodicDirectedLandscape} on approximating the periodic directed landscape by a full-space directed landscape.
    Now for fixed $y\in \R$ and $t \geq 0$, we find $c_2,C_2>0$, depending only on $t$, such that
    \begin{equation*}
    \begin{split}
        \P\left( \sup_{z \in \R \colon |z|\geq p/2}\Lp_p(z,0;y,t) \leq -p^{2}/300\right) &= \P\left( \sup_{z \in \R \colon |z|\geq \frac{1}{2}}\Lp(z,0;y/p,t/p^{3/2}) \leq -p^{3/2}/300 \right) \\
        &  \geq 1- C_2\exp(-c_2p^{3/4}) 
    \end{split}
    \end{equation*} for all $p>0$ large enough. Here, we use Definition~\ref{def:PeriodicExtension} for the first step, and Lemma~\ref{lem:ModeratePeriodic} for the inequality. Since $\hfix_0 \in \UC$, we find some $C_{\ast}>0$ so that 
    \begin{equation*}
      \sup_{|z| \geq p/2}  \hfix^{\textup{cutoff}}_p(z) \leq C_{\ast} p 
    \end{equation*}
for all $p>0$ large enough. 
Then we get that for some $c_3,C_3>0$ and all $p>0$ large enough
    \begin{equation*}
             \P\left( \hfixp(\hfix^{\textup{cutoff}}_p;y,t)=\sup_{z \in \R \colon |z|\leq p/2}\hfix^{\textup{cutoff}}_p(z)+ \Lp_p(z,0;y,t)\right) \geq 1- C_3\exp(-c_3p^{3/4})  . 
    \end{equation*}
Similarly, using Proposition \ref{pro:DLvalues}, note that the full-space KPZ fixed point satisfies 
\begin{equation*}
   \P\left( \hfix(\hfix_0;y,t)=\sup_{z \in \R \colon |z|\leq p/2}\hfix_0(z)+ \mathcal{L}(z,0;y,t)\right) 1- C_4\exp(-c_4p) 
\end{equation*} for some $c_4,C_4>0$ and all $p$ large enough. Together with the coupling in \eqref{eq:TwoDLsCoupling}, and the fact that $\hfix^{\textup{cutoff}}_p(z)=\hfix_0(z)$ for all $|z|\leq p/2$, we conclude.
\end{proof}

\section{Coupled KPZ fixed point convergence for the periodic ASEP}\label{sec:ASEPconvergence}
 
In the following, we establish convergence of the periodic asymmetric simple exclusion processes (ASEP) to periodic KPZ fixed points, coupled to the same periodic directed landscape. In order to present the main result, we first introduce some basic definitions and properties of the ASEP on the integers $\Z$ and of the periodic ASEP on a torus of length $N \in \N$.

\subsection{The asymmetric simple exclusion process and couplings}

Let us start by defining the asymmetric simple exclusion process on the torus and on the integer lattice. 

\begin{figure}
\centering
\begin{tikzpicture}[scale=0.9]

\def\x{1.6}
\def\y{-1.35}

\newcommand{\drawlattice}[2]{
  \node[shape=circle,scale=1.2,draw] (#1C1) at (0,#2*\y){} ;
  \node[shape=circle,scale=1.2,draw] (#1C2) at (\x,#2*\y){} ;
  \node[shape=circle,scale=1.2,draw] (#1C3) at (2*\x,#2*\y){} ;
  \node[shape=circle,scale=1.2,draw] (#1C4) at (3*\x,#2*\y){} ;
  \node[shape=circle,scale=1.2,draw] (#1C5) at (4*\x,#2*\y){} ;
  \node[shape=circle,scale=1.2,draw] (#1C6) at (5*\x,#2*\y){} ;
  \node[shape=circle,scale=1.2,draw] (#1C7) at (6*\x,#2*\y){} ;

  \draw[thick,densely dotted] (#1C1) -- (-0.6,#2*\y);
  \draw[thick] (#1C1) -- (#1C2);
  \draw[thick] (#1C2) -- (#1C3);
  \draw[thick] (#1C3) -- (#1C4);
  \draw[thick] (#1C4) -- (#1C5);
  \draw[thick] (#1C5) -- (#1C6);
  \draw[thick] (#1C6) -- (#1C7);
  \draw[thick,densely dotted] (#1C7) -- (10.2,#2*\y);
}

\drawlattice{A}{0}

\node[shape=circle,scale=0.9,fill=red] at (AC1) {};
\node[shape=circle,scale=0.9,fill=blue] at (AC2) {};
\node[shape=circle,scale=0.9,fill=red] at (AC4) {};
\node[shape=circle,scale=0.9,fill=red] at (AC5) {};

\node at (-1.4,0*\y) {$\eta_t$};
\node[anchor=west] at (10.6,0*\y) {attempt right jump};

\draw[->,thick,blue,bend left=35] (AC2) to node[above] {$\frac{1}{1+q}$} (AC3);

\drawlattice{B}{1}

\node[shape=circle,scale=0.9,fill=red] at (BC1) {};
\node[shape=circle,scale=0.9,fill=red] at (BC3) {};
\node[shape=circle,scale=0.9,fill=red] at (BC4) {};
\node[shape=circle,scale=0.9,fill=red] at (BC5) {};

\node at (-1.4,1*\y) {$\eta_{t+}$};
\node[anchor=west] at (10.6,1*\y) {jump accepted};

\drawlattice{C}{2}

\node[shape=circle,scale=0.9,fill=red] at (CC1) {};
\node[shape=circle,scale=0.9,fill=red] at (CC3) {};
\node[shape=circle,scale=0.9,fill=blue] at (CC4) {};
\node[shape=circle,scale=0.9,fill=red] at (CC5) {};

\node at (-1.4,2*\y) {$\eta_s$};
\node[anchor=west] at (10.6,2*\y) {attempt left jump};

\draw[->,thick,blue,bend right=35] (CC4) to node[above] {$\frac{q}{1+q}$} (CC3);
\node at (2.5*\x,2*\y+0.38) {$\times$};

\drawlattice{D}{3}

\node[shape=circle,scale=0.9,fill=red] at (DC1) {};
\node[shape=circle,scale=0.9,fill=red] at (DC3) {};
\node[shape=circle,scale=0.9,fill=red] at (DC4) {};
\node[shape=circle,scale=0.9,fill=red] at (DC5) {};

\node at (-1.4,3*\y) {$\eta_{s+}$};
\node[anchor=west] at (10.6,3*\y) {jump suppressed};

\end{tikzpicture}
\caption{\label{fig:ASEP} Visualization of the ASEP on $\mathbb{Z}$ with jump attempts at times $t$ and $s$. The particle attempting the jump is drawn in blue.}
\end{figure}
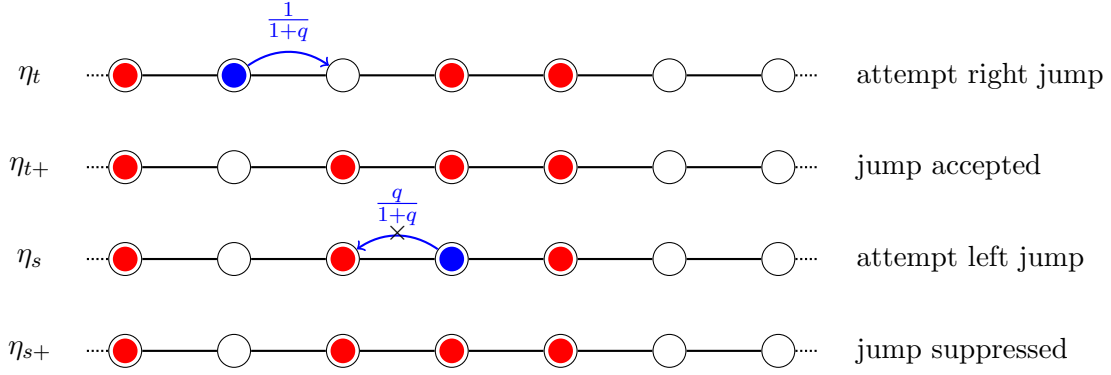

\begin{definition}\label{def:ASEP}
The \textbf{asymmetric simple exclusion process (ASEP)} on $\Z$ with parameter $q\in [0,1)$ is a continuous-time Markov chain $(\eta_t)_{t \geq 0}$ on $\{0,1\}^{\Z}$ with generator
\begin{equation}
L_{\Z}f(\eta) = \sum_{x \in \Z} \big(\eta(x)(1-\eta(x+1)) + q \eta(x+1)(1-\eta(x))\big)\big[ f(\eta^{x,x+1})- f(\eta) \big] 
\end{equation} acting on the core of cylinder functions $f \colon \{0,1\}^{\Z} \rightarrow \R$. 
Here, $\eta^{x,x+1}$ denotes the configuration, where we swap the values of $x$ and $x+1$, and leave all other values in $\eta$ unchanged, i.e.,
\begin{equation}
    \eta^{x,x+1}(y) := \begin{cases}
        \eta(x+1) & \text{ if } y=x \\
        \eta(x) & \text{ if } y=x+1 \\
        \eta(y) & \text{ otherwise. }
    \end{cases}
\end{equation}
 We define the \textbf{wrapped periodic asymmetric simple exclusion process (periodic ASEP)} as the continuous-time Markov chain $(\eta^{\per}_t)_{t \geq 0}$ on $\{0,1\}^{\Z/N\Z}$ with respect to the generator
\begin{equation}
L_{\per}f(\eta) = \sum_{x \in \Z/N\Z} \big(\eta(x)(1-\eta(x+1)) + q \eta(x+1)(1-\eta(x))\big)\big[ f(\eta^{x,x+1})- f(\eta) \big] 
\end{equation} for all cylinder functions $f \colon \{0,1\}^{\Z/N\Z} \rightarrow \R$. We extend the wrapped periodic ASEP $(\eta^{\per}_t)_{t \geq 0}$ to the \textbf{unwrapped periodic ASEP} $(\epb_t)_{t \geq 0}$ on $\{0,1\}^{\Z}$ by setting for all $t\geq 0$
\begin{equation}
    \epb_t(x) := \eta_t( x \textup{ mod } N) . 
\end{equation}
\end{definition}
In words, the processes $(\eta_t)_{t \geq 0}$ and $(\eta^{\per}_t)_{t \geq 0}$ have the following description.  Every edge $\{x,x+1\}$ on $\Z$, respectively the one-dimensional torus of length $N$, receives an independent rate $(1+q)$ clock. Whenever the clock rings at time $t$ for the edge $\{x,x+1\}$, and there is exactly one particle on either $x$ or $x+1$, then place the particle at $x+1$ with probability $(1+q)^{-1}$, and on $x$ with probability $q(1+q)^{-1}$. Otherwise, leave the configuration unchanged. We refer to Figure~\ref{fig:ASEP} for a visualization. 
The next statement is a classical result, see for example~\cite{L:interacting-particle}.
\begin{lemma}\label{lem:StationaryASEP} For all $\rho \in [0,1]$, the Bernoulli-$\rho$-product measures $\nu_{\rho}$, i.e., where every site is independently occupied with probability $\rho$, are invariant for the ASEP on $\mathbb{Z}$ and on $\mathbb{Z}/N\mathbb{Z}$. 
\end{lemma}
Let us stress that for the periodic ASEP, the number of particles is conserved. Thus, in order to obtain an irreducible Markov chain, we have to condition on the number of particles $k \in \{ 0,1,\dots,N\}$.  We now collect several couplings of ASEPs.

\begin{definition}\label{def:BasicCouplingASEP} Suppose that  $(\eta_t)_{t \geq 0}$ and $(\eta^{\prime}_t)_{t \geq 0}$ are defined on the same underlying graph $G=(V,E)$ (either the integer lattice or the torus of length $N$), but are allowed to have different initial data. 
We define the \textbf{basic coupling} $\mathbf{P}$ between the ASEPs $(\eta_t)_{t \geq 0}$ and $(\eta^{\prime}_t)_{t \geq 0}$ as follows.
For all edges in $E$, we place independent rate $(1+q)$ Poisson clocks. Whenever the clock of an edge $e=\{ x,x+1\}$ rings, and the exclusion processes are in states $\eta$ and $\eta^{\prime}$, respectively, we sample an independent Bernoulli-$(1+q)^{-1}$-distributed random variable $B$. If
\begin{itemize}
\item $B=1$: Then, if $\eta(x)=1-\eta(x+1)=1$ holds, we move the particle at site $x$ to site $x+1$ in configuration $\eta$. Similarly, if $\eta^{\prime}(x)=1-\eta^{\prime}(x+1)=1$ holds, we move the particle at site $x$ to site $x+1$ in configuration $\eta^{\prime}$. 
\item $B=0$. Then, if $\eta(x)=1-\eta(x+1)=0$ holds, we move the particle at site $x+1$ to site $x$ in configuration $\eta$. Similarly, if $\eta^{\prime}(x)=1-\eta^{\prime}(x+1)=0$ holds, we move the particle at site $x+1$ to site $x$ in configuration $\eta^{\prime}$. 
\end{itemize}
In all other cases, the configuration does not change.
\end{definition}

\subsection{Convergence of coupled full-space ASEPs to coupled full-space KPZ fixed points}

In the following, we recall the results from \cite{ACH:Fixed} on the convergence of the ASEP on the integers to a family of KPZ fixed points, coupled according to the same directed landscape. We start by defining a suitable height function for the ASEP on the integers $\Z$.
Consider a \textbf{Bernoulli path} $\lambda$ on $\Z$, i.e., a function $\lambda \colon \Z \rightarrow \Z$ such that
\begin{equation}\label{eq:Bernoulli}
\lambda(x+1)-\lambda(x) \in \{0,-1 \} 
\end{equation} for all $x\in \Z$. Note that we can obtain for every Bernoulli path $\lambda$ a particle configuration $\eta \in \{ 0,1 \}^{\Z}$ by setting for all $x\in \Z$
\begin{equation}\label{eq:PathConfiguration}
\eta(x) = \lambda(x-1)-\lambda(x) . 
\end{equation} 
\begin{definition}\label{def:HeightASEP}
For a given  Bernoulli path $h_0$, define the \textbf{ASEP height function} 
\begin{equation}
(y,t) \mapsto h^{\ASEP}(h_0,0;y,t)
\end{equation} as follows. For $t=0$, we set $h^{\ASEP}(h_0,0;y,0)=h_0(y)$ for all $y\in \Z$. Let $\eta_0$ denote the corresponding particle configuration for $h_0$ according \eqref{eq:PathConfiguration}, and let $(\eta_t)_{t \geq 0}$ be the corresponding ASEP. For times $t>0$, we define the height function $h^{\ASEP}(h_0,0;y,t)$ as follows. Whenever a particle jumps in $(\eta_t)_{t \geq 0}$ from site $y$ to site $y+1$ at time $r$, then we set
\begin{equation}
h^{\ASEP}(h_0,0;z,r) = \begin{cases}
h^{\ASEP}(h_0,0;z,r^-) + 1 & \text{ for } z=y ,  \\
h^{\ASEP}(h_0,0;z,r^-)  & \text{ for } z \neq y . 
\end{cases}
\end{equation} Similarly, whenever a particle jumps in $(\eta_t)_{t \geq 0}$ from site $y$ to site $y-1$ at time $r$, then we set
\begin{equation}
h^{\ASEP}(h_0,0;z,r) = \begin{cases}
h^{\ASEP}(h_0,0;z,r^-) - 1 & \text{ for } z=y-1 , \\
h^{\ASEP}(h_0,0;z,r^-)  & \text{ for } z \neq y-1 .
\end{cases}
\end{equation} For starting times $s\geq 0$ for $h^{\ASEP}(h_0,s;z,r)$, we consider the above construction, but with respect to $\eta_s(x)=h_0(x-1)-h_0(x)$ for all $x \in \Z$. 
\end{definition}
Let us mention that for $q=0$, this agrees up with the TASEP height function from Definition~\ref{def:HeightTASEP}.
Next, we aim to define a suitable scaling limit of the height function in a given direction, specified by a parameter $\alpha \in (-1,1)$. To this end, we define 
\begin{equation}\label{def:ParametersASEP}
\mu(\alpha):=\frac{1}{4}(1-\alpha)^{2}, \quad \mu^{\prime}(\alpha):= \frac{\alpha-1}{2} , \quad \sigma(\alpha):= \frac{1}{2}(1-\alpha^{2})^{\frac{2}{3}} , 
\end{equation}
as well as
\begin{equation}\label{def:betaASEP}
\beta(\alpha)= \frac{2\sigma(\alpha)^2}{|\mu^{\prime}(\alpha)|(1-|\mu^{\prime}(\alpha)|)} = 2 (1-\alpha^2)^{\frac{1}{3}} , 
\end{equation}
which allow us to define the rescaled ASEP height function. 
A detailed explanation for the choice of the exponents can be found in Section~2.2.4 in \cite{ACH:ScalingASEP}.

\begin{definition}[{\cite[Definition 2.11]{ACH:ScalingASEP}}]\label{def:ASEPheightRescaled}
For a given family of Bernoulli paths $h^{\varepsilon}_0 \colon \Z \rightarrow \Z$, and $\varepsilon>0$, we define the \textbf{rescaled ASEP height function} $(y,t) \mapsto \hAeps(\hfix^{\varepsilon}_0;y,t)$ by
\begin{align*}
\hAeps(\hfix^{\varepsilon}_0;y,t) := \sigma(\alpha)^{-1}\varepsilon^{\frac{1}{3}}\Big(2\mu(\alpha)\varepsilon^{-1}t &+\mu^{\prime}(\alpha)\beta(\alpha)y\varepsilon^{-\frac{2}{3}} \\
&- \hA(h^{\varepsilon}_0,0;2\alpha t \varepsilon^{-1} + \beta(\alpha)y\varepsilon^{-\frac{2}{3}},2(1-q)^{-1}\varepsilon^{-1}t) \Big) 
\end{align*} whenever the argument for $\hA$ is an integer, and linear interpolation otherwise. Here, we obtain $\hfix^{\varepsilon}_0$ from $h^{\varepsilon}_0$ by the rescaling
\begin{equation}\label{def:RescaledASEPheight}
    \hfix^{\varepsilon}_0(x) := \sigma(\alpha)^{-1}\varepsilon^{\frac{1}{3}}\Big(\mu^{\prime}(\alpha)\beta(\alpha)x\varepsilon^{-\frac{2}{3}} - h^{\varepsilon}_0(\beta(\alpha)x\varepsilon^{-\frac{2}{3}}) \Big).
\end{equation}
\end{definition}
The next statement 
provides the convergence of a collection of ASEP height functions on the integers to a family of coupled KPZ fixed points over 
an arbitrary countable subset $\mathcal{T}$  of $[0,\infty)$.

\begin{theorem}[{\cite[Theorem 2.10]{ACH:Fixed}}]\label{thm:ACHmain}
Let $\alpha \in (-1,1)$, and fix $M\in \N$. Consider a family of initial Bernoulli paths $h_0^{(i)}$ for $i\in \lbr M \rbr$. Let $\hfix_{0}^{(i),\varepsilon}$ denote the corresponding the rescaled ASEP height functions from \eqref{def:RescaledASEPheight} in Definition~\ref{def:ASEPheightRescaled}. Suppose that for all $i \in \lbr M \rbr$, the rescaled height functions  $\hfix_{0}^{(i),\varepsilon}$ converge as $\varepsilon \rightarrow 0$ in the topology of uniform convergence on compact sets to some function $\hfix_{0}^{(i)} \colon \R \rightarrow \R \cup \{-\infty \}$ in $\UC$ such that the following holds: There exists a constant $C>0$ such that, for all $\varepsilon>0$ and $i \in \lbr M \rbr$, we have 
\begin{equation}\label{eq:ExponentialGrowthInitial1}
\hfix_{0}^{(i),\varepsilon}(x) \leq C(1 + |x|) 
\end{equation} for all $x\in \R$, as well as 
\begin{equation}\label{eq:ExponentialGrowthInitial2}
\sup_{x\in [-C,C]}  \hfix_{0}^{(i),\varepsilon}(x) \geq - C . 
\end{equation}
Then we have weak convergence as $\varepsilon\to 0$ of
\begin{equation}
\left(\hAeps(\hfix_0^{(i),\varepsilon};\cdot,t)\right)_{i \in \lbr M \rbr, t \in \mathcal{T}} \rightarrow \left(\hfix(\hfix_0^{(i)};\cdot,t)\right)_{i \in \lbr M \rbr, t \in \mathcal{T}}
\end{equation} under the topology of uniform convergence on compact sets. Here, the functions on the left-hand side are coupled according to the basic coupling for the underlying ASEPs on $\Z$, while the functions on the right-hand side are coupled as in Definition~\ref{def:KPZfixedFullSpace} with respect to a common full-space directed landscape $\mathcal{L}$, taking values in $\mathcal{T}$ for the time component.
\end{theorem}

\subsection{Convergence of the periodic ASEP to coupled periodic KPZ fixed points}\label{sec:ASEPconvergencePeriodic}

In the following we will establish the convergence of the periodic ASEPs under the basic coupling to periodic KPZ fixed points, coupled according to the same periodic directed landscape $\Lp$. We  say that a configuration $\eta \in \{0,1\}^{\Z}$ is $\mathbf{N}$\textbf{-periodic} if
\begin{equation}
\eta(x) = \eta(x+N)
\end{equation} for all $x\in \Z$. Note that every configuration of the unwrapped periodic ASEP $(\bar{\eta}^{\per}_{t})_{t \geq 0}$ is $N$-periodic by construction. We call a Bernoulli path $\lambda$ $\mathbf{(N,k)}$\textbf{-periodic} if it maps  via \eqref{eq:PathConfiguration} to an $N$-periodic configuration with $k$ particles. We can now define the height function of the periodic ASEP, similarly to Definition~\ref{def:HeightASEP}. 

\begin{definition}\label{def:PeriodicASEPHeight}
Fix an $(N,k)$-periodic Bernoulli path $h^{\per}_0$ and define the \textbf{periodic ASEP height function} $h^{\pASEP}(h_0,0;y,t)$ for $t=0$ by $h^{\pASEP}(h^{\per}_0,0;y,0)=h^{\per}_0(y)$ for all $y\in \Z$. For times $t>0$, we define the height function $h^{\pASEP}(h^{\per}_0,0;y,t)$ in terms of the underlying periodic ASEP $(\eta^{\per}_t)_{t \geq 0}$ as follows. When a particle jumps from  $y$ to  $y+1$ at time $r$, set
\begin{equation}
h^{\pASEP}(h^{\per}_0,0;z,r) = \begin{cases}
h^{\pASEP}(h^{\per}_0,0;z,r^-) + 1 & \text{ if } z=y+jN \text{ for some } j \in \Z,  \\
h^{\pASEP}(h^{\per}_0,0;z,r^-)  & \text{ if } z \neq y+jN \text{ for some } j \in \Z . 
\end{cases}
\end{equation} Similarly, when a particle jumps in $(\eta^{\per}_t)_{t \geq 0}$ from site $y$ to site $y-1$ at time $r$, then we set
\begin{equation}
h^{\pASEP}(h^{\per}_0,0;z,r) = \begin{cases}
h^{\pASEP}(h^{\per}_0,0;z,r^-) - 1 & \text{ if } z=y-1+jN \text{ for some } j \in \Z , \\
h^{\pASEP}(h^{\per}_0,0;z,r^-)  & \text{ if } z \neq y-1+jN \text{ for some } j \in \Z .
\end{cases}
\end{equation} For starting times $s\geq 0$ for $h^{\pASEP}(h^{\per}_0,s;z,r)$, we consider the above construction, but with respect to $\eta^{\per}_s(x)=h^{\per}_0(x-1)-h^{\per}_0(x)$ for all $x \in \Z$. 
\end{definition}

In words, we define the height function for the periodic ASEP in the same way as in Definition~\ref{def:HeightASEP}, but with respect to periodic updates. Similarly to Definition \ref{def:ASEPheightRescaled}, we consider now a rescaled height function of the periodic ASEP in order to state our convergence results. 

\begin{definition}\label{def:ASEPheightRescaledPeriodic} Let $\varepsilon=\varepsilon(N,k):=N^{-3/2}\beta^{3/2}$ for $\beta=\beta(\alpha)$ from \eqref{def:betaASEP}. Moreover, define for all $N$ and $k$ the parameters $\rho=k/N$, and $\alpha=1-2\rho$.
For a given family of $(N,k)$-periodic Bernoulli paths $h^{\per,\varepsilon}_0 \colon \Z \rightarrow \Z$, define the \textbf{rescaled periodic ASEP height function} as the mapping $(y,t) \mapsto \hApeps(\hfix^{\per,\varepsilon}_0;y,t)$ by
\begin{align*}
\hApeps(\hfix^{\per,\varepsilon}_0;y,t) := \sigma(\alpha)^{-1}\varepsilon^{\frac{1}{3}}&\Big(2\mu(\alpha)\varepsilon^{-1}t +\mu^{\prime}(\alpha)\beta(\alpha)y\varepsilon^{-\frac{2}{3}} \\
&- h^{\pASEP}(h^{\per,\varepsilon}_0,0;2\alpha t \varepsilon^{-1} + \beta(\alpha)y\varepsilon^{-\frac{2}{3}},2(1-q)^{-1}\varepsilon^{-1}t) \Big) 
\end{align*} whenever the argument for $h^{\pASEP}$ is an integer, and linear interpolation otherwise. Here, as for \eqref{def:RescaledASEPheight}, we obtain $\hfix^{\per,\varepsilon}_0$ from $h^{\per,\varepsilon}_0$ by the rescaling
\begin{equation}\label{def:RescaledPeriodicASEPheight}
    \hfix^{\per,\varepsilon}_0(x) := \sigma(\alpha)^{-1}\varepsilon^{\frac{1}{3}}\Big(\mu^{\prime}(\alpha)\beta(\alpha)x\varepsilon^{-\frac{2}{3}} - h^{\per,\varepsilon}_0(\beta(\alpha)x\varepsilon^{-\frac{2}{3}}) \Big).
\end{equation}
\end{definition}
Observe that for any $t\geq 0$ and $(N,k)$-periodic Bernoulli path $h^{\per,\varepsilon}_0$, we can treat with a slight abuse of notation the map $y \mapsto \hApeps(\hfix^{\per,\varepsilon}_0;y,t)$ as an element of $\UCp$. 
Let  $\mathcal{T}$  be an arbitrary countable subset of $[0,1]$.
We are ready to state our main result on the convergence of coupled periodic ASEPs to coupled periodic KPZ fixed points, which formalizes Theorem \ref{thm:ConvergenceASEPMAIN}.

\begin{theorem}\label{thm:ASEPcirclePDL}
Fix $M\in \N$. For $N\in \N$ and $k=k(N) \in \lbr N \rbr$. Set $\rho=\frac{1}{2}(1-\alpha)=k/N$ for all $N$ large enough. Assume that $\rho \in [\mathfrak{a},1-\mathfrak{a}]$ for some $\mathfrak{a}>0$. 
Consider a family of $(N,k)$-periodic Bernoulli paths $h_0^{\per,(i)}$ for $i\in \lbr M \rbr$ and assume that the particle configurations $\eta_0^{\per,(i)}$, corresponding to $h_0^{\per,(i)}$, satisfy for all $i\in \lbr M \rbr$
\begin{equation}\label{eq:DominationAssumption}
\Ber_{\rho+CN^{-1/2}}^{N} \succeq \P( \eta_0^{\per,(i)} \in \, \cdot \, ) \succeq \Ber_{\rho-CN^{-1/2}}^{N}
\end{equation}
with some constant $C>0$. Let $\varepsilon=\varepsilon(N,k)>0$, and let $\hfix_{0}^{\per,(i),\varepsilon}$ denote the rescaled ASEP height functions from \eqref{def:RescaledPeriodicASEPheight} in Definition~\ref{def:ASEPheightRescaledPeriodic} at time $0$.
Suppose that for all $i \in \lbr M \rbr$, there exists a continuous function $\hfix_{0}^{\per,(i)} \colon \mathbb{T} \rightarrow \R$ 
such that as $\varepsilon \rightarrow 0$, $\hfix_{0}^{\per,(i),\varepsilon}$ converges weakly to $\hfix_{0}^{\per,(i)}$ with respect to the product topology induced by the uniform topology on $\mathcal{C}(\mathbb{T})$. Let $\hfixp$ denote a periodic KPZ fixed point. 
Then we have weak convergence as $\varepsilon\to 0$ of
\begin{equation}
\left(\hApeps(\hfix_0^{\per,(i),\varepsilon};\cdot,t)\right)_{i \in \lbr M \rbr, t \in \mathcal{T}} \rightarrow \left(\hfixp(\hfix_0^{\per,(i)};\cdot,t)\right)_{i \in \lbr M \rbr, t \in \mathcal{T}} . 
\end{equation}  Here, the functions on the left-hand side are coupled according to the basic coupling for the underlying periodic ASEPs, while the functions on the right-hand side are coupled as in Definition~\ref{def:PeriodicKPZfixed} with respect to the same periodic directed landscape $\Lp$, taking values in $\mathcal{T}$ for the time component.
\end{theorem}

Let us now describe the overall strategy for the proof of Theorem~\ref{thm:ASEPcirclePDL}. In Section~\ref{sec:ReferenceFrameCoupling}, we define the reference frame basic coupling between the unwrapped periodic ASEP and the ASEP on the integers. Under this coupling, the two processes will agree up with probability going to $1$ as $\delta \rightarrow 0$ on an interval of length $\frac{3}{4}N$ until time $\delta N^{3/2}$; see Proposition~\ref{pro:ModerateDeviationsSecondClass}. In Section~\ref{sec:PatchedASEP}, we define the patched ASEP as an approximation of the periodic ASEP, similar to the unwrapped patched directed landscape from Definition~\ref{def:PatchedDLExt}. In Section~\ref{sec:ConvergencePatchedASEP}, we define the patched KPZ fixed point via a variational characterization using the patched directed landscape. We show patched ASEPs under the basic coupling converge to patched KPZ fixed points, coupled according to the same patched directed landscape. Moreover, we show that the coupled patched KPZ fixed points are close to periodic KPZ fixed points, coupled according to a common periodic directed landscape. 
Combining the above steps, the proof of Theorem~\ref{thm:ASEPcirclePDL} is completed in Section~\ref{sec:ASEPfinalPart}. 

\begin{remark}\label{rem:ASEPtopDL}
There are several, likely accessible, ways to extend our ASEP results that we do not pursue here.
In Corollary 2.12(a) of \cite{ACH:ScalingASEP}, the authors establish the convergence of the ASEP on $\Z$ to the directed landscape under (colored) step initial data. We believe that a similar convergence can be shown for the periodic ASEP under step initial data to the periodic directed landscape, and leave this to future work.
We also believe that the assumption \eqref{eq:DominationAssumption} can be relaxed considerably, and that it should be possible to consider the scaling limit for multiple initial data with different filling fractions. Finally, rather than considering ASEP (or LPP) in a periodic domain, we could consider its unwrapped version in full-space with a periodic environment of Poisson clocks defining its dynamics, but now consider non-periodic initial data. The scaling limit should be accessible and related to the scaling limit of the unwrapped patched landscape that we work with in our construction in Section \ref{sec:PeriodicDL}.
\end{remark}

\subsection{The reference frame basic coupling for the ASEP}\label{sec:ReferenceFrameCoupling}

In preparation for proving Theorem~\ref{thm:ASEPcirclePDL}, we require an extension of the basic coupling. This allows to couple ASEP when the underlying graphs differ, i.e., we couple an ASEP on the integer with a periodic ASEP; see for example Section 2 in  \cite{GNS:MixingOpen} for a similar treatment for the open ASEP. 
For $a,b\in \R$ with $a<b$,  $\delta \in \R_{>0} \cup \{+\infty\}$, and  $N\in \N$ define the space-time $\boldsymbol{(\rho,a,b-a,N)}$\textbf{-reference frame}
\begin{align}\label{def:MovingFrame}
\mathcal{I}^{\rho}_{\delta,N}(a,b) &:= \left\{ (x,t) \in \R \times [0,\delta N^{3/2}) \colon x \in \big[ aN +f(t), bN + f(t)\big)  \right\}, 
\end{align} where we set for all $t\geq 0$
\begin{equation}\label{def:Functionf}
    f(t):=(1-2\rho)(1-q)t .
\end{equation} We write $\mathcal{I}^{\rho}_{N}(a,b)=  \mathcal{I}^{\rho}_{\infty,N}(a,b)$.

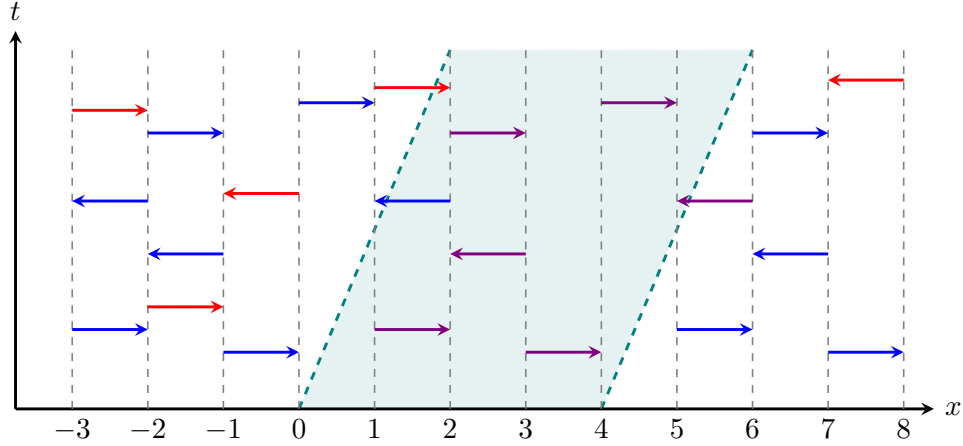
\begin{figure}
\centering
\begin{tikzpicture}[
    scale=1.0,
    >=stealth,
    particle/.style={circle, fill=black, inner sep=2.3pt},
    harrisarrow/.style={->, line width=1.1pt},
    pathline/.style={line width=1.1pt},
    siteline/.style={gray, dashed, line width=0.6pt}
]

\def\tmax{5}

\fill[teal!10] (0,0) -- (2,4.75) -- (6,4.75)--(4,0)--(0,0);

\draw[pathline,dashed,teal] (0,0) -- (2,4.75);
\draw[pathline,dashed,teal] (4,0) -- (6,4.75);

\draw[->, line width=1.1pt] (-3.75,0) -- (8.4,0) node[right] {$x$};
\draw[->, line width=1.1pt] (-3.75,0) -- (-3.75,\tmax) node[above] {$t$};

\foreach \x in {-3,-2,-1,0,1,2,3,4,5,6,7,8} {
    \draw[siteline] (\x,0) -- (\x,\tmax-0.2);
    \node[below] at (\x,0) {$\x$};
}







\draw[harrisarrow,red] (-2,1.35) -- (-1,1.35);

\draw[harrisarrow,violet] (1,1.05) -- (2,1.05);
\draw[harrisarrow,violet] (3,2.05) -- (2,2.05);
\draw[harrisarrow,violet] (2,3.65) -- (3,3.65);
\draw[harrisarrow,violet] (3,0.75) -- (4,0.75);
\draw[harrisarrow,violet] (4,4.05) -- (5,4.05);
\draw[harrisarrow,violet] (6,2.75) -- (5,2.75);

\draw[harrisarrow,blue] (1-4,1.05) -- (2-4,1.05);
\draw[harrisarrow,blue] (3-4,2.05) -- (2-4,2.05);
\draw[harrisarrow,blue] (2-4,3.65) -- (3-4,3.65);
\draw[harrisarrow,blue] (3-4,0.75) -- (4-4,0.75);
\draw[harrisarrow,blue] (4-4,4.05) -- (5-4,4.05);
\draw[harrisarrow,blue] (1+4,1.05) -- (2+4,1.05);
\draw[harrisarrow,blue] (3+4,2.05) -- (2+4,2.05);
\draw[harrisarrow,blue] (2+4,3.65) -- (3+4,3.65);
\draw[harrisarrow,blue] (3+4,0.75) -- (4+4,0.75);
\draw[harrisarrow,blue] (6-4,2.75) -- (5-4,2.75);
\draw[harrisarrow,blue] (6-8,2.75) -- (5-8,2.75);

\draw[harrisarrow,red] (1,4.25) -- (2,4.25);
\draw[harrisarrow,red] (0,2.85) -- (-1,2.85);
\draw[harrisarrow,red] (8,4.35) -- (7,4.35);
\draw[harrisarrow,red] (-3,3.95) -- (-2,3.95);


\end{tikzpicture}
\caption{\label{fig:Harris}Construction of the ASEP $(\eta_t)_{t \geq 0}$ on $\mathbb Z$ with blue arrows, and the unwrapped periodic ASEP $(\bar{\eta}^{\per}_t)_{t \geq 0}$ under the $(\rho,\theta,1,N)$-reference frame basic coupling with $\theta=0$, $N=4$ and some $\rho<\frac{1}{2}$ shaded in green. Arrows indicate the Poisson clocks along the edges, while the direction of the arrows indicate the outcome of the Bernoulli-$(1+q)^{-1}$-random variables. Purple arrows belong to both processes, while red arrows belong solely to the process $(\eta_t)_{t \geq 0}$, and blue arrows solely to the process $(\bar{\eta}^{\per}_t)_{t \geq 0}$. }
\end{figure}

\begin{definition}\label{def:BasicCouplingASEPFrame} Suppose that $(\eta_t)_{t \geq 0}$ is an ASEP on $\Z$, and that $(\bar{\eta}^{\per}_t)_{t \geq 0}$ is an unwrapped periodic ASEP. Fix some $\rho \in [0,1]$ and $\theta \in \R$, which are allowed to depend on $N$. Set
\begin{equation}
    E_t := \left\{ \{ x, x+1 \} \colon x \in \mathcal{I}^{\rho}_{N}(\theta,\theta+1) \right\}  
\end{equation} for all $t\geq 0$. We define the $\boldsymbol{(\rho,\theta,1,N)}$\textbf{-reference frame basic coupling} as follows. 
For all edges, we place independent rate $(1+q)$ Poisson clocks. 
Whenever the clock of an edge $e=\{ x,x+1\}$ rings, and the exclusion processes are in states $\eta$ and $\eta^{\prime}$, respectively, we sample an independent Bernoulli-$(1+q)^{-1}$-distributed random variable $B$. If
\begin{itemize}
\item $B=1$ and $e\in E_t$. Then, if $\eta(x)=1-\eta(x+1)=1$ holds, we move the particle at site $x$ to site $x+1$ in configuration $\eta$. Similarly, if $\bar{\eta}^{\per}(x)=1-\bar{\eta}^{\per}(x+1)=1$ holds, we move for all $j\in \Z$ the particle at site $x+jN$ to site $x+jN+1$ in $\bar{\eta}^{\per}$. 
\item $B=0$ and $e\in E_t$. Then, if $\eta(x)=1-\eta(x+1)=0$ holds, we move the particle at site $x+1$ to site $x$ in configuration $\eta$. Similarly, if $\bar{\eta}^{\per}(x)=1-\bar{\eta}^{\per}(x+1)=0$ holds, we move for all $j\in \Z$ the particle at site $x+jN+1$ to site $x+jN$ in  $\bar{\eta}^{\per}$. 
\item $B=1$ and $e \notin E_t$. Then, if $\eta(x)=1-\eta(x+1)=1$ holds, we move the particle at site $x$ to site $x+1$ in configuration $\eta$.
\item $B=0$ and $e\notin E_t$. Then, if $\eta(x)=1-\eta(x+1)=0$ holds, we move the particle at site $x+1$ to site $x$ in configuration $\eta$.
\end{itemize}
In all other cases, the configuration remains unchanged. 
\end{definition}
A visualization of the $(\rho,\theta,1,N)$-reference frame basic coupling is given in Figure~\ref{fig:Harris}.
Consider the basic coupling as in Definition~\ref{def:BasicCouplingASEP} with component-wise ordered initial states $\eta_0$ and $\eta^{\prime}_0$, i.e., 
\begin{equation}\label{def:ComponentwiseOrder}
\eta_0 \succeq \eta^{\prime}_0 \quad \Leftrightarrow \quad \eta_0(x) \geq  \eta^{\prime}_0(x) \text{ for all } x \in V . 
\end{equation} Then we get that under the basic coupling 
\begin{equation}\label{def:attactivity}
\mathbf{P}\big( \eta_t \succeq \eta^{\prime}_t \text{ for all } t \geq 0 \, \big| \,  \eta_0 \succeq \eta^{\prime}_0 \big) = 1 ,
\end{equation} We refer to \eqref{def:attactivity} as \textbf{attractivity}.
This motivates the following definition on the joint evolution of a collection of ASEPs under the basic coupling. 
To this end, we let $M \in \N$ be fixed, and consider the state space $\Omega^{M}_{N} = \lbr M \rbr^{N}$, where we refer to $1,\dots,M$ as colors of the particles.

\begin{definition}\label{def:ColoredASEP}
Suppose that we are given a family of $M-1$ ASEPs $(\eta^{(i)}_t)_{t \geq 0}$ for $i\in \lbr M-1\rbr$ on the same graph $G=(V,E)$ such that
\begin{equation}
 \eta^{(M-1)}_0 \succeq  \eta^{(M-2)}_0 \succeq \dots \succeq \eta^{(2)}_0 \succeq \eta^{(1)}_0 
\end{equation} according to the basic coupling $\mathbf{P}$. Then we define the \textbf{colored ASEP} $(\zeta_t)_{t \geq 0}$ for all $t \geq 0$ by
\begin{equation}
\zeta_t(i) := 1+\sum_{j=1}^{M-1} \big(1-\eta^{(j)}_t(i)\big)  
\end{equation}  for all $t \geq 0$ and $i \in  V$. We refer to the particles of colors $1,2,\dots,M-1$ as \textbf{first class particles}, \textbf{second class particles}, and so on, as well as to type $M$ particles as \textbf{empty sites}.
\end{definition}

In words, we assign to every edge $e$ a rate $1+q$ clock. When the clock rings at an edge $\{x,x+1\}$ at time $t$, we sample an independent Bernoulli-$(1+q)^{-1}$-random variable $B$. If $B=1$, then we sort the endpoints of $e$ in (weakly) increasing order, i.e., $
\eta_t(x+1) \geq \eta_t(x)$. If $B=0$, then we sort the endpoints of $e$   in decreasing order.
While the colored ASEP is suited to study the ASEP with ordered initial data, we will also require an extended reference frame coupling, following the ideas from \cite{FS:triple}, to compare ASEPs on different graphs and non-ordered initial data. For our needs, we state this extended coupling only when $M=3$. This is similar to Definition 4.1 in \cite{FS:triple}, which provides an extended coupling between the open ASEP and the ASEP on the integers with respect to the reference frame $\mathcal{I}^{1/2}_N(0,1)$.

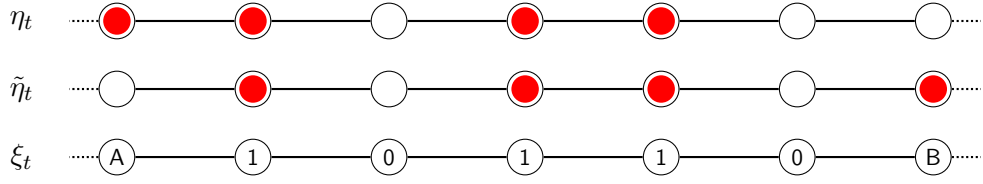
\begin{figure}
\centering
\begin{tikzpicture}[scale=0.9]

\def\x{2};
\def\y{-1};

 	\node[shape=circle,scale=1.3,draw] (C1) at (0,0*\y){} ; 
 	\node[shape=circle,scale=1.3,draw] (C2) at (\x,0*\y){} ; 
 	\node[shape=circle,scale=1.3,draw] (C3) at (2*\x,0*\y){} ; 
 	\node[shape=circle,scale=1.3,draw] (C4) at (3*\x,0*\y){} ; 
 	\node[shape=circle,scale=1.3,draw] (C5) at (4*\x,0*\y){} ; 
 	\node[shape=circle,scale=1.3,draw] (C6) at (5*\x,0*\y){} ; 
 	\node[shape=circle,scale=1.3,draw] (C7) at (6*\x,-0*\y){} ;

 	\draw[thick,densely dotted] (C1) -- (-0.7,0*\y);
 	\draw[thick] (C1) -- (C2);
 	\draw[thick] (C2) -- (C3);
 	\draw[thick] (C3) -- (C4);
 	\draw[thick] (C4) -- (C5);
 	\draw[thick] (C5) -- (C6);
 	\draw[thick] (C6) -- (C7);
 	\draw[thick,densely dotted] (C7) -- (12.7,0*\y);

	\node[shape=circle,scale=1,fill=red] (F5) at (C1) {};
	\node[shape=circle,scale=1,fill=red] (F5) at (C2) {}; 	
	\node[shape=circle,scale=1,fill=red] (F5) at (C4) {};
	\node[shape=circle,scale=1,fill=red] (F5) at (C5) {}; 	
 
	\node (H3) at (-1.4,0*\y) {$\eta_t$};

 	\node[shape=circle,scale=1.3,draw] (C1) at (0,1*\y){} ; 
 	\node[shape=circle,scale=1.3,draw] (C2) at (\x,1*\y){} ; 
 	\node[shape=circle,scale=1.3,draw] (C3) at (2*\x,1*\y){} ; 
 	\node[shape=circle,scale=1.3,draw] (C4) at (3*\x,1*\y){} ; 
 	\node[shape=circle,scale=1.3,draw] (C5) at (4*\x,1*\y){} ; 
 	\node[shape=circle,scale=1.3,draw] (C6) at (5*\x,1*\y){} ; 
 	\node[shape=circle,scale=1.3,draw] (C7) at (6*\x,1*\y){} ;  

 	\draw[thick,densely dotted] (C1) -- (-0.7,1*\y); 
 	\draw[thick] (C1) -- (C2);		
 	\draw[thick] (C2) -- (C3);
 	\draw[thick] (C3) -- (C4);
 	\draw[thick] (C4) -- (C5);
 	\draw[thick] (C5) -- (C6);
	\draw[thick] (C6) -- (C7);
 	\draw[thick,densely dotted] (C7) -- (12.7,1*\y);

	\node[shape=circle,scale=1,fill=red] (F5) at (C2) {}; 	
	\node[shape=circle,scale=1,fill=red] (F5) at (C7) {};	
	\node[shape=circle,scale=1,fill=red] (F5) at (C4) {}; 	
	\node[shape=circle,scale=1,fill=red] (F5) at (C5) {};

	\node (H3) at (-1.4,1*\y) {$\tilde{\eta}_t$};




 

 	\node[shape=circle,scale=1.3,draw] (C1) at (0,2*\y){} ; 
 	\node[shape=circle,scale=1.3,draw] (C2) at (\x,2*\y){} ; 
 	\node[shape=circle,scale=1.3,draw] (C3) at (2*\x,2*\y){} ; 
 	\node[shape=circle,scale=1.3,draw] (C4) at (3*\x,2*\y){} ; 
 	\node[shape=circle,scale=1.3,draw] (C5) at (4*\x,2*\y){} ; 
 	\node[shape=circle,scale=1.3,draw] (C6) at (5*\x,2*\y){} ; 
 	\node[shape=circle,scale=1.3,draw] (C7) at (6*\x,2*\y){} ;  
 
 	\draw[thick,densely dotted] (C1) -- (-0.7,2*\y); 
 	\draw[thick] (C1) -- (C2); 	
 	\draw[thick] (C2) -- (C3);
 	\draw[thick] (C3) -- (C4);
 	\draw[thick] (C4) -- (C5);
 	\draw[thick] (C5) -- (C6);
 	\draw[thick] (C6) -- (C7);
 	\draw[thick,densely dotted] (C7) -- (12.7,2*\y); 	
 	
	\node[scale=0.8] (H3) at (C1) {$\Aup$}; 	
	\node[scale=0.8] (H3) at (C2) {$\one$}; 
	\node[scale=0.8] (H3) at (C3) {$\zero$};  	
	\node[scale=0.8] (H3) at (C4) {$\one$}; 	
	\node[scale=0.8] (H3) at (C5) {$\one$}; 
	\node[scale=0.8] (H3) at (C6) {$\zero$};  	
	\node[scale=0.8] (H3) at (C7) {$\Bup$};  	
	

 
	\node (H3) at (-1.4,2*\y) {$\xi_t$};

	\end{tikzpicture}
\caption{\label{fig:ColoredASEP}Visualization of the extended colored ASEP $(\xi_t)_{t \geq 0}$ from Definition~\ref{def:ExtendedCouplingASEP} with respect to $(\eta_t)_{t \geq 0}$ and $(\tilde{\eta}_t)_{t \geq 0}$ on $\Z$. 
}
 \end{figure}

\begin{definition}\label{def:ExtendedCouplingASEP}
Let $(\eta_t)_{t \geq 0}$ be an ASEP on $\Z$ and $(\bar{\eta}^{\per}_t)_{t \geq 0}$ be an unwrapped periodic ASEP evolving  according to the $(\rho,\theta,1,N)$-reference frame basic coupling for some $\rho \in (0,1)$, $\theta \in \R$ and $N\in \N$. We set
\begin{equation}\label{eq:ProjectionExtension}
\xi_t(x) := \begin{cases}
\one & \text{ if } \eta_t(x)=1 \text{ and } \bar{\eta}^{\per}_t(x)=1, \\
\zero & \text{ if } \eta_t(x)=0 \text{ and }\bar{\eta}^{\per}_t(x)=0, \\
\Aup & \text{ if } \eta_t(x)=1 \text{ and } \bar{\eta}^{\per}_t(x)=0, \\
\Bup & \text{ if } \eta_t(x)=0 \text{ and } \bar{\eta}^{\per}_t(x)=1, \\
\end{cases}
\end{equation} for all $x\in \Z$ and $t \geq 0$.
 We refer to $(\xi_t)_{t \geq 0}$ on the state space $\{\one,\zero,\Aup,\Bup\}^{\Z}$ as the \textbf{extended colored ASEP with reference frame $\boldsymbol{\mathcal{I}_{N}^{\rho}(\theta,\theta+1)}$}. 
 
 Note that in the same way, replacing the process $(\bar{\eta}^{\per}_t)_{t \geq 0}$ by a full-space ASEP $(\eta'_t)_{t \geq 0}$, we can couple two full space ASEPs under the $(\rho,\theta,1,N)$-reference frame basic coupling, i.e., we use the same Poisson clocks for all $(x,t) \in \mathcal{I}^{\rho}_N(\theta,\theta+1)$ in the two ASEPs, and independent Poisson clocks for all for all $(x,t) \notin \mathcal{I}^{\rho}_N(\theta,\theta+1)$. We refer to the resulting process as \textbf{full space extended colored ASEP under the $(\rho,\theta,1,N)$-reference frame basic coupling}.
\end{definition}
A visualization of the extended colored ASEP is given in Figure~ \ref{fig:ColoredASEP}. 
 In words, we obtain $(\xi_t)_{t \geq 0}$  by the same update procedure as in the construction of the colored ASEP in Definition~\ref{def:ColoredASEP}, with the exception of an update of an $\{\Aup, \Bup\}$ pair, which gets replaced by a $\{\one,\zero\}$ pair. We write $\Pex$ for the law of $(\xi_t)_{t \geq 0}$ and refer to the law of $(\xi_t)_{t \geq 0}$ as the \textbf{extended $\boldsymbol{(\rho,\theta,1,N)}$-reference frame basic coupling} of $(\eta_t)_{t \geq 0}$ and $(\bar{\eta}^{\per}_t)_{t \geq 0}$. 

For two probability measures $\nu,\nu^{\prime}$ on a common probability space $\Omega$ with a partial order $\succeq$, we say that $\nu$ \textbf{stochastically dominates} $\nu^{\prime}$, and write with a slight abuse of notation $\nu \succeq \nu^{\prime}$, if there exists a coupling $\mathbf{P}_{\nu,\nu^{\prime}}$ such that
\begin{equation}
    \mathbf{P}_{\nu,\nu^{\prime}}( X \succeq X^{\prime}) = 1 
\end{equation} with $X \sim \nu$ and $X^{\prime} \sim \nu^{\prime}$. 
 We have the following comparison between $(\eta_t)_{t \geq 0}$ and $(\bar{\eta}^{\per}_t)_{t \geq 0}$ evolving according to the extended $(\rho,\theta,1,N)$-reference frame basic coupling.

\begin{proposition}\label{pro:ModerateDeviationsSecondClass}
Let $C>0$, $\theta \in \R$, $N\in \N$ and $\mathfrak{a}\in (0,1)$. Let $\rho \in [\mathfrak{a},1-\mathfrak{a}]$ and consider ASEPs $(\eta_t)_{t \geq 0}$ and $(\bar{\eta}^{\per}_t)_{t \geq 0}$ evolving according to the extended $(\rho,\theta,1,N)$-reference frame basic coupling. Consider (random) data $\eta^{N} \in \{0,1\}^{N}$ with 
\begin{equation}\label{eq:BernoulliDomination}
 \textup{Ber}_{\rho+CN^{-1/2}}^{N} \succeq \P( \eta^{N} \in \, \cdot \, ) \succeq  \textup{Ber}_{\rho-CN^{-1/2}}^{N} ,  
\end{equation} where $\textup{Ber}_{\rho}^{N}$ denotes the Bernoulli-$\rho$-product measure on $\{0,1\}^{N}$. Define initial data $\eta_0\in \{0,1\}^{\mathbb{Z}}$
\begin{equation}\label{eq:AgreeLocally}
\eta_0(x):=\eta^N(x- \lfloor \theta N \rfloor )
\end{equation} for all $x \in \lbr \theta N+1, (1+\theta)N \rbr$. For all
$x \notin \lbr \theta N+1, (1+\theta)N\rbr$, let $\eta_0(x)$ be occupied independently with probability $\rho$. For initial data $\tilde{\eta}_0$, assume that \eqref{eq:AgreeLocally} holds, but extend the configuration outside of $\lbr \theta N+1, (1+\theta)N \rbr$ periodically. Then there exist some $\tilde{c},\tilde{C}>0$, depending only on $\mathfrak{a}$, $q$ and $C$, and not $\theta,N,\rho$, such that the extended colored ASEP  $(\xi_t)_{t \geq 0}$ from Definition~\ref{def:ExtendedCouplingASEP} and  
\begin{equation*}
    \mathcal{A}^{\theta}_{\delta,N} := \left\{ \xi_t(x) \in \{\zero,\one \}  \, \forall (x,t) \in \mathcal{I}^{\rho}_{\delta,N}\Big(\frac{1}{8}+\theta,\frac{7}{8}+\theta\Big)  \right\}
\end{equation*} satisfy for all $\delta>0$ 
\begin{equation}\label{eq:GoodControlASEP}
\liminf_{N \rightarrow \infty}\Pex\big( \mathcal{A}^{\theta}_{\delta,N} \big)   \geq 1- \tilde{C} \exp(-\tilde{c} \delta^{-1} ) . 
\end{equation} 
Furthermore,  \eqref{eq:GoodControlASEP} continues to hold when $(\xi_t)_{t\geq 0}$ is defined the as a full-space extended colored ASEP with reference frame $\mathcal{I}^{\rho}_{N}\big(\theta,1+\theta\big)$ for two full-space ASEPs, i.e., we obtain $(\xi_t)_{t\geq 0}$ as in Definition~\ref{def:ExtendedCouplingASEP} with respect to two ASEPs on the integers, which evolve under the (extended) basic coupling along edges in $\mathcal{I}^{\rho}_{N}\big(\theta,1+\theta\big)$.
\end{proposition}
Since the arguments are similar to Proposition~4.10 in \cite{FS:triple}, using the (extended) colored ASEP, we outline the proof of Proposition~\ref{pro:ModerateDeviationsSecondClass} in Appendix~\ref{sec:AppendixASEP}.

\subsection{The patched asymmetric simple exclusion process}\label{sec:PatchedASEP}

In order to show Theorem \ref{thm:ASEPcirclePDL}, we utilize that the periodic ASEP can by Proposition~\ref{pro:ModerateDeviationsSecondClass} be approximated by  families of ASEPs on the integers. 
 Recall the sets $\mathcal{I}^{\rho}_{\delta,N}(a,b)$ from \eqref{def:MovingFrame}.
Fix $M\in \N$, and consider two families of ASEPs on the integers $(\eta^{N,i,1}_t)_{t \geq 0}$ and $(\eta^{N,i,2}_t)_{t \geq 0}$ in $N \in \N$ and $i\in \lbr M \rbr$. The processes $(\eta^{N,i,1}_t)_{t \geq 0}$ for $i\in \lbr M \rbr$  and $(\eta^{N,i,2}_t)^{i\in \lbr M \rbr}_{t \geq 0}$ for $i\in \lbr M \rbr$ evolve under the following modified version of the $(\rho,\theta,1,N)$-reference frame basic coupling until time $\delta N^{3/2}$: 
\begin{enumerate}
    \item whenever an edge $e=\{ x,x+1\}$ with $(x,t) \in \mathcal{I}^{\rho}_{\delta,N}(\frac{1}{4},\frac{3}{4})$ is updated for $\eta^{N,\cdot,1}_t$, then so is the corresponding edge $e$ in $\eta^{N,\cdot,2}_t$,
    \item whenever an edge $e=\{ x,x+1\}$ with $(x,t) \in \mathcal{I}^{\rho}_{\delta,N}(-\frac{1}{4},\frac{1}{4})$ is updated for $\eta^{N,\cdot,1}_t$, then so is the edge $e'=\{ x+N,x+1+N\}$ with $(x,t) \in \mathcal{I}^{\rho}_{\delta,N}(\frac{3}{4},\frac{5}{4})$ in $\eta^{N,\cdot,2}_t$.
    \item All other updates are performed according to two independent basic couplings for the families of ASEPs $(\eta^{N,i,1}_t)^{i \in \lbr M \rbr}_{t \geq 0}$ and $(\eta^{N,i,2}_t)^{i \in \lbr M \rbr}_{t \geq 0}$, respectively. 
\end{enumerate}
We refer to this as the \textbf{patched coupling} $\Ppatch$ between $(\eta^{N,i,1}_t)^{i\in \lbr M \rbr}_{t \geq 0}$ and $(\eta^{N,i,2}_t)^{i\in \lbr M \rbr}_{t \geq 0}$. 
\begin{lemma}\label{lem:patchedASEP} Fix some $M\in \N$. Assume that the initial data $ (\eta_0^{N,i,j})$ for $i \in \lbr M \rbr, j\in \{ 1,2\}$ of the above exclusion processes is $N$-periodic with $k=k(N)$ particles, satisfies $\eta_0^{N,i,1}=\eta_0^{N,i,2}$ for all $i\in \lbr M\rbr$, and for some constant $C>0$ 
    \begin{equation}\label{eq:DominationAssumptionPatch}
\Ber_{\rho+CN^{-1/2}}^{N} \succeq \P( \eta_0^{N,i,1} \in \, \cdot \, ) \succeq \Ber_{\rho-CN^{-1/2}}^{N}
\end{equation} for all $N \in \N$, $i\in \lbr M \rbr$, and  $\rho \in [\mathfrak{a},1-\mathfrak{a}]$ for some $\mathfrak{a}>0$. Then there exist  $c_0,C_0>0$, depending only on $\mathfrak{a},q$ and $C$, but not $\theta,\rho,N$, such that under the patched coupling,  for all $\delta>0$
\begin{equation}\label{eq:GoodCouplingPatched}
\begin{split}
    \Ppatch\left( \eta_t^{N,i,1}(x)=\eta_t^{N,i,2}(x) \ \forall (x,t)\in \mathcal{I}^{\rho}_{\delta,N}\Big(\frac{3}{8},\frac{5}{8}\Big)  \right) &\geq 1 - C_0 \exp(-c_0\delta^{-1}) , \\
    \Ppatch\left( \eta_t^{N,i,1}(x)=\eta_t^{N,i,2}(x+N) \ \forall (x,t)\in \mathcal{I}^{\rho}_{\delta,N}\Big(-\frac{1}{8},\frac{1}{8}\Big)\right) &\geq 1 - C_0 \exp(-c_0\delta^{-1}) .
\end{split}
\end{equation}
\end{lemma}
\begin{proof} For all $i\in \lbr M \rbr$ and $N\in \N$, let $(\bar{\eta}_t^{N,i})_{t \geq 0}$ denote an unwrapped periodic ASEP with initial data $\eta_0^{N,i,1}$. By Proposition \ref{pro:ModerateDeviationsSecondClass} under the $(\rho,0,1,N)$-reference frame basic coupling $\mathbf{P}_0$, 
\begin{equation}\label{eq:Extended1}
    \mathbf{P}_0\left( \eta_t^{N,i,1}(x)=\bar{\eta}_t^{N,i}(x) \ \forall (x,t)\in \mathcal{I}^{\rho}_{\delta,N}\Big(-\frac{1}{8},\frac{5}{8}\Big)  \right) \geq 1 - \frac{1}{2}C_0 \exp(-c_0\delta^{-1})  
\end{equation}
holds for some constants $c_0,C_0>0$ and all $\delta>0$. Similarly, by Proposition \ref{pro:ModerateDeviationsSecondClass} and the $(\rho,\frac{1}{2},1,N)$-reference frame basic coupling $\mathbf{P}_{\frac{1}{2}}$,
\begin{equation}\label{eq:Extended2}
    \mathbf{P}_{\frac{1}{2}}\left( \eta_t^{N,i,2}(x)=\bar{\eta}_t^{N,i}(x) \ \forall (x,t)\in \mathcal{I}^{\rho}_{\delta,N}\Big(\frac{3}{8},\frac{9}{8}\Big)  \right) \geq 1 - \frac{1}{2}C_0 \exp(-c_0\delta^{-1})
\end{equation} for all $\delta>0$.
Combining \eqref{eq:Extended1} and \eqref{eq:Extended2} with the fact that $\bar{\eta}_t^{N,i}$ is $N$-periodic for all $t\geq 0$ and $i \in \lbr M \rbr$ yields the desired result.
\end{proof}
Using the patched coupling, we define an approximation of the periodic ASEP.
\begin{definition}\label{def:PatchedASEPNew}
For each $j\in \N$, consider families of ASEPs on the integers $(\eta_t^{N,i,1,j})_{t \geq 0}$ and $(\eta_t^{N,i,2,j})_{t \geq 0}$ under the patched coupling with $i \in \lbr M \rbr$ and $N \in \N$. We define the $\boldsymbol{(\delta,N)}$\textbf{-patched ASEP} $(\eta^{\textup{patch},N,i}_t)_{t \geq 0}$, taking values in $\{0,1\}^{\Z/N\Z}$, as follows. First, let $j=1$, and define
\begin{equation}\label{def:PatchedProcess}
\eta^{\textup{patch},N,i}_t(x) := \begin{cases} \eta_t^{N,i,1,1}(x) & \text{ if } (x,t) \in \mathcal{I}_{\delta,N}^{\rho}\big(0,\frac{1}{2}\big) ,   \\
\eta_t^{N,i,2,1}(x) & \text{ if } (x,t) \in \mathcal{I}_{\delta,N}^{\rho}\big(\frac{1}{2},1\big) , 
\end{cases}
\end{equation} where the points $x$ are understood modulo $N$. 
In order to define $\eta^{\textup{patch},N,i}_t$ for $t>\delta N^{3/2}$, we proceed inductively. Suppose that we have constructed $(\eta^{\textup{patch},N,i}_t)_{t \in [0,j \delta N^{3/2}]}$ for some $j\in \N$ (the base case $j=1$ is covered above). Then consider two families of ASEPs $(\eta^{N,i,1,j+1}_t)_{t \geq 0}$ and $(\eta^{N,i,2,j+1}_t)_{t \geq 0}$, started at time $j\delta N^{3/2}$, where assign for all $i\in \lbr M \rbr$
\begin{equation}\label{eq:NewStart}
\begin{split}
        \eta^{N,i,1,j+1}_{j\delta N^{3/2}}(x) &= \bar{\eta}^{N,i}_{j\delta N^{3/2}}(x) \ \forall x \colon (x,j\delta N^{3/2}) \in \mathcal{I}_{\delta,N}^{\rho}\Big(-\frac{1}{8},\frac{5}{8}\Big) , \\
        \eta^{N,i,2,j+1}_{j\delta N^{3/2}}(x)&=  \bar{\eta}^{N,i}_{j\delta N^{3/2}}(x) \ \forall x \colon (x,j\delta N^{3/2}) \in \mathcal{I}_{\delta,N}^{\rho}\Big(\frac{3}{8},\frac{9}{8}\Big)  .
\end{split}
\end{equation} For all other $x$, we assign independent Bernoulli-$\rho$-random variables.
Then $(\eta^{N,i,1,j+1}_t)_{t \geq 0}$ and $(\eta^{N,i,2,j+1}_t)_{t \geq 0}$ evolve according to the patched coupling between time $j\delta  N^{3/2}$ and $(j+1)\delta N^{3/2}$, independently of the processes up to time $\delta j N^{3/2}$. Set for all $s \in [j \delta N^{3/2},(j+1)\delta N^{3/2}]$,
\begin{equation}\label{def:PatchedProcess2}
\eta^{\textup{patch},N,i}_{s}(x) := \begin{cases} \eta_{s}^{N,i,1,j+1}(x) & \text{ if } (x,s) \in \mathcal{I}_{\delta,N}^{\rho}\big(0,\frac{1}{2}\big) ,   \\
\eta_s^{N,i,2,j+1}(x) & \text{ if } (x,s) \in \mathcal{I}_{\delta,N}^{\rho}\big(\frac{1}{2},1\big) . 
\end{cases}
\end{equation} 
\end{definition}
We record the following coupling between the $(\delta,N)$-patched ASEP and the periodic ASEP.  

\begin{lemma}\label{lem:CoupleASEPperiodic}
Fix $M\in \N$ and consider a family of $(\delta,N)$-patched ASEPs $\big(\eta^{\textup{patch},N,i}_t\big)_{t \geq 0}$ for $i\in \lbr M \rbr$. Let $\big(\bar{\eta}^{N,i}_t\big)_{t \geq 0}$ be a family of unwrapped $N$-periodic ASEPs such that for all $i\in \lbr M \rbr$ and $N\in \N$
\begin{equation}
    \eta^{\textup{patch},N,i}_0 =  \bar{\eta}^{N,i}_0  
\end{equation} almost surely. Fix some $T>0$. Assume that \eqref{eq:DominationAssumptionPatch} holds for the processes in the definition of $\big(\bar{\eta}^{N,i}_t\big)_{t \geq 0}$. Then there exists a coupling $\mathbf{P}_{\ast}$ and  $c_1,C_1>0$ such that for all $\delta>0$
\begin{align}\label{eq:CouplingCircle}
\liminf_{N \rightarrow \infty}\mathbf{P}_{\ast}\left( \bar{\eta}^{N,i}_t = \eta^{\textup{patch},N,i}_t  \ \forall t \in [0,T N^{3/2}] , \,  i \in \lbr M \rbr\right) \geq 1- C_1\exp(-c_1\delta^{-1}) . 
\end{align}  
\end{lemma}
\begin{proof} We start by coupling the underlying full-space ASEPs in the definition of the $(\delta,N)$-patched ASEP
as in \eqref{eq:Extended1} and \eqref{eq:Extended2} with a periodic ASEP with respect to the basic coupling via a suitable reference frame. More precisely, \eqref{eq:Extended1} and \eqref{eq:Extended2} ensure that we find  $c_0,C_0>0$ and a coupling $\mathbf{P}_{\ast}$ such that for all $\delta>0$
\begin{align}\label{eq:CouplingCircleNew}
\liminf_{N \rightarrow \infty}\mathbf{P}_{\ast}\left( \bar{\eta}^{N,i}_t = \eta^{\textup{patch},N,i}_t  \ \forall t \in [0,\delta N^{3/2}] \, , i \in \lbr M \rbr\right) \geq 1- \frac{1}{\lceil T\delta^{-1}\rceil}C_0\exp(-c_0\delta^{-1}) . 
\end{align} 
Using attractivity \eqref{def:attactivity} of the basic coupling, we see that the configurations in the construction of $\bar{\eta}^{N,i}_{\delta N^{3/2}}$ satisfy again assumption \eqref{eq:DominationAssumptionPatch} for all $i\in \lbr M\rbr$. Iterating \eqref{eq:CouplingCircleNew} now $ \lceil T\delta^{-1} \rceil$ many times, and applying a union bound, we conclude \eqref{eq:CouplingCircle}.
\end{proof}

Let $(h_{j\delta N^{3/2}}^{i,1,j})_{i \in \lbr M \rbr}$ and $(h_{j\delta N^{3/2}}^{i,2,j})_{i \in \lbr M \rbr}$ denote the Bernoulli paths associated to the configurations $\eta_{j\delta N^{3/2}}^{N,i,1,j}$ and $\eta_{j\delta N^{3/2}}^{N,i,2,j}$ in the definition of the $(\delta,N)$-patched ASEP. 
Let $\hA_{1,j}$ and $\hA_{2,j}$ denote the height functions from Definition~\ref{def:HeightASEP} for the ASEPs $(\eta_t^{N,i,1,j})_{t \geq 0}$ and $(\eta_t^{N,i,2,j})_{t \geq 0}$ under $\mathbf{P}_{\textbf{pat}}$,  and by $\hApeps_{1,j}$ and $\hApeps_{2,j}$  their rescaled height functions from Definition~\ref{def:ASEPheightRescaledPeriodic}.

\begin{definition}\label{def:PatchedHeightFunctions}
Fix $\delta>0$ and $M\in \N$, and assume that the Bernoulli paths $(h_0^{i,1,1})_{i \in \lbr M \rbr}$ and $(h_0^{i,2,1})_{i \in \lbr M \rbr}$ are all $(N,k)$-periodic and satisfy  $h_0^{i}:=h_0^{i,1,1}=h_0^{i,2,1}$ for all $i \in \lbr M \rbr$. We define for all $t\in [0,\delta N^{3/2}]$ and $i \in \lbr M \rbr$ the $\boldsymbol{(\delta,N)}$\textbf{-patched ASEP height function} 
\begin{equation}\label{def:PatchedHeightFunction}
h^{\textup{patch},\delta}(h^{i}_0;y,t) := \begin{cases} \hA_{1,1}\big(h_0^{i,1,1},0;y,t\big) & \text{ if } (y,t) \in \mathcal{I}_{\delta,N}^{\rho}\big(0,\frac{1}{2}\big) ,   \\
\hA_{2,1}(h_0^{i,2,1},0;y,t) & \text{ if } (y,t) \in \mathcal{I}_{\delta,N}^{\rho}\big(\frac{1}{2},1\big) . 
\end{cases}
\end{equation} For all $t\geq 0$ and $i\in \lbr M \rbr$, recall $f(t)$ from \eqref{def:Functionf}, and let the number of particles between $f(t)+\frac{1}{4}N$ and $f(t)+\frac{5}{4}N$ for the patched ASEP with initial conditions $h_0^{i}$ be defined as
\begin{equation}
    K^{i}_t := h^{\textup{patch},\delta}\Big(h^{i}_0;f(t),t\Big) -  h^{\textup{patch},\delta}\Big(h^{i}_0;f(t)+N-1,t\Big) . 
\end{equation}
We now obtain $h^{\textup{patch},\delta}$ for all $(y,t) \notin \mathcal{I}_{\delta,N}^{\rho}\big(0,1\big)$ by extension so that $h^{\textup{patch},\delta}(h^{i}_0;\cdot,t)$ is an $(N,K_t^{i})$-periodic Bernoulli path. Now suppose that we have constructed the $(\delta,N)$-patched ASEP height function until time $j\delta N^{3/2}$ (the base case $j=1$ is covered above). Then we set 
\begin{equation}\label{def:PatchedHeightFunctionGeneral}
h^{\textup{patch},\delta}(h^{i}_0;y,t) := \begin{cases} \hA_{1,j}\big(h_{j\delta N^{3/2}}^{i,1,j},0;y,t\big) & \text{ if } (y,t) \in \mathcal{I}_{\delta,N}^{\rho}\big(0,\frac{1}{2}\big) ,   \\
\hA_{2,j}(h_{j\delta N^{3/2}}^{i,2,j},0;y,t) & \text{ if } (y,t) \in \mathcal{I}_{\delta,N}^{\rho}\big(\frac{1}{2},1\big) , 
\end{cases}
\end{equation}
for all $(y,t) \in \mathcal{I}_{\delta,N}^{\rho}\big(0,1\big)$ with $t \in [j\delta N^{3/2},(j+1) \delta N^{3/2}]$, and define $h^{\textup{patch},\delta}(h^{i}_0;y,t)$ outside of $\mathcal{I}_{\delta,N}^{\rho}\big(0,1\big)$ by extension to a periodic path as for the base case $j=1$.
\end{definition}

It will be convenient in the following to consider the sigma-algebra $\sigma(h^{\textup{patch},\delta})$ generated by the ASEP height functions $(\hA_{1,j},\hA_{2,j})$ in Definition~\ref{def:PatchedHeightFunctions}.

\begin{definition}\label{def:PatchedHeightFunctionsRescaled}
Recall that $\rho=k/N$, $\alpha=1-2\rho$, as well as $\varepsilon=\varepsilon(N,k)=N^{-3/2}\beta^{3/2}$ for $\beta=\beta(\alpha)$ from \eqref{def:betaASEP}. For all $i\in \lbr M \rbr$, $z\in \R$, $t \geq 0$, and periodic initial data $ h^{\per,\varepsilon}_0$, we define the \textbf{rescaled patched ASEP height function} as 
\begin{align*}
\hfix^{\textup{patch},\delta,\varepsilon}( \hfix^{\per,\varepsilon}_0;z,t) := \sigma(\alpha)^{-1}\varepsilon^{\frac{1}{3}}&\Big(2\mu(\alpha)\varepsilon^{-1}t +\mu^{\prime}(\alpha)\beta(\alpha)z\varepsilon^{-\frac{2}{3}} \\
&- h^{\textup{patch},\delta}( h^{\per,\varepsilon}_0,0;2\alpha t \varepsilon^{-1} + \beta(\alpha)z\varepsilon^{-\frac{2}{3}},2(1-q)^{-1}\varepsilon^{-1}t) \Big) 
\end{align*} whenever the argument for $h^{\textup{patch},\delta}$ is an integer, and linear interpolation otherwise. Here, we obtain $\hfix^{\per,\varepsilon}_0$ from $h^{\per,\varepsilon}_0$ as in \eqref{def:RescaledPeriodicASEPheight}.
\end{definition}

\begin{remark}\label{rem:PeriodicPatchedHeight}
    Note that if the event in \eqref{eq:CouplingCircle} occurs, then $K_t^{i}=K_0^{i}=k$ for all $i\in \lbr M\rbr$ and $t \in [0,cT]$, with $c=\beta^{\frac{3}{2}}(1-q)$. In particular, as $k = \rho N = \mu'(\alpha) N$, we get that
\begin{equation}
    \hfix^{\textup{patch},\delta,\varepsilon}(\hfix^{\per,\varepsilon}_0;z,t)=\hfix^{\textup{patch},\delta,\varepsilon}(\hfix^{\per,\varepsilon}_0;z+1,t)
\end{equation} for all $z\in \R$ and (rescaled) periodic initial data $\hfix^{\per,\varepsilon}_0$. 
\end{remark}

We record the following result on the total-variation distance between the rescaled patched ASEP height functions and the rescaled periodic ASEP height functions.

\begin{corollary}\label{cor:CouplePatchedPeriodicASEPs}
 Fix $T'\geq 0$, and initial data $h_0^{\per,(i),\varepsilon}$ for all $i\in \lbr M \rbr$. Then under the same assumptions as in Lemma~\ref{lem:CoupleASEPperiodic}, we find for every $\varepsilon'>0$ some $\delta_0=\delta_0(\varepsilon',T')$ so that 
\begin{equation}\label{eq:CloseToPeriodicASEPPatch}
     \bigg\lVert \Big(\hApeps(\hfix_0^{\per,(i),\varepsilon};\cdot,t)\Big)_{t \in [0,T'],i \in \lbr M \rbr} -  \Big(\hfix^{\textup{patch},\delta,\varepsilon}(\hfix_0^{\per,(i),\varepsilon};\cdot,t)\Big)_{t \in [0,T'],i \in \lbr M \rbr}\bigg\rVert_{\textup{TV}}\leq \varepsilon^{\prime} 
\end{equation} for all $\delta \leq \delta_0$, and $\varepsilon=\varepsilon(\delta)>0$ small enough. 
\end{corollary}
\begin{proof}
 This is immediate from Lemma~\ref{lem:CoupleASEPperiodic} with $T=T'\beta^{3/2}(1-q)$, and the coupling representation of the total-variation distance. 
\end{proof}

\subsection{Convergence of the patched ASEP to a patched KPZ fixed point}\label{sec:ConvergencePatchedASEP}

In the following, we establish a convergence result for the rescaled patched ASEP height function. 
Let $n\in \N$, and consider the $(n^{-1},N)$-patched ASEP on the enlarged probability space with $\sigma(h^{\textup{patch},\delta})$ given below Definition~\ref{def:PatchedHeightFunctions}. We denote the law of the rescaled height function $\hfix^{\textup{patch},n^{-1},\varepsilon}$ on this space by $\mathbb{P}^{n}_{\textup{pat}}$, and write $\mathbb{E}^{n}_{\textup{pat}}$ for its expectation. Let $\mathcal{T}$ denote an arbitrary countable subset of $[0,1]$, and recall that by Remark~\ref{rem:PeriodicPatchedHeight}, we can treat the functions $y \mapsto  \hfix^{\textup{patch},n^{-1},\varepsilon}(\hfix^{\per,(i),\varepsilon}_0;y,t)$ for some $t\geq 0$, as well as the initial data $\hfix^{\per,(i),\varepsilon}$ as elements of $\UCp$, provided that the event in \eqref{eq:CouplingCircle} occurs. We make the following definition. 


\begin{definition} \label{def:PatchedCoupledDLConvergence}
Fix $M\in \N$ and let $T>0$. Assume that for all $i\in \lbr M\rbr$, we have some $\hfix^{\per,(i)}_0 \in \UCp$. Consider the random function $(y,t) \mapsto \hfixpa(\mathfrak{h}^{\per,(i)}_0;y,t)$, where for all $y\in \mathbb{T}$, all $t \in [0,T]$, and all initial data $\mathfrak{h}^{\per,(i)}_0$,
    \begin{equation}\label{def:PatchedKPZFixedPoint}
         \hfixpa\big(\hfix^{\per,(i)}_0;y,t\big) :=  \sup_{x\in \mathbb{T}} \big(\hfix^{\per,(i)}_0(x) + \Lpa(x,0;y,t)\big). 
    \end{equation} Here, $\Lpa$ is a wrapped patched directed landscape from Definition~\ref{def:PatchedDLExt}. 
    We refer to  $(y,t) \mapsto \hfixpa(\mathfrak{h}^{\per}_0;y,t)$ as a \textbf{patched KPZ fixed point} of scale $n$ and with periodic initial data $\hfixp_0$.
\end{definition}

 Let us stress that the construction of the patched KPZ fixed point depends on the choice of the coupling to construct $\Lpa$.
In order to state the convergence of the rescaled patched ASEP height function to the patched KPZ fixed point, we recall and introduce some notation. As in  \eqref{eq:GoodCoupling}, we define for directed landscapes $\mathcal{L}^{1,j}$ and $\mathcal{L}^{2,j}$ with $j\in \N$ 
\begin{equation}\label{eq:GoodCouplingASEP}
 \mathcal{A}^{j}_n:=\left\{  \mathcal{L}^{1,j}_{\mathcal{O}^{1,j}_n} =   \mathcal{L}^{2,j}_{\mathcal{O}^{3,j}_n} \right\} \cap\left\{ \mathcal{L}^{1,j}_{\mathcal{O}^{2,j}_n} =   \mathcal{L}^{2,j}_{\mathcal{O}^{2,j}_n}  \right\} , 
\end{equation}
where we recall $\mathcal{O}^{i,j}_n$ from \eqref{def:RectanglesDL}.
For the $(n^{-1},N)$-patched ASEP, we denote for all $j\in \N$ by $\hAeps_{1,j}$ and $\hAeps_{2,j}$ the rescaled height functions of the families of ASEPs at the right-hand side \eqref{def:PatchedProcess} and \eqref{def:PatchedProcess2} of Definition~\ref{def:PatchedASEPNew} under the patched coupling.

\begin{lemma}\label{lem:patchedASEPfixed1} Fix some $M\in \N$ and $n\in \N$. 
Assume that for all $i\in \lbr M\rbr$, we have  $\hfix^{\per,(i),\varepsilon}_0$ in $\UCp$ converging to some $\hfix^{\per,(i)}_0 \colon \mathbb{T}\rightarrow \R$ as $\varepsilon \rightarrow 0$ with respect to the UC topology, and initial data $\hfix^{(i),1,\varepsilon,j}_0$ and $\hfix^{(i),2,\varepsilon,j}_0$ in $\UC$ with
\begin{equation}
    \begin{split}
        \hfix^{(i),1,\varepsilon,j}_0(x) &= \hfix^{\textup{patch},n^{-1},\varepsilon}\big(\hfix^{\per,(i),\varepsilon}_0;x,(j-1)n^{-1}\big) \ \forall x\in \Big[ -\frac{1}{8}, \frac{5}{8}\Big] \\
        \hfix^{(i),2,\varepsilon,j}_0(x) &= \hfix^{\textup{patch},n^{-1},\varepsilon}\big(\hfix^{\per,(i),\varepsilon}_0;x,(j-1)n^{-1}\big) \ \forall x\in \Big[ \frac{3}{8}, \frac{9}{8}\Big] , 
    \end{split}
\end{equation}
 such that $\hfix^{(i),1,\varepsilon,j}_0$ and $\hfix^{(i),2,\varepsilon,j}_0$ converge to some $\hfix^{(i),1,j}_0$ and $\hfix^{(i),2,j}_0$ in $\UC$. Then for all $j\in \N$, there exist full-space KPZ fixed points $\hfix^{1,j}$ and $\hfix^{2,j}$ such that as $\varepsilon \rightarrow 0$
\begin{equation}\label{eq:ASEPconvPatchtwice}
\begin{split}
    \hAeps_{1,j}\big( \hfix^{(i),1,\varepsilon,j}_0;\cdot,t\big) &\rightarrow \hfix^{1,j}\big(\hfix^{(i),1,j}_0;\cdot,t\big) ,  \\
    \hAeps_{2,j}\big( \hfix^{(i),2,\varepsilon,j}_0;\cdot,t\big) &\rightarrow \hfix^{2,j}\big(\hfix^{(i),2,j}_0;\cdot,t\big) , 
\end{split}
\end{equation} jointly in $i\in \lbr M \rbr$, and $t \in [0,n^{-1}] \cap \mathcal{T}$ with respect to the topology of uniform convergence on compact subsets. On the left sides of \eqref{eq:ASEPconvPatchtwice}, 
the functions $\hAeps_{1,j}\big(\hfix_0^{(i),1,\varepsilon,j};\cdot,t\big)$ and $\hAeps_{2,j}\big(\hfix_0^{(i),2,\varepsilon,j};\cdot,t\big)$ are coupled according to $\Ppatch$. On the right sides of \eqref{eq:ASEPconvPatchtwice}, the functions $\hfix^{1,j}(\hfix^{(i),1,j}_0;\cdot,t)$ are coupled through \eqref{eq:VariationFixedPeriodic} for all $i\in \lbr M \rbr$  under the same directed landscape $\mathcal{L}^{1,j}$; the functions $\hfix^{2,j}(\hfix^{(i),2,j}_0;\cdot,t)$ are coupled through \eqref{eq:VariationFixedPeriodic} for all $i\in \lbr M \rbr$  under the same directed landscape $\mathcal{L}^{2,j}$, and $\mathcal{L}^{1,j}$ and  $\mathcal{L}^{2,j}$  are coupled under some probability measure $\mathbf{P}$, so that
\begin{equation}\label{eq:QuantitativeDLs}
  \mathbf{P}\big(  \mathcal{A}^{j}_{n} \big) \geq 1 - C_0\exp(-c_0 n ) 
\end{equation}
holds for some $c_0,C_0>0$, which do not depend on $j$ and $n$, and all $n\in \N$. 
\end{lemma}

\begin{proof} We will in the following only consider $j=1$ as the case $j>1$ follows similarly. 
Note that we can extend $\hfix^{\per,(i),\varepsilon}_0$ periodically in order to obtain the initial data $\hfix^{(i),1,\varepsilon,1}_0$ and $\hfix^{(i),2,\varepsilon,1}_0$. 
The convergence in \eqref{eq:ASEPconvPatchtwice} to two families of KPZ fixed points is immediate from Theorem~\ref{thm:ACHmain}, noting that $\hfix^{(i),1,\varepsilon,1}_0$ and $\hfix^{(i),2,\varepsilon,1}_0$ satisfy assumptions \eqref{eq:ExponentialGrowthInitial1} and \eqref{eq:ExponentialGrowthInitial2}. 
To show that
the respective directed landscapes of these KPZ fixed points must satisfy \eqref{eq:QuantitativeDLs}, let $(\eta_t^{N,i,1})_{t \geq 0}$ and $(\eta_t^{N,i,2})_{t \geq 0}$ denote the ASEPs belonging to $\hAeps_{1,j}$ and $\hAeps_{2,j}$. We define an ASEP $(\eta_t^{N,i,3})_{t \geq 0}$ by
\begin{equation}
    \eta_t^{N,i,3} (x) := \eta_t^{N,i,2} (x+N)  
\end{equation} for all $x \in \Z$. Moreover, for all $i \in \lbr M\rbr$, consider the ASEPs $(\tilde{\eta}_t^{N,i,2})_{t \geq 0}$ and $(\tilde{\eta}_t^{N,i,3})_{t \geq 0}$ with initial states $\eta_0^{N,i,2}$ and $\eta_0^{N,i,3}$, respectively, such that $(\eta_t^{N,i,1})_{t \geq 0}$, $(\tilde{\eta}_t^{N,i,2})_{t \geq 0}$ and $(\tilde{\eta}_t^{N,i,3})_{t \geq 0}$ jointly evolve together according to the basic coupling. 
This construction allows us to apply the KPZ fixed point convergence in Theorem \ref{thm:ACHmain} jointly for all $3M$ many ASEPs $(\eta_t^{N,i,1})_{t \geq 0}$, $(\tilde{\eta}_t^{N,i,2})_{t \geq 0}$ and $(\tilde{\eta}_t^{N,i,3})_{t \geq 0}$ with $i\in \lbr M\rbr$ with respect to the same directed landscape. 
We will now transfer this result to the ASEPs $(\eta_t^{N,i,1})_{t \geq 0}$ and $(\eta_t^{N,i,2})_{t \geq 0}$. Note that under the above coupling, $(\eta_t^{N,i,2})_{t \geq 0}$ and $(\tilde{\eta}_t^{N,i,2})_{t \geq 0}$ evolve according to the $(\rho,0,1,N)$-reference frame basic coupling, and similarly for $(\eta_t^{N,i,3})_{t \geq 0}$ and $(\tilde{\eta}_t^{N,i,3})_{t \geq 0}$. This can be seen directly by verifying the marginal transition rates.  In view of Proposition~\ref{pro:ModerateDeviationsSecondClass} for two ASEPs on the integers, and the first statement of  \eqref{eq:GoodCouplingPatched} in Lemma~\ref{lem:patchedASEP}, we can apply Proposition~\ref{pro:ModerateDeviationsSecondClass} with $\theta=-\frac{1}{4}$ for the processes $(\eta_t^{N,i,2})_{t \geq 0}$ and $(\tilde{\eta}_t^{N,i,2})_{t \geq 0}$, as well as for the processes $(\eta_t^{N,i,1})_{t \geq 0}$ and $(\eta_t^{N,i,2})_{t \geq 0}$, to see that for some coupling $\mathbf{P}_{\ast}$
\begin{equation}\label{eq:ShiftPatch}
\mathbf{P}_{\ast} \left(\eta^{N,i,1}_t(x)= \eta_t^{N,i,2} (x) =\tilde{\eta}_t^{N,i,2} (x) \ \forall (x,t) \in \mathcal{I}^{\rho}_{n^{-1},N}\Big(\frac{3}{8},\frac{5}{8}\Big)\right) \geq 1- C_1\exp(-c_1 n)  
\end{equation}
with some $c_1,C_1>0$ and all $N$ large enough. Similarly, we get that 
\begin{equation}\label{eq:ShiftPatch2}
\mathbf{P}_{\ast} \left(\eta^{N,i,1}_t(x)= \eta_t^{N,i,3} (x) =\tilde{\eta}_t^{N,i,3} (x)\ \forall (x,t) \in \mathcal{I}^{\rho}_{n^{-1},N}\Big(-\frac{1}{8},\frac{1}{8}\Big)\right) \geq 1- C_1\exp(-c_1 n) 
\end{equation} for all $N$ large enough.
Now Theorem~\ref{thm:ACHmain} guarantees for the exclusion processes $(\eta_t^{N,i,1})_{t \geq 0}$, $(\tilde{\eta}_t^{N,i,2})_{t \geq 0}$ and $(\tilde{\eta}_t^{N,i,3})_{t \geq 0}$ with $i \in \lbr M \rbr$ that the respective rescaled height functions converge to KPZ fixed points, coupled with respect to the same directed landscape. Hence, when the events in  \eqref{eq:ShiftPatch} and \eqref{eq:ShiftPatch2} hold, the event $\mathcal{A}^1_n$ occurs with respect to the directed landscapes $\mathcal{L}^{1,1}$ and $\mathcal{L}^{2,1}$ in the definition of $\hfix^{1,1}$ and $\hfix^{2,1}$, respectively, allowing us to conclude \eqref{eq:QuantitativeDLs}.
\end{proof}


\begin{lemma}\label{lem:patchedASEPfixed2} 
Fix some $M,K\in \N$ and let  $n\in \N$. Assume that for all $i\in \lbr M\rbr$ and $\varepsilon>0$, we have some $\hfix^{\per,(i),\varepsilon}_0$ in $\UCp$ so that the associated particle configurations satisfy \eqref{eq:DominationAssumption} at time $0$, and such that the rescaled height functions $\hfix^{\textup{patch},n^{-1},\varepsilon}(\hfix^{\per,(i),\varepsilon}_0;\cdot,0)$ of the $(n^{-1},N)$-patched ASEP at time $0$ converge to $\hfix^{\per,(i)}_0 \colon \mathbb{T}\rightarrow \R$ as $\varepsilon \rightarrow 0$ with respect to the UC topology. 
Then there exist $c,C>0$, not depending on $n$ or $\varepsilon$,  such that for any $ t_1 <t_2<\dots<t_{K}$ in $\mathcal{T}$,
with probability $1-C \exp(-c n^{1/2})$, we find a coupling between $\hfix^{\textup{patch},n^{-1},\varepsilon}$ and $\hfixpa$ so that
\begin{equation}
  \max_{i \in \lbr M \rbr} \max_{\ell \in \lbr K \rbr}  \sup_{y\in \mathbb{T}} \left|\hfix^{\textup{patch},n^{-1},\varepsilon}(\mathfrak{h}^{\per,(i),\varepsilon}_0;y,t_\ell)- \hfixpa(\mathfrak{h}^{\per,(i)}_0;y,t_\ell) \right|   \leq (2n+1) 2^{-n}
\end{equation} for all $\varepsilon>0$ sufficiently small.
Here, the functions on the left-hand side are coupled according to the patched coupling $\Ppatch$ for the underlying ASEPs in the definition of the patched ASEPs, while the functions on the right-hand side are coupled as in Definition~\ref{def:PeriodicKPZfixed} with respect to the same patched directed landscape $\Lpa$, taking values in $\mathcal{T}$ for the time component.
\end{lemma}
\begin{proof}  
Let $K,M \in \N$ and note that the rescaled patched ASEP height function $\hfix^{\textup{patch},n^{-1},\varepsilon}$ from Definitions~\ref{def:PatchedHeightFunctions} and~\ref{def:PatchedHeightFunctionsRescaled} can be expressed for $t \in[0,n^{-1}]$ in terms of the ASEP height functions $\hAeps_{1,1}$ and $\hAeps_{2,1}$. We start with the case where  $t_1 <t_2<\dots<t_K$ are in $\mathcal{T} \cap [0,n^{-1}]$.  Lemma~\ref{lem:patchedASEPfixed1} together with Skorokhod's representation theorem yields that there exist $c_0,C_0>0$, not depending on $n$ or $\varepsilon$, and some $\hfixpa_{\ast}$ such that with probability $1-C_0 \exp(-c_0n^{1/2})$, we find a coupling so that  for all $t_1 <t_2<\dots<t_K$ in $\mathcal{T} \cap [0,n^{-1}]$,
\begin{equation}\label{eq:InitialConvergencePatchedVariaint}
 \max_{i \in \lbr M \rbr} \max_{\ell \in \lbr K \rbr}  \sup_{y\in \mathbb{T}} \left|\hfix^{\textup{patch},n^{-1},\varepsilon}(\mathfrak{h}^{\per,(i),\varepsilon}_0;\cdot,t_\ell)- \hfixpa_{\ast}(\mathfrak{h}^{\per,(i)}_0;\cdot,t_\ell) \right|   \leq 2^{-n}
\end{equation} for all $n\in N$, and all $\varepsilon>0$ sufficiently small.
More precisely, $\hfixpa_{\ast}$ is given for all $i \in \lbr M \rbr$ by the variational characterization
\begin{equation}\label{eq:VarCharHN}
    \hfixpa_{\ast}(\hfix^{\per,(i)}_0; y,t) = \begin{cases} \sup_{z\in \R} \hfix^{\per,(i)}_0(z)+ \mathcal{L}^{1,1}(z,0;y,t)
        & \text{ if } y \in \big[0, \frac{1}{2}\big) ,  t \in [0,n^{-1}] \cap \mathcal{T}, \\
    \sup_{z\in \R} \hfix^{\per,(i)}_0(z)+ \mathcal{L}^{2,1}(z,0;y,t)
        & \text{ if } y \in \big[\frac{1}{2}, 1\big) , t \in [0,n^{-1}] \cap \mathcal{T},
    \end{cases}
\end{equation} as we use the evolution in $\hAeps_{1,1}$ for all $y \in \big[0, \frac{1}{2}\big)$, and the evolution of  $\hAeps_{2,1}$ for all $y \in \big[\frac{1}{2}, 1\big)$ in the construction of $\hfix^{\textup{patch},n^{-1},\varepsilon}$, and $\mathcal{L}^{1,1}$ and $\mathcal{L}^{2,1}$ are coupled as in Definition~\ref{def:PeriodicKPZfixed}. We now show that for some $c_1,C_1>0$ and all $n\in \N$
\begin{equation}\label{eq:TVboundPatched}
 \left\lVert  \Big(\hfixpa_{\ast}\big(\hfix^{\per,(i)}_0;\cdot,\tau\big) \Big)^{i \in \lbr M \rbr}_{\tau \in \big[0,\frac{1}{n}\big] \cap \mathcal{T}}  - \Big(\hfixpa\big(\hfix^{\per,(i)}_0;\cdot,\tau\big) \Big)^{i \in \lbr M \rbr}_{\tau \in \big[0,\frac{1}{n}\big] \cap \mathcal{T}}  \right\rVert_{\textup{TV}} \leq C_1 \exp\big(-c_1 n^{\frac{1}{2}}\big) , 
\end{equation}
where $\hfixpa$ is a patched KPZ fixed point with respect to $\mathcal{L}^{1,1}$ and $\mathcal{L}^{2,1}$. Then \eqref{eq:InitialConvergencePatchedVariaint} together with \eqref{eq:TVboundPatched} yields that for some  $c_2,C_2>0$, not depending on $n$ or $\varepsilon$, with probability $1-C_2 \exp(-c_2n^{1/2})$, we find a coupling so that  for all $t_1 <t_2<\dots<t_K$ in $\mathcal{T} \cap [0,n^{-1}]$, 
\begin{equation}\label{eq:FirstIntervalBound}
 \max_{i \in \lbr M \rbr}  \max_{\ell \in \lbr K \rbr}  \sup_{y\in \mathbb{T}} \left|\hfix^{\textup{patch},n^{-1},\varepsilon}(\mathfrak{h}^{\per,(i),\varepsilon}_0;\cdot,t_\ell)- \hfixpa(\mathfrak{h}^{\per,(i),\varepsilon}_0;\cdot,t_\ell) \right|   \leq  2^{-n}
\end{equation} for all $n\in N$, and all $\varepsilon>0$ sufficiently small.
Note that for \eqref{eq:TVboundPatched}, it suffices to show 
\begin{equation*}
 \left\lVert  \Big(\hfixpa_{\ast}\big(\hfix^{\per,(i)}_0;y,\tau\big) \Big)^{i \in \lbr M \rbr}_{y\in \big[0,1\big),\tau \in \big[0,\frac{1}{n}\big] \cap \mathcal{T}}  - \Big(\hfixpa\big(\hfix^{\per,(i)}_0;y,\tau\big) \Big)^{i \in \lbr M \rbr}_{y\in \big[ 0,1\big),\tau \in \big[0,\frac{1}{n}\big] \cap \mathcal{T}}  \right\rVert_{\textup{TV}} \leq C_3 \exp\big(-c_3 n^{\frac{1}{2}}\big) , 
\end{equation*}
for some  $c_3,C_3>0$ and $n\in \N$. To see this, note that by Remark \ref{rem:PeriodicPatchedHeight}, the rescaled patched ASEP height function is $1$-periodic up to time $n^{-1}$ on an event of probability $1-C_4\exp(-c_4 n)$ for some  $c_4,C_4>0$, while the periodic KPZ fixed point is $1$-periodic by construction.  
Recall the set of maximizers $\Sigma_1(y,t)$ from \eqref{def:LocalExtremes} in the directed landscape $\mathcal{L}^{1,1}$ for all $(y,t) \in \R \times (0,\infty)$. Similarly, $\Sigma_2(y,t)$ denotes the set of maximizers in the directed landscape $\mathcal{L}^{2,1}$. We write $\pi^{1,1,-},\pi^{2,1,-}$ and $\pi^{1,1,+},\pi^{2,1,+}$ for the leftmost and rightmost geodesics in $\mathcal{L}^{1,1}$ and $\mathcal{L}^{2,1}$, respectively. Recalling Definition~\ref{def:PatchedDLExt}, $\hfixpa_{\ast}$ is a patched KPZ fixed point when $ \mathcal{A}^1_n$ from \eqref{eq:GoodCouplingASEP} and the events
\begin{equation*}
\begin{split}
     \tilde{\mathcal{B}}_{n}^{1} &:= \left\{   \pi^{1,1,-}_{(x,0;y,t)},\pi^{1,1,+}_{(x,0;y,t)} \subseteq \Big[-\frac{1}{8}, \frac{5}{8}\Big)\times [0,t]\ \forall x \in \Sigma_1(y,t), y \in \Big[0, \frac{1}{2}\Big] \ , t \in [0,n^{-1}] \cap \mathcal{T} \right\} \\
          \tilde{\mathcal{B}}_{n}^{2} &:= \left\{ \pi^{2,1,-}_{(x,0;y,t)},\pi^{2,1,+}_{(x,0;y,t)} \subseteq \Big[\frac{3}{8}, \frac{9}{8}\Big)\times [0,t]\ \forall x \in \Sigma_2(y,t), y \in \Big[\frac{1}{2},1\Big] \ , t \in [0,n^{-1}] \cap \mathcal{T} \right\}
\end{split}
\end{equation*}
occur. 
Note that by assumption \eqref{eq:DominationAssumption}, the initial data $\hfix^{\per,(i)}_0$ satisfies for all $i\in \lbr M \rbr$
\begin{equation}\label{eq:MinMax}
  -\infty <  \inf_{x \in \mathbb{T}} \hfix_0^{\per,(i)}(x) \leq \sup_{x \in \mathbb{T}} \hfix_0^{\per,(i)}(x) < \infty .
\end{equation}  
By Proposition~\ref{pro:ModulusOfContinuity} to get an upper bound on $\mathcal{L}^{1,1}(x,0;y,t)$ and $\mathcal{L}^{2,1}(x,0;y,t)$ when $|x-y| \geq \frac{1}{16}$, and a lower bound when $|x-y| \leq \frac{1}{32}$, we see that there exist $c_5,C_5>0$, depending only on $(\hfix^{\per,(i)}_0)_{i \in \lbr M \rbr}$, such that 
\begin{equation*}
\begin{split}
     \hat{\mathcal{B}}_{n}^{1} &:= \left\{   \Sigma_1(y,t) \subseteq \Big[-\frac{1}{16}, \frac{9}{16}\Big)\ \forall y \in \Big[0, \frac{1}{2}\Big] \, , t \in [0,n^{-1}] \cap \mathcal{T} \right\} \\
          \hat{\mathcal{B}}_{n}^{2} &:= \left\{  \Sigma_2(y,t) \subseteq  \Big[\frac{7}{16}, \frac{17}{16}\Big) \ \forall y \in \Big[\frac{1}{2}, 1\Big] \, , t \in [0,n^{-1}] \cap \mathcal{T} \right\}
\end{split}
\end{equation*}
satisfy for all $n \in \N$
\begin{equation}\label{eq:LocateStartPreFinal} 
       \P\left(\hat{\mathcal{B}}_{n}^{1} \cap \hat{\mathcal{B}}_{n}^{2} \right) \geq 1- C_5\exp(-c_5n^{1/2}) .  
\end{equation}
Together with Proposition~\ref{pro:ModerateDL} to bound the transversal fluctuations of geodesics, there exist 
$c_6,C_6>0$, depending only on $(\hfix^{\per,(i)}_0)_{i \in \lbr M \rbr}$, such that for all $n \in \N$
\begin{equation}\label{eq:LocateStartFinal} 
       \P\left(\tilde{\mathcal{B}}_{n}^{1} \cap \tilde{\mathcal{B}}_{n}^{2} \right) \geq 1- C_6\exp(-c_6n^{1/2}) ,  
\end{equation}
allowing us to conclude \eqref{eq:TVboundPatched}. 
It remains to extend the argument from $t_1 <t_2<\dots<t_K$ in $\mathcal{T} \cap [0,n^{-1}]$ to arbitrary $t_1 <t_2<\dots<t_K$ in $\mathcal{T}$. Observe that by the variational construction in \eqref{eq:VarCharHN}, the mapping
\begin{equation*}
  \hfix^{\per}_0 \mapsto   \hfixpa_{\ast}(\hfix^{\per}_0; y,t)
\end{equation*}
is $1$-Lipschitz in the uniform norm, i.e., for all $\hfix^{\per}_0,\tilde{\hfix}^{\per}_0 \in \UCp$ and $t\in [0,n^{-1}]$, 
\begin{equation*}
  \norm{ \hfixpa_{\ast}(\hfix^{\per}_0; \cdot,t) - \hfixpa_{\ast}(\tilde{\hfix}^{\per}_0; \cdot,t)}_{\infty} \leq \norm{ \hfix^{\per}_0 - \tilde{\hfix}^{\per}_0 }_{\infty} . 
\end{equation*} 
The same argument applies to the rescaled patched ASEP height function using attractivity of the underlying full-space ASEPs. In particular, note that assumption \eqref{eq:DominationAssumption} continues to hold at time $n^{-1}$.
Hence, using  \eqref{eq:FirstIntervalBound} on $[0,n^{-1}]$ and on $[n^{-1},2n^{-1}]$ with respect to the initial data from the patched ASEP and patched KPZ fixed point at times $n^{-1}$, we see that with probability $1-2 C_2 \exp(-c_2n^{1/2})$, we find a coupling so that for all $t_1 <t_2<\dots<t_K$ in $\mathcal{T} \cap [0,2n^{-1}]$, 
\begin{equation}\label{eq:SecondIntervalBound}
 \max_{i \in \lbr M \rbr}  \max_{\ell \in \lbr K \rbr}  \sup_{y\in \mathbb{T}} \left|\hfix^{\textup{patch},n^{-1},\varepsilon}(\mathfrak{h}^{\per,(i),\varepsilon}_0;\cdot,t_\ell)- \hfixpa(\mathfrak{h}^{\per,(i)}_0;\cdot,t_\ell) \right|   \leq 3 \cdot 2^{-n}
\end{equation} for all $n\in N$, and all $\varepsilon>0$ sufficiently small, where the constants $c_2,c_2>0$ are taken from \eqref{eq:FirstIntervalBound}, and the factor $3$ at the left-hand side of \eqref{eq:SecondIntervalBound} comes from replacing time $0$ by time $n^{-1}$ in the patched ASEP, replacing time $0$ by time $n^{-1}$ in the patched KPZ fixed point, and from the error in \eqref{eq:FirstIntervalBound}. Iterating this argument $n$ many times, we see that 
with probability $1- n C_2 \exp(-c_2n^{1/2})$, we find a coupling so that for all $t_1 <t_2<\dots<t_K$ in $\mathcal{T}$, 
\begin{equation}\label{eq:ThirdIntervalBound}
 \max_{i \in \lbr M \rbr}  \max_{\ell \in \lbr K \rbr}  \sup_{y\in \mathbb{T}} \left|\hfix^{\textup{patch},n^{-1},\varepsilon}(\mathfrak{h}^{\per,(i),\varepsilon}_0;\cdot,t_\ell)- \hfixpa(\mathfrak{h}^{\per,(i)}_0;\cdot,t_\ell) \right|   \leq (2n+1)  2^{-n}
\end{equation} for all $n\in N$, and all $\varepsilon>0$ small enough, allowing us to conclude. 
\end{proof}

\subsection{From the patched KPZ fixed point to the periodic KPZ fixed point}\label{sec:ASEPfinalPart}

As a final ingredient for Theorem~\ref{thm:ASEPcirclePDL}, we require an approximation of the periodic KPZ fixed point using patched KPZ fixed points.

\begin{lemma}\label{lem:PatchedKPZToPeriodicKPZ}
Fix some $K\in \N$ and $0 < t_1 < t_2 < \dots < t_K$. Then for any $\varepsilon'>0$ and initial data $\hfixp_0 \in \UCp$, we have
\begin{equation}\label{eq:CloseToPeriodicFixedPoint}
     \bigg\lVert \Big(\hfixpa(\hfixp_0;\cdot,t_\ell)\Big)_{\ell \in \lbr K \rbr} -  \Big(\hfixp(\hfixp_0;\cdot,t_\ell)\Big)_{\ell \in \lbr K \rbr}  \bigg\rVert_{\textup{TV}} \leq \varepsilon^{\prime} 
\end{equation} for all $n=n(\delta,\varepsilon^{\prime})$ large enough.
\end{lemma}

\begin{proof}
We will only consider the case $K=1$ and $t=t_K \in (0,1)$ as the result for general $K$ and $t_K<\infty$ follows by iterating \eqref{eq:CloseToPeriodicFixedPoint} and using the triangle inequality for the total-variation distance. 
We start by recalling from Proposition~\ref{pro:CauchySequence}  the periodic directed  landscape as the unique limit point of the wrapped patched directed landscapes $(\Lpa)_{n\in \N}$ under the metric $\Wast$, and the sets $\mathcal{D}_{\delta}$ for $\delta>0$ from \eqref{def:SlopeConti}. Let $\Lp$ be a periodic directed landscape. Then by iterating Lemma~\ref{lem:PatchedDLsConsistently} as in \eqref{eq:CauchyProperty2}, 
we see that for every fixed $t \in (0,1)$, and $n\in \N$ large enough, 
there exists a coupling $\mathbf{P}_n$ of $\Lpa$ and $\Lp$ such that 
 \begin{equation}\label{eq:PatchedPeriodicCouple}
  \tilde{\mathcal{A}}_{n} := \left\{ \Lpa(x,s;y,s^{\prime})= \Lp(x,s;y,s^{\prime}) \text{ for all } (x,s;y,s^{\prime}) \in \SetT \cap \mathcal{D}_{t^2}\right\} 
\end{equation} 
satisfies with some constants $c_1,C_1>0$ for all $n\in \N$ large enough
\begin{equation}\label{eq:PatchedPeriodicCoupleBound}
    \mathbf{P}_n( \tilde{\mathcal{A}}_{n} ) \geq 1- C_1\exp\big(-c_1 n^{-1/16}\big) . 
\end{equation}
As a consequence, we see that on the event $\tilde{\mathcal{A}}_{n}$,
\begin{equation}\label{eq:Finalperiodic}
    \Lpa(x,0;y,t) = \Lp(x,0;y,t) \text{ for all } x,y \in \mathbb{T} . 
\end{equation} In particular, whenever the event in \eqref{eq:Finalperiodic} occurs, we get $\hfixpa(\mathfrak{h}^{\per}_0;\cdot,t)=\hfixp(\mathfrak{h}^{\per}_0;\cdot,t)$. Combining \eqref{eq:PatchedPeriodicCoupleBound}
with the coupling representation of the total variation distance, we conclude.
\end{proof}

 \begin{proof}[Proof of Theorem~\ref{thm:ASEPcirclePDL}] 
 In following, for some fixed $K,M\in \N$, we consider the space $F_{M,K}$ of continuous bounded function $f$ which map from  
\begin{equation}\label{def:TestSpace}
 G_{M,K} := \bigotimes_{i \in \lbr M \rbr,j \in \lbr K \rbr} \left\{ f_{i,j} \colon \mathbb{T} \rightarrow \R \cup \{-\infty\} \text{ upper semi-continuous}\right\} 
\end{equation} to $\R$, where we endow $G_{M,K}$ with the topology coming from the supremum metric.
It suffices to show that for any initial data $\hfix^{\per,(i)}_0$, 
for any $f \in F_{M,K}$
with some $K\in \N$, and for any $0 < t_1 < t_2 < \dots < t_K$, we have
\begin{equation}\label{eq:WeakConvergenceASEP}
\begin{split}
     \lim_{\varepsilon \rightarrow 0}  \mathbb{E}\Big[ f(\hApeps(\hfix^{\per,(i),\varepsilon}_0&;\cdot,t_\ell), i\in \lbr M \rbr , \ell\in \lbr K \rbr )\Big] \\ &= \E^{\per}\left[ f( \hfixp(\hfix^{\per,(i)}_0;\cdot,t_\ell), i\in  \lbr M \rbr, \ell\in \lbr K \rbr) \right] .
\end{split}
\end{equation}
 Here, recall that we let $\mathbb{E},\mathbb{E}^{n}_{\textup{pat}},\mathbb{E}^{(n)}$ and $\mathbb{E}^{\per}$ denote the expectation under the laws of $\hApeps$, $\hfix^{\textup{patch},n^{-1},\varepsilon}$, $\hfixpa$ and $\hfixp$, respectively, where we use for $\hfix^{\textup{patch},\varepsilon,n^{-1}}$ the sigma-algebra $\sigma(\hfix^{\textup{patch},n^{-1}})$ defined below Definition~\ref{def:PatchedASEPNew}. We write in the following $f(\hfix)$ instead of $f(\hfix(\hfix^{\per,(i)}_0; \cdot , t_\ell), i\in  \lbr M \rbr, \ell\in \lbr K \rbr )$ to simplify notation, and
fix some $\varepsilon'>0$. 
Then Corollary~\ref{cor:CouplePatchedPeriodicASEPs} ensures that for every $\varepsilon'>0$, we find some $n_0 \in \N$ so that for all $n\geq n_0$, and all $\varepsilon>0$ small enough, 
\begin{equation}\label{eq:FinalASEP1}
      \left| \mathbb{E}\left[ f(\hApeps)\right] -  \mathbb{E}^n_{\textup{pat}}\left[ f( \hfix^{\textup{patch},n^{-1},\varepsilon} )  \right]  \right| \leq \frac{1}{3} \varepsilon' .
\end{equation} 
 Moreover, by Lemma~\ref{lem:patchedASEPfixed2}, using that $f$ is bounded, we find some $n_1=n_1(\varepsilon) \in \N$ so that
\begin{equation}\label{eq:FinalASEP2}
  \left|  \mathbb{E}^n_{\textup{pat}}\left[ f( \hfix^{\textup{patch},n^{-1},\varepsilon} )  \right] -  \mathbb{E}^{(n)}\left[ f( \hfixpa ) \right] \right| \leq \frac{1}{3} \varepsilon' .
\end{equation} for all $n\geq n_1$, and all $\varepsilon>0$ small enough. 
Then by Lemma~\ref{lem:PatchedKPZToPeriodicKPZ}, we find some $n_2(\varepsilon') \in \N$ so that for all $n\geq n_2$, and $\varepsilon>0$ small enough,  
\begin{equation}\label{eq:FinalASEP3}
   \left| \mathbb{E}^{(n)}\left[ f( \hfixpa )  \right] - \E^{\per}\left[ f( \hfixp ) \right] \right| \leq \frac{1}{3} \varepsilon' .
\end{equation}
Combining the bounds \eqref{eq:FinalASEP1}, \eqref{eq:FinalASEP2} and \eqref{eq:FinalASEP3}, we see that for all $n\geq \max(n_0,n_1,n_2)$ and $\varepsilon>0$ small enough, 
\begin{equation*}
     \left| \mathbb{E}\left[ f(\hApeps)\right] - \E^{\per}\left[ f( \hfixp ) \right] \right| \leq \varepsilon' . 
\end{equation*} Since $\varepsilon'>0$ was arbitrary, this allows us to conclude \eqref{eq:WeakConvergenceASEP}.
\end{proof}

\appendix

\section{Moderate deviations for exponential last passage percolation}\label{sec:AppendixLPP}

We will now give the proofs of the moderate deviation results on geodesics in periodic last passage percolation, which were deferred in Section \ref{sec:ModerateLPP}. 

\begin{proof}[Proof of Lemma \ref{lem:GeodesicsFullPeriodicAgree}]
We start by locating the endpoint in the periodic geodesic $\gper_{(0,0),v}$.
Proposition~3.12 in \cite{SS:TASEPcircle} states that for $m_{\ast}=k^2(N-k)^{-2}$, and some constants $c_1,\theta_1>0$, we get that 
\begin{equation}\label{eq:FixedTargetPeriodic}
 \Pper\big( \gper_{(0,0),(n,m_{\ast}n)}= \gamma_{(0,0),(n,m_{\ast}n)} \big) \geq 1-\exp({-c_1 \theta^{2/3}}) 
\end{equation} for all $\theta \geq \theta_1$, uniformly in $n \leq N$ from \eqref{eq:CorrectPeriodicScale}, provided that $N$ is large enough. We denote by $\slz(u;v)$ the slope between $u=(u_1,u_2)$ and $v=(v_1,v_2)$, i.e.,
\begin{equation}\label{def:Slope}
\slz(u;v):=(v_{2}-u_{2})/(v_1-u_1) . 
\end{equation}
Then for given $u,v$, we define two starting sites $u_-,u_+$ as
\begin{align*}
u_- &:= \textup{argmax}_{\{u \in \TR{(0,0)} \colon \slz(u;v)\leq m_{\ast} \}}\big\{ \slz(u,v)\big\} \\
u_+ &:= \textup{argmin}_{\{u \in \TR{(0,0)} \colon \slz(u;v)\geq m_{\ast} \}}\big\{ \slz(u,v) \big\} ,
\end{align*} 
i.e., $u_-$ and $u_+$ correspond to the sites in $\TR{(0,0)}$ with the slope closest to $m_{\ast}$. From Lemma~4.5 in~\cite{SS:TASEPcircle}, and assumption \eqref{eq:CorrectPeriodicScale}, we get that for any two sites $w,w^{\prime}$ on the geodesic $\gamma_{(0,0),(n,m_{\ast}n)}$, the geodesic between $w$ and $w^{\prime}$ will with high probability not cross over to a periodic translate, i.e., there exist constants $c_2,\theta_2>0$ such that for all $\theta \geq \theta_2$
\begin{equation}\label{eq:NonCrossingCondition}
 \Pper\left( T_{w,w^{\prime}} > T_{\TR(w)\setminus \{ w \},w^{\prime}} \ \forall w,w^{\prime} \in \gamma_{(0,0),(n,m_{\ast}n)} \ , w^{\prime} \succeq w\right) \geq 1- \exp(-c_2 \theta^{2/3}) .
\end{equation}
Then using the ordering of geodesics, together with \eqref{eq:FixedTargetPeriodic}, we get that for some constant $c_3>0$ and all $n=n(N)$ from \eqref{eq:CorrectPeriodicScale} with $\theta$ large enough that
\begin{equation}\label{eq:PeriodicToPoint}
\Pper(\gper_{(0,0),v} = \gamma_{u_{+},v}  \, \vee \,  \gper_{(0,0),v} = \gamma_{u_{-},v} ) \geq 1 - \exp(-c_3 \theta^{2/3}) .  
\end{equation}
In order to bound the transversal fluctuations of the geodesics $\gamma_{u_{+},v}$ and $\gamma_{u_-,v}$ under the measure $\Pper$, assume without loss of generality that $m \geq m_{\ast}$. Then consider $\gamma_1=\gamma_{u_{-},u_{-}+(n,m_{\ast}n)}$ as well as $\gamma_2=\gamma_{v-(n,m_{\ast}n),v}$, and note that the geodesic $\gamma_{u_-,v}$ is sandwiched between the two paths $\gamma_1$ and $\gamma_2$. A similar statement holds with respect to $\gamma_{u_+,v}$.
Now Corollary 3.15 in \cite{SS:TASEPcircle} states moderate deviations on the transversal fluctuations of $\gamma_1,\gamma_2$, i.e., there exists some $x_0>0$ such that for all $x \in [x_0,\theta^{2/3}]$ 
\begin{equation}\label{eq:PeriodicTF}
\Pper\left( \max(\TF(\gamma_1),\TF(\gamma_2)) \leq x m_{\ast}^{2/3}n^{2/3}\right) \geq 1-\exp(-c_4 x)
\end{equation}
with some constant $c_4>0$. We define the event 
\begin{align}\label{eq:IndividualA}
A=\left\{\max(\TF(\gamma_1),\TF(\gamma_2)) \leq \frac{k}{4}  \right\} , 
\end{align} and note that from \eqref{eq:FixedTargetPeriodic} and \eqref{eq:PeriodicTF} for some constant $c_5>0$, and all $N$ large enough 
\begin{equation*}
\begin{split}
\Pfull(A)&\geq 1-\exp\big(-c_5 \theta^{2/3}\big) ,  \\
\Pper(A)&\geq 1- \exp\big(-c_5 \theta^{2/3}\big) . 
\end{split}
\end{equation*}
For $i\in \N$, let $R_i$ denote the sites of $\ell_1$-distance at most $ik/4$ away from the line connecting $u_-$ to $v$. Note that $u_{-}+(n,m_{\ast}n),v-(n,m_{\ast}n) \in R_1$ for all $\theta$ large enough, while on the event $A$, we get that
$ \gamma_{u_-,v} \subseteq R_2$.
Using the same environment on the sites in $R_2$, there exist couplings $\mathbf{P}_1$ and $\mathbf{P}_2$ between $\Pfull$ and $\Pper$, and constants $c_6,C_6>0$ such that for all $\theta \geq \theta_0$, and all $N$ large enough
\begin{equation}\label{eq:AgreementPerIID}
\begin{split}
\mathbf{P}_1( \gamma_{u_-,v}= \tilde{\gamma}_{u_-,v} ) \geq 1- C_6\exp(-c_6 \theta^{2/3}) , \\
\mathbf{P}_2( \gamma_{u_+,v}= \tilde{\gamma}_{u_+,v} ) \geq 1- C_6\exp(-c_6 \theta^{2/3}) , 
\end{split}
\end{equation} where we write $\tilde{\gamma}$ for the geodesic in the periodic environment. Here, we use Lemma~\ref{lem:TransversalFluctuations} for the full-space model and \eqref{eq:PeriodicTF} for the periodic model. Controlling the transversal fluctuations of $\gamma_{u_-,v}$ and $\gamma_{u_+,v}$ by Lemma \ref{lem:TransversalFluctuations}, we see that for some constants $c_7,C_7>0$, and all $N$ large enough
\begin{equation}
\Pper\big( \max(\TF(\gamma_{u_+,v}),\TF(\gamma_{u_-,v})) \leq x \theta^{-2/3}N \big) \geq 1-C_7\exp({-c_7 x}) ,  
\end{equation} provided that $x \in [x_0, \theta^{2/3}]$. Together with \eqref{eq:PeriodicToPoint}, this yields the desired result. 
\end{proof}

Next, we establish a bound on the local transversal fluctuations of geodesics in periodic and full-space environments. 

\begin{proof}[Proof of Lemma \ref{lem:ModerateLocal}]
Since the case \eqref{eq:FullLocalModerate} is covered as Theorem 3 in \cite{BSS:Coalescence}, we will only show \eqref{eq:PeriodicLocalModerate}. Recall  \eqref{eq:AgreementPerIID} from the proof of Lemma \ref{lem:GeodesicsFullPeriodicAgree} for $u=(0,0)$ and $v=(\np,m\np)$. Together with \eqref{eq:FullLocalModerate},  for all $ x  \leq \theta^{1/3}$, and all $L \in \lbr \np \rbr$ with $L \geq L_0$, we find constants $c_0,C_0>0$ so that
\begin{equation}\label{eq:PeriodicLocalModerateProof}
\Pper\Big( \max\big(\TF_L(\gamma_{u_-,v}), \TF_L(\gamma_{u_+,v})\big) \geq x L^{2/3} \Big) \leq C_0\exp({-c_0 x^2}) 
\end{equation} for all $N$ and $\np=\np(N)$ large enough.
 Using now \eqref{eq:PeriodicToPoint} to bound the probability that $\gper_{u,v}$ agrees with either $\gper_{u_-,v}$ or $\gper_{u_+,v}$, we conclude.
\end{proof}

We now argue that the local transversal fluctuations of geodesics can be controlled uniformly in the choice of the starting and ending point, thereby establishing Proposition \ref{pro:TypicallyGoodLPP}.  We start by a rough estimate on the location of the endpoint of the lattice path $\gper_{u,v}$ in Lemma \ref{lem:EndpointConcentration}. We then show in Lemma \ref{lem:LocalFluctuationsLPP} that for lattice paths with endpoints in a given set $W^n$, we have that all paths with a sufficiently small slope are with high probability $(\delta,N)$-good. In a last step, we dominate all paths by the paths with their respective endpoints in $W^n$. 
For all $i\in \Z$, we define the discrete periodic lines in an $(N,k)$ periodic environment as
\begin{equation}
\mathbb{L}_\ell^n :=  \{ \exists s \in \R \, \colon \, (m_{\ast}^{-1} \ell n^{-1}N^{3/2},\ell n^{-1}N^{3/2}) + (\lfloor s (N-k) \rfloor , - \lfloor s k \rfloor ) \} .
\end{equation} and set for all $n\in \N$ and $\ell \in \lbr n \rbr$
\begin{equation}
    Y^{n}_{\ell} := \left\{ w \in \Z^2 \, \colon \exists   u\in \mathbb{L}_{\ell-1}^n ,  v\in \mathbb{L}_{\ell}^n , v \succeq w \succeq u \right\} . 
\end{equation}
Moreover, define the sets
\begin{equation*}
    X^{n}_{\ell} := \left\{   u=\big(u_1, \lfloor \ell n^{-1}N^{3/2} \rfloor \big) \in \Z^2  \, \colon \, u_1 - m_{\ast}^{-1}\ell n^{-1}N^{3/2} \in [-4nN,4nN]\right\} 
\end{equation*} for all $n\in \N$ and $\ell \in \lbr n\rbr \cup \{ 0 \}$, with the convention that $X^n := \bigcup_{\ell \in \lbr n+1\rbr} X^{n}_{\ell-1}$.
Define the collection of sites
\begin{equation*}
W^n := \left\{ \big( \big\lfloor i n^{-1} N \big\rfloor , \big\lfloor j n^{-1} N^{3/2} \big\rfloor   \big) \ \colon \,  i,j\in \Z \right\} , 
\end{equation*} and recall from \eqref{eq:VerticalLPPdistance} the horizontal coordinates $\gamma^{\per,\pm}_{u,v}(\vartheta)$ for all lattice paths $\gper_{u,v}$ from $u=(u_1,u_2)$ to $v=(v_1,v_2)$ with $u_2 \leq \vartheta \leq v_2$. 
\begin{lemma}\label{lem:EndpointConcentration}
For all $n\in \N$ and $\ell \in \lbr n \rbr$, we define the sets
\begin{equation}
    U_{\ell}^n := \left\{ (u;v)\in \Z^4_{\uparrow}  \,  \colon \,  \gper_{u,v} \cap X^{n}_{\ell-1} \neq \emptyset \, \wedge \, \gper_{u,v} \cap X^{n}_{\ell} \neq \emptyset \right\}
\end{equation}
and the events
\begin{equation*}
    A_n :=  \bigcap_{\ell \in \lbr n-1\rbr} \bigcap_{(u;v) \in U^{n}_{\ell}}  \Big\{ \Big| \gamma^{\per,\pm}_{u,v}( \lfloor \ell n^{-1}N^{3/2} \rfloor ) - \gamma^{\per,\pm}_{u,v}( \lfloor (\ell-1) n^{-1}N^{3/2} \rfloor) -  m_{\ast}^{-1}n^{-1}N^{3/2}\Big| \leq 4 N \Big\} . 
\end{equation*} 
 Then there exist constants $c,C>0$ such that for all $n \in \N$, we get that
\begin{equation}\label{eq:TFBound1}
    \P_{\per}(A_n) \geq 1 - C\exp(-c n^{2/3} )
\end{equation} for all $N$ large enough. Similarly, for all $n\in \N$, we define the events
\begin{equation}
    \tilde{A}_n :=  \bigcap_{\ell \in \lbr n\rbr} \bigcap_{(u;v) \in Y_{\ell}^{n} \times Y_{\ell}^{n}}  \Big\{ \TF(\gper_{u,v}) \leq 4 N \Big\} . 
\end{equation} 
 Then there exist constants $\tilde{c},\tilde{C}>0$ such that for all $n \in \N$, we get that
\begin{equation}\label{eq:TFBound2}
    \P_{\per}\big(\tilde{A}_n\big) \geq 1 - \tilde{C}\exp(-\tilde{c} n^{2/3} )
\end{equation} for all $N$ large enough.
\end{lemma}
\begin{proof} 
Throughout this proof, we omit the floors for notational convenience.
Let us start with  \eqref{eq:TFBound1}.
Note that on the event $A_n^{\complement}$, we get that
\begin{equation}\label{eq:LargeFluctuations}
   \gper_{u^{\prime},v^{\prime}}=\gamma_{u^{\prime},v^{\prime}} \text{ and } \big| v^{\prime}_1-u^{\prime}_1 - m_{\ast}^{-1}n^{-1}N^{3/2} \big| > 4N 
\end{equation}   
   for some $(u^{\prime};v^{\prime}) \in X^{n}_{\ell-1} \times X^{n}_{\ell}$ and $ \ell \in \lbr n \rbr$. Without loss of generality, we assume that $\ell=1$,  that $u^{\prime}=(u^{\prime}_1,u^{\prime}_2)$ satisfies $u^{\prime}_1 \in [0,N]$, and
\begin{equation}\label{eq:LargeTransversal}
    v^{\prime}_1 - u^{\prime}_1 - m_{\ast}^{-1}n^{-1}N^{3/2} > 4N
\end{equation}
as the remaining cases are treated similarly. Define the sites
\begin{equation*}
    u^{1}=(0,0) \quad \text{ and } \quad v^{1}=\big(m_{\ast}^{-1} n^{-1}N^{3/2}  + 4m_{\ast}^{-1}k,n^{-1}N^{3/2}+4k\big)  ,  
\end{equation*}
as well as their periodic translates
\begin{equation}
\begin{split}\label{eq:ShiftPoints}
    u^{2}=u^{1} + (2(N-k),-2k) \quad &\text{ and } \quad  u^{3}=u^{1} + (4(N-k),-4k) ,  \\
        v^{2}=v^{1} + (2(N-k),-2k) \quad &\text{ and } \quad  v^{3}=v^{1} + (4(N-k),-4k)  . 
    \end{split}
\end{equation} Note that $v^{3} \in X^{n}_{2}$. 
Since $v^{1}$ takes the form $(\tilde{n},m_{\ast}\tilde{n})$ for some $\tilde{n}\in \N$, 
Corollary 3.15 from \cite{SS:TASEPcircle} ensures that there exist some constants $c_1,C_1>0$ such that for all $n\in \N$ 
\begin{equation}\label{eq:PeriodicsNice}
    \P_{\per}\left( \gper_{u^{1},v^{1}}= \gamma_{u^{1},v^{1}} \, \wedge \,  \TF(\gamma_{u^{1},v^{1}}) \leq \frac{1}{2}N   \right) \geq 1 - C_1\exp(-c_1 n^{2/3} ) , 
\end{equation} provided $N$ is large enough. Since $u^{2},u^{3} \in \TR(u^{1})$ and $v^{2},v^{3} \in \TR(v^{1})$, the same holds for the periodic geodesics between $u^{2}$ and $v^{2}$, as well as $u^{3}$ and $v^{3}$, respectively. We claim that under the event in \eqref{eq:PeriodicsNice}, no pair $(u^{\prime},v^{\prime})$ with \eqref{eq:LargeTransversal} exists. 
Indeed, for a proof by contradiction, under \eqref{eq:LargeTransversal} and the event \eqref{eq:PeriodicsNice}, for all $N$ large enough, we can find some $w^- \in \gamma_{u^{2},v^{2}}$ and $w^+ \in \gamma_{u^{3},v^{3}}$ such that $w^-,w^+ \in \gamma_{u^{\prime},v^{\prime}}$. In this case, we note that 
\begin{equation}
    \tilde{w}^{+}=w^{+}+(-2(N-k),2k) \in \gamma_{u^{2},v^{2}} .  
\end{equation}
 Write $\Tper(\gamma)$ for the weight associated to a lattice path $\gamma$ in periodic last passage percolation. 
 Since $\gamma_{u^{\prime},v^{\prime}}$ is a periodic geodesic, and periodic geodesics are almost surely unique, this implies 
\begin{equation*}
     \Tper(\gamma_{w^{-},w^{+}}) >  \Tper(\gamma_{w^{-},\tilde{w}^{+}}) , 
\end{equation*}
 and hence
\begin{equation}
    \Tper(\gamma_{u^{2},v^{2}}) <  \Tper(\gamma_{u^{2},w^{-}}) +  \Tper(\gamma_{w^{-},w^{+}}) + \Tper(\gamma_{w^{+},v^{2}}) , 
\end{equation}
 contradicting the assumption that $\gamma_{u^{2},v^{2}}$ is a  geodesic, which satisfies \eqref{eq:PeriodicsNice}. In particular, 
\begin{equation*}
    \P_{\per}\left(  v^{\prime}_1 - u^{\prime}_1 - \frac{1}{m_{\ast}n}N^{3/2} \leq  4N \, \forall \, (u^{\prime};v^{\prime}) \in X^{n}_{0} \times X^{n}_{1} \colon    \gamma^{\per,\pm}_{u,v}(0) \in [0,N] \right) \geq 1 - C_1\exp(-c_1 n^{2/3}) , 
\end{equation*} provided $N$ is large enough, and we conclude \eqref{eq:TFBound1} by a union bound. For \eqref{eq:TFBound2}, we proceed similarly. Let $\tilde{v}^{1}\in \mathbb{L}_1^{n}$ be such that $\slz(u^{1},\tilde{v}^{1})=m_{\ast}$, and set
\begin{equation*}
        \tilde{v}^{2}=\tilde{v}^{1} + (2(N-k),-2k) \quad \text{ and } \quad  \tilde{v}^{3}=\tilde{v}^{1} + (4(N-k),-4k)  . 
\end{equation*}
Note that $\tilde{v}^{1},\tilde{v}^{2},\tilde{v}^{3} \in \mathbb{L}_1^{n}$ by construction, while $u^{1},u^{2},u^{3} \in \mathbb{L}_0^{n}$. Let $\textsf{R}(u^{1},\tilde{v}^{1};u^{2},\tilde{v}^{2})$ denote the area between the geodesics $\gper_{u^{1},\tilde{v}^1}$ and $\gper_{u^{2},\tilde{v}^2}$ and the lines $\mathbb{L}_0^{n}$ and $\mathbb{L}_1^{n}$.  Now as \eqref{eq:PeriodicsNice} continues to hold for the paths $\gper_{u^{1},\tilde{v}^1}$, and its periodic translates $\gper_{u^{2},\tilde{v}^2}$ and $\gper_{u^{3},\tilde{v}^3}$, we get by the same arguments as in the first case that for some constants $\tilde{c},\tilde{C}>0$ and all $n$ large enough
\begin{equation*}
    \P_{\per}\left(  \TF\big(\gper_{u^{\prime},v^{\prime}}\big) \leq  4N \ \forall \, (u^{\prime};v^{\prime}) \in Y_{1}^{n} \times Y_{1}^{n} \colon  u^{\prime} \in \textsf{R}(u^{1},\tilde{v}^{1};u^{2},\tilde{v}^{2})\right) \geq 1 - \tilde{C}\exp\big(-\tilde{c} n^{2/3}\big) , 
\end{equation*} provided $N$ is large enough. Since the environment is $(N,k)$-periodic, this implies \eqref{eq:TFBound2}.
\end{proof}

\begin{remark}\label{rem:ShiftBound}
    Note that the bounds from \eqref{eq:TFBound1} in Lemma \ref{lem:EndpointConcentration} remain valid when we translate the sets $X^{n}_{\ell}$ by some fixed $w \in \Z^2$, i.e., we set
    \begin{equation}
        \bar{X}^{n}_{\ell} := \left\{ v \in \Z^2 \colon v-w \in X^{n}_{\ell}\right\} . 
    \end{equation}
\end{remark}
Next, for all $(u;v) \in X^n \times X^n$ with $n\in \N$, all $\ell \in \N$ fixed, and all $i\in \lbr n \rbr$, we define 
\begin{equation*}
 C_{u,v}^{i,n,\ell} := \left\{  \Big| \gamma^{\per,\pm}_{u,v}\big( in^{-1} N^{3/2}\big) - \gamma^{\per,\pm}_{u,v}\big( (i-1)n^{-1} N^{3/2}\big)  - m_{\ast}^{-1} n^{-1}  N^{3/2}   \Big| \leq 5(\ell n)^{1/2}N  \right\}  , 
\end{equation*} provided that $(i-1)n^{-1}N^{3/2}, in^{-1}N^{3/2} \in [u_2,v_2]$, and we let $ C_{u,v}^{i,n,\ell}$ be trivially satisfied, otherwise. For all $n\in \N$, we set 
\begin{equation*}
    \tilde{X}^{n} = \bigcup_{\ell \in \lbr n \rbr \cup \{ 0\}}\left\{   u=\big(u_1, \lfloor \ell n^{-1}N^{3/2} \rfloor \big) \in \Z^2  \, \colon \, u_1 - m_{\ast}^{-1}\ell n^{-1}N^{3/2} \in [-2nN,2nN]\right\} 
\end{equation*}
and let
\begin{equation}
E_n = (W^n \times W^n) \cap (\tilde{X}^{n} \times \tilde{X}^{n}) \cap D_{n^{-1}}^N,
\end{equation}
where we recall $D_\delta^N$ for some $\delta>0$ and $N \in \N$ from \eqref{def:SlopeDn}.
We have the following result on the local fluctuations of paths with endpoints in $E_{\ell n}$ for some $n,\ell\in \N$.

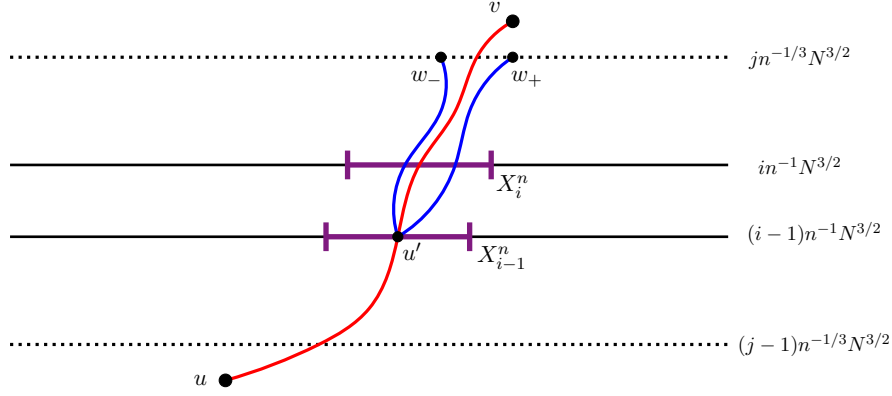
\begin{figure}
\centering
\begin{tikzpicture}[scale=.95]

   \draw[line width=1.2pt, dotted] (0,0) -- (10,0);
 \draw[line width=1pt] (0,1.5) -- (10,1.5);
 \draw[line width=1pt] (0,2.5) -- (10,2.5);   
   \draw[line width=1.2pt, dotted] (0,4) -- (10,4);


    \draw[darkblue, line width =2 pt] (4.4,1.5) -- (6.4,1.5);
    \draw[darkblue, line width =2 pt] (4.7,2.5) -- (6.7,2.5);
    
    \draw[darkblue, line width =2 pt] (4.4,1.3) -- (4.4,1.7);
    \draw[darkblue, line width =2 pt] (6.4,1.3) -- (6.4,1.7);

    \draw[darkblue, line width =2 pt] (4.7,2.3) -- (4.7,2.7);
    \draw[darkblue, line width =2 pt] (6.7,2.3) -- (6.7,2.7);
    
    
   
\draw[line width =1.2 pt,red] (3,-0.5) to[curve through={(4.3,0)..(5,0.5)..(5.4,1.5) .. (5.7,2.5) .. (6.2,3.25)..(6.5,4) }] (7,4.5);

\draw[line width =1.2 pt,blue] (5.4,1.5) to[curve through={(5.5,2.5) .. (6,3.25)}] (6,4);

\draw[line width =1.2 pt,blue] (5.4,1.5) to[curve through={(6.2,2.5) .. (6.4,3.25)}] (7,4);

\filldraw (3,-0.5) circle (2.5pt);

\filldraw (7,4.5) circle (2.5pt);

\filldraw (5.4,1.5) circle (2pt);
\filldraw (6,4) circle (2pt);
\filldraw (7,4) circle (2pt);

\node[scale=0.7] (site) at (11,4){$jn^{-1/3}N^{3/2}$};  
\node[scale=0.7] (site) at (11.2,0){$(j-1)n^{-1/3}N^{3/2}$};

\node[scale=0.7] (site) at (11.2,1.55){$(i-1)n^{-1}N^{3/2}$};  
\node[scale=0.7] (site) at (11,2.5){$in^{-1}N^{3/2}$};  

\node[scale=0.8] (site) at (2.65,-0.5){$u$};
\node[scale=0.8] (site) at (6.75,4.7){$v$};

\node[scale=0.8] (site) at (5.8,3.7){$w_-$};
\node[scale=0.8] (site) at (7.2,3.7){$w_+$};

\node[scale=0.8] (site) at (5.6,1.3){$u^{\prime}$};

\node[scale=0.8] (site) at (7,2.2){$X_{i}^{n}$};
\node[scale=0.8] (site) at (6.85,1.2){$X_{i-1}^{n}$};

	\end{tikzpicture}	
\caption{\label{fig:SplitAcrossAxis} Visualization of the path decomposition in the proof of Lemma \ref{lem:LocalFluctuationsLPP} in the case $v_2-u_2>4n^{-1/3}N^{3/2}$. The path $\gper_{u,v}$ is sandwiched between the paths $\gper_{u^{\prime},w_-}$ and $\gper_{u^{\prime},w_+}$ between the levels $(i-1)n^{-1}N^{3/2}$ and $in^{-1}N^{3/2}$.}
 \end{figure}

\begin{lemma}\label{lem:LocalFluctuationsLPP} Fix some $\ell \in \N$. 
Then there exist some constants $c,C>0$, depending only on $\ell$, such that for all $n\in \N$, we get that 
\begin{equation}
    \P_{\per}\left(C_{u,v}^{i,n,\ell} \, \forall (u;v)\in E_{\ell n} , i \in \lbr n\rbr \right) \geq 1-C\exp\big(-c n^{1/6} \big) . 
\end{equation}
\end{lemma}
\begin{proof}
We will in the following argue that for some constants $c_0,C_0>0$
\begin{equation*}
    \P_{\per}\left(C_{u,v}^{i,n,\ell} \right) \geq 1-C_0 \exp\big(-c_0 n^{1/6} \big) 
\end{equation*}
 for fixed $u=(u_1,u_2)$ and $v=(v_1,v_2)$ with $(u;v) \in E_{\ell n}$, and conclude by a union bound as $|E_{\ell n}| \leq C_0(\ell n)^4$ for all $n\in \N$ large enough with some absolute constant $C_0>0$. Depending on the vertical distance between $u$ and $v$, we distinguish two cases. First, suppose that $|v_2-u_2|\leq 4n^{-1/3}N^{3/2}$ and assert that $i\leq n/2$ as the arguments will by symmetry be similar for $i > n/2$. In this case, we get from Lemma \ref{lem:GeodesicsFullPeriodicAgree} with $x=n^{1/6}$ and $\theta=n^{1/3}/2$ that there exists some constants $c_1,C_1>0$ such that
\begin{equation}\label{eq:GlobalTFGood}
  \P_{\per}\left(  \TF(\gper_{u,v}) \leq  n^{1/6} (v_2-u_2)^{2/3} \right) \geq 1 - C_1\exp(-c_1 n^{1/6}) . 
\end{equation} Since $(u,v)\in (X^{n} \cap X^{n}) \cap D^N_{(\ell n)^{-1}}$ by our assumptions, the event in \eqref{eq:GlobalTFGood} ensures that
\begin{align*}
    \slz\big( (\gamma^{\per,\pm}_{u,v}((i-1)n^{-1}N^{3/2}),(i-1)n^{-1}N^{3/2}) ; (v_1,v_2) \big)^{-1} &\leq m_{\ast}^{-1} + N^{-1/2}(\ell n)^{1/2} + \frac{2(v_2-u_2)^{2/3}n^{1/6}}{v_2-u_2} \\
    &\leq m_{\ast}^{-1} +  N^{-1/2} (4 \ell n)^{1/2} , 
\end{align*}
and similarly for a lower bound. Hence, we get that for all $N$ large enough
\begin{equation*}
    \P_{\per}\left(   \big( \gamma^{\per,\pm}_{u,v}((i-1)n^{-1}N^{3/2}),(i-1)n^{-1}N^{3/2} ; v_1,v_2 \big) \in D^N_{(4\ell n)^{-1} } \right) \geq 1 - C_1\exp(-c_1 n^{1/6}) .
\end{equation*}  The desired bound on the local fluctuations from level $(i-1)n^{-1}N$ to level $in^{-1}N$ now follows from Lemma \ref{lem:ModerateLocal} with $L=n^{-1}N^{3/2}$ and $x=n^{1/6}$. 
Next, assume that $|v_2-u_2| > 4n^{-1/3}N^{3/2}$. Then by splitting the geodesic $\gper_{u,v}$ according to the intersection points with the lines $(\mathbb{L}_j^{n^{1/3}})_{j\in \Z}$, \eqref{eq:TFBound2} in Lemma~\ref{lem:EndpointConcentration} yields that for some constants $c_2,C_2>0$, and all $n\in  \N$,
\begin{equation}\label{eq:RoughTF}
    \Pper\left( \TF(\gper_{u,v}) \leq 4Nn^{1/3} \right) \geq 1- C_2 \exp(-c_2 n^{2/9}) , 
\end{equation} provided that $N$ is sufficiently large. In particular, on the event in \eqref{eq:RoughTF}, we get that
\begin{align*}
   u^{\prime}=(u^{\prime}_1,u^{\prime}_2) := \big( \gamma^{\per,\pm}_{u,v}((i-1)n^{-1}N^{3/2}),(i-1)n^{-1}N^{3/2}  \big) \in X_{i-1}^n .  
\end{align*} 
Moreover, we find some $w_-,w_+ \in W^n \cap X^{n}_j$ of the form
\begin{equation*}
    w_{-}=(z ,j n^{-1} N^{3/2}) \quad \text{and} \quad w_{+}=( z+n^{-1}N  ,j n^{-1} N^{3/2})
\end{equation*}
with some $j,z\in \N$  such that $j-(i-1)= n^{2/3}$ and $v_2 \geq jn^{-1}N^{3/2}$, as well as
\begin{equation*}
    z \leq  \gamma^{\per,\pm}_{u,v}\big( j n^{-1/3} N^{3/2}\big) \leq z+n^{-1}N . 
\end{equation*} 
A visualization of this construction in given in Figure~\ref{fig:SplitAcrossAxis}. 
Then we get from \eqref{eq:TFBound1} in Lemma~\ref{lem:EndpointConcentration} and Remark~\ref{rem:ShiftBound} that on $A_{n^{1/3}}$ for the sets $X^{n}_{\ell}$ shifted by $(m_{\ast}^{-1}(i-1)n^{-1}N^{3/2}, (i-1)n^{-1}N^{3/2})$, for all $n$ large enough,
\begin{equation*}
     (u^{\prime};w_-),(u^{\prime};w_+)\in D_{4n^{-4/3}}^N .  
\end{equation*} 
Thus, by Lemma \ref{lem:ModerateLocal} we obtain the desired bound on the transversal fluctuations from level $(i-1)n^{-1}N^{3/2}$ to level $in^{-1}N^{3/2}$, allowing us to conclude. 
\end{proof}
Next, we define for all $(u;v) \in W^n \times W^n$ with $u=(u_1,u_2)$ and $v=(v_1,v_2)$ the events 
\begin{equation*}
 B_{u,v}^{j,n,\ell} := \left\{  \sup_{\varphi \in [0,n^{-1}]} \Big| \gamma^{\per,\pm}_{u,v}( (jn^{-1}+\varphi) N^{3/2}) - \gamma^{\per,\pm}_{u,v}(jn^{-1} N^{3/2})  - m_{\ast}^{-1} \varphi N^{3/2}   \Big| \leq 11(\ell n)^{1/2} N \right\} , 
\end{equation*} provided that $\gper_{u,v}$ intersects levels  $jn^{-1}N^{3/2}$ and $(j+1)n^{-1}N^{3/2}$ for $j\in \N_0$, and let $B_{u,v}^{j,n}$ be trivially satisfied, otherwise. 
The next lemma allows us to control the local fluctuations of paths in $E_{\ell n}$. 

\begin{lemma}\label{lem:GridIsGood}
Fix some $\ell \in \N$ and consider the event 
\begin{equation}
    B^{\ell}_n := \bigcap_{(u;v) \in E_{\ell n}} \bigcap_{j \in \lbr n\rbr}  B_{u,v}^{j-1,n,\ell} . 
\end{equation} Then there exist constants $c_0,C_0>0$, depending only on $\ell$, such that for all $n\in \N$,
\begin{equation}\label{eq:FirstGood}
    \P(B^{\ell}_n) \geq 1- C_0 \exp(-c_0 n^{1/6}) . 
\end{equation} 
\end{lemma}
\begin{proof}  
Fix some $i\in \N$, and assume that $\gper_{u,v}$ for $(u;v) \in E_{\ell n}$ intersects levels  $(i-1)n^{-1}N^{3/2}$ and $in^{-1}N^{3/2}$. Then we can find some random sites $u_-,u_+,v_-,v_+\in W^{102\ell n}$ with
\begin{align*}
     u_{-}=(z,(i-1)n^{-1} N^{3/2}) \quad &\text{and} \quad u_{+}=(z+n^{-1}N,(i-1)n^{-1} N^{3/2}) , \\
      v_{-}=(z^{\prime},in^{-1} N^{3/2}) \quad &\text{and} \quad v_{+}=(z^{\prime}+n^{-1}N,in^{-1} N^{3/2}) , 
\end{align*}
as well as 
\begin{align*}
    \gamma^{\per,\pm}_{u,v}\big( (i-1)n^{-1} N^{3/2}\big) \in [z, z+n^{-1}N] \quad &\text{and} \quad
    \gamma^{\per,\pm}_{u,v}\big( in^{-1} N^{3/2}\big)  \in [z^{\prime}, z^{\prime}+n^{-1}N]. 
\end{align*}
Now by Lemma \ref{lem:EndpointConcentration}, splitting the periodic geodesic $\gper_{u,v}$ for $(u,v)\in \tilde{X}_{\ell n}\times \tilde{X}_{\ell n}$ according to the intersection with the lines $\mathbb{L}_{\ell}^n$ for some $\ell \in \lbr n\rbr$, the event $\tilde{A}_n$ in \eqref{eq:TFBound2} guarantees that
\begin{equation}
        (u_-,v_{-}), (u_+,v_{+}) \in \tilde{X}_{32\ell n} \times \tilde{X}_{32\ell n}
\end{equation} for all $n$ large enough. 
Moreover, on the event $C^{i,n,\ell}_{u,v}$ from Lemma \ref{lem:LocalFluctuationsLPP} with $(u,v)\in \tilde{X}_n\times \tilde{X}_n$, we get that for all $n\in \N$ large enough
\begin{equation*}
    (u_-,v_{-}), (u_+,v_{+}) \in E_{102\ell n} ,  
\end{equation*} where in order to guarantee that $(u_-,v_{-}), (u_+,v_{+}) \in D^{N}_{(102\ell n)^{-1}}$, we account for horizontal fluctuations of $10(\ell n)^{1/2}N$ from the set $C^{i,n,\ell}_{u,v}$, and the remaining term from the choice of $u_{-},u_{+},v_{-},v_{+}$.  
Thus, together with Lemma \ref{lem:GeodesicsFullPeriodicAgree} for $x=n^{1/6}$ on the transversal fluctuations between $u^{\pm}$ and $v^{\pm}$, we see that with probability at least $1- C_1 \exp(-c_1 n^{1/6})$ with some $c_1,C_1>0$, 
\begin{equation*}
 \sup_{\varphi \in [0,\delta]} \Big| \gamma^{\per,\pm}_{u,v}\big( ((i-1)n^{-1}+\varphi) N^{3/2}\big) - \gamma^{\per,\pm}_{u,v}\big((i-1)n^{-1} N^{3/2}\big)  - m_{\ast}^{-1} \varphi N^{3/2}   \Big| \leq 11(\ell n)^{1/2} N 
\end{equation*} holds. Since $j\in \lbr n\rbr$ was arbitrary, we conclude by a union bound.  
\end{proof}

We have now all tools in order to show Proposition \ref{pro:TypicallyGoodLPP}.

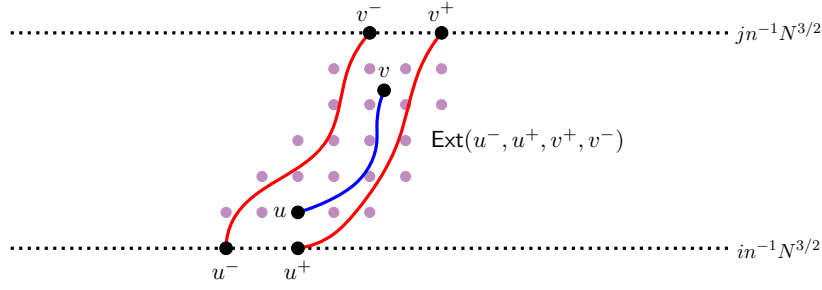
\begin{figure}
\centering
\begin{tikzpicture}[scale=.95]

   \draw[line width=1.2pt, dotted] (0,0) -- (10,0);
   \draw[line width=1.2pt, dotted] (0,3) -- (10,3);


    



\filldraw[darkblue!50] (3,0.5) circle (2pt);
\filldraw[darkblue!50] (3.5,0.5) circle (2pt);
\filldraw[darkblue!50] (4,0.5) circle (2pt);
\filldraw[darkblue!50] (4.5,0.5) circle (2pt);
\filldraw[darkblue!50] (5,0.5) circle (2pt);

\filldraw[darkblue!50] (3.5,1) circle (2pt);
\filldraw[darkblue!50] (4,1) circle (2pt);
\filldraw[darkblue!50] (4.5,1) circle (2pt);
\filldraw[darkblue!50] (5,1) circle (2pt);
\filldraw[darkblue!50] (5.5,1) circle (2pt);

\filldraw[darkblue!50] (4,1.5) circle (2pt);
\filldraw[darkblue!50] (4.5,1.5) circle (2pt);
\filldraw[darkblue!50] (5,1.5) circle (2pt);
\filldraw[darkblue!50] (5.5,1.5) circle (2pt);

\filldraw[darkblue!50] (4.5,2) circle (2pt);
\filldraw[darkblue!50] (5,2) circle (2pt);
\filldraw[darkblue!50] (5.5,2) circle (2pt);
\filldraw[darkblue!50] (6,2) circle (2pt);

\filldraw[darkblue!50] (4.5,2.5) circle (2pt);
\filldraw[darkblue!50] (5,2.5) circle (2pt);
\filldraw[darkblue!50] (5.5,2.5) circle (2pt);
\filldraw[darkblue!50] (6,2.5) circle (2pt);


    
   
\draw[line width =1.2 pt,red] (3,0) to[curve through={(3.1,0.4)..(4.4,1.5) .. (4.7,2.5)}] (5,3);

\draw[line width =1.2 pt,red] (4,0) to[curve through={(4.5,0.2)..(4.8,0.5)..(5.4,1.5) .. (5.7,2.5)}] (6,3);

\draw[line width =1.2 pt,blue] (4,0.5) to[curve through={(4.5,0.7)..(4.8,0.9)..(5.1,1.5) .. (5.1,1.7)}] (5.2,2.2);



\filldraw (3,0) circle (2.5pt);
\filldraw (5,3) circle (2.5pt);

\filldraw (4,0) circle (2.5pt);
\filldraw (6,3) circle (2.5pt);


\node[scale=0.7] (site) at (10.7,3){$jn^{-1}N^{3/2}$};  
\node[scale=0.7] (site) at (10.7,0){$in^{-1}N^{3/2}$};

\node[scale=0.8] (site) at (3.75,0.5){$u$};
\node[scale=0.8] (site) at (5.2,2.45){$v$};

\filldraw (4,0.5) circle (2.5pt);
\filldraw (5.2,2.2) circle (2.5pt);

\node[scale=0.8] (site) at (5,3.3){$v^-$};
\node[scale=0.8] (site) at (6,3.3){$v^+$};

\node[scale=0.8] (site) at (3,-0.3){$u^-$};
\node[scale=0.8] (site) at (4,-0.3){$u^+$};

\node[scale=0.8] (site) at (7.2,1.5){$\textsf{Ext}(u^{-},u^{+},v^{+},v^{-} )$};

	\end{tikzpicture}	
\caption{\label{fig:FinalLPP} Visualization of the path sandwiching in the proof of Proposition~\ref{pro:TypicallyGoodLPP}. The set $\textsf{Ext}(u^{-},u^{+},v^{+},v^{-} )$ is depicted as purple dots.}
 \end{figure}

\begin{proof}[Proof of Proposition \ref{pro:TypicallyGoodLPP}]
Note that for all $u,v \in D^N_{n^{-1}} \cap \mathbb{X}^{N}_{n^{-1}}$, Lemma \ref{lem:LocalFluctuationsLPP} and Lemma \ref{lem:GridIsGood} with $\ell=2$ guarantee that 
there exists  $u^+,u^-,v^+,v^- \in W^n$ such that $u^{-},u^{+}\in X^{n}_{i}$ and $v^{-},v^{+}\in X^{n}_{j}$ with some $i\leq j$, and
\begin{equation*}
\begin{split}
    \gamma^{\per,\pm}_{u^{-},v^{-}}(u_2) \leq u_1 \leq     \gamma^{\per,\pm}_{u^{+},v^{+}}(u_2)  \quad \text{and} \quad 
    \gamma^{\per,\pm}_{u^{-},v^{-}}(v_2) \leq v_1 \leq     \gamma^{\per,\pm}_{u^{+},v^{+}}(v_2) , 
\end{split} 
\end{equation*} while with respect to some absolute constants $c_0,C_0>0$, we get that 
\begin{equation*}
\begin{split}
   \Big| \gamma^{\per,\pm}_{u^{+},v^{+}}(u_2) -  \gamma^{\per,\pm}_{u^{-},v^{-}}(u_2)  \Big| \leq C_0n^{-1/2} N  \quad \text{and} \quad 
    \Big| \gamma^{\per,\pm}_{u^{+},v^{+}}(v_2) -   \gamma^{\per,\pm}_{u^{-},v^{-}}(v_2) \Big| \leq C_0n^{-1/2} N 
\end{split} 
\end{equation*} on an event of probability at least $1-\exp(-c_0n^{1/6})$  for all $n\in \N$ large enough. A visualization of this construction in given in Figure~\ref{fig:FinalLPP}. 
We define for all $(u^{-},u^{+},v^{+},v^{-})$ the set 
\begin{equation*}
    \textsf{Int}(u^{-},u^{+},v^{+},v^{-}) := \left\{ w=(w_1,w_2) \in W^n \, \colon \, \gamma^{\per,-}_{u^{-},v^{-}}(w_2) \leq w_1 \leq  \gamma^{\per,+}_{u^{+},v^{+}}(w_2) , w_2 \notin \{ u_2^{-},u_2^{+},v_2^{-},v_2^{+} \}\right\} 
\end{equation*} of points in $W^n$ sandwiched between the paths $\gper_{u^{-},v^{-}}$ and $\gper_{u^{+},v^{+}}$. Moreover, define 
\begin{equation*}
    \textsf{Ext}(u^{-},u^{+},v^{+},v^{-} ) := \left\{ (x,0)+w \colon \exists w\in \textsf{Int}(u^{-},u^{+},v^{+},v^{-} ),  |x| \leq  2n^{-1}N^{3/2}   \right\} \cap W^n  
\end{equation*}  as the 
extended hull of the sites $(u^{-},u^{+},v^{+},v^{-})$. 
Observe that by Lemma \ref{lem:GridIsGood} for the paths $\gper_{u^{-};v^{-}}$ and $\gper_{u^{+};v^{+}}$, the event $B^{\ell}_n$ with $\ell=2C_0$ implies that we find some universal constant $C_1>0$ such that for all $n\in \N$
\begin{equation*}
    (u^{\pm},w),(w,v^{\pm}) \in E_{C_1n} \text{ for all } w=(w_1,w_2)\in \textsf{Ext}(u^{-},u^{+},v^{+},v^{-} ) . 
\end{equation*}
Now observe that we find some $w=(w_1,w_2)$ and $w^{\prime}=(w_1+n^{-1}N,w_2)$ in $\textsf{Ext}(u^{-},u^{+},v^{+},v^{-} )$ such that 
\begin{equation*}
    w_1 \leq \gamma^{\per}_{u,v}(w_2)  \leq  w_1+n^{-1}N
\end{equation*} and $w_2=(i-1)n^{-1}N^{3/2}$ with $i\in \lbr n\rbr$. Hence, Lemma \ref{lem:GridIsGood} with $\ell=C_1$ implies the desired bound on the transversal fluctuations of $\gamma_{w,v^{-}}$ and $\gamma_{w^{\prime},v^{+}}$, and hence $\gamma_{u,v}$, between levels $(i-1)n^{-1}N^{3/2}$ and $in^{-1}N^{3/2}$, allowing us to conclude.
\end{proof}

In a similar way as in Lemma~\ref{lem:GridIsGood}, we deduce Lemma~\ref{lem:ConsistentLPPLocally} on the local transversal fluctuations of geodesics in periodic last passage percolation.

 \begin{proof}[Proof of Lemma~\ref{lem:ConsistentLPPLocally}] We will only show \eqref{eq:Ball1Local} as the same arguments apply for \eqref{eq:Ball2Local} by symmetry.
With a slight abuse of notation, set
\begin{equation*}
\begin{split}
       u^-&=(u_1^-,u_2^-)=u - (3\delta^{3/5}N + 2\delta m_{\ast}^{-1}N^{3/2},2\delta N^{3/2} ) , \\
       u^+&=(u_1^+,u_2^+)=u + (3\delta^{3/5}N - 2\delta m_{\ast}^{-1}N^{3/2},-2\delta N^{3/2} ) . 
\end{split}
\end{equation*} We refer to Figure~\ref{fig:FinalLPP2} for a visualization. 
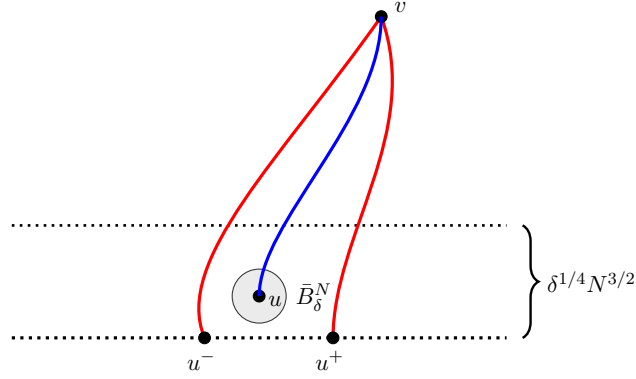
\begin{figure}
\centering
\begin{tikzpicture}[scale=0.85]

    \coordinate (uminus) at (3,0);
    \coordinate (uplus)  at (5,0);
    \coordinate (u)      at (3.85,0.65);
    \coordinate (v)      at (5.75,5.00);

    \draw[line width=1.2pt, dotted] (0,0) -- (7.7,0);

    \draw[line width=1.0pt, dotted] (0,1.75) -- (7.7,1.75);

    \draw[line width=1.2pt, red]
        (uminus)
        .. controls (2.5,1.2) and (4.2,2.8) ..
        (v);

    \draw[line width=1.2pt, red]
        (uplus)
        .. controls (5.0,1.6) and (6.4,3.2) ..
        (v);

    \filldraw[fill=gray!20, draw=black, opacity=0.8] (u) circle (0.42);

    \filldraw (uminus) circle (2.6pt);
    \filldraw (uplus)  circle (2.6pt);
    \filldraw (u)      circle (2.6pt);
    \filldraw (v)      circle (2.6pt);

   \draw[line width=1.2pt, blue]
        (u)
        .. controls (3.85,1.8) and (5.75,3.4) ..
        (v);

    \node[scale=0.8] at (2.95,-0.35) {$u^{-}$};
    \node[scale=0.8] at (4.95,-0.35) {$u^{+}$};
    \node[scale=0.8] at (4.10,0.58) {$u$};
    \node[scale=0.8] at (6.05,5.15) {$v$};

    \node[scale=0.8] at (4.7,0.63) {$\bar{B}_{\delta}^{N}$};

    \draw[decorate,decoration={brace,amplitude=6pt}, line width=0.9pt]
        (7.95,1.75) -- (7.95,0)
        node[midway,right=6pt,scale=0.8] {$\delta^{1/4}N^{3/2}$};

\end{tikzpicture}
\caption{\label{fig:FinalLPP2} Visualization of the site $u^-$ and $u^+$ and the ball $\bar{B}_{\delta}^{N}$ around $u$ used in the proof of Lemma~\ref{lem:ConsistentLPPLocally}, where we see that the geodesic from $u$ to $v$ is sandwiched between the geodesics from $u^-$ to $v$ and from $u^+$ to $v$, respectively.} 
\end{figure}
By Proposition~\ref{pro:TypicallyGoodLPP}, we get that for some $C_1,c_1>0$ and all $N$ large enough
\begin{equation}\label{eq:GoodLocal}
    \P\left( \text{$\gper_{u^{-},v}$ and $\gper_{u^{+},v}$ are $(\delta^{1/4},N)$-good} \right) \geq 1- C\exp(-c \delta^{-1/24}) . 
\end{equation} In particular, on the event in \eqref{eq:GoodLocal}, we see  that for all $z\in [u_2^-,u_2^-+\delta^{1/4} N^{3/2}]$ 
\begin{equation}\label{eq:SplitPoint}
\begin{split}
    \big| \gamma_{u^{-},v}(z)-\gamma_{u^{-},v}(u^-_2) - m_{\ast}^{-1} (z-u^{-}_2) \big| &\leq \frac{N}{8} ,  \\
    \big| \gamma_{u^{+},v}(z)-\gamma_{u^{+},v}(u^+_2)  - m_{\ast}^{-1} (z-u^{-}_2)\big| &\leq \frac{N}{8} .
\end{split}
\end{equation} 
Define the points
\begin{equation*}
    w^{-} := (\gamma_{u^{-},v}(u_2^{-}+\delta^{1/4}N^{3/2}),u_2^{-}+\delta^{1/4}N^{3/2}) , \quad     w^{+} := (\gamma_{u^{+},v}(u_2^{+}+\delta^{1/4}N^{3/2}),u_2^{-}+\delta^{1/4}N^{3/2})
\end{equation*}
lying on the geodesics $\gamma_{u^{-},v}$ and $\gamma_{u^{+},v}$, respectively. Then by \eqref{eq:SplitPoint},  Lemma~\ref{lem:LPPCoupleSmallSlaps} with
$n=\delta^{-1/4}$ for the geodesics $\gamma_{u^{-},w^-}$ and $\gamma_{u^{+},w^+}$, and the almost sure existence of unique maximizers for the last passage times $T$ and $\Tper$, we find a coupling $\mathbf{P}$ between $(N,k)$-periodic last passage percolation and the full-space last passage percolation on $(R^{i,j}_n)^{2}_{\uparrow}$ for some $i \in \Z$ and $j\in \N_0$ so that for some $c_2,C_2>0$ and all $N$ large enough
\begin{equation*}
\begin{split}
     \mathbf{P}&\left( \gamma_{u^{-},v}(z)= \gper_{u^{-},v}(z) \wedge \gamma_{u^{+},v}(z)= \gper_{u^{+},v}(z) \, \forall z\in [u_2-2\delta m_{\ast}^{-1}N^{3/2},u_2+\delta m_{\ast}^{-1}N^{3/2}]\right) \\
     &
     \geq 1- C_2\exp(-C_2\delta^{-1/8}) . 
     \end{split}
\end{equation*}
Recall the local transversal fluctuations $\TF_L(\gamma)$ from \eqref{def:TransversalFluctuationsLocal}, and note that on the event in \eqref{eq:SplitPoint} on the transversal fluctuations at height $u_2^{-}+\delta^{1/4}N^{3/2}$, we find a constant $C_{\ast}>0$ so that
\begin{equation}
\begin{split}
    \sup_{L\in [0,3\delta N^{3/2}]} \big| \TF_{L}(\gper_{u^-,v}) - \TF_{L}(\gper_{u^-,w^{-}})\big|   &\leq  C_{\ast}\delta^{3/4} N \\
    \sup_{L\in [0,3\delta N^{3/2}]} \big| \TF_{L}(\gper_{u^+,v}) - \TF_{L}(\gper_{u^+,w^{+}})\big|   &\leq  C_{\ast}\delta^{3/4} N 
\end{split}
\end{equation} as the geodesics $\gper_{u^-,w^{-}}$ and $\gper_{u^+,w^{+}}$ are subpaths of $\gper_{u^-,v}$ and  $\gper_{u^+,v}$ by construction, but the reference line in the definition of $\TF$ may change. 
By \eqref{eq:SplitPoint} for $\gper_{u^-,w^{-}}$ and $\gper_{u^+,w^{+}}$, and Lemma~\ref{lem:ModerateLocal} with respect to some $n$ of order $\delta^{1/4}N^{3/2}$ and $x=\delta^{1/15}/3$, we find $c_3,C_3>0$ so 
\begin{equation}
\begin{split}
       \P\left( \big|\TF_{3\delta N^{3/2}}(\gper_{u^-,v})\big|\leq \delta^{3/5}N \right) &\geq 1-C_3 \exp(-c_3\delta^{-1/24}) , \\
         \P\left( \big|\TF_{3\delta N^{3/2}}(\gper_{u^+,v})\big|\leq \delta^{3/5}N \right) &\geq 1-C_3 \exp(-c_3\delta^{-1/24}) , 
\end{split}
\end{equation} for all $\delta>0$, and $N$ large enough.
Moreover, Lemma~\ref{lem:TransversalFluctuations} with $n$ of order $\delta N^{3/2}$ and $x=\delta^{1/15}/3$ for the geodesic between $u^-$ and $(\gamma_{u^-,v}(u^-_2+3\delta N^{3/2}),u^-_2+3\delta N^{3/2})$ as well as the geodesic between $u^+$ and $(\gamma_{u^+,v}(u^-_2+3\delta N^{3/2}),u^-_2+3\delta N^{3/2})$ guarantees that we find $c_4,C_4>0$ so that
\begin{equation}
\begin{split}
       \P\left(  \sup_{L\in [0,3\delta N^{3/2}]}\big|\TF_L(\gper_{u^-,v})\big|\leq 2\delta^{3/5}N \right) &\geq 1-C_4 \exp(-c_4\delta^{-1/24}) , \\
          \P\left(  \sup_{L\in [0,3\delta N^{3/2}]}\big|\TF_L(\gper_{u^+,v})\big|\leq 2\delta^{3/5}N \right) &\geq 1-C_4 \exp(-c_4\delta^{-1/24}) , 
\end{split}
\end{equation} for all $\delta>0$, and $N$ large enough.
The ordering of geodesics now ensures that 
\begin{equation}
 \P\left(  \sup_{u'\in \bar{B}^{N}_{\delta}(u)}\sup_{z\in [\delta N^{3/2},3\delta N^{3/2}]}\big|\TF_{u_2+z}(\gper_{u';v})\big|\leq 4\delta^{3/5}N \right) \geq 1-C_5 \exp(-c_5\delta^{-1/24})
\end{equation} for some $c_5,C_5>0$, all $\delta>0$ and $N$ large enough. From the definition of the local transversal fluctuations, we deduce \eqref{eq:Ball1Local}, allowing us to conclude.
\end{proof}

\section{Moderate deviations for asymmetric simple exclusion processes}\label{sec:AppendixASEP}

In this section, we provide the deferred proof of Proposition \ref{pro:ModerateDeviationsSecondClass}.
Since the arguments are analogous to Proposition~4.11 in~\cite{FS:triple}, we focus on describing the necessary changes in the arguments rather than providing full details. Moreover, we will only consider the case of the extended colored ASEP defined between an ASEP on the integers and an unwrapped periodic ASEP. The case of two ASEPs on the integers according to the $(\rho,\theta,1,N)$-reference frame basic coupling is treated analogously. 

 By shift-invariance, we can in the following assume without loss of generality that $\theta=0$. Note that by assumption \eqref{eq:AgreeLocally}, we see that $\xi_0(x)\in \{ \zero,\one\}$ for all $x\in \lbr N \rbr$. We start by introducing three extra types of particles ($\onep,\zerop$, and $\Cup$). Moreover, we define an order of types as
\begin{equation}\label{def:ExtendedPartialOrder}
\one  \succeq \onep \succeq \Aup \succeq  \Cup \succeq \Bup \succeq \zerop \succeq \zero . 
\end{equation} The type $\Cup$ particles will only serve as auxiliary types and will be treated later on. 
Recall that a type $\Aup$ particle  in the extended colored ASEP at site $x$ means that the ASEP on the integers $(\eta_t)_{t \geq 0}$ has a particle at site $x$ while the periodic ASEP $(\bar{\eta}^{\per}_t)_{t \geq 0}$ has an empty site, and vice versa for type $\Bup$ particles. Note that whenever an edge $e$, containing a type $\Aup$ and a type $\Bup$ particle, receives an update, we see that under the basic coupling, the endpoints of $e$ turn by this update into a type $\one$ and a type $\zero$ particle. We aim now to trace back the collisions of type $\Aup$ and type $\Bup$ particles via the following extended disagreement process; see Definition 4.1 in \cite{FS:triple} for a similar construction. 

\begin{definition}\label{def:BasicCouplingExtended}
We define the extended disagreement process $(\xi^{\modi}_t)_{t \geq 0}$ as a Markov process on the state space
\begin{equation}
\{ \one , \onep , \Aup , \Bup , \zerop, \zero \}^{\Z} .
\end{equation}
For all edges $e=\{ x,x+1\}$ with $x\in \Z$, assign a rate $1+q$ Poisson clock. Whenever a clock rings at time $t$, we sample a Uniform-$[0, 1]$-random variable $U$ independently of all previous samples. First, consider $x\in \mathcal{I}^{\rho}_{N}\big(\theta,1+\theta\big)$. We distinguish two cases:
\begin{itemize}
\item[(1)] Assume that $U\leq (1+q)^{-1}$. If $\{ \xi^{\modi}_{t_-}(x), \xi^{\modi}_{t_-}(x+1) \}  \neq \{\Aup,\Bup\}$, then sort the endpoints in $\xi^{\modi}_t$ with respect to  \eqref{def:ExtendedPartialOrder} in increasing order. If $\{ \xi^{\modi}_{t_-}(x), \xi^{\modi}_{t_-}(x+1) \}  = \{\Aup,\Bup\}$, then set $\xi^{\modi}_t(x+1)=\onep$ and $\xi^{\modi}_t(x)=\zerop$. Moreover, for all $y=x \mod N$ with $y \neq x$, we distinguish the following four sub-cases:
\begin{itemize}
    \item[(1a)] If $\xi^{\modi}_{t_-}(y) \in \{\one,\onep\}$ and $\xi^{\modi}_{t_-}(y+1) \in \{\zero,\zerop\}$, set $(\xi^{\modi}_{t}(y),\xi^{\modi}_{t}(y+1))=(\Aup,\Bup)$,
    \item[(1b)] If $\xi^{\modi}_{t_-}(y) \in \{\one,\onep\}$ and $\xi^{\modi}_{t_-}(y+1) = \Aup$, then swap the endpoints.
    \item[(1c)] If $\xi^{\modi}_{t_-}(y) =\Bup$ and $\xi^{\modi}_{t_-}(y+1) = \Aup$, then set $(\xi^{\modi}_{t}(y),\xi^{\modi}_{t}(y+1))=(\zerop,\onep)$,  
    \item[(1d)] If $\xi^{\modi}_{t_-}(y) =\Bup$ and $\xi^{\modi}_{t_-}(y+1) \in \{\zero,\zerop\}$, then swap the endpoints. 
\end{itemize}
In all other cases, leave the configuration unchanged. 
\item[(2)] Assume that  $U> (1+q)^{-1}$. If $\{ \xi^{\modi}_{t_-}(x), \xi^{\modi}_{t_-}(x+1) \}  \neq \{\Aup,\Bup\}$, then sort the endpoints in $\xi^{\modi}_t$ with respect to  \eqref{def:ExtendedPartialOrder} in decreasing order. If $\{ \xi^{\modi}_{t_-}(x), \xi^{\modi}_{t_-}(x+1) \}  = \{\Aup,\Bup\}$, then set $\xi^{\modi}_t(x)=\onep$ and $\xi^{\modi}_t(x+1)=\zerop$. Moreover, for all $y=x \mod N$ with $y \neq x$, we distinguish the following four sub-cases:
\begin{itemize}
    \item[(2a)] If $\xi^{\modi}_{t_-}(y) \in \{\zero,\zerop\}$ and $\xi^{\modi}_{t_-}(y+1) \in \{\one,\onep\}$, set $(\xi^{\modi}_{t}(y),\xi^{\modi}_{t}(y+1))=(\Bup,\Aup)$,
    \item[(2b)] If $\xi^{\modi}_{t_-}(y) = \Aup$ and $\xi^{\modi}_{t_-}(y+1) \in \{\one,\onep\}$, then swap the endpoints.
    \item[(2c)] If $\xi^{\modi}_{t_-}(y) =\Aup$ and $\xi^{\modi}_{t_-}(y+1) = \Bup$, then set $(\xi^{\modi}_{t}(y),\xi^{\modi}_{t}(y+1))=(\onep,\zerop)$,  
    \item[(2d)] If $\xi^{\modi}_{t_-}(y)\in \{\zero,\zerop\}$ and $\xi^{\modi}_{t_-}(y+1)=\Bup$, then swap the endpoints. 
\end{itemize}
In all other cases, leave the configuration unchanged. 
\end{itemize}
Now assume that $x\notin \mathcal{I}^{\rho}_{N}\big(\theta,1+\theta\big)$. Depending on the outcome of the random variable $U$, we proceed with an update of the edge $\{x,x+1\}$ as follows:
\begin{itemize}
\item[(3)] Assume that $U\leq (1+q)^{-1}$. We distinguish the following four sub-cases:
\begin{itemize}
    \item[(3a)] If $\xi^{\modi}_{t_-}(x) \in \{\one,\onep\}$ and $\xi^{\modi}_{t_-}(x+1) \in \{\zero,\zerop\}$, set $(\xi^{\modi}_{t}(x),\xi^{\modi}_{t}(x+1))=(\Bup,\Aup)$,
    \item[(3b)] If $\xi^{\modi}_{t_-}(x) \in \{\one,\onep\}$ and $\xi^{\modi}_{t_-}(x+1) = \Bup$, then swap the endpoints.
    \item[(3c)] If $\xi^{\modi}_{t_-}(x) =\Aup$ and $\xi^{\modi}_{t_-}(x+1) = \Bup$, then set $(\xi^{\modi}_{t}(x),\xi^{\modi}_{t}(x+1))=(\zerop,\onep)$,  
    \item[(3d)] If $\xi^{\modi}_{t_-}(x) =\Aup$ and $\xi^{\modi}_{t_-}(x+1) \in \{\zero,\zerop\}$, then swap the endpoints. 
\end{itemize}
In all other cases, leave the configuration unchanged. 
\item[(4)] Assume that  $U> (1+q)^{-1}$. We distinguish the following four sub-cases:
\begin{itemize}
    \item[(4a)] If $\xi^{\modi}_{t_-}(x) \in \{\zero,\zerop\}$ and $\xi^{\modi}_{t_-}(x+1) \in \{\one,\onep\}$, set $(\xi^{\modi}_{t}(x),\xi^{\modi}_{t}(x+1))=(\Aup,\Bup)$,
    \item[(4b)] If $\xi^{\modi}_{t_-}(x) = \Bup$ and $\xi^{\modi}_{t_-}(x+1) \in \{\one,\onep\}$, then swap the endpoints.
    \item[(4c)] If $\xi^{\modi}_{t_-}(x) =\Bup$ and $\xi^{\modi}_{t_-}(x+1) = \Aup$, then set $(\xi^{\modi}_{t}(x),\xi^{\modi}_{t}(x+1))=(\onep,\zerop)$,  
    \item[(4d)] If $\xi^{\modi}_{t_-}(x)\in \{\zero,\zerop\}$ and $\xi^{\modi}_{t_-}(x+1)=\Aup$, then swap the endpoints. 
\end{itemize}
In all other cases, leave the configuration unchanged. 
\end{itemize}
\end{definition}

The extended disagreement process $(\xi^{\modi}_t)_{t \geq 0}$ has the following intuitive description. For all edges, we perform the same updates as under the $(\rho,\theta,1,N)$-reference frame basic coupling. However, whenever a pair $(\Aup,\Bup)$ or $(\Bup,\Aup)$ receives an update, we 
we immediately replace it by $(\onep,\zerop)$ or $(\zerop,\onep)$ (instead of $(\one,\zero)$ or $(\zero,\one)$ as in the basic coupling), respectively, and then perform the update. 
In the following, let $\xi_0^{\modi}$ agree with the initial configuration $\xi_0$ of the extended colored ASEP from  Proposition~\ref{pro:ModerateDeviationsSecondClass}.
Let $(\kappa_i)_{i \in \lbr N\rbr}$ be a family of independent Bernoulli-$N^{-1/2}$-random variables. We define a new configuration $\tilde{\xi}_0$ by
\begin{equation}\label{eq:xiProcess}
\tilde{\xi}_0(x) := \begin{cases}
\Cup & \text{ if } x\in \big[\frac{N}{11},\frac{N}{11}+N^{3/4}\big], \kappa_x=1 \text{ and } \xi^{\modi}_0(x) = \one  , \\
\xi_0 & \text{ otherwise. }
\end{cases}
\end{equation} In words, we replace every first class particle in $\xi_0$ on the interval $\big[\frac{N}{11},\frac{N}{11}+N^{3/4}\big]$ by a type $\Cup$ particle with probability $N^{-1/2}$. The dynamics $(\tilde{\xi}_t)_{t \geq 0}$ has now the same evolution as $(\xi^{\modi}_t)_{t \geq 0}$ with respect to the partial order \eqref{def:ExtendedPartialOrder}. More precisely, for any update along an edge $e=\{ x,x+1 \}$ at time $s$ such that $(x,s) \in \mathcal{I}^{\rho}_{\delta,N}\big(\frac{1}{8},\frac{7}{8}\big)$, we sort the endpoints in increasing order with probability $(1+q)^{-1}$, and in decreasing order with probability $q(1+q)^{-1}$ (recalling the exception of creating type $\onep$ and type $\zerop$ particles after an update of a type $\Aup$ and type $\Bup$ pair). 
For all updates outside of $\mathcal{I}^{\rho}_{\delta,N}\big(\frac{1}{8},\frac{7}{8}\big)$, we treat the type $\Cup$ particles like type $\one$ particles. Let us stress that by the construction, type $\Aup$ and type $\Bup$ particles can in $(\tilde{\xi}_t)_{t \geq 0}$ only be created at sites to the left of $f(t)$ and to the right of $f(t)+N$ for $t \geq 0$ (recall the function $f$ from \eqref{def:Functionf}). 

We will in the following argue that the type $\Aup$ and type $\Bup$ particles created at sites to the left of $f(t)$ will with high probability not reach $f(s)+N/8$ by time $s$ for all $0 \leq t \leq s \leq \delta N^{3/2}$. The argument that the type $\Aup$ and type $\Bup$ particles created at sites to the right of $f(t)+N$ will with high probability not reach $f(s)+7N/8$ by time $s$ for all $0 \leq t \leq s \leq \delta N^{3/2}$ is similar. 
To this end, for all $t \in [0,\delta N^{3/2}]$, let $X^{\Bup}_t$ denote the location of the rightmost type $\Bup$, which has entered at sites to the left of $f(s)$ for some $s\in [0,t]$. Here, we use the convention that $X^{\Bup}_t= -\infty$ if there is no such particle at time $t$. For all $t \in [0,\delta N^{3/2}]$, let $X^{\Cup}_t$ denote the location of the rightmost type $\Cup$ particle at time $t$, with the convention that $X^{\Cup}_t=-\infty$ if there exists some $s\leq t$ such that a type $\Cup$ is at site $x$ at time $s$ with $(x,s) \notin \mathcal{I}^{\rho}_{\delta,N}\big(\frac{1}{12},\frac{1}{10}\big)$. We define the events
\begin{align*}
A^{\rho}_{\delta,N} &:= \left\{  \tilde{\xi}_t(x) \neq \Cup \, \forall  (x,t) \notin \mathcal{I}^{\rho}_{\delta,N}\Big(\frac{1}{12},\frac{1}{10}\Big) \right\} ,\\
B^{\rho}_{\delta,N} &:= \left\{  X^{\Cup}_t > X^{\Bup}_t \, \forall t \in [0,\delta N^{3/2}] \right\} .
\end{align*}
In words, the event $A^{\rho}_{\delta,N}$ ensures that the type $\Cup$ particles remain with high probability in a given interval, while the event $B^{\rho}_{\delta,N}$ requires that all type $\Bup$ particles created on sites to the left of $f(t)$ stay to the left of at least one type $\Cup$ particle. 
Our goal is to show that there exist constants $C_0,c_0>0$ such that for all $\delta>0$ sufficiently small
\begin{equation}\label{eq:TwoConditionsASEP}
\P\left(A^{\rho}_{\delta,N} \cap B^{\rho}_{\delta,N}\right) \geq 1- C_0\exp(-c_0 \delta^{-1}) .  
\end{equation} 
In particular, this ensures that no type $\Bup$ particle will reach sites to the right of $f(t)+\frac{N}{10}$ at times $t \in [0,\delta N^{3/2}]$. We start by bounding the probability of the event $A^{\rho}_{\delta,N}$, following the arguments from Lemma~4.8 in \cite{FS:triple}.

\begin{lemma}\label{lem:BoundAn}
    There exist constants $c_1,C_1>0$, depending only on $\mathfrak{a}>0$, such that 
\begin{equation}\label{eq:ProbASEP1}
\P\big(A^{\rho}_{\delta,N}\big) \geq 1- C_1\exp(-c_1 \delta^{-1})  
\end{equation} for all $\delta>0$, and $N$ large enough.
\end{lemma}

\begin{proof}[Sketch of proof]
We will in the following only argue that
\begin{equation*}
   \P \left(  \tilde{\xi}_t(x) \neq \Cup \, \forall  (x,t) \, \colon \, x\leq f(t)+\frac{N}{12} ,  t\in [0,\delta N^{3/2}] \right) \geq 1- C_0\exp(-c_0 \delta^{-1}) 
\end{equation*}
as the corresponding upper bound on the location of the type $\Cup$ particles follows similarly.  Consider the process $(\zeta_t)_{t \geq 0}$ given by
\begin{equation}\label{eq:zetaProcess}
\zeta_t(x) := \begin{cases}
1 & \text{ if } \tilde{\xi}_t(x) \in \{\one,\onep,\Aup\}, \\
2 & \text{ if } \tilde{\xi}_t(x)=\Cup, \\
5 & \text{ if } \tilde{\xi}_t(x) \in \{\Bup,\zerop,\zero\}, \\
\end{cases}
\end{equation} for all $x\in \Z$ and $t\geq 0$. Let $\tau_{\Cup}$ denote the first time a type $\Cup$ particle reaches site $f(t)+N/12$ or $f(t)+N/10$. By verifying the marginal transition rates, we see that the process $(\zeta_t)_{t \geq 0}$ has until time $\tau_{\Cup}$  the law of a colored ASEP on the integers. Now let $(\kappa^{\prime}_i)_{i \in \lbr N\rbr}$ be a family of independent Bernoulli-$N^{-1/2}$-random variables, and set
\begin{equation}
    x_{\ast} := \inf\left\{ x\geq \frac{N}{11}-N^{3/4} \, \colon \zeta_0(x)=5 \right\} . 
\end{equation}
 Consider a colored ASEP $(\zeta^{\prime}_t)_{t \geq 0}$ on $\Z$ with initial data
\begin{equation}\label{eq:xiProcess2}
\zeta^{\prime}_0(x) := \begin{cases}
3 & \text{ if } x\in \big[\frac{N}{11}-N^{3/4},\frac{N}{11}\big] \setminus \{x_{\ast}\} \text{ and }\kappa^{\prime}_x=1 \text{ and } \zeta_0(x) = 5  , \\
4 & \text{ if } x=x_{\ast}  ,  \\
\zeta_0 & \text{ otherwise. }
\end{cases}
\end{equation} 
In words, the colored ASEP $(\zeta^{\prime}_t)_{t \geq 0}$ has the same initial data as $(\zeta_t)_{t \geq 0}$, but where we change the some of the type $5$ particles to types $3$ and $4$. 
Let $\tau_{2}$ denote the first time a type $2$ particle reaches site $f(t)+N/12$ or $f(t)+N/10$ in $(\zeta^{\prime}_t)_{t \geq 0}$.  Note that we can evolve $(\zeta_t)_{t \geq 0}$ and $(\zeta^{\prime}_t)_{t \geq 0}$ together under the basic coupling on the integers until time $\tau_{\Cup}$, and that in this case $\tau_2=\tau_{\Cup}$. Hence, it suffices to show that for all $\delta>0$ and $N$ large enough
\begin{equation*}
   \P \left(  \zeta^{\prime}_t(x) \neq 2 \ \forall   (x,t) \, \colon \,  x\leq f(t)+\frac{N}{12} \, \forall t\in [0,\delta N^{3/2}] \right) \geq 1- C_0\exp(-c_0 \delta^{-1}) .  
\end{equation*}
This follows by similar arguments as Lemma~4.8 in~\cite{FS:triple} with $\kappa=0$ and $\theta=\delta^{-1}$ and $y=\delta^{-1/3}$ in their notation. More precisely, let $(X_t^{1})_{t \geq 0}, (X_t^{2})_{t \geq 0},\dots$ denote the positions of the type $2$ particles in $(\zeta^{\prime}_t)_{t \geq 0}$ from left to right. Similarly, let $(Y_t^{1})_{t \geq 0}, (Y_t^{2})_{t \geq 0},\dots$ denote the positions of the type $3$ particles in $(\zeta^{\prime}_t)_{t \geq 0}$ from left to right. The location of the type $4$ particle is denoted by $(Z_t)_{t \geq 0}$. For all $N \in \N$, consider the events 
\begin{equation}
    \begin{split}
        \mathcal{A}_1 &:= \Big\{ \sum_{x\in \Z} \mathds{1}_{ \{ \zeta^{\prime}_0(x)=3 \} } \geq 2N^{\frac{1}{5}} \Big\} \\
\mathcal{A}_2 &:= \Big\{ Y_t^{N^{1/5}} < X_t^{1} \ \forall 0 \leq t \leq  N^{3} \Big\} \\
\mathcal{A}_3 &:= \Big\{ Z_t < Y_t^{N^{1/5}} \ \forall  0 \leq t \leq  N^3 \Big\} . 
    \end{split}
\end{equation}Then by the same arguments as for equation~(4.40) in \cite{FS:triple}, we get that 
\begin{equation}\label{eq:CensoringClaim}
 \lim_{N \rightarrow \infty}  \P\left( \mathcal{A}_1 \cap \mathcal{A}_2 \cap \mathcal{A}_{3}\right) =1 . 
\end{equation} 
Here, a lower bound on the probability of the event $\mathcal{A}_1$ follows from a Chernoff bound, recalling that $\rho \in [\mathfrak{a},1-\mathfrak{a}]$ for some $\mathfrak{a}>0$, while a lower bound on the probability of the events $\mathcal{A}_2$ and $\mathcal{A}_3$ uses the censoring inequality from \cite{PW:Censoring} for the ASEP on $\Z$; see also Lemma~2.9 in \cite{FS:triple}, and Figure~9 in \cite{FS:triple} for an illustration. 
From the definition of $\mathcal{A}_1,\mathcal{A}_2$ and $\mathcal{A}_3$, we see that
\begin{equation}\label{eq:OrderingClaim}
     \P\big( Z_{t} \leq  X_t^{1} \ \forall  t \in [0,N^{3}] \, \big| \, \mathcal{A}_1 \cap \mathcal{A}_2 \cap \mathcal{A}_{3}\big) = 1 . 
\end{equation}
Combining \eqref{eq:CensoringClaim} and \eqref{eq:OrderingClaim}, it remains to show that there exist constants $c_2,C_2>0$ such that
\begin{equation}\label{eq:SecondMax} 
      \P \left(  Z_t \geq f(t)+\frac{N}{12} \ \forall  t\in [0,\delta N^{3/2}] \right) \geq 1- C_2\exp(-c_2 \delta^{-1})  
\end{equation} for all $\delta>0$ and $N$ large enough. Projecting all type $2$ and type $3$ second class particles in  $(\zeta^{\prime}_t)_{t \geq 0}$ to type $1$ particles, \eqref{eq:SecondMax} follows from Lemma~4.6 in \cite{FS:triple} (with the $(\kappa, \theta, y)=(0,\delta^{-1},\frac{1}{200}\delta^{-1/3})$ in their notation).  
Here, let us remark that Lemma~4.6 and Lemma~4.8 in~\cite{FS:triple} asserts that the colored ASEP $(\zeta^{\prime}_t)_{t \geq 0}$ is started from a Bernoulli product measure. However, as noted in Remark~4.10 in~\cite{FS:triple}, this condition can be relaxed only requiring that the laws of all underlying single-species exclusion processes in the colored ASEP $(\zeta^{\prime}_t)_{t \geq 0}$ are stochastically dominated from above and below by Bernoulli-$\big( \rho \pm C^{\prime}N^{-1/2}\big)$-product measures for some constant $C^{\prime}>0$.
\end{proof}

Next, we bound the probability of the event $B^{\rho}_{\delta,N}$  following Lemma~4.12 in \cite{FS:triple}.

\begin{lemma}\label{lem:BoundBn}
    There exist constants $C_3,c_3>0$ such that for all $N$ large enough 
\begin{equation}\label{eq:ProbASEP2}
\P\big(B^{\rho}_{\delta,N} \, \big| \, A^{\rho}_{\delta,N} \big) \geq 1- C_3\exp(-c_3 N^{1/5}) . 
\end{equation}
\end{lemma}

\begin{proof}[Sketch of proof] We assume in the following the same notation as in Lemma~\ref{lem:BoundAn}.
Note that by an application of the exponential Markov inequality, we see that
\begin{equation}\label{eq:ProbASEPmiddle}
\P\Big( \Big| \{ x \in \Big[\frac{N}{12}, \frac{N}{10} \Big] \colon \tilde{\xi}_0(x) = \Cup \Big| \geq  N^{1/5} \Big) \geq 1-C_4\exp(-c_4N^{1/5})
\end{equation} with some constants $c_4,C_4>0$ and all $N$ large enough. Hence, it suffices to show that \begin{equation}\label{eq:ProbASEPNumber}
\P\bigg( B^{\rho}_{\delta,N} \, \bigg| \, \Big| \{ x \in \Big[\frac{N}{12}, \frac{N}{10} \Big] \colon \tilde{\xi}_0(x) = \Cup \Big| \geq  N^{1/5} , A^{\rho}_{\delta,N} \bigg)
 \geq 1- C_5\exp(-c_5 N^{1/5}) 
\end{equation} for some constants $c_5,C_5>0$ and all $N$ large enough. 
To this end, consider the process $(\zeta^{\prime\prime}_t)_{t \geq 0}$ where we set for all $x\in \Z$ and $t\geq 0$
\begin{equation}\label{eq:zetaProcess2}
\zeta^{\prime\prime}_t(x) := \begin{cases}
1 & \text{ if } \tilde{\xi}_t(x) \in \{\one,\onep,\Aup\}, \\
2 & \text{ if } \tilde{\xi}_t(x)=\Cup, \\
3 & \text{ if } \tilde{\xi}_t(x) \in \{\zerop,\Bup\}, \\
4 & \text{ if } \tilde{\xi}_t(x) =\zero . 
\end{cases}
\end{equation} Note that we can interpret the process $(\zeta^{\prime\prime}_t)_{t \geq 0}$  for all $t \in [0,\tau_{\Cup})$ on  $[f(t),f(t)+N]$ as a colored ASEP, i.e.~for any update along an edge $e=\{ x,x+1 \}$ at time $s$ such that $(x,s) \in \mathcal{I}^{\rho}_{\delta,N}\big(0,1\big)$, we sort the endpoints in increasing order with probability $(1+q)^{-1}$, and in decreasing order with probability $q(1+q)^{-1}$. Let us stress that at this point, it is crucial to distinguish between type $\zero$ and type $\zerop$ particles in the definition of the extended disagreement process. 

In order to control the positions of the type $3$ particles in $(\zeta''_t)_{t \geq 0}$ with respect to the type $2$ particles until time $\tau_{\Cup}$, we proceed as follows. For all $t\geq 0$, let $M_2(t)$  be the number of type $2$ particles to the right of $f(t)+N/12$ in $\zeta^{\prime\prime}_t$. Note that on the event $A^{\rho}_{\delta,N}$, we get that $M_2(t)=M_2(0)$ for all $t\in [0,\delta N^{3/2}]$. Similarly, for all $t\geq 0$, let $M_3(t)$ denote the number of type  $3$ particles to the right of $f(t)+N/12$ and to the left of $X_t^{\Bup}$ in $\zeta^{\prime\prime}_t$.
With a slight abuse of notation, we denote the respective positions of the type $2$ particles in $(\zeta^{\prime\prime}_t)_{t \geq 0}$ from left to right by $(X_t^{1})_{t \geq 0,}, (X_t^{2})_{t \geq 0},\dots,(X_t^{M_2(t)})_{t \geq 0}$, and the positions of the type $3$ particles in $(\zeta^{\prime}_t)_{t \geq 0}$ from left to right by $(Y_t^{1})_{t \geq 0}, (Y_t^{2})_{t \geq 0},\dots,(Y_t^{M_3(t)})_{t \geq 0}$. Let 
\begin{equation}
    \tau := \inf\left\{ t\geq 0 \, \colon \,  Y_t^{M_3(t)} > X_t^{N^{1/5}} \right\} 
\end{equation} be the first time a type $3$ particle is to the right of $N^{1/5}$ many type $2$ particles. We claim that 
\begin{equation} \label{eq:ProbASEPmiddleEnd}
\P\bigg( \tau  > \delta N^{\frac{3}{2}} \, \bigg| \, \Big| \{ x \in \Big[\frac{N}{12}, \frac{N}{10} \Big] \colon \tilde{\xi}_0(x) = \Cup \Big| \geq  N^{1/5} , A^{\rho}_{\delta,N} \bigg)
 \geq 1- C_5\exp(-c_5 N^{1/5}) 
\end{equation} for some constants $c_5,C_5>0$, for all $\delta>0$, and all $N$ large enough. 

To see that \eqref{eq:ProbASEPmiddleEnd} holds, we follow the arguments of Lemma~4.12 of \cite{FS:triple}. 
For each $t\geq 0$, we assign a configuration $\bar{\eta}_t \in \{0,1\}^{M_2(t)+M_3(t)}$ by setting $\bar{\eta}_t(x)=1$ if the $x^{\text{th}}$ particle at the locations $X_t^{1},X_t^{2},\dots,X_t^{M_2(t)}$ and $Y_t^{1},Y_t^{2},\dots,Y_t^{M_3(t)}$, counted from left to right, is of type $2$, and $\bar{\eta}_t(x)=0$ otherwise. Consider the set of configurations
\begin{equation}
    \mathfrak{A}_0 := \left\{  \eta \in \{0,1\}^{\Z} \, \colon \,  \sum_{x \leq 0} \eta(x) = \sum_{x > 0} \big(1-\eta_t(x) \big) < \infty  \right\} .  
\end{equation} For all $t\geq 0$, we define configuration $\eta^{\Z}_t \in \mathfrak{A}_0$ by setting
\begin{equation}
    \eta^{\Z}_t(x) := \begin{cases} \bar{\eta}_t(x+M_3(t)) &\text{ if } x \in [-M_3(t)+1,M_2(t)] , \\
    1 &\text{ if } x > M_2(t) , \\
    0 &\text{ otherwise, }    
    \end{cases}
\end{equation}
for all $x\in \Z$. Moreover, we let 
\begin{equation}
    \tau_{\ast} := \inf\left\{ t\geq 0 \, \colon \, \exists  x \geq N^{1/5} \, \colon \, \eta^{\Z}_t(x) = 0  \right\} . 
\end{equation} Note that since $M_2(t)=M_2(0) \geq N^{1/5}$ for all $t\in [0,\tau_{\Cup}]$, we see that 
\begin{equation}
    \min(\tau_{\ast},\tau_{\Cup}) =  \min(\tau,\tau_{\Cup}) . 
\end{equation}
 As in the proof of Lemma 4.12 of \cite{FS:triple}, 
 the process $(\eta^{\Z}_t)_{t \geq 0}$ can until time $\min\big(\tau_{\ast}, \tau_{\Cup} \big)$ be interpreted as an ASEP on the integers with censoring. Moreover, at time $0$, as all type $2$ particles are to the right of the type $3$ particles, we see that $\eta_0^{\Z}(x)=\mathds{1}_{x>0}$ for all $x\in \Z$. Then from the censoring inequality for the ASEP on the integers to bound the location of the rightmost empty site in $(\eta^{\Z}_t)_{t \geq 0}$ -- see \cite[Lemma~2.9]{FS:triple} -- we deduce that  for some $c_6,C_6>0$,
\begin{equation}\label{eq:ProbASEPmiddleFinal}
    \P\left( \tau_{\ast} > \delta N^{\frac{3}{2}} \right) \geq 1- C_6\exp(-c_6 N^{1/5}) 
\end{equation} for all $\delta>0$, and all $N$ large enough. This allows us to conclude  \eqref{eq:ProbASEPmiddleEnd}, and hence~\eqref{eq:ProbASEPNumber}.
\end{proof}

We have now all tools to derive Proposition \ref{pro:ModerateDeviationsSecondClass} following Proposition~4.11 in \cite{FS:triple}.  

\begin{proof}[Proof of Proposition \ref{pro:ModerateDeviationsSecondClass} (Outline)]
Recall for all $t\geq 0$ that $X^{\Bup}_t$ denotes the location of the rightmost type $\Bup$ particles created at sites to the left of $f(t)$. Set
\begin{equation}
    \tau_{\Bup} := \inf\left\{ t \geq 0 \, \colon \, X^{\Bup}_t >  f(t)+\frac{N}{10} \right\}
\end{equation} and note that $\tau_{\Bup}$ is a stopping time for the process $(\xi_t)_{t \geq 0}$. From Lemma~\ref{lem:BoundAn} and Lemma~\ref{lem:BoundBn}, we see that the event
\begin{equation}
    \mathcal{B}_{\Bup} := \left\{ \tau_{\Bup} \geq \delta N^{\frac{3}{2}} \right\}
\end{equation} satisfies for some  $c_1,C_1>0$, for all $\delta>0$, and all $N$ large enough
\begin{equation}\label{eq:BBound}
    \Pex(\mathcal{B}_{\Bup}) \geq 1 - C_1\exp(-c_1 \delta^{-1}) . 
\end{equation}
In words, on $\mathcal{B}$, no type $\Bup$ particle will reach site $N/10+f(t)$ by time $t\leq \delta N^{3/2}$.

It remains to show a similar result for the location $(X^{\Aup}_t)_{t \geq 0}$ of the rightmost type $\Aup$ particle created to the left of $f(t)$. To this end, set 
\begin{equation}
    \tau_{\Aup} := \inf\left\{ t \geq 0 \, \colon \, X^{\Aup}_t >  f(t)+\frac{N}{8} \right\}
\end{equation} and note that $\tau_{\Aup}$ is a stopping time for the process $(\xi_t)_{t \geq 0}$. Recall the ASEP on the integers $(\eta_t)_{t \geq 0}$ and the unwrapped periodic ASEP $(\bar{\eta}^{\per}_t)_{t \geq 0}$ in the definition of $(\xi_t)_{t \geq 0}$. Let $(\kappa^{\prime\prime}_i)_{i \in \lbr N\rbr}$ be a family of independent Bernoulli-$(\delta^{-1/2}N^{-1/2})$-random variables.
Set
\begin{equation}
    \zeta_0(x) := \begin{cases} 2 &\text{ if } x \in \big[\frac{N}{9},\frac{N}{9}+\frac{N}{100}\big] \text{ and } \eta_t(x)=\kappa^{\prime\prime}_x=1 , \\
     \eta_t(x)    &\text{ otherwise, }
    \end{cases}
\end{equation} for all $x\in \Z$, and define $(\zeta_t)_{t \geq 0}$ as a multi-species ASEP on the integers with initial data $\zeta_0$. Let $(\tilde{\eta}_t)_{t \geq 0}$ denote the single-species ASEP, which we obtain as a projection from $(\zeta_t)_{t \geq 0}$ by merging colors $2$ and $0$ into empty sites. 
We let the processes $(\eta_t)_{t \geq 0}$, $(\bar{\eta}^{\per}_t)_{t \geq 0}$ and $(\tilde{\eta}_t)_{t \geq 0}$ 
evolve together under the $(\rho,\theta,1,N)$-reference frame coupling from Definition~\ref{def:BasicCouplingASEPFrame}, and use $(\eta_t)_{t \geq 0}$ and $(\tilde{\eta}_t)_{t \geq 0}$ to obtain the evolution of $(\zeta_t)_{t \geq 0}$. With a slight abuse of notation again write $\Pex$ for this coupling.  

Next, let $\tilde{M}_2$ denote the number of second class particles in $\zeta_0$. Recalling that $\rho \geq \mathfrak{a}>0$ by our assumptions, a Chernoff bound yields  
\begin{equation}
    \P\left(  \tilde{M}_2 \geq \frac{\mathfrak{a}}{200}\delta^{-1/2} N^{1/2} \right) \geq 1 - C_3\exp(-c_3 \delta^{-1/2} N^{1/2})
\end{equation}
for some $c_3,C_3>0$, and all $N$ large enough. Define the stopping time
\begin{equation}
    \tau_{2} := \inf\left\{ t \geq 0 \, \colon \, \exists x\in \Z \setminus \Big[f(t)+\frac{N}{9}-\frac{N}{100},f(t)+\frac{N}{9}+\frac{N}{80}\Big] \colon \zeta_t(x) =2 \right\} 
\end{equation} for the process $(\zeta_t)_{t \geq 0}$. Then by the same arguments as in Lemma 4.8 of \cite{FS:triple}, as also used below \eqref{eq:SecondMax} in the proof of Lemma~\ref{lem:BoundAn}, to bound the probability of the event $\{ \tau_2 > \delta N^{3/2} \}$,  
\begin{equation}
    \mathcal{D}^{\rho}_{\delta,N} := \{ \tau_2 > \delta N^{3/2} \} \cap \left\{ \tilde{M}_2 \geq \frac{\mathfrak{a}}{200}\delta^{-1/2} N^{1/2} \right\}
\end{equation} satisfies for some $c_2,C_2>0$ and any $\delta>0$,
\begin{equation}\label{eq:ProbASEP3}
\Pex\big(\mathcal{D}^{\rho}_{\delta,N} \big) \geq 1- C_2\exp\big(-c_2 \delta^{-1}\big) , 
\end{equation} 
for all $N$ large enough. Define the event
\begin{equation}
    \mathcal{B}_{\Aup} := \{ \tau_{\Aup} \geq \delta N^{3/2} \} . 
\end{equation}
We claim that when the events $\mathcal{D}^{\rho}_{\delta,N}$, $ \mathcal{B}_{\Aup}^{\complement}$ and $\mathcal{B}_{\Bup}$ occur, 
\begin{equation}
   \mathcal{G}^{\rho}_{\delta,N} := \left\{ \sum_{x\in [\frac{1}{10}N,\frac{1}{8}N]} \big(\eta_{T}(x)-\bar{\eta}^{\per}_T(x)\big) \geq \frac{\mathfrak{a}}{200}\delta^{-1/2} N^{1/2} \right\}  
\end{equation}
occurs as well. To see this, note that whenever an edge $\{x,x+1\}$ satisfying 
\begin{equation*}
(\zeta_{t_-}(x),\eta_{t_-}(x),\bar{\eta}^{\per}_{t_-}(x))= (1,1,0) \quad \text{ and } \quad (\zeta_{t_-}(x+1),\eta_{t_-}(x+1),\bar{\eta}^{\per}_{t_-}(x+1))= (2,1,1)
\end{equation*}
 is updated at time $t\geq 0$, then this results with probability $(1+q)^{-1}$ in 
\begin{equation*}
    (\zeta_{t}(x),\eta_{t}(x),\bar{\eta}^{\per}_{t}(x))= (2,1,0) \quad \text{ and } \quad (\zeta_{t}(x+1),\eta_{t}(x+1),\bar{\eta}^{\per}_t(x+1))= (1,1,1) , 
\end{equation*} and with probability $q(1+q)^{-1}$ in 
\begin{equation*}
    (\zeta_{t}(x),\eta_{t}(x),\bar{\eta}^{\per}_{t}(x))= (1,1,1) \quad \text{ and } \quad (\zeta_{t}(x+1),\eta_{t}(x+1),\bar{\eta}^{\per}_t(x+1))= (2,1,0) . 
\end{equation*}
In words, whenever a site containing a second class particle is reached for the first time by a type $\Aup$ particle in $(\xi_t)_{t \geq 0}$, then the type $\Aup$ gets coupled to a second class particle. Moreover, verifying the marginal transitions, this pairing is preserved under the events $\mathcal{B}_{\Bup}$ and $\mathcal{D}^{\rho}_{\delta,N}$ until time  $T:=\delta N^{3/2}$.
 Thus, provided that the events $\mathcal{D}^{\rho}_{\delta,N}$ and  $ \mathcal{B}_{\Aup}^{\complement}$ occur (which ensures that all second class particles remain in the interval $[N/10,N/8]$ until time $T$ and that we see at least one type $\Aup$ particle to the right of $N/8$ by time $T$, while due to the event $\mathcal{B}_{\Bup}$ no type $\Bup$ particles have entered the interval $[N/10,N/8]$ by time $T$), every second class particle in $\zeta_{T}$ must be paired with a site that is empty in $\bar{\eta}^{\per}_{T}$ and occupied in $\eta_T$, that is, a type $\Aup$ particle. 
 Since $\xi_T$ contains, due to the event $\mathcal{B}_{\Bup}$, only particles of types $\Aup,\one$ in the interval $[N/10,N/8]$, we see 
 \begin{equation}\label{eq:GImply}
     \sum_{x\in [\frac{1}{10}N,\frac{1}{8}N]} \big(\eta_{T}(x)-\bar{\eta}^{\per}_T(x)\big) = \sum_{x\in [\frac{1}{10}N,\frac{1}{8}N]} \mathds{1}_{\{\xi_T(x)=\Aup\}} .
 \end{equation}
 Noting that the number of type $\Aup$ particles in $[N/10, N/8]$ is at least equal to the number of second class particles in $\zeta_T$ in this interval, we see that the event in \eqref{eq:GImply} indeed implies that $\mathcal{G}^{\rho}_{\delta,N}$ occurs. This gives the desired claim. 
 
 Since $\eta_{T}$ and $\bar{\eta}^{\per}_T$ are by assumption \eqref{eq:BernoulliDomination} and attractivity stochastically dominated by Bernoulli-$(\rho \pm C N^{-1/2})$ product measures, we get from a standard moderate deviation estimate for the number of particles contained in $[\frac{1}{10}N,\frac{1}{8}N]$ for the configurations $\eta_{T}$ and $\bar{\eta}^{\per}_T$ that
\begin{equation}
    \Pex(\mathcal{G}^{\rho}_{\delta,N}) \leq  C_3\exp(-c_3\delta^{-1} )
\end{equation}
for some $c_3,C_3>0$, depending only on $\mathfrak{a},C$ and $q$.
Thus, together with \eqref{eq:BBound} and \eqref{eq:ProbASEP3}, we see that for some $c_4,C_4>0$, for all $\delta>0$, and all $N$ large enough
\begin{equation}\label{eq:BBoundnew}
\begin{split}
     \Pex(\mathcal{B}_{\Aup}^{\complement}) &\leq \Pex(\mathcal{B}_{\Aup}^{\complement} \cap \mathcal{B}_{\Bup} \cap \mathcal{D}^{\rho}_{\delta,N}) + \Pex(\mathcal{B}_{\Bup}^{\complement}) + \Pex\big(\big(\mathcal{D}^{\rho}_{\delta,N} \big)^{\complement}\big) \\
     & \leq \Pex(\mathcal{G}^{\rho}_{\delta,N})  +  \Pex(\mathcal{B}_{\Bup}^{\complement}) + \Pex\big(\big(\mathcal{D}^{\rho}_{\delta,N} \big)^{\complement}\big) \\  
     &\leq C_4\exp(-c_4 \delta^{-1}) . 
\end{split}
\end{equation}
Combining \eqref{eq:BBound} and \eqref{eq:BBoundnew}, this gives for some $C_5,c_5>0$ 
\begin{equation}\label{eq:ProbASEP5}
\Pex\left( \max\big(X^{\Aup}_t, X^{\Bup}_t\big) \leq f(t) + \frac{N}{8} \ \forall t \in [0,\delta N^{3/2}] \right) \geq 1- C_5\exp\big(-c_5 \delta^{-1}\big) 
\end{equation} for all $\delta>0$ and $N$ large enough, allowing us to conclude. 
\end{proof}

 \bibliographystyle{alpha}

\bibliography{ASEPCircle}

\end{document}